\pdfoutput=1
\documentclass[11pt]{article}
\usepackage{amsmath,amssymb,amsthm,mathtools}
\usepackage[margin=1in]{geometry}
\usepackage{booktabs}
\usepackage{graphicx}
\usepackage{algorithm}
\usepackage{algpseudocode}
\usepackage[round,authoryear]{natbib}
\usepackage{xcolor}
\usepackage[colorlinks=true,linkcolor=blue,citecolor=blue,urlcolor=blue]{hyperref}
\usepackage{caption}
\usepackage{enumitem}
\newcommand{\bb}{\mathbb}

\newcommand{\cD}{\mathcal{D}}

\newcommand{\cS}{\mathcal{S}}
\newcommand{\cT}{\mathcal{T}}

\newcommand{\ind}[1]{\mathbf{1}\{#1\}}
\newcommand{\E}{\bb{E}}
\newcommand{\Pbb}{\bb{P}}
\newcommand{\Var}{\operatorname{Var}}
\newcommand{\kl}{\operatorname{kl}}

\newcommand{\Bern}{\operatorname{Bern}}
\newcommand{\Pois}{\operatorname{Pois}}
\newcommand{\Cat}{\operatorname{Cat}}

\newcommand{\argmax}{\operatorname*{arg\,max}}
\newcommand{\Bfun}{\mathrm{B}}
\newcommand{\ER}{\mathrm{ER}}
\newcommand{\DC}{\mathrm{DC}}

\newcommand{\suppname}{appendix}
\newcommand{\suppref}[1]{Appendix~\ref{#1}}

\theoremstyle{plain}
\newtheorem{theorem}{Theorem}[section]
\newtheorem{lemma}[theorem]{Lemma}
\newtheorem{proposition}[theorem]{Proposition}
\newtheorem{corollary}[theorem]{Corollary}
\theoremstyle{definition}
\newtheorem{assumption}[theorem]{Assumption}

\theoremstyle{remark}
\newtheorem{remark}[theorem]{Remark}

\title{Evidence thresholds of collapsed structured block models}
\author{Marios Papamichalis\thanks{\raggedright Human Nature Lab, Yale University, New Haven, CT 06511, USA. \href{mailto:marios.papamichalis@yale.edu}{\nolinkurl{marios.papamichalis@yale.edu}}}
\and Regina Ruane\thanks{\raggedright Department of Statistics and Data Science, The Wharton School, University of Pennsylvania, Philadelphia, PA, USA. \href{mailto:ruanej@wharton.upenn.edu}{\nolinkurl{ruanej@wharton.upenn.edu}}}}
\date{}

\begin{document}
\maketitle

\begin{abstract}
A stochastic block model explains a network by splitting its nodes into communities and a matrix of tie probabilities between them. Collapsing those probabilities out scores each split by one number, its marginal likelihood, and every decision the model takes, how many communities, whether a weak one survives, what generates the ties inside it, is a comparison of two such numbers. The comparisons fail in known ways: weak structure is dropped, small communities are merged, a group of hubs is returned in place of the real groups. Each failure has a threshold, and the thresholds have not been computed. We compute them. We prove that the best single partition keeps a planted structure only above the first-moment bound, which sits above the Kesten–Stigum threshold for up to ten communities, while the evidence summed over all partitions switches at Kesten–Stigum itself: between the two, the posterior favours structure that every partition a sampler can report rejects. We prove three further thresholds, one for merging small communities, one for the density at which the representation alone decides, one for the degree prior that cuts a community along its degrees. The empirical finding confirms the theory: on planted partitions and on thirteen labelled networks, from the karate club to citation and co-purchase networks of up to 19,717 nodes, each threshold falls where it is predicted, against likelihood, spectral, penalised and nonparametric alternatives. 
\end{abstract}

\section{Introduction}\label{sec:intro}

A stochastic block model \citep{Holland1983} explains the ties of a network by a partition of its
nodes into communities and a matrix of tie probabilities between communities, and assumes that inside
a community every pair of nodes is equally likely to be tied. In a network of political blogs each
camp contains a few blogs that link to hundreds of others and many that link to a handful, and a
homogeneous block model fitted to it returns a community of hubs and a community of peripheral
blogs in place of the two camps, because splitting a camp by degree raises the likelihood more than
separating the camps does. Community structure is therefore confounded with the mechanism that
generates ties inside communities and, as this paper shows, with the prior placed on the degrees
when the mechanism is corrected for.

Bayesian block models answer such questions with one number, the marginal likelihood of a partition
with the block parameters integrated out, which we call the evidence. Allocation samplers and the
integrated classification likelihood optimise it, and every decision they take, how many communities,
whether a weak partition is kept, which mechanism explains a community, is a comparison of evidences.
The thesis of this paper is that each of these decisions is taken at a threshold that can be
computed, that the three objects which share the objective, the comparison of two labellings, the
maximum a posteriori labelling and the evidence summed over labellings, switch at different
thresholds, and that the failures observed in practice, the loss of weak partitions, the merging of
small communities, the merging of dense classes and the splitting of hub-dominated communities, are
properties of the objective and of its priors and persist under any optimiser of the same objective.
An analyst
who knows the thresholds knows, before any sampler is run, in which region of signal, number of
communities, density and degree heterogeneity a fitted partition is a reliable estimate.

Three lines of work come close to these questions and leave them open. \citet{Banks2016},
\citet{Abbe2015} and \citet{CojaOghlan2018} locate the information-theoretic threshold of the planted
partition between the Kesten--Stigum threshold of \citet{Decelle2011} and a first-moment bound, and
\citet{Peixoto2013} derives the resolution limit of the microcanonical description length, but no
result places the maximum a posteriori criterion of a conjugate block model, the objective of
\citet{Nowicki2001}, \citet{McDaid2013} and \citet{Come2015}, relative to either. \citet{Yan2014}
compare the homogeneous and the degree-corrected block model with the partition held fixed and the
block parameters not integrated, and \citet{Peixoto2017} corrects degrees for every node at once
under a prior on the degree sequence whose effect on the partition has not been analysed.
\citet{Zhang2020} and \citet{Amongero2025} enforce assortative structure through the prior on
between-community connectivity without quantifying what the constraint costs when the data violate
it.

This paper studies a collapsed structured block model (CSBM), a collapsed block model in which every
community carries a latent mechanism type, an Erd\H{o}s--R\'enyi type or a degree-corrected type,
whose parameters are integrated out, with degrees counted inside the community or over the whole
network and with one of three priors on the degree sequence, and in which the prior on
between-community connectivity may be truncated at a known cap or shrunk towards it. The evidence of
a partition remains a product of closed-form block terms, so the number of communities, the mechanism
inside each community and the admissible between-community connectivity are decided by one
evidence, and a single-site Gibbs update of a node costs special-function evaluations linear in its
degree plus the square of the number of communities. The homogeneous instance is the conjugate block
model of \citet{Nowicki2001} and the degree-corrected instance with total degrees is the collapsed
model of \citet{Peixoto2017}, so what we prove about them applies to that literature. The technical
move behind the results is a comparison of evidences through profile likelihoods: a Stirling
expansion places every conjugate block evidence within a logarithmic term of the profile
log-likelihood of its block, and each result reduces to a population divergence between the
generating rates and their best approximation in the competing class, and to the cost that the block,
partition and degree priors charge.

The contributions are the following. First, Theorem~\ref{thm:threshold} shows that the maximum a
posteriori criterion keeps a planted partition over the one-community partition only when its fit
gain per node exceeds the logarithm of the number of communities, and Corollary~\ref{cor:icl}
identifies this threshold with the first-moment bound of \citet{Banks2016}: the criterion, and with it
the exact integrated classification likelihood of \citet{Come2015}, pays the entropy of the union
bound. Theorem~\ref{thm:estimator} extends the statement from the criterion to the estimator, the
maximiser over all labellings, with an overlap guarantee, and Theorem~\ref{thm:two} shows that the
evidence summed over labellings switches at the Kesten--Stigum threshold instead, so that between
the two thresholds the evidence detects structure that every single labelling, the planted one
included, rejects. Section~\ref{sec:Kexp} measures both on the same networks by thermodynamic
integration. Second, Theorem~\ref{thm:resolution} and Corollary~\ref{cor:res} give the resolution
limit of the conjugate evidence, of order $\sqrt{n/\log n}$ with a per-block cost $(a+\tfrac12)\log\nu$
set by the block prior, and Section~\ref{sec:Kres} shows the flat evidence saturating at the
predicted $K_{\mathrm{res}}$ while the microcanonical prior of \citet{Peixoto2017} resolves on. Third,
Proposition~\ref{prop:family} and Corollary~\ref{cor:window} give the exact price that the Poisson
representation behind degree correction charges a binary block, convex in its density, and the window
of community sizes in which that price alone decides whether two classes are merged, which is where
the two representations disagree on the football and school networks. Fourth,
Theorem~\ref{thm:split} shows that the typed evidence selects the generating mechanism of a large
community, does not split a heterogeneous community into homogeneous blocks and does not absorb a
genuine two-community structure, Proposition~\ref{prop:finite} gives its finite-size form, and
Propositions~\ref{prop:degsplit} and~\ref{prop:gtyped} show that the uniform degree prior of
\citet{Peixoto2017} pays a gain linear in the community size for cutting a community along its
degrees, that a partition-invariant degree prior pays nothing, and that the global typed instance
under such a prior selects the mechanism from total degrees with the margins of the theorem, which
halves the over-splitting of degree-corrected communities with $K$ inferred and keeps their type. Fifth,
Proposition~\ref{thm:cap} shows that a hard cap on between-community connectivity adds at most a
bounded number of nats for separation and penalises a violated cap at a rate linear in the block
size, so that it does not protect boundaries against contamination, while shrinkage of fixed strength
has an effect bounded uniformly in the block size (Proposition~\ref{prop:soft}). Sixth, the
experiments of Sections~\ref{sec:simulations} and~\ref{sec:real} test every result with paired seeds
and bootstrap intervals against the degree-corrected likelihood of \citet{Karrer2011}, the test of
\citet{Yan2014}, the estimators of $K$ of \citet{Hu2020}, \citet{Wang2017}, \citet{Le2022} and
\citet{Li2020}, the partition priors of \citet{Legramanti2022} and \citet{Miller2018} and the
spectral and modularity methods a practitioner would use, on planted partitions and on labelled
networks, and locate the regime where the typed evidence improves on every alternative, heavy
within-community degree tails at moderate or strong signal and mixtures of mechanisms, the regimes
where it ties, and the regimes where it loses. Table~\ref{tab:rule} states the decision rule that
results.

\begin{table}[tbp]\centering\footnotesize
\caption{Decision rule for the collapsed evidence. Each row names a regime, the instance or method to
use, the result that justifies the choice and the section whose experiment confirms it.}
\label{tab:rule}
\begin{tabular}{p{4.3cm}p{4.3cm}p{3.9cm}}\toprule
Regime & Use & Result, experiment\\\midrule
Signal between the Kesten--Stigum threshold and the first-moment bound & the evidence over labellings or a spectral method with $K$ given, in place of the MAP of a conjugate evidence & Theorems~\ref{thm:threshold} to~\ref{thm:two}, Section~\ref{sec:Kexp}\\
Many small communities ($K$ near $K_{\mathrm{res}}$) & Jeffreys block priors or the uniform prior on the edge-count matrix & Theorem~\ref{thm:resolution}, Corollary~\ref{cor:res}, Section~\ref{sec:Kres}\\
Dense binary classes inside the window of Corollary~\ref{cor:window} & the Bernoulli representation in place of a Poisson degree-corrected evidence & Proposition~\ref{prop:family}, Section~\ref{sec:real-family}\\
Heavy within-community degree tail, moderate or strong signal, or mixed mechanisms & the typed evidence & Theorem~\ref{thm:split}, Section~\ref{sec:phase}\\
Heterogeneity on between-community ties (hub-dominated networks) & the global typed instance with the global degree prior & Propositions~\ref{prop:degsplit} and~\ref{prop:gtyped}, Section~\ref{sec:real-hubs}\\
$K$ inferred under degree correction & a partition-invariant degree prior in place of the uniform or histogram prior & Proposition~\ref{prop:degsplit}, Sections~\ref{sec:phase} and~\ref{sec:real-hubs}\\
Known assortativity, contamination possible & soft shrinkage in place of a hard cap & Propositions~\ref{thm:cap} and~\ref{prop:soft}, Section~\ref{sec:capexp}\\
\bottomrule\end{tabular}
\end{table}

The paper makes no comparison with neural methods, does not run the agglomerative sampler of
\citet{Peixoto2014} (its baselines are the microcanonical evidence with the uniform edge-count prior
inside our sampler and a script for graph-tool in the released code), confines mechanism types to
diagonal blocks, and reports results only for the two exact types and, on fixed partitions, the
Bernoulli degree-corrected type. Proofs are in the \suppname.

\section{Related work}\label{sec:related}

Conjugate collapsing of the block parameters of a stochastic block model and Gibbs sampling of the
allocations are due to \citet{Nowicki2001}, allocation samplers with an unknown number of communities
to \citet{Nobile2007} and \citet{McDaid2013}, the latent block model to \citet{Wyse2012}, and the
exact integrated classification likelihood, the same objective optimised greedily with Jeffreys
priors, to \citet{Come2015} and \citet{Come2021}. Variational alternatives are \citet{Daudin2008},
\citet{Mariadassou2010} and \citet{Latouche2012}, and Bayesian nonparametric formulations are given
by \citet{Kemp2006}, by \citet{Geng2019} with the mixture of finite mixtures of \citet{Miller2018},
and by \citet{Legramanti2022} with the Pitman--Yor and Gnedin priors, which enter the experiments as
alternatives to the Dirichlet--multinomial and CRP priors of the theorems. The sampler of this paper
belongs to this family and adds the per-community mechanism type, the degree priors, the constrained
off-diagonal prior and the analysis of the evidence.

Detectability thresholds are due to \citet{Decelle2011}, proved by \citet{Mossel2015,Mossel2018} and
\citet{Massoulie2014} and reviewed by \citet{Abbe2018}. \citet{Banks2016} bound the
information-theoretic threshold from above by a first-moment argument, whose upper bound
Corollary~\ref{cor:icl} identifies with the threshold of the MAP criterion, \citet{Abbe2015} show
detection below the Kesten--Stigum threshold for four or more communities, \citet{CojaOghlan2018}
give the exact threshold and the condensation phase, which Theorem~\ref{thm:two} uses, and
\citet{Chin2025} study the same first-moment quantities for the maximum-likelihood labelling. None of
these places the conjugate evidence of an allocation sampler, in which the block parameters are
integrated and a partition prior is paid, relative to the thresholds, which is what
Theorems~\ref{thm:threshold} to~\ref{thm:two} do. The resolution limit of the microcanonical
description length is derived by \citet{Peixoto2013} and removed by the nested prior of
\citet{Peixoto2014}. Theorem~\ref{thm:resolution} is the conjugate-evidence form of that limit with
the per-block cost set by the block prior.

Degree-corrected block models are due to \citet{Karrer2011}. \citet{Peixoto2017} gives exact
microcanonical marginals for the degree-corrected model with the uniform and the degree-histogram
priors on the degree sequence, whose effect on the partition Proposition~\ref{prop:degsplit}
analyses, and \citet{Peixoto2020} gives the merge-split sampler that the released graph-tool script
runs. The decision between a homogeneous and a degree-corrected block model is the model selection
problem of \citet{Yan2014}, whose likelihood-ratio test with a bootstrap null enters the experiments
per community, and which Theorem~\ref{thm:split}(a) takes inside a single community with the
parameters integrated. \citet{Peixoto2018}, \citet{Aicher2015} and \citet{Ludkin2020} extend the collapsed treatment to
edge weights, and the multilayer, signed, dyad-state and zero-inflated block models of
\citet{Han2015}, \citet{Peixoto2015}, \citet{Stanley2016}, \citet{Barbillon2017}, \citet{Jiang2015},
\citet{Holland1981}, \citet{Wang1987} and \citet{Lu2025} are instances of the conjugate template of
\suppref{app:catalogue}. \citet{Schweinberger2015}, \citet{Vu2013}, \citet{Fosdick2019} and
\citet{Peixoto2022} model dependence inside blocks, and the type library of Section~\ref{sec:types}
is a collapsed, per-community version of this programme. Assortativity constraints are studied by
\citet{Zhang2020} and \citet{Amongero2025}, and the truncated conjugate cap of Section~\ref{sec:cap}
is the simplest such constraint.

Model selection for the number of communities is treated by \citet{Yan2014}, \citet{Newman2016},
\citet{Saldana2017}, \citet{Wang2017}, \citet{Riolo2017} and \citet{Hu2020}, by network
cross-validation \citep{Chen2018,Li2020} and spectrally by \citet{Le2022}. The penalised likelihood
of \citet{Wang2017}, the corrected BIC of \citet{Hu2020}, the edge cross-validation of
\citet{Li2020}, the spectral estimator of \citet{Le2022} and the microcanonical evidence of
\citet{Peixoto2017} enter Section~\ref{sec:Kres} as alternatives to the collapsed evidence.
Posterior consistency for the stochastic block model is established by \citet{vanderPas2018},
\citet{Geng2019}, \citet{Ghosh2020} and \citet{Jiang2021}, profile-likelihood consistency by
\citet{Bickel2009}, \citet{Zhao2012} and \citet{Choi2012}, minimax rates by \citet{Gao2017} and
spectral consistency by \citet{Lei2015}, and Proposition~\ref{thm:concentration} in the \suppname{}
adds nothing to these results in generality. Collapsed Gibbs sampling is analysed by \citet{Liu1994}
and \citet{LiuWongKong1994}, split and merge moves are due to \citet{Jain2004}, thermodynamic
integration to \citet{Friel2008} and parallel tempering to \citet{Geyer1991}. The spectral, modularity and flow methods of the experiments are due to
\citet{Saade2014}, \citet{Jin2015}, \citet{Qin2013}, \citet{Amini2013}, \citet{Blondel2008},
\citet{Traag2019} and \citet{Rosvall2008}, point estimation of partitions and credible balls in the
variation of information to \citet{Wade2018}, and the distinction between metadata and ground truth
to \citet{Peel2017}.

\section{Model and inference}\label{sec:model}

This section defines the collapsed structured block model (CSBM), its two mechanism types with three
degree priors, the constraints on between-community ties and the sampler. Throughout, a \emph{network} is the observed object, a
\emph{community} is a set of nodes, and a \emph{block} is a pair of communities together with the dyads
between them.

\subsection{Notation}\label{sec:data}
Let $V=[n]=\{1,\dots,n\}$ be the nodes, let $\cD=\{(i,j):1\le i<j\le n\}$ for undirected data and
$\cD=\{(i,j):i\neq j\}$ for directed data, and let $Y=\{y_{ij}:(i,j)\in\cD\}$ be the observed dyads,
binary or count-valued in the main text and of the other types of \suppref{app:catalogue} there. For
binary data we also write $A$ for the adjacency matrix. A partition is a vector
$z\in[K_{\max}]^n$, where $K_{\max}$ is the number of available labels, with communities
$C_k=\{i:z_i=k\}$ of sizes $n_k=|C_k|$ and $K(z)$ nonempty communities. For undirected data the block
containing dyad $(i,j)$ is the unordered pair $\{z_i,z_j\}$,
\[
\cD_{rs}(z)=\{(i,j)\in\cD:\ \{z_i,z_j\}=\{r,s\}\},\qquad N_{rs}(z)=|\cD_{rs}(z)|,\qquad r\le s,
\]
so that $N_{rr}=\binom{n_r}2$ and $N_{rs}=n_rn_s$ for $r<s$. For directed data blocks are ordered pairs
and $N_{rr}=n_r(n_r-1)$. The block sufficient statistic is $S_{rs}(z)=\sum_{(i,j)\in\cD_{rs}(z)}T(y_{ij})$
for a family-specific statistic $T$. For binary data $S_{rs}=e_{rs}$, the number of edges in the block.
Blocks with $r=s$ are \emph{diagonal} and blocks with $r<s$ are \emph{off-diagonal}. The
\emph{evidence} of a partition is the marginal likelihood $p(Y\mid z)$, and the evidence of a block or
of a community is the corresponding factor of $p(Y\mid z)$. The symbol $\Bfun(\cdot,\cdot)$ denotes the
beta function, $I_x(\alpha,\beta)$ the regularised incomplete beta function, and $(a,b)$ the two
hyperparameters of a conjugate block prior. ARI is the adjusted Rand index, MAP the maximum a
posteriori partition and AUC the area under the receiver operating characteristic curve. $O$, $o$
and $\Theta$ refer to $n\to\infty$ or to the community size $\nu\to\infty$ as stated, $O_p$ and $o_p$
are their in-probability versions, and $\asymp$ denotes two-sided bounds up to constants. Unless
stated otherwise, every block prior is Beta$(1,1)$ or Gamma$(1,1)$, the Jeffreys instances use
Beta$(\tfrac12,\tfrac12)$, the partition priors use concentration $\alpha=1$, and the type prior is
$\pi(\ER)=\pi(\DC)=\tfrac12$.

\subsection{Conjugate block evidences}\label{sec:conjugate}
Every off-diagonal block is an exponential-family block with a conjugate prior, so that
\[
p(Y_{rs}\mid z)=\int\prod_{(i,j)\in\cD_{rs}(z)}p(y_{ij}\mid\theta)\,\pi(\theta)\,d\theta ,
\]
where $Y_{rs}$ is the restriction of $Y$ to $\cD_{rs}(z)$, is available in closed form as a function of
$(S_{rs},N_{rs})$. The instances used in the main text are the Beta--Bernoulli
block $\Bfun(a+e,b+N-e)/\Bfun(a,b)$ for a block with $e$ edges among $N$ dyads and the Gamma--Poisson
block, and Table~\ref{tab:catalogue} of \suppref{app:catalogue} collects these and the other conjugate
blocks of the released code.

\subsection{Mechanism types and type collapsing}\label{sec:types}
Each community $C_k$ carries a latent type $t_k\in\cT$ with prior $\pi(t)$ and a within-community
evidence $p_t(Y_{C_k})=\int p_t(Y_{C_k}\mid\Theta)\pi_t(\Theta)\,d\Theta$, where $Y_{C_k}$ is the
restriction of $Y$ to the dyads inside $C_k$. The type-collapsed community evidence and the evidence
of a partition are
\begin{equation}\label{eq:csbm}
p(Y_{C_k})=\sum_{t\in\cT}\pi(t)\,p_t(Y_{C_k}),\qquad
p(Y\mid z)=\prod_{k=1}^{K(z)}p(Y_{C_k})\prod_{r<s}p(Y_{rs}\mid z),
\end{equation}
and the posterior is $p(z\mid Y)\propto p(z)p(Y\mid z)$ for a partition prior $p(z)$, which is the
Dirichlet--multinomial prior with $K_{\max}$ labels, the Chinese restaurant process (CRP) prior with
$K_{\max}$ labels, or the uniform prior of \citet{Peixoto2017}. The posterior type of a community is
$p(t_k=t\mid Y,z)\propto\pi(t)p_t(Y_{C_k})$.

We use the multigraph (Poisson) representation of the within-community network with $E$ edges among
$N=\binom{n_k}2$ dyads and within-community degrees $\mathbf d=(d_i)_{i\in C_k}$, and two exact types.
\begin{itemize}[leftmargin=*,itemsep=1pt]
\item \emph{Erd\H{o}s--R\'enyi type} $t=\ER$: $y_{ij}\mid\lambda\sim\Pois(\lambda)$ and
$\lambda\sim\mathrm{Gamma}(a,b)$, so that $p_{\ER}(Y_{C_k})=b^a\Gamma(a+E)/\{\Gamma(a)(b+N)^{a+E}\}$
up to the factor $\prod_{i<j}(y_{ij}!)^{-1}$. This factor is common to both types and to every
partition and is dropped throughout. All evidences of a given data set are taken with respect to the
same dominating measure, so type and partition comparisons are well defined, while evidences of
different data types, binary against count-valued for instance, are not comparable with each other.
\item \emph{Degree-corrected type} $t=\DC$: write $p_{\ER}=P(E)\,P(Y_{C_k}\mid E)$ with $P(Y_{C_k}\mid E)$
the multinomial placement of $E$ edges among $N$ dyads. The degree-corrected type keeps the same $P(E)$
and replaces the placement by the configuration-model count of \citet{Peixoto2017},
\begin{equation}\label{eq:dc}
P_{\DC}(Y_{C_k}\mid \mathbf d,E)=\frac{\prod_{i}d_i!\;(2E)!!}{(2E)!\,\prod_{i<j}y_{ij}!},\qquad
P(\mathbf d\mid E)=\binom{n_k+2E-1}{2E}^{-1},
\end{equation}
so that $p_{\DC}(Y_{C_k})=P(E)P_{\DC}(Y_{C_k}\mid \mathbf d,E)P(\mathbf d\mid E)$. The factor $P(\mathbf d\mid E)$ is the uniform
prior on degree sequences, whose removal gives the model conditioned on the degrees of
\citet{Karrer2011}.
\end{itemize}
Both evidences depend on the community only through $(E,N,n_k,\sum_i\log d_i!)$, so a single-site
move changes them through the within-community degree of the moving node and those of its neighbours,
at a cost of $O(\deg(i)+K^2)$ special-function evaluations (\suppref{app:catalogue}). The same
construction with node degrees taken over the whole network gives the collapsed degree-corrected
block model of \citet{Peixoto2017}, which we call the \emph{global-DC} instance. We call the
homogeneous Beta--Bernoulli instance the \emph{flat} evidence and the instance with
$\cT=\{\ER,\DC\}$ the \emph{typed} evidence.

The typed evidence corrects degrees inside communities only. A second exact pair of types makes the
per-community choice with total degrees: in the multigraph counting of \citet{Peixoto2017} the
$e_r=\sum_{i\in C_r}k_i$ half-edges of community $r$ are placed uniformly over its $n_r$ nodes, with
factor $n_r^{-e_r}$, or according to the total degrees $(k_i)_{i\in C_r}$, with factor
$\prod_{i\in C_r}k_i!/e_r!$ times a prior on the degree sequence $P(\mathbf k_r\mid e_r)$, and the
block terms are common to both. Collapsing the choice per community with prior $\pi(t)$ gives the
\emph{global typed} instance, whose limits $\pi(\DC)\to1$ and $\pi(\DC)\to0$ are the global-DC
instance and the microcanonical homogeneous model, and which is the instance for networks on which
heterogeneity acts on the between-community ties as well. Three degree priors are used: the uniform
prior $1/\binom{n_r+e_r-1}{e_r}$ of \citet{Peixoto2017}, his degree-histogram hyperprior
$P(\mathbf k_r\mid e_r)=\prod_k\eta^r_k!/(n_r!\,q(e_r,n_r))$ with $\eta^r_k$ the number of nodes of degree
$k$ in the community and $q(m,n)$ the number of partitions of the integer $m$ into at most $n$ parts,
and a \emph{global degree prior} $P(\mathbf k\mid z)=\prod_ip(k_i)$ with $p$ the empirical degree
law of the network, which does not depend on the partition (Proposition~\ref{prop:degsplit}). Mechanisms without a closed-form evidence can enter \eqref{eq:csbm}
through the Laplace form of Theorem~\ref{thm:occam} in the \suppname, used here only for the
Bernoulli degree-corrected type of Section~\ref{sec:real-family}.

\subsection{Constraints between communities}\label{sec:cap}
For binary data with Beta$(a,b)$ off-diagonal priors, a hard cap truncates the prior to $[0,x]$ for a
known $x\in(0,1]$, and the block evidence becomes
\begin{equation}\label{eq:cap}
p(Y_{rs}\mid z)=\frac{\Bfun(a+e,b+N-e)}{\Bfun(a,b)}\cdot\frac{I_{x}(a+e,\,b+N-e)}{I_{x}(a,b)} .
\end{equation}
The Poisson analogue uses the regularised incomplete gamma function (\suppref{app:catalogue}).
Proposition~\ref{thm:cap} shows what the truncation does under contamination, and
Proposition~\ref{prop:soft} that \emph{soft shrinkage}, an off-diagonal Beta prior with mean $x$ and
fixed strength $a+b=m$, here $m=50$, has a bounded effect on every block. Both are compared with a
\emph{plug-in relative cap} $p_{rs}\le c\min(\hat p_{rr},\hat p_{ss})$, $c=0.1$, with $\hat p_{rr}$
the empirical within density under the current partition, which is not a coherent joint prior, and
with the \emph{planted-partition} instance, in which all between-community dyads share one block, the
assortative structure of \citet{Zhang2020} in conjugate form.

\subsection{Sampler, initialisation and summaries}\label{sec:inference}
A single-site update of node $i$ removes it from its community, scores every candidate community with
the ratios of the block evidences that change, and re-inserts it (Algorithm~\ref{alg:sweep} in
\suppref{app:catalogue}). Single-site kernels create and destroy communities slowly, so we add the
restricted-Gibbs split and merge move of \citet{Jain2004} with a random or spectral launch state,
whose acceptance probability is in \suppref{app:catalogue}, and both kernels are validated against
an exact enumerated posterior on a 12-node network.

Before sampling, a set of cheap partitions is scored by the collapsed posterior and the best one starts
the chain: Bethe Hessian spectral clustering \citep{Saade2014} and Leiden \citep{Traag2019} on the
binary network, and the layer-wise and category-wise spectral partitions of \suppref{app:zero} for
the other data types (the model-based initialisation). Every reported partition is the MAP sample
after a deterministic polishing sweep, pooled over chains where more than one chain is run and
compared with the one-community labelling, which the sampler rarely reaches from an informed start
and which Theorem~\ref{thm:threshold} makes the reference candidate, because
Theorems~\ref{thm:threshold} to~\ref{thm:resolution} are statements about the objective that this
estimator maximises. The greedy optimiser of the same objective \citep{Come2015} reaches the same
partitions at a small fraction of the cost on unimodal problems (Section~\ref{sec:mixing}) and is the
recommended point estimator. The sampler supplies what an optimiser cannot: posterior type
probabilities, the credible ball of \citet{Wade2018} around the reported partition, and the detection
of multimodality on real networks. The number of communities is chosen either with a fixed number of
labels under the Dirichlet--multinomial prior and the collapsed posterior of the MAP partition
compared over a grid of $K$, the integrated classification likelihood (ICL) criterion of
\citet{Biernacki2000}, or with the CRP prior and $K$ read off the occupied labels, with
$K_{\max}=\max(2K^\star,10)$ in the simulations unless stated otherwise, where $K^\star$ is the number
of planted communities, written $K$ in the theorems. The Pitman--Yor and Gnedin priors of
\citet{Legramanti2022} and the mixture of finite mixtures of \citet{Miller2018} are implemented as
alternative partition priors of the flat evidence and run as baselines.
Theorems~\ref{thm:threshold} to~\ref{thm:resolution} are stated for the Dirichlet--multinomial and
CRP priors and apply to both criteria. Collapsing is a computational device
(Remark~\ref{prop:collapse} in the \suppname), evaluated as one in Section~\ref{sec:mixing}.

\section{Theory}\label{sec:theory}

This section states the thresholds at which the collapsed evidence of Section~\ref{sec:model} decides
whether a planted partition survives, how many communities to keep, the representation of a dense
block and the mechanism inside a community, and what a constraint between communities does to these
decisions. The first subsection separates three objects that share one objective: the comparison of
the planted partition with the one-community partition, the maximum a posteriori labelling, and the
evidence summed over labellings. The next three give the resolution limit with a per-block cost fixed
by the block prior, the gap between the Bernoulli and the Poisson evidences of a binary network with
the window in which it decides, and mechanism selection inside a community with its finite-size form
and the effect of the degree prior. The last gives the effect of a hard cap on between-community ties.
Proofs are in the \suppname, which also contains posterior concentration under primitive conditions on
the block matrix (Proposition~\ref{thm:concentration}), the zero power of reduced data
(Theorem~\ref{thm:zero}) and the bounds for shrinkage of fixed strength (Proposition~\ref{prop:soft}).
Throughout, $A$ is an undirected binary network on $n$ nodes unless stated otherwise,
$\kl(p\,\|\,q)=p\log\frac pq+(1-p)\log\frac{1-p}{1-q}$ is the Bernoulli divergence,
$\kl_{\Pois}(\lambda\,\|\,\mu)=\lambda\log\frac\lambda\mu-\lambda+\mu$ the Poisson divergence,
$h(p)=p\log p+(1-p)\log(1-p)$, and $\ell(z)=\sum_{r\le s}N_{rs}(z)\,h(\hat p_{rs}(z))$ with
$\hat p_{rs}=e_{rs}/N_{rs}$ is the profile log-likelihood of the homogeneous block model. The function
\begin{equation}\label{eq:J}
J(u,v)=u\log\frac{2u}{u+v}+v\log\frac{2v}{u+v},\qquad u,v>0,
\end{equation}
is the sum of the two Poisson divergences $\kl_{\Pois}(u\,\|\,\frac{u+v}2)+\kl_{\Pois}(v\,\|\,\frac{u+v}2)$
and satisfies $J(u,v)\ge(u-v)^2/(4\max(u,v))$.

\subsection{The MAP criterion, the estimator and the evidence over partitions}\label{sec:thm-threshold}
Let $A$ be a symmetric planted partition $z^\star$ with $K$ equal communities, within probability
$p_{\mathrm{in}}=c_{\mathrm{in}}/n$ and between probability $p_{\mathrm{out}}=c_{\mathrm{out}}/n$,
$c_{\mathrm{in}}>c_{\mathrm{out}}$, mean degree $d=(c_{\mathrm{in}}+(K-1)c_{\mathrm{out}})/K$ fixed as
$n\to\infty$, with the flat evidence under Beta$(a,b)$ block priors, $a,b\ge\frac12$, and a CRP or
Dirichlet--multinomial partition prior with $K_{\max}\ge K$ labels. Let $z_1$ be the one-community
partition, and let
$\mathrm{SNR}_{\mathrm{KS}}=(c_{\mathrm{in}}-c_{\mathrm{out}})^2/(K(c_{\mathrm{in}}+(K-1)c_{\mathrm{out}}))$
be the Kesten--Stigum signal-to-noise ratio, which equals $1$ at the Kesten--Stigum threshold
\citep{Decelle2011,Mossel2015}. Three objects are compared: the log posterior of $z^\star$ against that
of $z_1$, which is the criterion of Theorem~\ref{thm:threshold}, the maximiser $\hat z$ of the log
posterior over all labellings, which is the estimator, and the posterior mass of all labellings against
the one-block model, which is the Bayes factor for structure.
\begin{theorem}[MAP criterion]\label{thm:threshold}
Under the setting just described,
\[
\frac1n\Big[\log p(z^\star\mid A)-\log p(z_1\mid A)\Big]\xrightarrow{p}\frac{J_K}{2}-\log K,
\quad J_K=\frac{c_{\mathrm{in}}}K\log\frac{c_{\mathrm{in}}}d+\frac{(K-1)c_{\mathrm{out}}}{K}\log\frac{c_{\mathrm{out}}}d,
\]
so the planted partition has higher collapsed posterior than the one-community partition, with
probability tending to one, if $J_K>2\log K$, and lower if $J_K<2\log K$. In the second case $z^\star$
is not the maximum a posteriori labelling. Moreover
$J_K=\frac{K-1}{2}\,\mathrm{SNR}_{\mathrm{KS}}\,(1+O((c_{\mathrm{in}}-c_{\mathrm{out}})/d))$, and the block
hyperparameters enter the left side only at order $K^2\log n/n$.
\end{theorem}
The mechanism is the partition prior: a specific balanced $K$-labelling costs $n\log K$ nats under the
CRP or Dirichlet--multinomial prior while the fit gain of the planted partition is $\tfrac n2J_K$. The
same two quantities appear in the first-moment argument of \citet{Banks2016}, where $n\log K$ is the
entropy of the union bound over the $K^n$ labellings and $\tfrac n2J_K$ is the large-deviation cost of
a labelling with the edge profile of $z^\star$. The threshold of the theorem is therefore the
first-moment bound on the information-theoretic threshold, and the corollary records the identity.
\begin{corollary}[The threshold is the first-moment bound of Banks et al.]\label{cor:icl}
Write $c_{\mathrm{in}}=d(1+(K-1)\lambda)$ and $c_{\mathrm{out}}=d(1-\lambda)$ with $0<\lambda<1$.
\begin{enumerate}[label=(\roman*),leftmargin=*]
\item $J_K>2\log K$ if and only if
\[
d\ >\ d^{\mathrm{upper}}_{\mathrm c}(K,\lambda)=\frac{2K\log K}{(1+(K-1)\lambda)\log(1+(K-1)\lambda)+(K-1)(1-\lambda)\log(1-\lambda)},
\]
the upper bound on the detectability threshold in Theorem~1 of \citet{Banks2016}. Let
$\mathrm{SNR}^\star(d,K)$ be the value of $\mathrm{SNR}_{\mathrm{KS}}=d\lambda^2$ at which $J_K=2\log K$.
\item As $d\to\infty$, $\mathrm{SNR}^\star(d,K)\to4\log K/(K-1)$, which exceeds the Kesten--Stigum
value $1$ if and only if $K\le10$. For $K\ge11$ the value $\mathrm{SNR}^\star(d,K)$ exceeds $1$ below a mean degree
$d^\star(K)$ and falls below it above, with $d^\star(11)\approx4.7\cdot10^3$, $d^\star(12)\approx850$,
$d^\star(15)\approx160$ and $d^\star(20)\approx61$ (Table~\ref{tab:iclthr} of the \suppname). The
crossing at $K=11$ is the statement of \citet{Banks2016} that detection below the Kesten--Stigum
threshold is possible for eleven or more assortative communities.
\item The exact integrated classification likelihood of \citet{Come2015}, which maximises the same
objective with Jeffreys priors, does not select the planted partition over the one-community partition
when $\mathrm{SNR}_{\mathrm{KS}}<\mathrm{SNR}^\star(d,K)$, with probability tending to one.
\end{enumerate}
\end{corollary}
Part (i) is an algebraic identity, part (ii) is the linearisation $J_K=\frac{K-1}2\mathrm{SNR}_{\mathrm{KS}}(1+o(1))$
and a computation, and part (iii) restates the theorem for the objective of \citet{Come2015} through
Remark~\ref{prop:collapse}(i). What the corollary adds to \citet{Banks2016} is the identification: the
maximum a posteriori criterion of the conjugate block model, and with it the exact integrated
classification likelihood, pays the entropy of the union bound and cannot keep a planted partition
between the information-theoretic threshold, which \citet{CojaOghlan2018} locate between the two
bounds of \citet{Banks2016}, and the first-moment bound. For $K=2$ the threshold in the units
$s=(c_{\mathrm{in}}-c_{\mathrm{out}})^2/(c_{\mathrm{in}}+c_{\mathrm{out}})=2\,\mathrm{SNR}_{\mathrm{KS}}$
is $5.3$ at $d=10$, $5.4$ at $d=16$ and $8\log2=5.55$ in the large-degree limit, while the
Kesten--Stigum threshold is $s=2$.

Theorem~\ref{thm:threshold} compares two labellings. The estimator is the maximiser over all of them,
and the first-moment argument extends to it. For a labelling $z$ with at most $K_{\max}$ labels let
$\omega(z)$ be the $K\times K_{\max}$ matrix with entries $\omega_{ab}=n^{-1}\#\{i:z^\star_i=a,\,z_i=b\}$,
let $\pi_b=\sum_a\omega_{ab}$ be its column sums, and let the overlap be
$\mathrm{ov}(z)=\max_\phi\sum_b\omega_{\phi(b)b}$ over all maps $\phi$ from the labels of $z$ to those
of $z^\star$, so that $\mathrm{ov}(z^\star)=1$ and $\mathrm{ov}(z)=1/K$ when every block of $z$ contains
the planted communities in equal proportions. The blocks of $z$ are the pairs $(b,b')$ with $b\le b'$,
with dyad fractions $\nu_{bb'}=\pi_b\pi_{b'}$ for $b<b'$ and $\pi_b^2/2$ for $b=b'$, and expected
per-node edge counts $\mu_{bb'}(\omega)=\sum_{a,a'}\omega_{ab}\omega_{a'b'}c_{aa'}$ for $b<b'$ and half
that sum for $b=b'$, where $c_{aa}=c_{\mathrm{in}}$ and $c_{aa'}=c_{\mathrm{out}}$. Put
$\Psi_\omega(x)=\sum_{b\le b'}x_{bb'}\log(x_{bb'}/\nu_{bb'})$,
$\Psi^\star=\frac{c_{\mathrm{in}}}{2K}\log c_{\mathrm{in}}+\frac{(K-1)c_{\mathrm{out}}}{2K}\log c_{\mathrm{out}}$,
$H(\omega)=\frac1K\sum_a H(K\omega_{a\cdot})$ with $H$ the entropy of a probability vector, and
\begin{equation}\label{eq:F}
F(\omega)=H(\omega)-\inf_{x\in\mathcal X_\omega}\ \sum_{b\le b'}\kl_{\Pois}\big(x_{bb'}\,\|\,\mu_{bb'}(\omega)\big),
\end{equation}
where $\mathcal X_\omega=\{x:\ \sum_{b\le b'}x_{bb'}=\tfrac d2,\ \Psi_\omega(x)\ge\Psi^\star-\log K+H(\pi)\}$.
$F(\omega)$ is the exponential rate, per node, of the expected number of labellings with overlap
matrix $\omega$ whose log posterior is at least that of $z^\star$: $H(\omega)$ counts them, the infimum
is the large-deviation cost of the edge profile they need, and $\log K-H(\pi)$ is the prior slack of a
labelling with column sums $\pi$ relative to $z^\star$.
\begin{theorem}[The estimator]\label{thm:estimator}
In the setting of Theorem~\ref{thm:threshold}, let $\hat z$ maximise the collapsed log posterior over
all labellings with at most $K_{\max}$ labels, and let
$\beta^\star=\sup\{v\in[1/K,1]:\ \sup_{\omega:\,\mathrm{ov}(\omega)\le v}F(\omega)<0\}$, with
$\beta^\star=1/K$ when the set is empty.
\begin{enumerate}[label=(\alph*),leftmargin=*]
\item For every $\delta>0$, $P(\mathrm{ov}(\hat z)\ge\beta^\star-\delta)\to1$, and the probability
of the complement is exponentially small in $n$.
\item On the set of overlap matrices whose columns are proportional to the uniform row,
$F(\omega)=\min\{H(\pi),\ \log K-J_K/2\}$, and $F$ is continuous, so $\beta^\star>1/K$ if and only
if $J_K>2\log K$.
\end{enumerate}
\end{theorem}
Part (b) says that the threshold of the estimator is the threshold of the criterion: above
$J_K=2\log K$ the maximum a posteriori labelling is correlated with the planted partition, below it the
planted partition loses to the one-community partition. Part (a) is quantitative: for $K=2$ and the
design of Section~\ref{sec:Kexp} ($d=16$, first-moment bound $\mathrm{SNR}_{\mathrm{KS}}=2.70$),
$\beta^\star$ is $0.62$ at $\mathrm{SNR}_{\mathrm{KS}}=3$, $0.76$ at $4$, $0.86$ at $6$ and $0.93$ at
$10$ (Table~\ref{tab:beta_star} of the \suppname, computed from~\eqref{eq:F} on a grid), while between
$1$ and $2.7$ the bound gives nothing, because $F$ is positive on the uncorrelated labellings.
Section~\ref{sec:Kexp} shows what the estimator returns below the threshold: the one-community
partition once that partition is compared explicitly, and over-fitted states whose posterior is below
it when a sampler does not compare it.

The third object behaves differently. With the rates $c_{\mathrm{in}},c_{\mathrm{out}}$ known, let
$L(A\mid z)$ be the likelihood of $A$ under the labelling $z$, $L_0(A)$ the likelihood under the
Erd\H{o}s--R\'enyi model with edge probability $d/n$, and
$Z_n(A)=\sum_zK^{-n}L(A\mid z)/L_0(A)$ the likelihood ratio of the block model with uniformly random
labels against no structure, which is the Bayes factor for structure when the rates are known. With the
rates integrated, let $\mathrm{BF}_n(A)=\sum_zp(z)\,p(A\mid z)/p_1(A)$, where $p(A\mid z)$ is the flat
conjugate evidence and $p_1(A)$ the evidence of the one-block model under the same Beta prior.
\begin{theorem}[Two thresholds]\label{thm:two}
Let $K=2$ and $c_{\mathrm{in}}>c_{\mathrm{out}}$ be fixed.
\begin{enumerate}[label=(\alph*),leftmargin=*]
\item (\emph{Criterion.}) $n^{-1}\log\big[2^{-n}L(A\mid z^\star)/L_0(A)\big]\to J_2/2-\log2$ in
probability, so the planted labelling beats no structure if and only if $J_2>2\log2$, the bound of
Corollary~\ref{cor:icl}(i).
\item (\emph{Evidence, rates known.}) If $\mathrm{SNR}_{\mathrm{KS}}<1$ then $Z_n(A)=O_p(1)$ under the
planted model. If $\mathrm{SNR}_{\mathrm{KS}}>1$ then $Z_n(A)\to\infty$ in probability.
\item (\emph{Evidence, rates integrated.}) If $\mathrm{SNR}_{\mathrm{KS}}<1$ then
$\log\mathrm{BF}_n(A)\le\tfrac32\log n+O_p(1)$. If $\mathrm{SNR}_{\mathrm{KS}}>1$ then
$\log\mathrm{BF}_n(A)\ge\log Z_n(A)-10\log n-O_p(1)$, so $\mathrm{BF}_n(A)\to\infty$ in probability
whenever $\log Z_n(A)/\log n\to\infty$, which holds when $n^{-1}\log Z_n(A)$ has a positive limit.
\end{enumerate}
\end{theorem}
Part (b) is the contiguity and the consistent parameter estimation of \citet{Mossel2015}, and the
positive limit in (c) is the free energy above the condensation threshold, which for two assortative
communities is the Kesten--Stigum threshold \citep{CojaOghlan2018}. Together the two theorems place the
three objects: the criterion and the estimator switch at the first-moment bound $J_K=2\log K$, which for
$K\le10$ lies above the Kesten--Stigum threshold, and the evidence summed over labellings switches at the
Kesten--Stigum threshold itself, so that between the two the evidence detects structure that the
maximum a posteriori labelling rejects. The price the estimator pays is the entropy of the labellings,
$n\log K$, which the sum does not pay because its mass is spread over exponentially many labellings
correlated with $z^\star$. Section~\ref{sec:Kexp} measures both quantities on the same networks. The
threshold is a property of the objective and not of collapsing, since the collapsed posterior is the
$z$-marginal of the joint posterior (Remark~\ref{prop:collapse}(i) in the \suppname), so any
estimator that maximises it, sampled or not, has it.

\subsection{Resolution limit of the flat evidence}\label{sec:thm-resolution}
Consider a planted partition $z^\star$ with $K$ (even) communities of size $\nu$, $n=K\nu$, in which every
node has expected degree $d_{\mathrm{in}}$ inside its community and $d_b>0$ to each other community, in
the sparse regime $\nu\to\infty$ with $d_{\mathrm{in}},d_b$ fixed, with the flat evidence under Beta$(a,b)$
priors, $a,b\ge\tfrac12$, on every block and a CRP or Dirichlet--multinomial partition prior with
$K_{\max}\ge K$ labels. Let $z'$ merge the communities in pairs.
\begin{theorem}[Resolution limit]\label{thm:resolution}
Under the setting just described,
\begin{align*}
\log p(z^\star\mid A)-\log p(z'\mid A)
&=\tfrac{n}{2}\,\big(J(d_{\mathrm{in}},d_b)-2\log2\big)-\tfrac{3K^2+2K}{8}\,\big(a+\tfrac12\big)\log\nu\\
&\qquad+O(K^2+K\log\nu)+o_p(n),
\end{align*}
with $J$ as in \eqref{eq:J}.
\end{theorem}
In the expansion, $-n\log2$ is the partition-prior cost of keeping $K$ instead of $K/2$ communities and
$(a+\frac12)\log\nu+O(1)$ is the cost of every block that carries $\Theta(\nu)$ expected edges among
$\Theta(\nu^2)$ dyads, by equation~\eqref{eq:blockexp} of the \suppname{} (a block expected to be
empty costs $2a\log\nu+O(1)$). The fit gain from separating all $K/2$ pairs is $\Theta(n)$ and does
not grow with $K$, while the number of block terms grows like $K^2$. The Jeffreys prior, $a=\tfrac12$,
charges $\log\nu$ per block against $\tfrac32\log\nu$ for the uniform prior, which is why the Jeffreys
instance merges fewer cliques in Table~\ref{tab:R5} of the \suppname.
\begin{corollary}[Two-sided resolution limit]\label{cor:res}
In the setting of Theorem~\ref{thm:resolution} with $J=J(d_{\mathrm{in}},d_b)>2\log2$, put
\[
K_{\mathrm{res}}(\nu)=\frac{8\,\nu\,(J-2\log2)}{3(2a+1)\log\nu},
\qquad\text{equivalently}\qquad K_{\mathrm{res}}=\sqrt{\frac{8n\,(J-2\log2)}{3(2a+1)\log\nu}} .
\]
For every $\delta\in(0,1)$, if $K\le(1-\delta)K_{\mathrm{res}}(\nu)$ then
$\log p(z^\star\mid A)-\log p(z'\mid A)>0$ with probability tending to one, and if
$(1+\delta)K_{\mathrm{res}}(\nu)\le K=O(\nu/\log\nu)$ then it is negative with probability tending to one. If
$J<2\log2$ the merger is preferred for every $K$, which is the $K=2$ case of Theorem~\ref{thm:threshold}.
\end{corollary}
The corollary follows because both leading terms of the expansion are of order $n$ at
$K=K_{\mathrm{res}}$ while the remainder $O(K^2+K\log\nu)+o_p(n)$ is $o(n)$ for $K=O(\nu/\log\nu)$.
At a fixed total between degree $d_{\mathrm{out}}$, $d_b=d_{\mathrm{out}}/(K-1)\to0$ and
$J\to d_{\mathrm{in}}\log2$, so the number of communities that the flat evidence resolves grows like
$\nu/\log\nu$ at community size $\nu$, equivalently like $\sqrt{n/\log n}$, the order that
\citet{Peixoto2013} derives for the microcanonical model at fixed degree. The theorem is the
conjugate-evidence form of that limit, and what it adds is the per-block cost $(a+\tfrac12)\log\nu$,
which the block prior sets: the Jeffreys prior resolves $3/2$ times as many communities as the uniform
prior. At a fixed total between degree the proof goes through provided $\nu d_b/\log\nu\to\infty$,
which forces $K=o(\nu/\log\nu)$, so the merge statement of the corollary is not covered there, and
the crossing at $K\approx20$ in Figure~\ref{fig:resolution} of the \suppname{} is an exact
finite-$\nu$ computation. For the design of Section~\ref{sec:Kres}, $\nu=50$ nodes per community,
within degree $4.9$ and between degree $0.125$ to each other community, the corollary gives
$K_{\mathrm{res}}(50)=17$ under the uniform prior and $26$ under the Jeffreys prior, and
$K_{\mathrm{res}}(100)=86$. Section~\ref{sec:Kres} tests it against $K$ and reports the exact
description lengths of a ring of cliques, on which the uniform prior on the block edge-count matrix
of \citet{Peixoto2017} and the nested prior of \citet{Peixoto2014} keep the cliques that independent
block priors merge.

\subsection{Bernoulli and Poisson evidences}\label{sec:thm-family}
The degree-corrected type of Section~\ref{sec:types} and the global-DC instance are built on the
multigraph representation, in which a binary block is scored by the Gamma--Poisson evidence. The next
proposition gives the exact price of that representation on a binary network. Let
$g(p)=p+(1-p)\log(1-p)$ for $p\in[0,1)$ and $g(1)=1$. Both evidences below are probabilities of the
same binary array $A$, the Poisson one through the counts $y_{ij}\in\{0,1\}$, so they are taken with
respect to the same dominating measure and their ratio is a Bayes factor.
\begin{proposition}[Family gap]\label{prop:family}
Let $A$ be a binary network on $n$ nodes and $z$ a partition with $K$ communities and block densities
$\hat p_{rs}=e_{rs}/N_{rs}$. Write $p_{\mathrm{BB}}$ and $p_{\mathrm{GP}}$ for the flat evidences under
Beta$(a,b)$ and Gamma$(a',b')$ block priors, all fixed, and
$G(z)=\sum_{r\le s}N_{rs}\,g(\hat p_{rs})$.
\begin{enumerate}[label=(\roman*),leftmargin=*]
\item $\log p_{\mathrm{GP}}(A\mid z)-\log p_{\mathrm{BB}}(A\mid z)=-G(z)+R(z)$ with
$|R(z)|\le C\,K(K+1)\log(n^2+c)$ for constants $C,c$ that depend only on the hyperparameters.
\item $g$ is convex and increasing on $[0,1)$, $g(p)\ge p^2/2$ and $g''(p)=1/(1-p)\ge1$. If $z'$
refines $z$, so that every block of $z$ is a union of blocks of $z'$, then
\[
G(z')-G(z)\ \ge\ \tfrac12\sum_{\text{blocks }(r,s)\text{ of }z}\ \sum_{\text{sub-blocks }i}N_i\,(\hat p_i-\hat p_{rs})^2\ \ge\ 0 .
\]
\end{enumerate}
\end{proposition}
Part (i) is Lemmas~\ref{lem:stirling} and~\ref{lem:poisson}, which place each evidence within a
logarithmic term of its profile log-likelihood, together with the identity
$N h(\hat p)-N(\hat p\log\hat p-\hat p)=N g(\hat p)$ between the Bernoulli and the Poisson profiles of
a block, and part (ii) is Jensen's inequality with the strong convexity of $g$. The charge grows with
the density ($g(0.04)=0.0008$, $g(0.75)=0.40$ and $g(0.92)=0.72$ nats per dyad) and is larger for the
finer of two nested partitions by the Jensen gap, so relative to the Beta--Bernoulli evidence the
Gamma--Poisson evidence prefers to merge dense blocks of unequal density, and every degree-corrected
evidence built on the multigraph representation, including the microcanonical one of
\citet{Peixoto2017}, carries the effect on dense binary blocks. The corollary turns the proposition
into a threshold.
\begin{corollary}[Density window of the family effect]\label{cor:window}
Let two classes of $\nu$ nodes each have within density $p$ and between density $q<p$, and compare
the partition that keeps them apart with the one that merges them under the CRP or
Dirichlet--multinomial prior. Let $\varphi(p)=p\log p-p$ and let
$\mathrm{JG}_f(p,q)=\tfrac12[f(p)+f(q)]-f(\tfrac{p+q}2)$ be the Jensen gap of a convex $f$. As
$\nu\to\infty$ with $p,q$ fixed, the Bernoulli evidence keeps the classes apart if
$\nu\,\mathrm{JG}_h(p,q)>\log2$ and merges them if the reverse inequality holds, and the Poisson
evidence does the same with $\mathrm{JG}_\varphi$ in place of $\mathrm{JG}_h$. Since
$\mathrm{JG}_h=\mathrm{JG}_\varphi+\mathrm{JG}_g$ with $\mathrm{JG}_g\ge(p-q)^2/8$, the Poisson evidence
merges classes that the Bernoulli evidence keeps apart exactly when
\[
\frac{\log2}{\mathrm{JG}_h(p,q)}\ <\ \nu\ <\ \frac{\log2}{\mathrm{JG}_\varphi(p,q)},
\]
a window whose width in $\nu$ is the factor $\mathrm{JG}_h/\mathrm{JG}_\varphi=1/(1-\bar p)+O(p-q)$
with $\bar p=(p+q)/2$: a factor $2$ at density one half and a factor $10$ at density $0.9$.
\end{corollary}
The window is where the representation, and nothing else, decides. On the football and school
networks the classes have $5$ to $25$ nodes, within densities from $0.1$ to $1.0$ and between densities
up to $0.52$, and with the pair's own densities four of the $66$ pairs of conferences and one of the
$55$ pairs of classes lie inside their windows (script \texttt{check\_window.py} of the released
code), the school pair being the two classes of the fifth grade, with between density $0.52$, which
every instance mixes (Section~\ref{sec:real-school}).

\subsection{Type selection and false splits}\label{sec:thm-split}
Let $C$ be a community of $\nu$ nodes whose within-community network is generated by the Poisson model
$y_{ij}\sim\Pois(\theta_i\theta_j\omega_0)$ independently for $i<j$, with propensities
$\theta_i\stackrel{iid}{\sim}F$, $\E_F\theta=1$, $F$ supported in a compact subset of $(0,\infty)$, and
$\omega_0>0$ fixed while $\nu\to\infty$, which is the dense regime. The expected degree is
$\bar d=\omega_0(\nu-1)$, and $F$ degenerate at $1$ is the Erd\H{o}s--R\'enyi case. The two types of
Section~\ref{sec:types} give the evidences $p_{\ER}(Y_C)$ and $p_{\DC}(Y_C)$, with the uniform degree
prior. For a labelling $z$ of $C$ into two communities, $p_{\ER\text{-}2}(Y_C\mid z)$ denotes the homogeneous
Gamma--Poisson evidence of the two diagonal blocks and the one off-diagonal block that $z$ creates. All
three are taken with respect to the same dominating measure, and the Gamma hyperparameters are fixed.
\begin{theorem}[Type selection, false splits and false merges]\label{thm:split}
\begin{enumerate}[label=(\alph*),leftmargin=*]
\item (\emph{Mechanism selection.}) If $F$ is degenerate,
\[
\log p_{\DC}(Y_C)-\log p_{\ER}(Y_C)=-\tfrac\nu2\log\bar d+c_1\nu+O_p(\sqrt\nu),\qquad
c_1=\tfrac{\log(2\pi)-1}{2}-\tfrac{\bar d}{2\nu}+O(\bar d^{-1}),
\]
so the typed evidence selects $\ER$ with a margin of order $\nu\log\bar d$. If $F$ is not degenerate,
\[
\log p_{\DC}(Y_C)-\log p_{\ER}(Y_C)=\nu\bar d\,\E_F[\theta\log\theta]\,\big(1+o_p(1)\big),
\]
and $\E_F[\theta\log\theta]>0$, so the typed evidence selects $\DC$ with probability tending to one.
\item (\emph{No false split.}) If the support of $F$ contains at least three points, there is
$\kappa=\kappa(F,\omega_0)>0$ such that, with probability tending to one,
\[
\min_{z}\ \big[\log p_{\DC}(Y_C)-\log p_{\ER\text{-}2}(Y_C\mid z)\big]\ \ge\ \kappa\,\nu^2,
\]
the minimum being over all labellings of $C$ into two communities.
\item (\emph{No false merge.}) Let instead $C$ consist of two Erd\H{o}s--R\'enyi communities of $\nu/2$
nodes each, $\nu$ even, with Poisson rates $\lambda_{\mathrm{in}}$ within and
$\lambda_{\mathrm{out}}\neq\lambda_{\mathrm{in}}$ between them, both fixed as $\nu\to\infty$, and let $z^\star$
be the planted labelling. Then
\[
\log p_{\ER\text{-}2}(Y_C\mid z^\star)-\log p_{\DC}(Y_C)\ \ge\ \tfrac{\nu^2}{4}\,J(\lambda_{\mathrm{in}},\lambda_{\mathrm{out}})\,\big(1-o_p(1)\big).
\]
\end{enumerate}
\end{theorem}
Part (a) is a Stirling computation on $\sum_i\log d_i!$: under a degenerate $F$ the configuration
count gains only $\tfrac\nu2\log(2\pi\bar d)+O(\nu)$ nats over the multinomial placement of the edges
while the uniform degree prior costs $\nu\log\bar d+O(\nu)$, and under a non-degenerate $F$ the
configuration count gains $\bar d\sum_i\theta_i\log\theta_i$ nats, with $\E_F[\theta\log\theta]>0$ by
Jensen's inequality. Parts (b) and (c) compare evidences through profile likelihoods
(Lemmas~\ref{lem:poisson} and~\ref{lem:dcprofile} in the \suppname): the margin in (b) is the
population Poisson divergence between the truth and its best two-block homogeneous approximation,
of order $\nu^2$ as soon as $F$ has three points in its support, and the margin in (c) is the
divergence between a two-block truth and its best rank-one approximation. Relative to the
likelihood-ratio test of \citet{Yan2014}, which compares the two mechanisms for the whole network
with the partition fixed and the block parameters not integrated, the theorem makes the decision
inside a single community with the parameters integrated out and adds the two split and merge
statements, and Section~\ref{sec:phase} runs the test of \citet{Yan2014} per community as the
baseline of the theorem. Two qualifications apply. Part (b) covers splits into two homogeneous
blocks, and if $F$ has exactly two atoms the split by atom is a correctly specified two-block model
with no $\Theta(\nu^2)$ margin in either direction. And the theorem is asymptotic in the dense
regime: at $\nu=80$ and $\bar d=10$, the sizes of Section~\ref{sec:phase}, the $O(\nu)$ terms of (a)
are of the same order as $\tfrac\nu2\log\bar d$. The next proposition gives the constants that the
theorem leaves implicit.
\begin{proposition}[Finite-size form of Theorem~\ref{thm:split}(b)]\label{prop:finite}
In the setting of Theorem~\ref{thm:split}(b), let $z_s$ split $C$ by atom into two blocks, let
$\bar\lambda_b$ be the mean Poisson rate of the ordered atom pairs that fall in block $b$ of $z_s$, and put
\[
\kappa_F(z_s)=\tfrac12\sum_{t,t'}F(t)F(t')\,\kl_{\Pois}\big(tt'\omega_0\,\big\|\,\bar\lambda_{b(t,t')}\big),
\]
the Jensen gap of the Poisson rates over the blocks of $z_s$, which is positive when $F$ has at least
three atoms. Let $\hat\ell_{\DC}$ be the Poisson profile log-likelihood of the degree-corrected model
at its maximiser and $R_\nu=\hat\ell_{\DC}-\log p_{\DC}(Y_C)$ its degree-sequence cost. Then
\begin{align*}
\log p_{\DC}(Y_C)-\log p_{\ER\text{-}2}(Y_C\mid z_s)&=\kappa_F(z_s)\,\nu^2-R_\nu+O_p(\nu),\\
R_\nu&\le\tfrac{\omega_0}2\sum_i\theta_i^2+\nu\big(1+\log(1+2\bar d)\big)+L_P+1 ,
\end{align*}
with $L_P$ the constant of Lemma~\ref{lem:dcprofile}, so the degree-corrected type beats the split once
$\kappa_F(z_s)\nu$ exceeds $1+\log(1+2\bar d)+\tfrac{\omega_0}2\E_F[\theta^2]+o(1)$.
\end{proposition}
The proposition replaces the constant of the proof of Theorem~\ref{thm:split}(b), of order $10^{-6}$,
by the leading coefficient. In the design of Table~\ref{tab:splitb} ($F$ uniform on $\{0.5,1,1.5\}$,
$\omega_0=0.2$) the coefficients are $\kappa_F=0.0034$ and $0.0057$ for the two by-atom splits, the
measured cost is $R_\nu/\nu=\tfrac12\log\bar d+0.25\pm0.04$ over $300\le\nu\le2400$ (script
\texttt{check\_split\_finite.py}), and the crossing that solves
$\kappa_F\nu=\tfrac12\log\bar d(\nu)+0.25$, at $\nu\approx830$ and $430$, is within a factor of $1.4$
of the exact crossings of the evidences, between $600$ and $800$ and between $300$ and $450$.
Section~\ref{sec:phase} runs the sampler in this design.

The split that Theorem~\ref{thm:split}(b) does not cover, into two degree-corrected blocks, is
governed by the degree prior of the degree-corrected evidence, and the next proposition shows that the
prior can create a gain linear in the community size for a split that carries no block structure.
\begin{proposition}[Splitting by degree under a degree prior]\label{prop:degsplit}
Let $C$ be a community of $\nu$ nodes with total degrees $k_i$, $e=\sum_ik_i$, and let $z_s$ split it
into two blocks of $\nu_1,\nu_2$ nodes with degree sums $e_1,e_2$ and mean degrees
$\bar d_1=e_1/\nu_1$, $\bar d_2=e_2/\nu_2$, $\bar d=e/\nu$. Let $\phi(x)=(1+x)\log(1+x)-x\log x$, which is
concave on $(0,\infty)$.
\begin{enumerate}[label=(\roman*),leftmargin=*]
\item Under the microcanonical degree-corrected evidence with the uniform degree prior
$P(\mathbf k_r\mid e_r)=\binom{e_r+\nu_r-1}{e_r}^{-1}$ (the global-DC and global typed instances), the
degree prior contributes to $\log p(A\mid z_s)-\log p(A\mid z_1)$ the amount
\[
\nu\Big[\phi(\bar d)-\tfrac{\nu_1}\nu\phi(\bar d_1)-\tfrac{\nu_2}\nu\phi(\bar d_2)\Big]+O(\log\nu)\ \ge\ 0,
\]
which is of order $\nu$ whenever $\bar d_1\ne\bar d_2$, while the configuration count and the prior on
the block edge counts contribute $O(\log e)$ when the split carries no block structure. A degree-corrected
community is therefore split along its degrees at a gain linear in its size.
\item Under a prior on the degree sequence that factorises over nodes with a law that does not depend
on the partition, $P(\mathbf k\mid z)=\prod_ip(k_i)$, the degree prior contributes nothing to the
comparison, and the gain of a structureless split is $O(\log e)$.
\end{enumerate}
\end{proposition}
Part (i) is the mechanism behind the degree-corrected communities that the global typed instance
calls homogeneous with $K$ inferred (Section~\ref{sec:phase}): the instance first cuts the community
along its degrees and then finds each part homogeneous. Part (ii) is satisfied by the empirical
degree law of the network, $p(k)=\#\{i:k_i=k\}/n$, the global degree prior of
Section~\ref{sec:types}. The degree-histogram hyperprior of \citet{Peixoto2017} does not satisfy
(ii), and on a heavy-tailed community its incentive to split by degree exceeds that of the uniform
prior (Table~\ref{tab:degsplit} of the \suppname), while the within-degree type of the typed
instance has no such incentive. The last proposition carries part (a) of the theorem over to the
global typed instance, whose relevant quantity is the total degree.
\begin{proposition}[Type selection with total degrees]\label{prop:gtyped}
Let $C$ be a community of $\nu$ nodes whose total degrees $k_i$ are, conditionally on their sum $e$,
distributed as the multinomial allocation of $e$ half-edges with node weights $\theta_i$ drawn from
$F$ with mean one and compact support in $(0,\infty)$, $\bar d=e/\nu\to\infty$ with $\bar d=o(\nu)$.
Under the global typed instance with the uniform degree prior, the row factors satisfy
\[
\log\frac{\prod_ik_i!\,\binom{e+\nu-1}{e}^{-1}/e!}{\nu^{-e}}
=\begin{cases}-\tfrac\nu2\log\bar d+\nu\,\tfrac{\log(2\pi)-1}2+O_p(\sqrt\nu)+O(\nu/\bar d) & F\ \text{degenerate},\\[2pt]
\nu\bar d\,\E_F[\theta\log\theta]\,(1+o_p(1)) & F\ \text{not degenerate},\end{cases}
\]
so the global typed instance selects the mechanism of a community from its total degrees with the
margins of Theorem~\ref{thm:split}(a). With the global degree prior of Proposition~\ref{prop:degsplit}(ii)
the uniform-prior term $\log\binom{e+\nu-1}{e}$ is replaced by $-\sum_i\log p(k_i)$, and the comparison
becomes one between two laws for the degrees of the community, the multinomial allocation and the
empirical law of the network: its margin is $\nu$ times the divergence between the two when the
community is homogeneous and the network is not, and it vanishes for a community whose degrees follow
the law of the network.
\end{proposition}

\subsection{Constraints between communities}\label{sec:thm-cap}
Let two communities of sizes $n_1,n_2$ have $e_b$ between edges among $N_b=n_1n_2$ dyads. Write
$\mathrm{BF}_x$ for the Bayes factor of the separated against the merged partition under Beta$(a,b)$
block priors, $a,b\ge\tfrac12$, with the off-diagonal prior truncated to $[0,x]$, so that $x=1$ is the
uncapped model, and let $\Delta_x(e,N)=\log I_x(a+e,b+N-e)-\log I_x(a,b)$ with $I_x$ the regularised
incomplete beta function. Truncation divides the block evidence by the prior mass of $[0,x]$ and
multiplies it by the posterior mass of the same set, which is the identity behind the proposition.
\begin{proposition}[The cap]\label{thm:cap}
\begin{enumerate}[label=(\roman*),leftmargin=*]
\item $\log\mathrm{BF}_x=\log\mathrm{BF}_1+\Delta_x(e_b,N_b)$, the function $\Delta_x(e,N)$ is
non-increasing in $e$ for fixed $N$, and $\Delta_x(e,N)\le\log(1/I_x(a,b))$ for all $e,N$, which is
$\log(1/x)$ for $a=b=1$: the cap adds a bounded amount of evidence for separation, whatever the data.
\item For positive integers $a,b$, put $n'=a+b+N-1$ and $\hat\mu=(a+e)/n'$, which is within $1/N$ of
the posterior mean of the between density. If $\hat\mu>x$ then
$\Delta_x(e,N)\le-\,n'\,\kl(\hat\mu\,\|\,x)+\log(1/I_x(a,b))$: a block whose density exceeds the cap
is penalised at a rate linear in its number of dyads.
\item (\emph{Contamination.}) Let $e^\star_x=\min\{e_b:\mathrm{BF}_x<1\}$ be the smallest number of
between edges at which the merged partition is preferred. If $e^\star_1<\infty$ and
$\log\mathrm{BF}_1(e+1)-\log\mathrm{BF}_1(e)\le-\zeta<0$ for $e^\star_1\le e<e^\star_1+j_x$, where
$j_x=\lceil\log(1/I_x(a,b))/\zeta\rceil$ and $e^\star_1+j_x\le N_b$, then $e^\star_x\le e^\star_1+j_x$:
the cap postpones the merger under contamination by at most $j_x$ cross edges. For positive integers
$a,b$, whenever $\hat\mu_b=(a+e_b)/(a+b+N_b-1)>x$ and
$\log\mathrm{BF}_1(e_b)<n'\kl(\hat\mu_b\|x)-\log(1/I_x(a,b))$, the capped model prefers the merged
partition while the uncapped one need not.
\end{enumerate}
\end{proposition}
Part (i) is the product form of the evidence with $I_x\le1$, part (ii) is a Chernoff bound on the
posterior mass of $[0,x]$, and part (iii) follows from (i). The over-splitting statement, that the
cap penalises a split of a community of density above the cap at rate $N'\kl(p_C\|x)$, is
Proposition~\ref{prop:capsplit} of the \suppname. In the exact example of Figure~\ref{fig:cap}
there, two communities of $50$ nodes with eleven cross edges among $2500$ between dyads, $x=0.02$
and $a=b=1$, the cap adds $3.9$ nats to $149.8$, and under uniform contamination the merged partition
becomes preferred at $114$ cross edges with the cap and at $197$ without one. The cap is a hard
assortativity constraint that does not protect boundaries against contamination, an effect that the
assortative models of \citet{Zhang2020} and \citet{Amongero2025} do not analyse, whereas shrinkage of
fixed strength has an effect on the same Bayes factor bounded uniformly in the block size
(Proposition~\ref{prop:soft} of the \suppname): for the prior of Section~\ref{sec:cap} the cost at a
between density of $0.04$, $0.08$ and $0.20$ is $+1.9$, $-0.1$ and $-6.8$ nats whatever the block
size, where the truncation costs $-12$, $-85$ and $-479$ nats at $N=1600$. Section~\ref{sec:capexp}
tests both devices under contamination.

\section{Simulation study}\label{sec:simulations}

The simulations test the results of Section~\ref{sec:theory} in their order and locate the regimes in
which the typed evidence is ahead of the alternatives, ties them and loses. Every network is generated
by the released code from a planted partition with the block matrix stated in the caption, every
method receives the same networks (paired seeds), the same $K$ when $K$ is given and its default
tuning, and every cell reports the mean ARI against the planted partition with a $95\%$ percentile
bootstrap interval and the mean $\hat K$. Head-to-head claims rest on the mean paired difference over
seeds, its $95\%$ bootstrap interval and a Wilcoxon signed-rank test with Holm adjustment over all
cells and comparisons of the experiment at level $0.05$. The alternatives have four roles. For the
mechanism question they are the degree-corrected likelihood of \citet{Karrer2011}, the per-community
likelihood-ratio test of \citet{Yan2014} with a bootstrap null, the global-DC instance and the global
typed instance under its three degree priors. For the choice of $K$ they are the corrected BIC of
\citet{Hu2020}, the penalised likelihood of \citet{Wang2017}, the microcanonical evidence of
\citet{Peixoto2017}, variational EM with the ICL criterion, edge cross-validation \citep{Li2020}, the
spectral estimate of \citet{Le2022} and the flat evidence under the Pitman--Yor and Gnedin priors of
\citet{Legramanti2022} and the mixture-of-finite-mixtures prior of \citet{Miller2018}. The Bethe
Hessian \citep{Saade2014}, regularised spectral clustering \citep{Qin2013,Lei2015}, Leiden and Infomap
are the standards a practitioner would use, and SCORE \citep{Jin2015}, designed for degree
heterogeneity, is the adversarial case. The collapsed samplers use the model-based initialisation of
Section~\ref{sec:inference}, one split and merge proposal per sweep, $300$ sweeps with the first third
discarded, the polished MAP partition and the explicit comparison with the one-community labelling.
Bold entries mark the best mean in a column at two decimals. All numbers come from
\texttt{experiments/bench.py} of the released code, the four reduced-data designs that test
Theorem~\ref{thm:zero} are in the \suppname, and the sampler of \citet{Peixoto2014} was not run
(Section~\ref{sec:discussion}).

\subsection{The criterion, the estimator and the evidence over labellings}\label{sec:Kexp}
This experiment tests Theorems~\ref{thm:threshold}, \ref{thm:estimator} and~\ref{thm:two} on the same
networks: $K=2$, $n=1000$, mean degree $16$, twenty paired seeds per value of the signal
$s=(c_{\mathrm{in}}-c_{\mathrm{out}})^2/(c_{\mathrm{in}}+c_{\mathrm{out}})=2\,\mathrm{SNR}_{\mathrm{KS}}$,
so that the Kesten--Stigum threshold is $s=2$ and the first-moment bound of Corollary~\ref{cor:icl} is
$s=5.4$.
\begin{table}[tbp]\centering\scriptsize
\caption{The MAP criterion and the estimator of Theorems~\ref{thm:threshold} and~\ref{thm:estimator} with the sampler ($K=2$, $n=1000$, mean degree $16$, 20 paired seeds per row, $s=(c_{\mathrm{in}}-c_{\mathrm{out}})^2/(c_{\mathrm{in}}+c_{\mathrm{out}})=2\,\mathrm{SNR}_{\mathrm{KS}}$, Kesten--Stigum threshold $s=2$, first-moment bound $s=5.4$). Flat CRP evidence with $10$ labels: fraction of seeds in which the planted partition has higher collapsed posterior than the one-community partition (exact), fraction of seeds in which the MAP is the one-community labelling (the sampler's best state compared with it), fraction with $\hat K=2$, mean ARI ($95\%$ bootstrap interval, $0$ when the MAP is the one-community labelling) and mean $\hat K$. Jeffreys: fraction with $\hat K=2$ for the Jeffreys instance. Bethe Hessian with its own $\hat K$ and Leiden: ARI [mean $\hat K$]. The Bethe Hessian with $K=2$ given has the same ARI from $s=4$ on and $0.04$ at $s=2$. Last three columns, on the 12 seeds per row with stored samples: median radius of the $95\%$ credible ball of \citet{Wade2018} around the reported partition (nats of variation of information), fraction of seeds with the planted partition inside the ball, and mean posterior probability of two communities. Higher ARI is better.}
\label{tab:R17}
\resizebox{\linewidth}{!}{\begin{tabular}{lccccccccccc}\toprule
$s$ & keeps $z^\star$ & MAP $=z_1$ & P($\hat K=2$) & ARI & $\hat K$ & Jeffreys P($\hat K=2$) & Bethe, $K$ inferred & Leiden & ball radius & $z^\star$ in ball & post.\ P($K=2$)\\\midrule
2 & 0.00 & 0.75 & 0.00 & 0.00 (0.00, 0.00) & 1.0 & 0.00 & 0.01 (0.00, 0.02) [1.1] & 0.01 (0.01, 0.01) [8.5] & 1.11 & 0.67 & 0.27\\
4 & 0.00 & 1.00 & 0.00 & 0.00 (0.00, 0.00) & 1.0 & 0.00 & 0.50 (0.47, 0.52) [2.0] & 0.04 (0.03, 0.05) [8.0] & 1.19 & 1.00 & 0.40\\
6 & 0.95 & 0.00 & 0.75 & 0.76 (0.73, 0.79) & 2.4 & 0.45 & 0.75 (0.73, 0.76) [2.0] & 0.12 (0.10, 0.14) [7.4] & 0.79 & 1.00 & 0.53\\
8 & 1.00 & 0.00 & 0.90 & 0.88 (0.84, 0.91) & 2.1 & 0.75 & 0.87 (0.86, 0.88) [2.0] & 0.34 (0.30, 0.39) [6.0] & 0.57 & 0.92 & 0.47\\
10 & 1.00 & 0.00 & 0.95 & 0.95 (0.93, 0.96) & 2.0 & 0.80 & 0.93 (0.92, 0.94) [2.0] & 0.52 (0.47, 0.57) [4.5] & 0.42 & 1.00 & 0.63\\
11 & 1.00 & 0.00 & 0.95 & 0.97 (0.96, 0.97) & 2.0 & 0.90 & 0.95 (0.94, 0.96) [2.0] & 0.63 (0.57, 0.71) [3.6] & 0.15 & 0.92 & 0.77\\
12 & 1.00 & 0.00 & 0.95 & 0.98 (0.97, 0.98) & 2.0 & 0.90 & 0.97 (0.96, 0.97) [2.0] & 0.84 (0.79, 0.90) [3.0] & 0.08 & 0.75 & 0.91\\
14 & 1.00 & 0.00 & 1.00 & 0.99 (0.99, 0.99) & 2.0 & 1.00 & 0.98 (0.98, 0.99) [2.0] & 0.98 (0.95, 0.99) [2.4] & 0.05 & 0.92 & 1.00\\
16 & 1.00 & 0.00 & 1.00 & 0.99 (0.99, 1.00) & 2.0 & 1.00 & 0.99 (0.99, 0.99) [2.0] & 0.99 (0.99, 1.00) [2.2] & 0.03 & 0.83 & 0.99\\
20 & 1.00 & 0.00 & 1.00 & 1.00 (1.00, 1.00) & 2.0 & 1.00 & 1.00 (1.00, 1.00) [2.0] & 1.00 (1.00, 1.00) [2.0] & 0.01 & 0.92 & 0.99\\
\bottomrule\end{tabular}}\end{table}

\begin{figure}[t]
\centering
\includegraphics[width=0.42\linewidth]{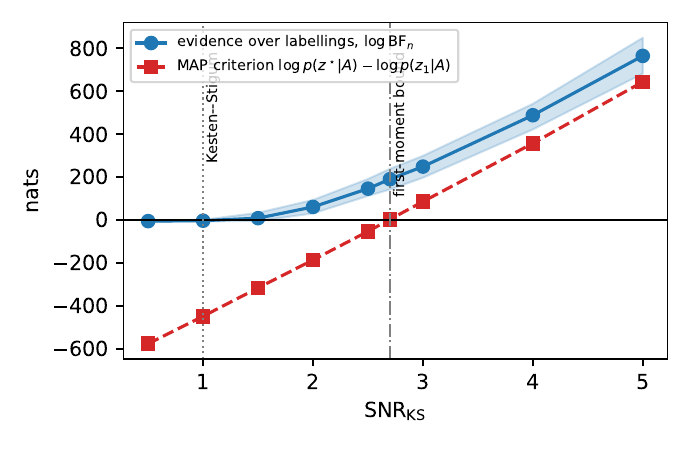}
\caption{The two thresholds of Theorem~\ref{thm:two} on the design of Table~\ref{tab:R17}: the log
Bayes factor of the flat evidence summed over all labellings against the one-block model, computed by
thermodynamic integration with parallel tempering (ten seeds, mean and range), and the MAP criterion
$\log p(z^\star\mid A)-\log p(z_1\mid A)$ of Theorem~\ref{thm:threshold}, exact on the same networks,
against $\mathrm{SNR}_{\mathrm{KS}}=s/2$. The dotted lines mark the Kesten--Stigum threshold and the
first-moment bound. Table~\ref{tab:ti} of the \suppname{} lists the values.}
\label{fig:two}
\end{figure}
The first two columns of Table~\ref{tab:R17} are exact. The planted partition has higher collapsed
posterior than the one-community labelling on no seed at $s\le4$, on $19$ of $20$ seeds at $s=6$ and
on every seed from $s=8$, the two sides of the first-moment bound $s=5.4$, and the criterion itself
follows the limit $n(J_2/2-\log2)$ of Theorem~\ref{thm:threshold} within $12$ nats at every value of
$s$ (Table~\ref{tab:ti}: $-186$ against $-182$ at $s=4$, $+3$ against $+2$ at $s=5.4$). The
estimator behaves as Theorem~\ref{thm:estimator}(b) says. At $s=4$ the best labelling that the
sampler finds has a posterior $39$ to $285$ nats below the one-community labelling on every seed, and
at $s=2$ on the $15$ seeds on which the chain did not reach one community itself, so once the two are
compared the MAP is $z_1$, $\hat K=1$ and the ARI is $0$, for the uniform and the Jeffreys instance
alike. A sampler that does not make this comparison stays in over-fitted states (two to five
communities with ARI $0.42$ at $s=4$), because the merge of two large blocks is never accepted at
weak signal. From $s=6$ the MAP is correlated with the planted partition (ARI $0.76$, then $0.95$ at
$s=10$ and $1.00$ at $s=20$). The overlap guarantee $\beta^\star$ of
Theorem~\ref{thm:estimator}(a) is conservative ($0.62$ at $s=6$, where the observed ARI of $0.76$
needs an overlap of about $0.94$), because the first-moment bound counts every labelling with a
given overlap as if it could beat $z^\star$.

Figure~\ref{fig:two} shows the third object. The log Bayes factor of the flat evidence summed over all
labellings, computed by thermodynamic integration along $26$ power posteriors with parallel tempering
\citep{Friel2008,Geyer1991}, is $-5$ and $-3$ at $\mathrm{SNR}_{\mathrm{KS}}=0.5$ and $1$, positive on
seven of ten seeds at $1.5$ ($+9$, range $-5$ to $+34$), $+61$ at $2$ ($s=4$), $+146$ at $2.5$ and
$+190$ at the first-moment bound, where the MAP criterion is $+3$. Between the Kesten--Stigum
threshold and the first-moment bound the evidence over labellings favours structure by factors
between $10^{26}$ and $10^{82}$ while every single labelling the sampler finds, and the planted one,
are below the one-block model. This is Theorem~\ref{thm:two}: the mass of the sum sits on exponentially many labellings
correlated with $z^\star$, none of which pays for itself. The spectral estimate of \citet{Le2022}
switches with the evidence and returns $\hat K=2$ from $s=4$ on with ARI $0.50$ there, where the MAP
returns one community, and the two are tied from $s=6$ on (paired differences of $+0.01$ to $+0.02$
in favour of the MAP for $8\le s\le14$, on $13$ to $19$ of the $20$ seeds). Above the
first-moment bound the sampler supplies what an optimiser does not: the median radius of the $95\%$
credible ball of \citet{Wade2018} around the reported partition falls from $0.8$ nats of variation of
information at $s=6$ to $0.4$ at $s=10$ and $0.01$ at $s=20$, the planted partition lies inside the
ball on every seed at $s=6$ and $10$ and on $9$ to $11$ of the $12$ seeds at the other values, and the
posterior probability of two communities crosses one half at $s=6$ and reaches $0.91$ at $s=12$, with
the rest of the mass on three or more communities.

\subsection{Number of communities}\label{sec:Kres}
This experiment tests Theorem~\ref{thm:resolution} and Corollary~\ref{cor:res}: how many communities
the flat evidence resolves, how the answer depends on the block prior, and how the estimators of $K$
of the literature behave on the same networks.
\begin{figure}[t]
\centering
\includegraphics[width=0.42\linewidth]{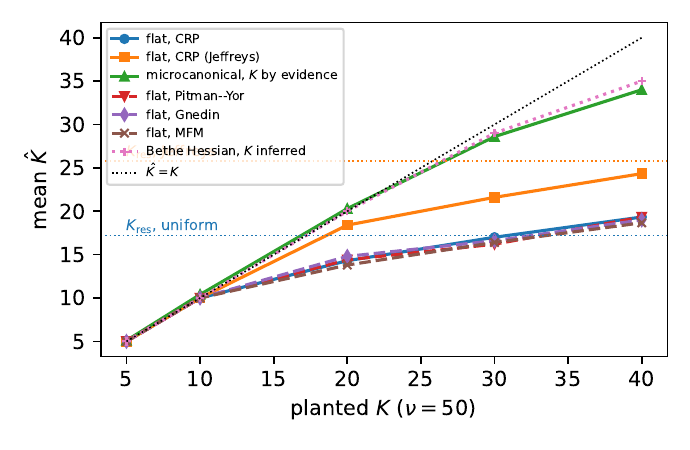}
\caption{Corollary~\ref{cor:res} against $K$: mean $\hat K$ on a planted partition with $\nu=50$ nodes
per community, within density $0.1$ and between density $0.0025$ (the rates of the $K=20$ design of
Section~\ref{sec:phase}), for $K$ from $5$ to $40$, with the flat CRP evidence under the uniform and
the Jeffreys block priors, the microcanonical evidence with the uniform prior on the edge-count
matrix, the Bethe Hessian with its own $\hat K$ and Leiden (ten paired seeds for $K\le20$, five
for $K=30$ and three for $K=40$). The horizontal lines are $K_{\mathrm{res}}(50)=17$ and $26$ of the
corollary for the two block priors. Table~\ref{tab:R21} of the \suppname{} lists the values with
the ARIs and the Pitman--Yor, Gnedin, mixture-of-finite-mixtures and corrected-BIC estimates.}
\label{fig:kres}
\end{figure}
Figure~\ref{fig:kres} varies $K$ at fixed community size $\nu=50$, within degree $4.9$ and between
degree $0.125$ to each other community, for which the corollary gives $K_{\mathrm{res}}(50)=17$ under
the uniform prior and $26$ under the Jeffreys prior. The flat CRP evidence returns every community up
to $K=10$, then $14.3$ of $20$, $17.0$ of $30$ and $19.3$ of $40$ under the uniform prior, and $18.4$,
$21.6$ and $24.3$ under the Jeffreys prior. The limit of the uniform prior is reached where the
corollary puts it, the Jeffreys prior resolves more, as the per-block cost $(a+\tfrac12)\log\nu$ says,
and both saturate below the planted $K$ while the microcanonical evidence with the uniform prior on
the edge-count matrix ($20.3$, $28.6$ and $34.0$) and the Bethe Hessian, which has no partition
prior ($20$, $29$ and $35$), keep resolving. The Pitman--Yor, Gnedin and mixture-of-finite-mixtures
priors follow the CRP within one community at every $K$, because they change the partition prior by
$O(K\log n)$ nats and leave the per-block cost, and the corrected BIC returns $17.2$, $20$ and $24$.

Table~\ref{tab:R5} of the \suppname{} compares the estimators of $K$ for $K^\star=5$ at $n=500$ and
mean degree $10$ over sixteen paired seeds, where the first-moment bound is $\mathrm{SNR}_{\mathrm{KS}}=2.08$. Below it,
at $1.5$, every method is near chance (ARI $\le0.23$), the flat CRP evidence returns $2.9$ communities
on average, the global instances with their uniform degree prior $7.4$ and $8.3$, the microcanonical
evidence $8.6$ and edge cross-validation one. Above it, at $2.5$, the flat CRP instance with either
block prior finds $\hat K=5$ on every seed and has the highest ARI ($0.81$, against $0.79$ for the
typed CRP instance and the microcanonical evidence, $0.74$ for the Bethe Hessian and $0.73$ for edge
cross-validation and variational EM with the ICL criterion), and the evidence-grid selection with
fixed labels is more conservative than the CRP ($\hat K=4.6$), as Theorem~\ref{thm:threshold} says
of the integrated classification criterion. The other estimators behave as their penalties say (eight paired seeds): the Gnedin and
mixture-of-finite-mixtures priors find $\hat K=5$ on every seed at $2.5$ (ARI $0.80$), the
Pitman--Yor prior $4.9$, the corrected BIC of \citet{Hu2020}, whose penalty $n\log k$ is the
partition-prior cost of Theorem~\ref{thm:threshold} on a fixed spectral partition, returns one
community at $1.5$ and $4.5$ at $2.5$, the criterion of \citet{Wang2017} returns $3.2$ and $5.0$, and
the global typed instance with the global degree prior returns one community at $1.5$, where the
uniform degree prior over-splits to $8.3$, and $5.5$ at $2.5$. At $4$
and $8$ every method with the right $K$ agrees. The ring of thirty six-cliques
(last column) is the classical form of the resolution limit, outside the regime of
Theorem~\ref{thm:resolution} because the community size is fixed at six. Every instance with
independent conjugate blocks merges cliques ($\hat K$ from $8$ to $14.5$, ARI $0.33$ to $0.50$), the
Jeffreys prior least, the Pitman--Yor, Gnedin and mixture-of-finite-mixtures priors like the CRP
($9.9$ to $10.6$) and the corrected BIC and the Wang--Bickel criterion most ($8.0$ and $7.0$), while
the microcanonical evidence returns all thirty cliques on every seed, as the exact description lengths of Section~\ref{sec:thm-resolution} say, and so do the Bethe
Hessian, Leiden, Infomap and edge cross-validation. For modularity, thirty cliques with $480$ edges
sit inside its own limit of about $\sqrt{2\cdot480}=31$ \citep{Fortunato2007}.

\subsection{Mechanism selection}\label{sec:phase}
This experiment tests Theorem~\ref{thm:split} and Propositions~\ref{prop:finite}
to~\ref{prop:gtyped}: whether the typed evidence recovers the generating mechanism of each community,
whether it avoids the splits of a heterogeneous community that the flat evidence makes, whether its
posterior type probabilities are calibrated, and which degree prior makes the global typed instance
calibrated.
\begin{figure}[t]
\centering
\includegraphics[width=0.62\linewidth]{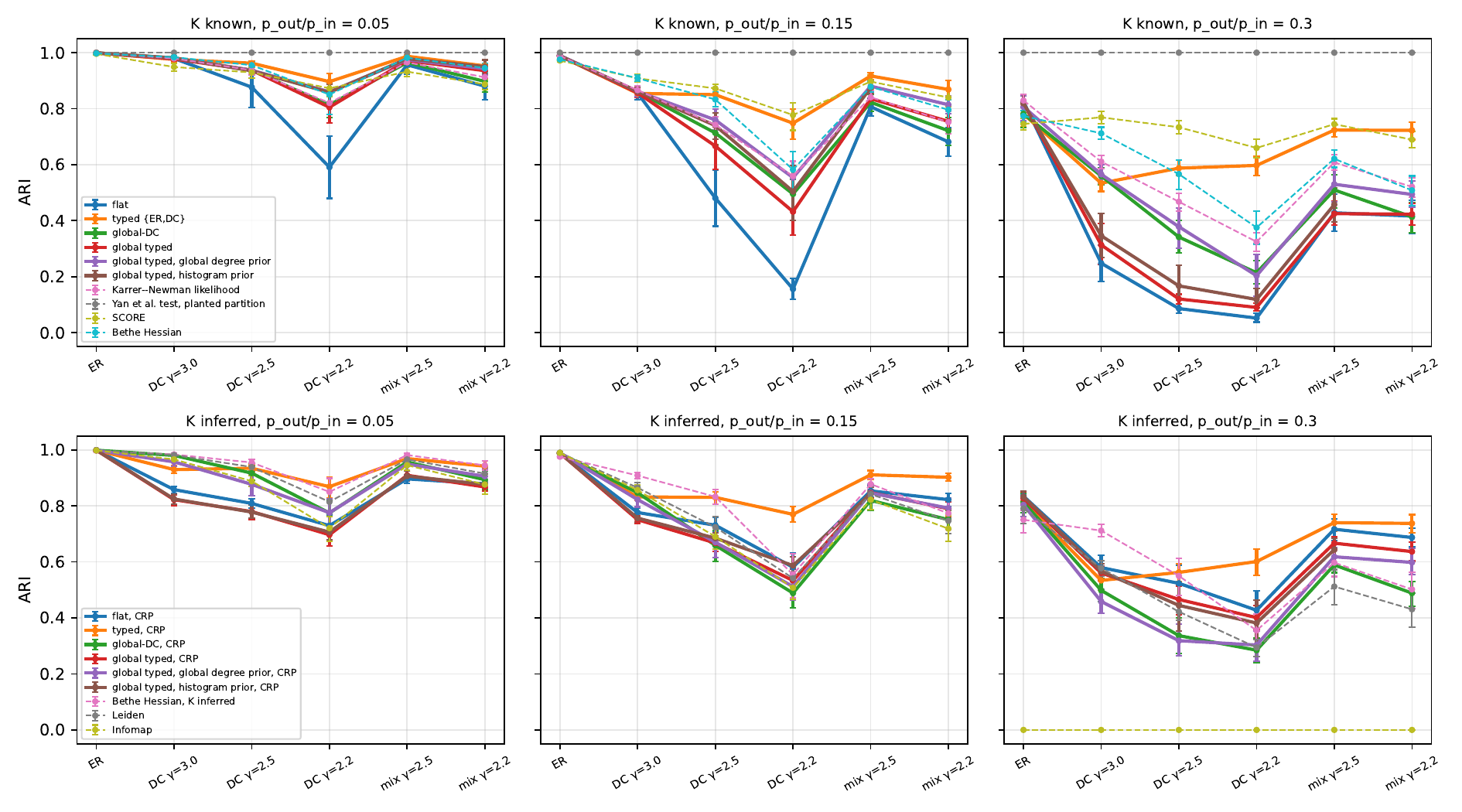}
\caption{Within-community heterogeneity against community signal ($n=320$, $K=4$, mean within
degree 10, twenty paired seeds per cell, error bars $95\%$ bootstrap intervals of the mean). The
horizontal axis of every panel lists the six designs: all communities Erd\H{o}s--R\'enyi (ER), all
communities degree-corrected (DC) with Pareto propensities of exponent $\gamma\in\{3,2.5,2.2\}$,
unbounded and normalised to mean one and therefore outside the compact-support condition of
Theorem~\ref{thm:split}, and the mixture of two ER and two DC communities at $\gamma\in\{2.5,2.2\}$.
The columns are $p_{\mathrm{out}}/p_{\mathrm{in}}\in\{0.05,0.15,0.30\}$, with
$\mathrm{SNR}_{\mathrm{KS}}=8.0$, $5.0$ and $2.6$. The top row reports the ARI with $K$ known for the
flat evidence, the typed evidence, the global typed instance with the global degree prior, the
global-DC instance, the Karrer--Newman likelihood, SCORE and the Bethe Hessian. The bottom row reports
the ARI with $K$ inferred for the CRP instances, Leiden, Infomap and the Bethe Hessian with its own
estimate of $K$. A higher value is better, and the false-split and false-merge rates of the bottom row
are in Table~\ref{tab:R15} of the \suppname.}
\label{fig:phase}
\end{figure}
Three regimes are visible in Figure~\ref{fig:phase}, in Figure~\ref{fig:winner} of the \suppname{}
and in the paired comparisons of Table~\ref{tab:R15pairedfull} of the \suppname. In the homogeneous design every
method recovers the planted partition at $p_{\mathrm{out}}/p_{\mathrm{in}}\le0.15$ (ARI $\ge0.97$),
and the typed evidence pays no measurable price for carrying the type indicator: its ARI equals that
of the flat evidence on $39$ of the $40$ seeds at the two stronger signals, and its type posterior
selects ER for every community, the sign of the leading term $-\tfrac\nu2\log\bar d$ of
Theorem~\ref{thm:split}(a). At the weakest signal the flat evidence and the Karrer--Newman likelihood
($0.83$ each) are ahead of the typed and global-DC instances ($0.78$ each), which empty one of their
four labels on three seeds of twenty, a difference that is not significant after adjustment.

In the heterogeneous designs the ARI of the flat evidence falls as the tail becomes heavier ($0.48$ and
$0.16$ at $\gamma=2.5$ and $2.2$ for $p_{\mathrm{out}}/p_{\mathrm{in}}=0.15$, $0.09$ and $0.05$ at
$0.30$), because a flat diagonal block can explain a hub-dominated community only by splitting it.
In all eight heavy-tail and mixture cells at the two weaker signals the typed evidence is ahead of the
flat evidence ($+0.11$ to $+0.59$), of the global-DC instance ($+0.09$ to $+0.38$) and of the
Karrer--Newman likelihood ($+0.08$ to $+0.27$), every difference significant after Holm adjustment.
Over the $18$ cells it is significantly ahead of the flat evidence in $12$, of the global-DC instance
in $11$ and of the Karrer--Newman likelihood in $11$, and significantly behind the Karrer--Newman
likelihood in one, the $\gamma=3$ design at the weakest signal ($0.53$ against $0.61$), where its type
posterior assigns the ER type to $17\%$ of the communities because the margin
$\nu\bar d\,\E_F[\theta\log\theta]$ of Theorem~\ref{thm:split}(a) is small for a $\gamma=3$ tail.
The global typed instance with the global degree prior, which decides from total degrees, is within
$0.04$ of the typed evidence at $\gamma=3$, behind it by $0.09$ to $0.40$ in the heavy-tail designs at
the two weaker signals ($0.76$ against $0.85$ and $0.38$ against $0.59$ at $\gamma=2.5$), because
the planted heterogeneity acts inside the communities and a diagonal type absorbs it exactly, and
ahead of the same instance with the uniform degree prior by $0.04$ to $0.26$ in the eight heavy-tail
and mixture cells at the two weaker signals, with type accuracy $0.96$ to $1.00$ in the
degree-corrected designs against $0.55$ to $0.85$ for the uniform prior and $0.70$ to $1.00$ for the
histogram prior. SCORE is the strongest alternative at weak signal: at
$p_{\mathrm{out}}/p_{\mathrm{in}}=0.30$ it is ahead of every likelihood-based method in the three
all-DC designs ($0.77$, $0.73$ and $0.66$ against $0.53$, $0.59$ and $0.60$ for the typed evidence,
the first two differences significant). Against the Bethe Hessian with $K$ given the typed evidence is significantly ahead in
seven cells and behind in the two $\gamma=3$ cells at the weaker signals. Figure~\ref{fig:winner}
summarises the map: the typed evidence leads every alternative for a heavy within-community tail
($\gamma\le2.5$) at moderate or strong signal and in every mechanism mixture at
$p_{\mathrm{out}}/p_{\mathrm{in}}\le0.15$, where it alone recovers both the partition and the
mechanism ($0.99$, $0.92$ and $0.72$ at $\gamma=2.5$ across the three signals with type accuracy
$1.00$, against $0.97$, $0.84$ and $0.61$ for the Karrer--Newman likelihood), a spectral method with
$K$ given leads near the threshold of Theorem~\ref{thm:threshold}, and every method ties on
homogeneous communities. The test of \citet{Yan2014}, run per community on the planted partition with a bootstrap
null at level $0.05$, names the mechanism of $99\%$ of the $1440$ communities of the experiment
(every one of its $474$ homogeneous calls and $99\%$ of its $966$ degree-corrected calls is right),
against $97\%$ for the typed evidence on its own fitted partition with $K$ known, whose errors are
all in the homogeneous direction: the test sees the planted communities and the parameters at their
maximum, and on a real network it needs the partition it is given, which is what the typed evidence
supplies.

With $K$ inferred, the flat CRP evidence splits heterogeneous communities in $83\%$ of its fits
($\hat K$ from $7.6$ to $9.9$ for $K=4$, Tables~\ref{tab:R15_auto_0.05} to~\ref{tab:R15_auto_0.3} of
the \suppname), the global typed CRP evidence with the uniform degree prior in $84\%$ and the
global-DC CRP evidence in $52\%$, which is Proposition~\ref{prop:degsplit}(i), and the typed CRP
evidence in $79\%$, mildly ($\hat K=5.8$ to $6.7$ in the all-DC designs and $5.0$ to $5.4$ in the
mixtures), which Theorem~\ref{thm:split}(b) does not forbid since it does not cover a split into two
degree-corrected blocks. The global degree prior removes the part of the splitting that the degree
prior pays: in the $\gamma=3$ designs the global typed CRP evidence returns $\hat K=4.4$, $5.0$ and
$5.4$ under the global prior against $7.9$, $8.0$ and $7.6$ under the uniform prior and $7.8$, $7.9$
and $7.6$ under the histogram prior, at $\gamma=2.5$ and strong signal $5.2$ against $7.8$ and $7.7$,
and its type accuracy with $K$ inferred is $0.79$, $0.72$ and $0.63$ at $\gamma=3$, $2.5$ and $2.2$
against $0.24$, $0.32$ and $0.44$ under the uniform prior, whose pieces are called homogeneous, as
the proposition says. Among the methods that infer $K$ (Table~\ref{tab:R15pairedauto} of the
\suppname), the typed CRP evidence is significantly ahead of the flat CRP evidence in $11$ cells and
behind in none, ahead of Leiden in nine and behind in one, and ahead of the Bethe Hessian with its
own $\hat K$ in six cells and behind it in the three $\gamma=3$ cells, where it over-splits.
Table~\ref{tab:R15calib} of the \suppname{} reports the calibration of the type posteriors. The
typed evidence with $K$ known is decisive and wrong in one direction only: $99\%$ of its $1007$
communities have posterior probability of the degree-corrected type below $0.05$ or above $0.95$,
every community above $0.95$ was generated degree-corrected, and $8\%$ of the $362$ below $0.05$
were degree-corrected, $23$ of them at $\gamma=3$. With $K$ inferred the fraction of wrong
homogeneous calls rises to $34\%$ for the typed CRP evidence and to $64\%$ and $65\%$ for the global
typed CRP evidence under the uniform and the histogram prior, because the pieces into which these
priors cut a heterogeneous community are closer to homogeneous. Under the global degree prior the
degree-corrected calls stay right ($98\%$ of the $641$ communities above $0.95$ with $K$ known,
$99\%$ of $592$ with $K$ inferred) and the homogeneous communities spread over the lower bins ($186$
below $0.05$ and $165$ between $0.05$ and $0.95$ with $K$ known, $96\%$ of the former and all of the
latter homogeneous), the vanishing margin of Proposition~\ref{prop:gtyped} for a community whose
degree law is that of the network: the global prior trades the decisiveness of the homogeneous call
for the removal of the degree split.

Table~\ref{tab:R22} of the \suppname{} runs the sampler in the design of Proposition~\ref{prop:finite},
two degree-corrected communities with propensities uniform on $\{0.5,1,1.5\}$ and $\omega_0=0.2$, at
$\nu\in\{300,600,1200\}$ with five paired seeds. The typed CRP evidence returns the two communities
with the degree-corrected type on every seed at every $\nu$, and so does the global typed instance
with the global degree prior up to $\nu=600$ (three blocks on two seeds at $1200$), while the flat CRP
evidence cuts each community into three or four blocks at every size ($\hat K=7.4$ to $8.0$, ARI
$0.31$ to $0.34$), the false split that Theorem~\ref{thm:split}(b) forbids to the degree-corrected type
and that a flat block makes. The by-atom split, which the exact evidences of Table~\ref{tab:splitb}
prefer to the degree-corrected type at $\nu\le600$, carries no block structure, and the sampler does
not reach it from the model-based start at any size, so the crossing of the proposition is a property
of the objective that the search does not expose in this design.

Table~\ref{tab:R18} of the \suppname{} takes the mechanism question to larger networks. Its first three columns repeat
the three heterogeneous designs of Figure~\ref{fig:phase} at $n=3000$ and
$p_{\mathrm{out}}/p_{\mathrm{in}}=0.15$ over six paired seeds. The ordering with $K$ known is
unchanged and the gaps widen: the typed evidence reaches $0.79$, $0.74$ and $0.89$ with type accuracy
$1.00$, ahead of the flat evidence, the global-DC instance, the Karrer--Newman likelihood ($0.68$,
$0.47$ and $0.46$) and the Bethe Hessian, which falls to $0.26$ and $0.07$ in the two all-DC designs
at this size, on every seed of every design, while SCORE is within $0.07$ of the typed evidence in
the all-DC designs and behind by $0.52$ in the mixture. The last three columns are the regime of
Theorem~\ref{thm:split} itself, $K=20$ communities of $\nu$ nodes, ten homogeneous and ten
degree-corrected with the propensities of Proposition~\ref{prop:finite}, at expected within degree
$0.1(\nu-1)$. At $\nu=250$ every method recovers the partition and the typed evidence recovers the
mechanism of every community. At $\nu=100$ the partition is recovered by every collapsed instance
and the Karrer--Newman likelihood ($0.96$ to $0.98$) but the type posteriors are near chance
($0.62$), because with within degree $10$ the margin of Theorem~\ref{thm:split}(a), about $86$ nats,
is of the order of the $O(\nu)$ terms, the finite-size effect that Proposition~\ref{prop:finite}
quantifies. At $\nu=50$ the typed and the flat evidences merge communities ($\hat K=13.7$ and $14.3$
of $20$) while the global-DC and Karrer--Newman instances keep all twenty ($0.81$ to $0.84$) and the
flat evidence under the Jeffreys prior keeps $17.8$ (ARI $0.75$ against $0.63$), which is the
resolution limit of Corollary~\ref{cor:res} at $K=20>K_{\mathrm{res}}(50)=17$, with the ordering of
the two block priors that the corollary gives, and not the mechanism.

\subsection{Contamination of constrained boundaries}\label{sec:capexp}
This experiment tests Proposition~\ref{thm:cap}(iii) and Proposition~\ref{prop:soft}: what a violated
cap or contamination by cross edges costs, and whether shrinkage of fixed strength protects the
boundaries that the truncation does not.
Figure~\ref{fig:capexp} of the \suppname{} ($n=120$, $K=3$, $p_{\mathrm{in}}=0.2$, cap $x=0.02$,
planted between density $p_{\mathrm{out}}/x\in\{0.25,0.5,1,2,4\}$, $0$, $60$ and $120$ contaminating
cross edges, thirty paired seeds per cell) and the paired differences of Table~\ref{tab:R16paired}
there show three things. First, the cap never helps measurably: in the six cells where it is valid and at
most $60$ edges are added, the capped and uncapped evidences have the same ARI on at least $27$ of
the $30$ seeds and the mean paired difference is at most $+0.002$, the bound $\log(1/x)=3.9$ nats of
Proposition~\ref{thm:cap}(i) against fit differences of hundreds of nats. Second, the cap lowers the
ARI as soon as it is wrong or contamination makes it wrong: at $p_{\mathrm{out}}=2x$ the hard cap
lowers the mean ARI from $0.94$ to $0.83$ without contamination, from $0.80$ to $0.20$ with $60$
contaminating edges and from $0.44$ to $0.04$ with $120$ (paired differences $-0.11$, $-0.60$ and
$-0.40$, the last two significant after Holm adjustment), with the CRP prior it merges communities
($\hat K$ from $2.9$ to $1.3$), the rate of Proposition~\ref{thm:cap}(ii), and at a valid cap $120$
contaminating edges push the empirical between density over it and the capped model falls from
$0.91$ to $0.56$. Third, soft shrinkage is never significantly worse than the uncapped evidence (largest
loss $0.005$) and is ahead in the contaminated cells at $p_{\mathrm{out}}=2x$ ($+0.05$ and $+0.07$),
while the plug-in relative cap behaves like the hard cap under contamination. At $p_{\mathrm{out}}=4x$
($\mathrm{SNR}_{\mathrm{KS}}=1.6$) the flat evidence merges the communities even with $K=3$ labels
fixed (ARI $0.05$), while the Bethe Hessian and variational EM with $K$ given still recover part of
the structure ($0.43$ and $0.31$): this is the threshold of Theorem~\ref{thm:threshold}, since for
$K=3$ the criterion $J_3>2\log3=2.20$ gives $J_3=3.39$ at $p_{\mathrm{out}}=2x$ and $J_3=1.49$ at
$4x$. The cap encodes a known assortativity constraint, useful against over-splitting
(Proposition~\ref{prop:capsplit} of the \suppname) and harmful against contamination, and shrinkage of
fixed strength is the device to use when contamination is possible.

\subsection{Computation}\label{sec:mixing}
The full account is in \suppref{app:computation}: with $n$ from $300$ to $10^4$ at
$\mathrm{SNR}_{\mathrm{KS}}=3$ the collapsed instances scale linearly in $n$ at fixed degree ($3.6$ to
$118$ seconds for the flat evidence with $200$ sweeps), the greedy optimiser of the same objective is
two orders of magnitude faster at no measurable cost in ARI, the split and merge move and an informed
start account for most of what a random-start chain lacks, and the sampler's advantage is confined
to the multimodal cases of Section~\ref{sec:real} and to the type posteriors and credible balls that
an optimiser does not supply.

\section{Real networks}\label{sec:real}

This section tests whether the ordering of Section~\ref{sec:simulations} holds on networks with
external labels, and it isolates the two effects that the simulations did not contain, the
under-fitting of dense binary blocks by the Poisson representation and the splitting of hub-dominated
communities by the degree prior. The eight networks of the main comparison, public, simple, undirected and reduced to the largest
connected component, are Zachary's karate club \citep[$34$
nodes, $78$ edges, two factions,][]{Zachary1977}, the Doubtful Sound dolphins \citep[$62$, $159$, two
groups,][]{Lusseau2003}, American college football \citep[$115$, $613$, twelve
conferences,][]{Girvan2002}, political books \citep[$105$, $441$, three labels,][]{Krebs2004},
political blogs \citep[$1222$, $16{,}714$, two camps,][]{Adamic2005}, the first day of the
primary-school contact network binarised \citep[$236$, $5899$, ten classes and the
teachers,][]{Stehle2011}, the email-Eu-core network \citep[$986$, $16{,}064$, $42$
departments,][]{Leskovec2007} and the Cora citation network \citep[$2485$, $5069$, seven subject
areas,][]{Sen2008}. Section~\ref{sec:real-large} adds the second day of the school network and
larger networks with up to $19{,}717$ nodes.
Metadata are external labels and not ground truth \citep{Peel2017}, so every method is also scored by
held-out link prediction on ten percent of the dyads, masked uniformly at random with one mask per
network (a balanced set of $5\%$ of the edges and as many non-edges on the larger networks): the
collapsed samplers remove the masked dyads from the likelihood and predict them by the posterior
predictive, every other method predicts with the block-density plug-in of its partition, and the AUC
is compared with the global density and the Chung--Lu degree model. Methods that need $K$ receive the
number of metadata classes, the instances that infer $K$ use $\max(3K,20)$ labels, capped at $60$,
and every collapsed instance is run with four chains of $300$ sweeps from over-dispersed starts (the
model-based candidate of Section~\ref{sec:inference}, Bethe Hessian spectral, uniformly random and
Leiden). Table~\ref{tab:real_main} reports the pooled MAP, the chain with the highest log posterior,
on the six networks on which the findings rest, and the \suppname{} reports all networks and methods,
the held-out AUC, the convergence diagnostics, four chains of $1500$ sweeps on the four larger
networks and what the sampler adds to its initialisation.
\begin{table}[tbp]\centering\small
\caption{Real networks: ARI against the metadata partition, with the number of metadata classes in parentheses after each network name and $K$ fixed to it for the methods that need $K$, and $\hat K$ in brackets for the methods that infer it. The collapsed samplers report the pooled MAP over four chains of $300$ sweeps. The dolphins and political-book networks, the Jeffreys, Pitman--Yor, Gnedin and MFM priors, the greedy optimiser, variational EM, the Bickel--Chen profile likelihood, the test of \citet{Yan2014}, Louvain, walktrap, Infomap and the held-out AUC of every method are in the \suppname. Higher is better.}
\label{tab:real_main}
\resizebox{\linewidth}{!}{\begin{tabular}{lcccccc}\toprule
Method & karate (2) & football (12) & polblogs (2) & school (11) & eu-core (42) & Cora (7)\\\midrule
flat & -0.01 & 0.86 & 0.00 & 0.44 & 0.16 & 0.16\\
typed \{ER,DC\} & -0.01 & 0.58 & 0.40 & 0.43 & 0.18 & 0.22\\
global typed & -0.01 & 0.81 & 0.77 & 0.64 & 0.48 & 0.19\\
global typed, global degree prior & 0.00 & 0.81 & 0.77 & 0.68 & 0.49 & 0.26\\
global-DC & 0.00 & 0.81 & 0.77 & 0.70 & 0.46 & 0.26\\
flat, CRP & 0.78 [5] & 0.86 [11] & 0.20 [20] & 0.59 [21] & 0.16 [42] & 0.28 [19]\\
typed, CRP & -0.01 [2] & 0.63 [8] & 0.19 [20] & 0.52 [17] & 0.18 [35] & 0.29 [18]\\
global typed, CRP & -0.01 [2] & 0.81 [10] & 0.15 [20] & 0.52 [17] & 0.47 [25] & 0.25 [21]\\
global typed, global degree prior, CRP & 0.00 [1] & 0.81 [10] & 0.14 [20] & 0.52 [17] & 0.43 [23] & 0.24 [21]\\
global-DC, CRP & 0.00 [1] & 0.81 [10] & 0.19 [20] & 0.53 [17] & 0.46 [25] & 0.25 [21]\\
Karrer--Newman likelihood & 0.77 & 0.89 & 0.79 & 0.74 & 0.42 & 0.27\\
corrected BIC (Hu et al.) & 0.00 [1] & 0.73 [9] & 0.32 [5] & \textbf{0.90 [11]} & \textbf{0.52 [17]} & \textbf{0.33 [8]}\\
Bethe Hessian & \textbf{0.88} & 0.90 & 0.11 & 0.89 & 0.14 & 0.10\\
regularised spectral & \textbf{0.88} & \textbf{0.91} & \textbf{0.81} & \textbf{0.90} & 0.39 & 0.29\\
SCORE & 0.67 & 0.89 & 0.80 & 0.25 & 0.07 & 0.22\\
Bethe Hessian, $K$ inferred & 0.00 [1] & 0.72 [9] & 0.08 [7] & 0.76 [8] & 0.14 [17] & 0.09 [50]\\
Leiden & 0.46 [4] & 0.71 [9] & 0.76 [23] & 0.56 [6] & 0.29 [17] & 0.25 [85]\\
\bottomrule\end{tabular}}\end{table}

\subsection{Dense classes and the family effect}\label{sec:real-family}
On the football and school networks, which are dense and assortative, the two block families order
the same partitions differently, and Table~\ref{tab:familygap} of the \suppname{} isolates the
reason with exact evidences on the full networks. On the football network the Gamma--Poisson evidence prefers the typed MAP,
which merges conferences ($\hat K=7$), to the flat MAP, which keeps eleven of the twelve, by $22$
nats, whereas the Beta--Bernoulli evidence prefers the flat MAP by $79$ nats, and on the school
network the Gamma--Poisson evidence ranks the global-DC MAP, the partition closest to the classes,
$244$ nats below the flat MAP, whereas the Beta--Bernoulli evidence ranks it first. The column $G(z)$
attributes the change to Proposition~\ref{prop:family}: the difference between the two log evidences
of every partition equals $G(z)=\sum Ng(\hat p)$ within $6$ nats, and $G(z)$ is $171$ nats for the
typed MAP of the football network against $277$ for the flat MAP and $2001$ nats for the typed MAP of
the school network against $2330$ for the global-DC MAP, so the Poisson representation subtracts less
from the coarser partition and prefers it. Corollary~\ref{cor:window} says where this can happen:
with the pairs' own densities four of the $66$ pairs of conferences and one of the $55$ pairs of
classes lie inside their windows, the school pair being the two fifth-grade classes with between
density $0.52$, which every instance mixes (Section~\ref{sec:real-school}). The last column of
Table~\ref{tab:familygap} evaluates the same partitions under a Bernoulli type library, the
Erd\H{o}s--R\'enyi type collapsed exactly and the Bernoulli degree-corrected type through the Laplace
form of Theorem~\ref{thm:occam} with $\alpha_i\sim N(0,4)$ on the logit scale. That library orders the
partitions as the Beta--Bernoulli evidence does, so the effect belongs to the representation and a
Bernoulli degree-corrected type removes it on fixed partitions.

\subsection{Hub-dominated networks and the degree prior}\label{sec:real-hubs}
On the political blogs, where degree heterogeneity acts on the between-camp ties as well as inside
the camps, degree correction must be global. The flat evidence returns the hub and periphery split, a
community of $332$ nodes of mean degree $75$ with half of them from each camp, with ARI $0.00$. The
global-DC and global typed instances reach $0.77$, the same partition, close to the Karrer--Newman
likelihood ($0.79$), SCORE ($0.80$) and regularised spectral clustering ($0.81$). The typed instance
sits in between ($0.40$), because a type confined to the diagonal cannot absorb heterogeneity that
acts on the between-camp ties, and it prefers its own partition to the camps by about $1100$ nats
(Table~\ref{tab:family} of the \suppname). With $K$ inferred every collapsed instance uses its $20$ labels on the blogs, the global typed
instance with the global degree prior included ($0.77$ with $K=2$, $0.14$ with $20$ labels), and so
does the flat evidence under the Pitman--Yor, Gnedin and mixture-of-finite-mixtures priors ($0.12$
to $0.20$), while the corrected BIC of \citet{Hu2020} stops at five communities ($0.32$) and Leiden at
$23$ ($0.76$). The degree prior is not the reason here: the Karrer--Newman profile log-likelihood,
which absorbs the degrees exactly, rises by $3500$ to $5100$ nats from the camps to the $19$ or $20$
blocks that the CRP instances find, against $122$ nats for a cut of each camp into its top quarter by
degree and the rest (\texttt{experiments/blogs\_profile\_gain.py}), so the finer partitions carry
block structure beyond the degrees, which every evidence sees, whereas Proposition~\ref{prop:degsplit}
concerns the gain of a structureless cut, isolated on planted communities in Section~\ref{sec:phase}.
The global typed instances and the test of \citet{Yan2014} call both camps degree-corrected.

On the karate club the degree split of \citet{Karrer2011} has higher flat evidence than the factions
by $22$ to $25$ nats, so every flat or typed instance with $K=2$ returns the leader and periphery split
(ARI $-0.01$), moving away from the spectral initialisation at the factions ($0.88$), and the type
follows the partition: at the degree split both communities are homogeneous and at the factions both
are degree-corrected. The global instances leave their second label empty ($\hat K=1$), because
after the cost of the degree sequence the two factions do not raise the evidence of a $78$-edge
network, whereas the Karrer--Newman likelihood, which charges nothing for $K=2$, recovers them
($0.77$). On the email network the global typed instances have the highest ARI of the
methods with $K$ given ($0.48$ and $0.49$ under the uniform and the global degree prior, against
$0.46$ for the global-DC instance and $0.42$ for Karrer--Newman) and the highest AUC ($0.95$), behind
only the corrected BIC with its own $17$ communities ($0.52$), and on Cora the global instances reach
$0.26$ against $0.33$ for the corrected BIC, with every collapsed instance's chains in different modes.

\subsection{Larger networks}\label{sec:real-large}
\begin{table}[tbp]\centering\scriptsize
\caption{Larger networks: ARI against the metadata classes. Node counts and metadata classes in the header, $K$ fixed to the number of classes for the methods that need it, $\hat K$ in brackets for the methods that infer it (label budget $\min(3K,60)$). The Pitman--Yor, Gnedin and MFM priors are in the held-out table of the supplement. Collapsed samplers: pooled MAP over four chains of $300$ sweeps, two chains of $200$ sweeps on the Amazon network and one chain on PubMed, whose size allowed only the instances shown. Higher is better.}
\label{tab:realv6_ari}
\resizebox{\linewidth}{!}{\begin{tabular}{lccccc}\toprule
Method & school day 2 (238, 11) & CiteSeer (1331, 6) & EuroSiS (1272, 13) & Amazon photo (7487, 8) & PubMed (19717, 3)\\\midrule
flat & 0.46 & 0.06 & 0.21 & 0.39 & 0.00\\
typed \{ER,DC\} & 0.53 & 0.09 & 0.58 & -- & --\\
global typed & 0.50 & 0.08 & 0.83 & -- & --\\
global typed, global degree prior & 0.79 & 0.10 & 0.71 & 0.49 & 0.01\\
global-DC & 0.83 & 0.09 & 0.71 & 0.49 & 0.01\\
Karrer--Newman likelihood & 0.78 & 0.14 & 0.58 & 0.43 & 0.11\\
regularised spectral & 0.92 & 0.24 & 0.63 & 0.44 & 0.29\\
Bethe Hessian & 0.90 & 0.00 & 0.07 & -- & --\\
SCORE & 0.41 & 0.12 & 0.40 & -- & --\\
flat, CRP & 0.62 [20] & 0.05 [14] & 0.28 [23] & 0.25 [24] & 0.05 [20]\\
typed, CRP & 0.57 [18] & 0.04 [11] & 0.58 [17] & -- & --\\
global typed, CRP & 0.57 [18] & 0.08 [18] & 0.56 [33] & -- & --\\
global typed, global degree prior, CRP & 0.57 [18] & 0.09 [18] & 0.55 [30] & 0.35 [24] & --\\
global-DC, CRP & 0.61 [17] & 0.12 [15] & 0.57 [29] & 0.35 [24] & --\\
corrected BIC & 0.78 [13] & 0.24 [6] & 0.66 [15] & 0.35 [24] & 0.08 [20]\\
Bethe Hessian, $K$ inferred & 0.85 [9] & 0.01 [28] & 0.12 [17] & 0.05 [50] & 0.00 [50]\\
Leiden & 0.44 [5] & 0.14 [85] & 0.86 [36] & 0.59 [22] & 0.09 [527]\\
\bottomrule\end{tabular}}\end{table}

Table~\ref{tab:realv6_ari} takes the comparison to the second day of the school network and to
four larger networks with their published labels: the CiteSeer citation network \citep[$1331$
nodes in the largest component, six subject areas,][]{Sen2008}, the EuroSiS web network
\citep[$1272$, $13$ countries,][]{Boguna2010}, the Amazon photo co-purchase network \citep[$7487$,
eight product categories,][]{Shchur2018} with two chains of $200$ sweeps, and the PubMed citation
network \citep[$19{,}717$, three classes,][]{Sen2008} with one chain per instance and the instances
that its size allowed. On PubMed every method is far from the three classes (regularised spectral clustering $0.29$,
the Karrer--Newman likelihood $0.11$, the collapsed instances with $K=3$ at most $0.01$), and the
flat CRP evidence with $20$ labels has the highest held-out AUC ($0.93$, against $0.90$ for Leiden's
$527$ communities).
On the other three, with $K$ given, the ordering of
Section~\ref{sec:real-hubs} repeats where degrees are heterogeneous: on EuroSiS the global typed
instance reaches $0.83$, the global-DC instance and the global typed instance with the global degree
prior $0.71$, the typed instance and the Karrer--Newman likelihood $0.58$, regularised spectral
clustering $0.63$ and the flat evidence $0.21$, and on the Amazon network the two global instances
reach $0.49$ against $0.43$ for Karrer--Newman, $0.44$ for regularised spectral clustering and $0.39$
for the flat evidence, while on CiteSeer every method is below $0.25$. With $K$ inferred every
collapsed instance finds more communities than the metadata (from $11$ on CiteSeer to the full
budget of $24$ on the Amazon network), so does Leiden ($22$ to $527$), and the finer partitions
predict held-out ties better (Table~\ref{tab:realv6_auc} of the \suppname): on the Amazon network
the global-DC and global typed CRP instances have the highest AUC of all methods ($0.98$, against
$0.96$ for the flat CRP evidence and $0.94$ for Leiden), on EuroSiS the CRP instances reach $0.93$ to
$0.95$ against $0.91$ for Leiden, and on CiteSeer Leiden's $85$ communities predict best ($0.88$
against $0.82$ for the typed CRP evidence).

\subsection{Modes of the evidence on the school network}\label{sec:real-school}
On the school network the flat evidence has two modes, and the one with the higher log posterior is
far from the class partition: the chains started from the model-based and spectral candidates stop
at the class partition (ARI $0.92$), and the chain started at random finds a partition $294$ nats
higher that keeps three classes whole, merges pupils of two or three classes of the same or adjacent
grades and cuts each merged group into two (\texttt{results/report/school\_crosstab.txt}), predicts
the held-out contacts better (AUC $0.88$ against $0.84$) and has ARI $0.44$, the number in
Table~\ref{tab:real_main}. Four chains of $1500$ sweeps leave the two modes in place
(split-$\hat R$ $133$). On the second day of the same school (Table~\ref{tab:realv6_ari}) the ranking
repeats: the flat, typed and global typed instances with $K=11$ stop at $0.46$ to $0.53$, the
global-DC instance and the global typed instance with the global degree prior reach $0.83$ and
$0.79$, regularised spectral clustering $0.92$, and every CRP instance finds $17$ to $20$ blocks with
the highest held-out AUC of all methods ($0.94$, against $0.85$ to $0.88$ for the spectral estimates
and the corrected BIC), so on both days the mode of the evidence is a partition finer than the classes
that predicts contacts better than the classes do.

\subsection{Semi-synthetic networks}\label{sec:real-semi}
The blogs, the football network and the school carry no ground truth for the mechanism inside a
community, so we plant one. On the blogs the within-camp ties of one camp are replaced by an
Erd\H{o}s--R\'enyi network with the same number of edges, with the other camp and the between-camp
ties kept, which tests Theorem~\ref{thm:split}(a) on real between-camp ties.
In Table~\ref{tab:R19} of the \suppname{} both typed instances assign the Erd\H{o}s--R\'enyi type to
the rewired camp and the degree-corrected type to the other on every one of the twenty rewirings of
each camp, and recover the camps with ARI $0.93$ to $0.94$ whichever camp is rewired, while the flat
evidence stays at the hub split ($0.27$ and $0.14$), the degree-corrected alternatives lose part of
the rewired camp, whose ties no longer follow its degrees (typed evidence ahead by $+0.03$ to $+0.09$
of the global-DC instance and Karrer--Newman on every rewiring, by $+0.12$ and $+0.14$ of SCORE), and
the Bethe Hessian is ahead by $0.03$ for one rewired camp and fails for the other ($0.12$).
Table~\ref{tab:R23} of the \suppname{} does the same on the football and school networks, where the within ties of one
of the three largest classes are replaced either by an Erd\H{o}s--R\'enyi network or by a
degree-corrected network with Pareto$(2.5)$ propensities, and the result is a boundary of the
mechanism question. The classes have $12$ to $25$ nodes and within densities from $0.64$ to $0.98$,
so a simple graph bounds every degree by the class size and the degree-corrected rewiring cannot
create hubs: the coefficient of variation of the within-class degrees is $0.22$ to $0.26$ after the
degree-corrected rewiring against $0.17$ to $0.20$ after the homogeneous one on the football
conferences, and $0.04$ to $0.11$ against $0.03$ to $0.07$ in the school classes. No instance
therefore separates the two rewirings: the typed evidence calls the Erd\H{o}s--R\'enyi type on both
networks for both, the global typed instances call the Erd\H{o}s--R\'enyi type on the football
network for both and the degree-corrected type for $37\%$ to $43\%$ of the rewired school classes
whichever mechanism was planted, since the total degrees of a class include its contacts with the
other classes, and the partitions are those of the unrewired networks (regularised spectral
clustering $0.90$, the typed evidence $0.71$ and $0.55$). In a near-complete class the two mechanisms coincide and
Theorem~\ref{thm:split}(a), whose margin $\nu\bar d\,\E_F[\theta\log\theta]$ presumes propensities that
the network can express, has nothing to select.
Table~\ref{tab:R24} of the \suppname{} runs the count-valued instances of the catalogue on the school
network with the number of $20$-second contacts as the tie weight and $K=11$ given (four chains). The
Gamma--Poisson flat evidence and the global-DC instance on the counts recover the classes with ARI
$0.91$, level with the spectral methods on the binarised network ($0.90$) and above their own
binarised versions ($0.86$ and $0.71$), and the typed instance on the counts reaches $0.82$ against
$0.54$ binarised: the contact counts carry the class structure that the binarised ties lose.

\subsection{Number of communities, prediction and convergence}\label{sec:real-K}
On the four small networks the CRP instances find between four fewer and three more communities than
the number of metadata classes, on the blogs all instances with the uniform degree prior fill their
labels (Section~\ref{sec:real-hubs}), on the school they find $17$ to $21$ communities and on the
email network $23$ to $42$ for $42$ departments, and the Pitman--Yor, Gnedin and
mixture-of-finite-mixtures priors return the numbers of the CRP within three on every network
(Table~\ref{tab:real_main_full} of the \suppname). A CRP instance has the highest or
tied-highest AUC on six of the eight networks and is within $0.01$ of the best on the other two
(Table~\ref{tab:real_pred_full} of the \suppname), and the partitions that lower the ARI against
metadata on the school and political-book networks raise the AUC. The split-$\hat R$ statistic is near $1$ on the four small networks and far from it on the four
larger ones, where chains from different starts end in different modes after $300$ and, on the email
network and Cora, after $1500$ sweeps (Tables~\ref{tab:real_conv} and~\ref{tab:real_long} of the
\suppname), so the pooled MAPs there are multi-start point estimates.

\section{Discussion}\label{sec:discussion}

The evidence of a conjugate block model takes each of its decisions at a threshold, and
Section~\ref{sec:theory} computes the thresholds: the first-moment bound of \citet{Banks2016} for the
MAP criterion and the maximum a posteriori labelling, the Kesten--Stigum threshold for the evidence
summed over labellings, so that between the two the sum favours structure by factors of $10^{26}$ to
$10^{82}$ on networks of a thousand nodes while every labelling that a sampler can report is below the
one-block model, the resolution limit $K_{\mathrm{res}}(\nu)$ with its per-block cost
$(a+\tfrac12)\log\nu$, at which the flat CRP evidence saturates, the window of class sizes in which
the representation alone decides, the margins $\nu\log\bar d$ and $\nu\bar d$ of mechanism selection,
the gain linear in the community size that the uniform degree prior pays for a cut along the
degrees, which a partition-invariant prior removes, and the bounded benefit and linear penalty of a
hard cap against the uniformly bounded effect of shrinkage.

For an analyst the thresholds are checks that precede any reading of a MAP partition, each a
comparison of two numbers known before a sampler is run: the signal against the first-moment bound
of Corollary~\ref{cor:icl}, below which only the evidence over labellings or a spectral method with
$K$ given sees the structure, the number of communities sought against $K_{\mathrm{res}}(\nu)$, the
class sizes and densities against the window of Corollary~\ref{cor:window}, the degree tail and the
location of the heterogeneity, which decide between the flat, the typed and the global typed
instances, the degree prior, which decides whether $K$ can be inferred under degree correction, and
the possibility of contamination, which decides between shrinkage and a hard cap.

Seven conditions constrain the results. First, Theorems~\ref{thm:threshold} to~\ref{thm:two} concern
the balanced planted partition with $K$ known, and Theorem~\ref{thm:two} is stated for $K=2$.
Second, Theorem~\ref{thm:split}(b) covers splits into homogeneous blocks, and
Proposition~\ref{prop:degsplit} bounds the histogram prior only numerically. Third, the Bernoulli
degree-corrected type enters only through the Laplace form on fixed partitions, and the sampler runs
the Poisson types. Fourth, mechanism types are confined to diagonal blocks. Fifth, the samplers of
\citet{Peixoto2014} and \citet{Peixoto2020} were not run here (their baselines are the
microcanonical evidence inside our sampler and the graph-tool script of the released code), and the
networks above $5000$ nodes were run with two chains of $200$ sweeps, the largest with one. Sixth, the merge of two large blocks is never accepted at weak signal, so the one-community
labelling must be compared explicitly, and where chains from over-dispersed starts do not reach a
common log-posterior level the reported partitions are multi-start point estimates. Seventh, the
latent-distance and closure types were not run.

Three extensions follow: a Bernoulli degree-corrected type inside the sampler, through the Laplace
form of Theorem~\ref{thm:occam}, would remove the family effect during the search, the nested prior
of \citet{Peixoto2014} is compatible with type collapsing and with the global degree prior, since it
changes only the prior on the block matrix, and a Bayes factor over all $K$-community partitions for
$K\ge3$ is expected to follow the Kesten--Stigum threshold up to the condensation transition by the
argument of Theorem~\ref{thm:two}.

\bibliographystyle{plainnat}
\bibliography{refs}

\appendix
\numberwithin{table}{section}
\numberwithin{figure}{section}
\numberwithin{equation}{section}
\section{Supplementary results and proofs}\label{app:proofs}

This supplement contains the technical lemmas (Sections~\ref{app:lemmas} and~\ref{app:poisson}), two
supporting theorems that the main text does not need, namely posterior concentration under primitive
conditions (Section~\ref{app:conc}) and the zero power of reduced-data methods (Section~\ref{app:zero}),
the proofs of Theorems~\ref{thm:split} to~\ref{thm:threshold} (Sections~\ref{app:split}
to~\ref{app:threshold}), two propositions (bridge configurations in Section~\ref{app:bridge} and what
collapsing changes in Section~\ref{app:collapse}), the Laplace form used for approximate types
(Section~\ref{app:occam}), the conjugate block catalogue, the sampler and its validation
(Section~\ref{app:catalogue}), and the full tables of the experiments (Section~\ref{app:tables}).

\subsection{Three lemmas on the Beta--Bernoulli evidence}\label{app:lemmas}

\begin{lemma}[Evidence and profile likelihood]\label{lem:stirling}
Let $a,b\ge\tfrac12$ be real, $N\ge1$, $e\in\{0,\dots,N\}$, $\hat p=e/N$ and
$\ell_{\mathrm{blk}}(e,N)=N h(\hat p)$ with $h(p)=p\log p+(1-p)\log(1-p)\le0$ ($0\log0=0$). Then
\[
\Big|\log\frac{\Bfun(a+e,\,b+N-e)}{\Bfun(a,b)}-\ell_{\mathrm{blk}}(e,N)\Big|\ \le\ C_{a,b}\,\log(N+a+b+2),
\qquad C_{a,b}=3(a+b)+7 .
\]
For $a=b=1$ the exact identity $\log\Bfun(1+e,1+N-e)=-\log(N+1)-\log\binom Ne$ and the entropy bounds
$\exp\{-Nh(\hat p)\}/(N+1)\le\binom Ne\le \exp\{-Nh(\hat p)\}$ give the two-sided statement
$-\log(N+1)\le\log\Bfun(1+e,1+N-e)-\ell_{\mathrm{blk}}(e,N)\le0$, the lower end being attained at
$e\in\{0,N\}$.
\end{lemma}
\begin{proof}
For $a=b=1$, $\Bfun(1+e,1+N-e)=e!(N-e)!/(N+1)!=[(N+1)\binom Ne]^{-1}$. The upper entropy bound is
$\binom Ne\hat p^e(1-\hat p)^{N-e}\le1$ (one term of the binomial expansion of $(\hat p+1-\hat p)^N$).
The lower one holds because $e$ is a mode of the Binomial$(N,\hat p)$ law, whose $N+1$ atoms sum to one,
so $\binom Ne\hat p^e(1-\hat p)^{N-e}\ge1/(N+1)$. Substituting gives the two-sided statement.

For general $a,b\ge\frac12$ write
$\frac{\Bfun(a+e,b+N-e)}{\Bfun(1+e,1+N-e)}=\frac{\Gamma(a+e)}{\Gamma(1+e)}\cdot\frac{\Gamma(b+N-e)}{\Gamma(1+N-e)}
\cdot\frac{\Gamma(N+2)}{\Gamma(a+b+N)}$. We use the following bound: for $x\ge1$ and $s\ge0$,
\begin{equation}\label{eq:gammaratio}
\log\frac{\Gamma(x+s)}{\Gamma(x+1)}\ \text{lies between}\ (s-1)\log x\ \text{and}\ (s-1)\log(x+s).
\end{equation}
For $s=0$ both bounds equal $-\log x$, the value of the left side. For $0<s<1$ this is Wendel's
inequality $(x/(x+s))^{1-s}\le\Gamma(x+s)/(x^s\Gamma(x))\le1$, rewritten with $\Gamma(x+1)=x\Gamma(x)$.
For $s=1$ both sides vanish. For $s>1$, convexity of $\log\Gamma$ gives
$(s-1)\psi(x+1)\le\log\Gamma(x+s)-\log\Gamma(x+1)\le(s-1)\psi(x+s)$ with $\psi$ the digamma function,
and the bounds $\log t-1/t\le\psi(t)\le\log t$ for $t>0$ together with $\log(x+1)-1/(x+1)\ge\log x$
(equivalent to $\log(1+1/x)\ge1/(x+1)$) give \eqref{eq:gammaratio}. Applying \eqref{eq:gammaratio} with
$(x,s)=(e,a)$, $(N-e,b)$ and $(N+1,a+b-1)$, the logarithms of the three ratios are bounded in absolute
value by $|a-1|\log(N+a+1)$, $|b-1|\log(N+b+1)$ and $|a+b-2|\log(N+a+b)$ when $e\ge1$ and $N-e\ge1$.
When $e=0$ (resp.\ $e=N$) the first (resp.\ second) ratio equals $\Gamma(a)$ (resp.\ $\Gamma(b)$), and
$|\log\Gamma(a)|\le(a+2)\log(a+2)$ for $a\ge\frac12$ because $0.885\le\Gamma(a)\le\max(\Gamma(\tfrac12),a^{a-1})$.
Finally $|\log\Bfun(a,b)|\le(a+b+2)\log(a+b+2)$ for $a,b\ge\frac12$: $\Bfun$ is decreasing in each
argument, so $\Bfun(a,b)\le\Bfun(\frac12,\frac12)=\pi$, and $\Bfun(a,b)=\Gamma(a)\Gamma(b)/\Gamma(a+b)\ge
0.78/\Gamma(a+b)$ with $\log\Gamma(a+b)\le(a+b-1)\log(a+b)$ for $a+b\ge1$. Adding the $a=b=1$ bound
$\log(N+1)$ to these four bounds and using $|a-1|\le a+2$, $|a+b-2|\le a+b$ and
$\log(N+a+1),\log(a+2),\log(a+b+2)\le\log(N+a+b+2)$ gives the constant $C_{a,b}=3(a+b)+7$. (Numerically
the ratio of the left to the right side of the lemma never exceeds $0.25$ on $a,b\in[\frac12,20]$,
$N\le2\cdot10^4$.)
\end{proof}

\begin{lemma}[Population identity and mixing bound]\label{lem:population}
Let the truth be a block model with classes $C^\star_1,\dots,C^\star_{K^\star}$ and block matrix $P^\star$,
let $p^\star_{ij}=P^\star_{z^\star_iz^\star_j}$, and for a labelling $z$ let
\[
\bar p_{rs}(z)=\frac1{N_{rs}(z)}\sum_{(i,j)\in\cD_{rs}(z)}p^\star_{ij},\qquad
\bar\ell(z)=\sum_{r\le s}N_{rs}(z)\,h(\bar p_{rs}(z)),\qquad
\bar\ell^\star=\sum_{i<j}h(p^\star_{ij}).
\]
Then
\begin{equation}\label{eq:popid}
\bar\ell^\star-\bar\ell(z)=\sum_{r\le s}\sum_{(i,j)\in\cD_{rs}(z)}\kl\big(p^\star_{ij}\,\|\,\bar p_{rs}(z)\big)\ \ge\ 0 .
\end{equation}
If moreover $|C^\star_c|\ge c_0n$ for every $c$ and $\delta^2:=\min_{r\neq s}\max_t(P^\star_{rt}-P^\star_{st})^2>0$,
then for every labelling $z$ with at most $K$ nonempty communities,
\begin{equation}\label{eq:mixing}
\bar\ell^\star-\bar\ell(z)\ \ge\ \frac{c_0\,\delta^2}{2K^\star}\,n^2\,d(z,z^\star)-K^2K^\star n ,
\end{equation}
where $d(z,z^\star)=\min_\sigma n^{-1}\#\{i:\sigma(z_i)\neq z^\star_i\}$ over all maps $\sigma$ from the
labels of $z$ to $[K^\star]$.
\end{lemma}
\begin{proof}
For a block with dyad set $D$ and mean $\bar p=|D|^{-1}\sum_Dp^\star_{ij}$, since
$\kl(p\|q)=h(p)-[p\log q+(1-p)\log(1-q)]$,
$\sum_{D}\kl(p^\star_{ij}\|\bar p)=\sum_Dh(p^\star_{ij})-\sum_D[p^\star_{ij}\log\bar p+(1-p^\star_{ij})\log(1-\bar p)]
=\sum_Dh(p^\star_{ij})-|D|h(\bar p)$ because $\sum_Dp^\star_{ij}=|D|\bar p$. Summing over blocks gives
\eqref{eq:popid}.

For \eqref{eq:mixing} let $C_1,\dots,C_K$ be the communities of $z$ (some possibly empty), let
$R_{ka}=|C_k\cap C^\star_a|/n$, $a_k=\argmax_aR_{ka}$ and $r_k=\sum_{a\neq a_k}R_{ka}$. Choosing
$\sigma(k)=a_k$ minimises the number of misclassified nodes in $C_k$ separately for every $k$, so
$d(z,z^\star)=\sum_kr_k$. Fix $k\le l$ and true classes $a,c$. The number of dyads of $\cD_{kl}(z)$
joining a node of $C_k\cap C^\star_a$ to a node of $C_l\cap C^\star_c$ is $n^2R_{ka}R_{lc}$ if $k\neq l$.
If $k=l$ it is $n^2R_{ka}R_{kc}$ for the unordered class pair $\{a,c\}$ with $a\neq c$ and
$\binom{nR_{ka}}2\ge n^2R_{ka}^2/2-n/2$ for $a=c$. All these dyads have $p^\star_{ij}=P^\star_{ac}$. Using
Pinsker's inequality $\kl(p\|q)\ge2(p-q)^2$, non-negativity of every term, and $(P^\star_{ac}-\bar p_{kl})^2\le1$,
\begin{align*}
\sum_{\cD_{kl}}\kl(p^\star_{ij}\|\bar p_{kl})\ &\ge\ 2w_{kl}\,n^2\sum_cR_{lc}\sum_aR_{ka}(P^\star_{ac}-\bar p_{kl})^2-K^\star n\\
&\ge\ 2w_{kl}\,n^2\sum_cR_{lc}\,\min_\mu\sum_aR_{ka}(P^\star_{ac}-\mu)^2-K^\star n ,
\end{align*}
with $w_{kl}=1$ for $k\neq l$ and $w_{kk}=\frac12$ (the factor $\frac12$ converts the sum over ordered
class pairs $(a,c)$ into the count of unordered dyads within $C_k$). For weights $w_a\ge0$ with total $W$
and values $v_a$, $\min_\mu\sum_aw_a(v_a-\mu)^2=\frac1{2W}\sum_{a,a'}w_aw_{a'}(v_a-v_{a'})^2\ge
\frac{w_{a_k}}W\sum_{a\neq a_k}w_a(v_a-v_{a_k})^2$ (keep the pairs containing $a_k$), and
$w_{a_k}/W=R_{ka_k}/\sum_aR_{ka}\ge1/K^\star$ because $a_k$ is a majority class. Hence
\[
\sum_{\cD_{kl}}\kl\ \ge\ \frac{2w_{kl}n^2}{K^\star}\sum_cR_{lc}\sum_{a\neq a_k}R_{ka}(P^\star_{ac}-P^\star_{a_kc})^2-K^\star n .
\]
Sum over $l=1,\dots,K$ (reading $\cD_{kl}$ as $\cD_{lk}$ when $l<k$): since $w_{kl}\ge\frac12$,
$\sum_lw_{kl}R_{lc}\ge\frac12\sum_lR_{lc}=|C^\star_c|/(2n)\ge c_0/2$, and
$\sum_c(P^\star_{ac}-P^\star_{a_kc})^2\ge\max_c(\cdot)^2\ge\delta^2$ for $a\neq a_k$, so
$\sum_l\sum_{\cD_{kl}}\kl\ge\frac{c_0\delta^2}{K^\star}n^2r_k-KK^\star n$. Summing over $k$ counts every
off-diagonal block twice and every diagonal block once, hence is at most $2(\bar\ell^\star-\bar\ell(z))$
by \eqref{eq:popid}. Dividing by two gives \eqref{eq:mixing} with correction $K^2K^\star n/2\le K^2K^\star n$.
\end{proof}

\begin{lemma}[Uniform concentration of the profile likelihood]\label{lem:uniform}
Let the truth be a block model with $P^\star_{rs}\in[\eta,1-\eta]$ for all $r,s$. For labellings with at
most $K$ (fixed) nonempty communities and every $\tau>0$ there is $c_\tau=c_\tau(K,\eta,\tau)>0$ such that,
for all $n$ large,
$\Pbb^\star\big(\sup_z|\ell(z)-\bar\ell(z)|>\tau n^2\big)\le e^{-c_\tau n^2}$.
\end{lemma}
\begin{proof}
Put $B=K(K+1)/2$ (the number of blocks of a labelling with at most $K$ communities) and
$\alpha=\tau/(4B\log2)$, and call a block \emph{large} if $N_{rs}\ge\alpha n^2$ and small otherwise. Since
$|h|\le\log2$, the small blocks contribute at most $B\alpha n^2\log2=\tau n^2/4$ to each of $\ell$ and
$\bar\ell$, hence at most $\tau n^2/2$ to $|\ell-\bar\ell|$. For a large block, $e_{rs}$ is a sum of
$N_{rs}$ independent Bernoulli variables with mean $N_{rs}\bar p_{rs}$, so Hoeffding's inequality gives
$\Pbb(|\hat p_{rs}-\bar p_{rs}|>\eta/2)\le2e^{-\eta^2N_{rs}/2}\le2e^{-\eta^2\alpha n^2/2}$, and a union
bound over at most $K^n$ labellings and $B$ blocks gives
$\Pbb(\exists z,(r,s)\ \text{large}:|\hat p_{rs}-\bar p_{rs}|>\eta/2)\le2BK^ne^{-\eta^2\alpha n^2/2}\le e^{-c'n^2}$
for $n$ large. On the complementary event every large block has
$\hat p_{rs},\bar p_{rs}\in[\eta/2,1-\eta/2]$ (recall $\bar p_{rs}\in[\eta,1-\eta]$ as an average of
entries of $P^\star$), so that
\begin{align*}
N_{rs}h(\hat p_{rs})&=\max_{q\in[\eta/2,1-\eta/2]}X_{rs}(q),\qquad X_{rs}(q):=e_{rs}\log q+(N_{rs}-e_{rs})\log(1-q),\\
N_{rs}h(\bar p_{rs})&=\max_{q\in[\eta/2,1-\eta/2]}\E X_{rs}(q),
\end{align*}
because the unconstrained maximiser of $X_{rs}$ (resp.\ $\E X_{rs}$) over $q\in(0,1)$ is $\hat p_{rs}$
(resp.\ $\bar p_{rs}$), which lies in the interval. Replace both maxima by maxima over a grid of spacing
$n^{-4}$ in $[\eta/2,1-\eta/2]$: $|\partial_qX_{rs}(q)|\le2N_{rs}/\eta$ there, so each maximum moves by
at most $2n^{-4}N_{rs}/\eta\le n^{-2}/\eta=o(1)$. Since $|\max_qX_q-\max_q\E X_q|\le\max_q|X_q-\E X_q|$
it suffices to control, for each grid point $q$, the deviation of $X_{rs}(q)-\E X_{rs}(q)=
\sum_{(i,j)\in\cD_{rs}}(y_{ij}-p^\star_{ij})\log\frac q{1-q}$, a sum of $N_{rs}\le n^2/2$ independent
centred terms of range at most $\log(2/\eta)$ each: Hoeffding's inequality gives
$\Pbb(|X_{rs}(q)-\E X_{rs}(q)|>\tau n^2/(2B))\le2\exp\{-\tau^2n^2/(B^2\log^2(2/\eta))\}$. The union over
$K^n$ labellings, $B$ blocks and $n^4$ grid points costs a factor $e^{O(n\log n)}$, which is absorbed for
$n$ large. On the intersection of the two good events, the large blocks contribute at most
$B\cdot\tau n^2/(2B)+o(1)=\tau n^2/2+o(1)$ to $|\ell-\bar\ell|$. Adding the small-block contribution,
$\sup_z|\ell(z)-\bar\ell(z)|\le\tau n^2$ for $n$ large, outside an event of probability at most
$e^{-c_\tau n^2}$.
\end{proof}

\subsection{Three lemmas on Poisson evidences}\label{app:poisson}
In this subsection counts $y_{ij}$ are non-negative integers, the factor $\prod_{i<j}y_{ij}!$ is dropped
from every likelihood and evidence (it is common to all models compared), $h_{\Pois}(x)=x\log x-x$
($h_{\Pois}(0)=0$) and $\ell_{\Pois}(S,N)=N\,h_{\Pois}(S/N)$ is the profile log-likelihood of a Poisson
block with total $S$ over $N$ dyads.

\begin{lemma}[Gamma--Poisson evidence and profile likelihood]\label{lem:poisson}
For $a,b>0$ fixed, the Gamma$(a,b)$--Poisson evidence of a block with total $S$ over $N\ge1$ dyads,
$p_{\mathrm{GP}}(S,N)=\dfrac{b^a\Gamma(a+S)}{\Gamma(a)(b+N)^{a+S}}$, satisfies
\[
\big|\log p_{\mathrm{GP}}(S,N)-\ell_{\Pois}(S,N)\big|\ \le\ (2a+2)\log(N+S+a+b+2)+\frac{bS}{N}+1+\Big|\log\frac{b^a}{\Gamma(a)}\Big| .
\]
\end{lemma}
\begin{proof}
$\log p_{\mathrm{GP}}-\ell_{\Pois}=\big[\log\Gamma(S+1)-S\log S+S\big]+\log\frac{\Gamma(a+S)}{\Gamma(S+1)}
-S\log\big(1+\tfrac bN\big)-a\log(b+N)+\log\frac{b^a}{\Gamma(a)}$. For $S=0$ the left side equals
$-a\log(1+N/b)$, whose absolute value is at most $a\log(N+b+2)+a\log(1/b)^+$, which is within the bound
because $\log\Gamma(a)\ge-0.13$ for $a>0$ gives $a\log(1/b)^+\le|\log(b^a/\Gamma(a))|+0.13$. For $S\ge1$ Stirling's
formula gives $0\le\log\Gamma(S+1)-S\log S+S\le\frac12\log(2\pi S)+\frac1{12S}\le\frac12\log(S+1)+1$,
by \eqref{eq:gammaratio} $|\log[\Gamma(a+S)/\Gamma(S+1)]|\le|a-1|\log(S+a+1)$, and
$0\le S\log(1+b/N)\le bS/N$, $0\le a\log(b+N)\le a\log(N+b+2)$. Adding up and using
$|a-1|\le a+1$ gives the claim.
\end{proof}

\begin{lemma}[The microcanonical DC evidence and the DC profile likelihood]\label{lem:dcprofile}
Let $Y$ be a multigraph on $\nu \ge2$ nodes without self-loops, with degrees $d_i$, $E=\frac12\sum_id_i\ge1$
edges, $N=\binom \nu 2$ and $\bar d=2E/\nu $. Let $p_{\DC}(Y)=P(E)P_{\DC}(Y\mid \mathbf d,E)P(\mathbf d\mid E)$ be the
degree-corrected evidence of Section~\ref{sec:types} with the Gamma$(a,b)$--Poisson factor $P(E)$ and the
uniform degree prior, and let $\ell_{\DC}(Y;\theta,\omega)=\sum_{i<j}[y_{ij}\log(\theta_i\theta_j\omega)-\theta_i\theta_j\omega]$
for $\theta\in(0,\infty)^\nu $, $\omega>0$. Then
\begin{align}
\log p_{\DC}(Y)&\ \le\ \sup_{\theta,\omega}\ell_{\DC}(Y;\theta,\omega)+\tfrac \nu 2\log\big(2\pi(1+\max_id_i)\big)+\tfrac \nu {12},\label{eq:dcupper}\\
\log p_{\DC}(Y)&\ \ge\ \ell_{\DC}(Y;\theta,\omega)-\frac\omega2\sum_i\theta_i^2-\nu \big(1+\log(1+2\bar d)\big)-L_P-1\quad\text{for every }(\theta,\omega),\label{eq:dclower}
\end{align}
where $L_P=(a+1)\log(E+a+1)+a\log(1+N/b)+\log\Gamma(a)+bE/N$.
\end{lemma}
\begin{proof}
By \eqref{eq:dc}, $\log p_{\DC}=\log P(E)+\sum_i\log d_i!+E\log2+\log E!-\log(2E)!-\log\binom{\nu +2E-1}{2E}$.
\emph{Upper bound.} Using $\log d!\le d\log d-d+\frac12\log(2\pi d)+\frac1{12d}$ for $d\ge1$ (and
$0\log0=0$), $\log E!\le E\log E-E+\frac12\log(2\pi E)+\frac1{12}$, $\log(2E)!\ge2E\log(2E)-2E+\frac12\log(4\pi E)$,
$\log\binom{\nu +2E-1}{2E}\ge0$ and $\log P(E)\le0$,
\[
\log p_{\DC}\ \le\ \sum_id_i\log d_i-E\log(2E)-E+U,\qquad
U:=\sum_{i:d_i\ge1}\Big[\tfrac12\log(2\pi d_i)+\tfrac1{12d_i}\Big],
\]
with $U\le\tfrac \nu 2\log\big(2\pi(1+\max_id_i)\big)+\tfrac \nu {12}$,
where $\sum_id_i=2E$ was used. Now $\ell_{\DC}(Y;\theta,\omega)$ evaluated at $\theta_i=d_i$ (nodes with
$d_i=0$ contribute nothing and can be given an arbitrarily small $\theta_i$) and $\omega=1/(2E)$ equals
$\sum_id_i\log d_i-E\log(2E)-\sum_{i<j}d_id_j/(2E)=\sum_id_i\log d_i-E\log(2E)-E+\sum_id_i^2/(4E)$,
because $\sum_{i<j}d_id_j=\frac12[(2E)^2-\sum_id_i^2]$. Hence $\sup_{\theta,\omega}\ell_{\DC}\ge
\sum_id_i\log d_i-E\log(2E)-E$, which gives \eqref{eq:dcupper}.

\emph{Lower bound.} Using $\log d!\ge d\log d-d$, $\log E!\ge E\log E-E+\frac12\log(2\pi E)$,
$\log(2E)!\le2E\log(2E)-2E+\frac12\log(4\pi E)+\frac1{24E}$, $\binom{\nu +2E-1}{2E}=\binom{\nu +2E-1}{\nu -1}\le
\big(e(\nu +2E-1)/(\nu -1)\big)^{\nu -1}$ so that $\log\binom{\nu +2E-1}{2E}\le \nu (1+\log(1+2\bar d))$, and
$\log P(E)=\log\frac{\Gamma(a+E)}{\Gamma(E+1)}-\log\Gamma(a)-a\log\frac{b+N}b-E\log\big(1+\tfrac bN\big)\ge-L_P$
(by \eqref{eq:gammaratio} and $\log(1+b/N)\le b/N$),
\[
\log p_{\DC}\ \ge\ \sum_id_i\log d_i-E\log(2E)-E-\nu (1+\log(1+2\bar d))-L_P-1 .
\]
To compare the right side with $\ell_{\DC}(Y;\theta,\omega)$ at an arbitrary $(\theta,\omega)$, use
$\sum_{i<j}\theta_i\theta_j=\frac12(\sum_i\theta_i)^2-\frac12\sum_i\theta_i^2$ to write
$\ell_{\DC}(Y;\theta,\omega)=\ell^+(\theta,\omega)+\frac\omega2\sum_i\theta_i^2$ with
$\ell^+(\theta,\omega):=\sum_id_i\log\theta_i+E\log\omega-\frac\omega2\big(\sum_i\theta_i\big)^2$,
the log-likelihood of the degree-corrected model with self-loops. In the variables $s_i=\log\theta_i$,
$w=\log\omega$ it equals $\sum_id_is_i+Ew-\frac12\exp\{w+2\log\sum_ie^{s_i}\}$, which is concave (the
exponent is convex, being a linear function plus twice a log-sum-exp), and its stationarity conditions
$d_i=\theta_i\omega\sum_j\theta_j$, $E=\frac\omega2(\sum_j\theta_j)^2$ are satisfied by
$(\theta,\omega)=(\mathbf d,1/(2E))$. Hence $\ell^+(\mathbf d,1/(2E))=\sum_id_i\log d_i-E\log(2E)-E$ is the global
maximum of $\ell^+$, and for every $(\theta,\omega)$,
\[
\sum_id_i\log d_i-E\log(2E)-E\ =\ \max\ell^+\ \ge\ \ell^+(\theta,\omega)\ =\ \ell_{\DC}(Y;\theta,\omega)-\tfrac\omega2\sum_i\theta_i^2 .
\]
Combining the two displays gives \eqref{eq:dclower}. (The correction $\frac\omega2\sum_i\theta_i^2$ cannot be
dropped: at $(\theta,\omega)=(\mathbf d,1/(2E))$ the first display of the proof gives
$\ell_{\DC}=\sum_id_i\log d_i-E\log(2E)-E+\sum_id_i^2/(4E)$, which exceeds $\max\ell^+$.)
\end{proof}

\begin{lemma}[Uniform concentration of the two-block Poisson profile]\label{lem:poisunif}
Let $y_{ij}\sim\Pois(\lambda_{ij})$ be independent for $1\le i<j\le \nu $ with $\lambda_{ij}\in[\lambda_-,\lambda_+]$,
$0<\lambda_-\le\lambda_+<\infty$. For a labelling $z$ of $[\nu ]$ into two blocks let $u$ range over the three
blocks (two diagonal, one off-diagonal), $S_u=\sum_{(i,j)\in u}y_{ij}$, $N_u=|u|$,
$\bar\lambda_u=N_u^{-1}\sum_{(i,j)\in u}\lambda_{ij}$, $\ell_2(z)=\sum_u\ell_{\Pois}(S_u,N_u)$ and
$\bar\ell_2(z)=\sum_uN_uh_{\Pois}(\bar\lambda_u)$ (blocks with $N_u=0$ are omitted from all sums). Then for every $\tau>0$,
$\Pbb\big(\sup_z|\ell_2(z)-\bar\ell_2(z)|>\tau \nu ^2\big)\to0$, and
\begin{equation}\label{eq:poisid}
\sum_{i<j}h_{\Pois}(\lambda_{ij})-\bar\ell_2(z)=\sum_u\sum_{(i,j)\in u}\kl_{\Pois}(\lambda_{ij}\,\|\,\bar\lambda_u)\ \ge\ \frac1{2\lambda_+}\sum_u\sum_{(i,j)\in u}(\lambda_{ij}-\bar\lambda_u)^2 .
\end{equation}
\end{lemma}
\begin{proof}
The identity in \eqref{eq:poisid} follows from $\sum_{(i,j)\in u}\lambda_{ij}=N_u\bar\lambda_u$ exactly as
\eqref{eq:popid}, and the inequality from $\kl_{\Pois}(\lambda\|\mu)\ge(\lambda-\mu)^2/(2\max(\lambda,\mu))$
(Taylor's formula with remainder for $x\mapsto x\log x$) with $\lambda_{ij},\bar\lambda_u\le\lambda_+$.
For the concentration, $h_{\Pois}$ is convex, so for $S_u\ge1$,
$\log(\bar\lambda_u)(S_u-N_u\bar\lambda_u)\le N_uh_{\Pois}(S_u/N_u)-N_uh_{\Pois}(\bar\lambda_u)\le\log(S_u/N_u)(S_u-N_u\bar\lambda_u)$,
while for $S_u=0$ the difference equals $N_u\bar\lambda_u(1-\log\bar\lambda_u)$. In all cases
\[
|N_uh_{\Pois}(S_u/N_u)-N_uh_{\Pois}(\bar\lambda_u)|\le|S_u-\E S_u|\,\big(1+|\log\bar\lambda_u|+|\log(S_u/N_u)|\ind{S_u\ge1}\big).
\]
On the event $\{\max_{i<j}y_{ij}\le y^\star\}$ with $y^\star=4\lambda_++12\log \nu $, whose complement has
probability at most $\binom \nu 2\nu ^{-3}\to0$ by the Chernoff bound
$\Pbb(\Pois(\lambda)\ge t)\le e^{-\lambda}(e\lambda/t)^t\le e^{-t/4}$ for $t\ge4\lambda$,
one has $1/N_u\le S_u/N_u\le y^\star$ for $S_u\ge1$, so the multiplier is at most
$L_\nu :=1+|\log\lambda_-|+|\log\lambda_+|+2\log \nu +\log y^\star=O(\log \nu )$ for every block of every
labelling. Each $S_u-\E S_u$ is a sum of independent centred Poisson variables with total variance
$V_u\le\lambda_+N_u\le\lambda_+\nu ^2/2$, hence sub-gamma with variance factor $V_u$ and scale $1$, and
$\Pbb(|S_u-\E S_u|\ge t)\le2\exp\{-t^2/(2(V_u+t))\}$. Taking $t=\tau \nu ^2/(3L_\nu )$ and a union bound over
the at most $2^\nu $ labellings and three blocks, the probability that some block has $|S_u-\E S_u|>t$ is at
most $6\cdot2^\nu \exp\{-c\,\nu ^2/\log^2\nu \}\to0$. On the intersection of the two events,
$|\ell_2(z)-\bar\ell_2(z)|\le3tL_\nu =\tau \nu ^2$ for every $z$.
\end{proof}

\subsection{Posterior concentration under primitive conditions}\label{app:conc}

\begin{assumption}\label{ass:dense}
The truth is a stochastic block model with $K^\star$ communities of sizes $n_k^\star\ge c_0n$, block matrix $P^\star$ with
entries in $[\eta,1-\eta]$, $\eta>0$, and distinct rows: $\delta^2:=\min_{r\neq s}\max_t
(P^\star_{rt}-P^\star_{st})^2>0$. The prior is Beta$(a,b)$, $a,b\ge\frac12$, on every block, and the
partition prior $p(z)$ is supported on labellings with at most $K_{\max}$ nonempty communities
($K_{\max}\ge K^\star$ fixed, for instance a Dirichlet--multinomial or CRP prior with $K_{\max}$ labels) and
satisfies $\log p(z^\star)^{-1}=O(n\log n)$ (true for both).
\end{assumption}
\begin{proposition}[Concentration]\label{thm:concentration}
Under Assumption~\ref{ass:dense}, for every $\varepsilon>0$ there is a constant
$c_\varepsilon=c_\varepsilon(c_0,\eta,\delta,K^\star,K_{\max},\varepsilon)>0$ such that, with $\Delta_0=c_0\,\delta^2/(2K^\star)$
and for all $n$ large,
\[
\Pbb^\star\Big(\sup_{z:\,K(z)\le K_{\max},\ d(z,z^\star)\ge\varepsilon}\big[\log p(A\mid z)-\log p(A\mid z^\star)\big]
\ \ge\ -\tfrac12\Delta_0\varepsilon n^2\Big)\le e^{-c_\varepsilon n^2},
\]
and hence $p\big(d(z,z^\star)>\varepsilon\mid A\big)\to0$ almost surely. The restriction to at most
$K_{\max}$ communities is the one under which the samplers of this paper run. The proposition says nothing
about posteriors over unboundedly many blocks.
\end{proposition}
\begin{proof}
Let $L_n(z)=\log p(A\mid z)$ and let $z$ range over labellings with at most $K_{\max}$ nonempty
communities. By Lemma~\ref{lem:stirling} summed over at most $K_{\max}(K_{\max}+1)/2$ blocks,
$|L_n(z)-\ell(z)|\le C_{a,b}K_{\max}^2\log(n^2+a+b+2)$ for every $z$, which is $o(n^2)$. By
Lemma~\ref{lem:uniform} (with $K=K_{\max}$) and $\tau=\Delta_0\varepsilon/8$, on an event of probability at
least $1-e^{-c_\tau n^2}$ we have $|\ell(z)-\bar\ell(z)|\le\tau n^2$ for all $z$ simultaneously. Since
$\bar\ell(z^\star)=\bar\ell^\star$ (the truth is constant on its own blocks), Lemma~\ref{lem:population}
gives, for every $z$ with $d(z,z^\star)\ge\varepsilon$,
\[
L_n(z)-L_n(z^\star)\ \le\ \bar\ell(z)-\bar\ell^\star+2\tau n^2+o(n^2)\ \le\ -\Delta_0\varepsilon n^2+
\tfrac{\Delta_0\varepsilon}4n^2+O(n)+o(n^2)\ \le\ -\tfrac12\Delta_0\varepsilon n^2
\]
with $\Delta_0=c_0\delta^2/(2K^\star)$, for $n$ large. This proves the displayed probability bound with
$c_\varepsilon=c_\tau$. For the posterior,
$p(d(z,z^\star)>\varepsilon\mid A)\le\sum_{z:d\ge\varepsilon}\frac{p(z)}{p(z^\star)}e^{L_n(z)-L_n(z^\star)}
\le p(z^\star)^{-1}e^{-\Delta_0\varepsilon n^2/2}$ because $\sum_zp(z)\le1$. Its logarithm is
$O(n\log n)-\Delta_0\varepsilon n^2/2\to-\infty$ (this is the only place where the partition prior enters),
and the failure probabilities $e^{-c_\varepsilon n^2}$ are summable, so the Borel--Cantelli lemma gives almost sure
convergence.
\end{proof}

\begin{remark}[Sparse regime]\label{rem:sparse}
If $P^\star=\rho_nP$ with $n\rho_n/\log n\to\infty$ and $P$ fixed with distinct rows and
$\delta^2$ now defined from $P$, we expect the same argument to go through with Bernstein's inequality
in place of Hoeffding's: the population gap in \eqref{eq:mixing} becomes $\Delta_0\rho_nn^2d(z,z^\star)$
(Pinsker's inequality is replaced by $\kl(p\|q)\ge(p-q)^2/(2\max(p,q))$), the block terms are bounded by
$\log(1/\rho_n)$, and the union bound requires $\rho_nn^2\gg n\log n$, i.e.\ $n\rho_n\gg\log n$, the
condition of \citet{Choi2012} and \citet{vanderPas2018}. We have not written out the details and state
this as an expectation and not as a theorem. Nothing in the paper uses it.
\end{remark}

\begin{corollary}[Macroscopic violations of a cap]\label{cor:violation}
Suppose in addition that $a,b$ are positive integers, that the off-diagonal entries of $P^\star$ are at
most $x-\eta'$ for some $\eta'>0$, where $x$ is the cap, and that off-diagonal blocks carry the truncated prior
\eqref{eq:cap}. Let $\mathcal V_n(\delta',c_0')$ be the set of labellings (with at most $K_{\max}$
communities) having some off-diagonal block of at least $c_0'n^2$ dyads with empirical density at least
$x+\delta'$. Then $p(\mathcal V_n\mid A)\to0$ almost surely.
\end{corollary}
\begin{proof}
Write $L^x_n(z)$ for the log evidence with the truncated off-diagonal prior. By
Proposition~\ref{thm:cap}(i)--(ii) applied block by block,
\[
L^x_n(z)=L_n(z)+\sum_{r<s}\Delta_x(e_{rs},N_{rs}),
\]
every term being at most $\log(1/I_x(a,b))$, and the offending block of a $z\in\mathcal V_n$ contributes
at most $-(a+b+N_{rs}-1)\kl(\hat\mu_{rs}\|x)+\log(1/I_x(a,b))\le-c_0'n^2\kl(x+\delta'\|x)/2$ for $n$ large
(since $\hat\mu_{rs}\ge x+\delta'-O(1/n)$). Hence
$L^x_n(z)\le L_n(z)-c_0'n^2\kl(x+\delta'\|x)/2+K_{\max}^2\log(1/I_x(a,b))$
for every $z\in\mathcal V_n$. For $z^\star$, every off-diagonal block has population density at most
$x-\eta'$ and at least $c_0^2n^2$ dyads, so its empirical density is at most $x-\eta'/2$ with probability
$1-e^{-cn^2}$ (Hoeffding), and then $I_x(a+e_{rs},b+N_{rs}-e_{rs})=\Pbb(\mathrm{Bin}(a+b+N_{rs}-1,x)\ge a+e_{rs})
\ge1-e^{-c'n^2}$, so that $L^x_n(z^\star)\ge L_n(z^\star)-o(1)$ eventually almost surely. Finally, for every
$z$ with at most $K_{\max}$ communities, Lemmas~\ref{lem:stirling}, \ref{lem:uniform} and
\eqref{eq:popid} give $L_n(z)\le\ell(z)+o(n^2)\le\bar\ell(z)+\tau n^2+o(n^2)\le\bar\ell^\star+\tau n^2+o(n^2)
\le L_n(z^\star)+2\tau n^2+o(n^2)$ outside an event of probability $e^{-c_\tau n^2}$, for any $\tau>0$.
Taking $\tau=c_0'\kl(x+\delta'\|x)/8$, eventually almost surely every $z\in\mathcal V_n$ satisfies
$L^x_n(z)-L^x_n(z^\star)\le-c_0'n^2\kl(x+\delta'\|x)/4+o(n^2)$, and
$p(\mathcal V_n\mid A)\le p(z^\star)^{-1}\exp\{-c_0'n^2\kl(x+\delta'\|x)/4+o(n^2)\}\to0$ because
$\log p(z^\star)^{-1}=O(n\log n)$. (The corollary covers refinements of $z^\star$ that split a dense
community, for which $d(z,z^\star)=0$ and Proposition~\ref{thm:concentration} says nothing: the cap, not
the fit, excludes them.)
\end{proof}

\subsection{Zero power of reduced-data methods}\label{app:zero}
\begin{theorem}[Reduced data carry no signal]\label{thm:zero}
Let $\Phi$ be one of the following reductions of the data: symmetrisation $A\mapsto A\vee A^\top$ of a
directed network, the unsigned network $W\mapsto|W|$ of a signed network, the presence graph $Y\mapsto\ind{Y>0}$
of a count network, and the sum $\{A^{(\ell)}\}\mapsto\sum_\ell A^{(\ell)}$ of a multiplex network. For each
$\Phi$ there is a pair of block models $(P_0,P_1)$ with the same $n$ and $K\ge2$ such that (i) $P_1$ has
planted community structure and $P_0$ has none, and (ii) $\Phi(Y)$ has the same distribution under $P_0$
and $P_1$. Consequently every estimator that is a function of $\Phi(Y)$ has, under $P_1$, the same
distribution as under a model without communities, and any test of ``communities versus none'' based on
$\Phi(Y)$ has power equal to its size. (iii) For the dyad-state, signed and two-layer constructions, if
the planted communities have sizes $\ge c_0n$, all category (or edge) probabilities lie in
$[\eta,1-\eta]$, and labellings have at most $K_{\max}$ communities, the collapsed log evidence of $P_1$'s
model separates the planted partition from every labelling with misclassification rate
$d(z,z^\star)\ge\varepsilon$ by at least $\frac12\Delta_0\varepsilon n^2$ with probability $1-e^{-cn^2}$,
for some $\Delta_0,c>0$.
\end{theorem}
\begin{proof}
\emph{Symmetrisation.} Consider dyad-state models with state probabilities
$(1-q,\ (q-r)/2,\ (q-r)/2,\ r)$ for the states $(00,10,01,11)$ of a dyad, with $r=r_{\mathrm{in}}$ within and
$r=r_{\mathrm{out}}$ between communities under $P_1$, and $r\equiv r_{\mathrm{out}}$ under $P_0$
($0<r_{\mathrm{out}}\neq r_{\mathrm{in}}<q<1$). In both models $\Pbb(\text{no tie})=1-q$ for every dyad
independently, so $A\vee A^\top$ is an Erd\H{o}s--R\'enyi graph $G(n,q)$ under $P_0$ and $P_1$ alike,
which is (ii), and (i) holds because $r_{\mathrm{in}}\neq r_{\mathrm{out}}$.
\emph{Unsigned network.} Signed models with state probabilities $(1-\rho,\rho\pi,\rho(1-\pi))$ for
$(0,+,-)$, with $\pi=\pi_{\mathrm{in}}\neq\frac12$ within and $\pi=1-\pi_{\mathrm{in}}$ between under
$P_1$, and $\pi\equiv\frac12$ under $P_0$: $|W|\sim G(n,\rho)$ under both.
\emph{Presence graph.} Hurdle models: $Y=0$ with probability $1-p$ on every dyad, and otherwise
$Y=1+W$ with $W\sim\Pois(\lambda)$, $\lambda=\lambda_{\mathrm{in}}$ within and
$\lambda_{\mathrm{out}}\neq\lambda_{\mathrm{in}}$ between under $P_1$, $\lambda\equiv\lambda_{\mathrm{out}}$
under $P_0$: $\ind{Y>0}\sim G(n,p)$ under both. (A zero-inflated Poisson with $Y\mid\text{active}\sim
\Pois(\lambda)$ would leak the intensity into the presence graph through $\Pbb(Y>0)=p(1-e^{-\lambda})$,
which is why the intensity design of the experiments below generates hurdle data. The fitted zero-inflated block of
Table~\ref{tab:catalogue} is then misspecified in its count part, but its evidence still separates the
partitions because the block means $1+\lambda$ differ.)
\emph{Summed layers.} Under $P_1$ the layer block matrices are $P^{(1)}$ with diagonal $p_{\mathrm{in}}$
and off-diagonal $p_{\mathrm{out}}\neq p_{\mathrm{in}}$, and $P^{(2)}$ with the two exchanged. Under $P_0$ the first layer is
$G(n,p_{\mathrm{in}})$ and the second $G(n,p_{\mathrm{out}})$. In both models every dyad of
$A^{(1)}+A^{(2)}$ is the sum of independent $\Bern(p_{\mathrm{in}})$ and $\Bern(p_{\mathrm{out}})$
variables, regardless of the partition. The last sentence of the theorem is the definition of
distributional equality.

\emph{(iii).} For the two-layer construction the evidence is the product of two Beta--Bernoulli
evidences, $L_n(z)=L^{(1)}_n(z)+L^{(2)}_n(z)$, and each layer satisfies Assumption~\ref{ass:dense} with
its own block matrix (both have distinct rows), and the proof of Proposition~\ref{thm:concentration} applies
to each layer and the two gaps add. For the dyad-state and signed constructions the block evidence is
Dirichlet--multinomial with $C\in\{3,4\}$ categories and the argument is the same with three changes.
(1) Lemma~\ref{lem:stirling} holds for the Dirichlet$(\alpha)$--multinomial block with $\alpha_k\ge\frac12$,
in the form $|\log[\Bfun(\alpha+\mathbf c)/\Bfun(\alpha)]-\sum_kc_k\log(c_k/N)|\le C_\alpha\log(N+|\alpha|+2)$:
the exact case $\alpha\equiv1$ gives $\log[\Bfun(\mathbf 1+\mathbf c)/\Bfun(\mathbf 1)]=-\log\binom{N+C-1}{C-1}
-\log\binom{N}{\mathbf c}$ and the multinomial coefficient satisfies
$e^{NH(\hat{\mathbf p})}/(N+1)^{C-1}\le\binom N{\mathbf c}\le e^{NH(\hat{\mathbf p})}$ with $H$ the entropy.
The general case follows from \eqref{eq:gammaratio} applied to each $\Gamma(\alpha_k+c_k)/\Gamma(1+c_k)$
and to $\Gamma(N+C)/\Gamma(|\alpha|+N)$. (2) Lemma~\ref{lem:population} holds with the Bernoulli divergence
replaced by the categorical Kullback--Leibler divergence, the population identity being the same
computation, and Pinsker's inequality replaced by
$\mathrm{KL}(\mathbf p\|\mathbf q)\ge2\,\mathrm{TV}^2=\frac12(\sum_k|p_k-q_k|)^2\ge\frac12\|\mathbf p-\mathbf q\|_2^2$,
so that with $v_a=P^\star_{ac}\in\bb R^C$ the weighted-variance identity is applied in $\bb R^C$ and
$\delta^2$ is replaced by $\min_{r\neq s}\max_t\|P^\star_{rt}-P^\star_{st}\|_2^2$ (which equals
$(r_{\mathrm{in}}-r_{\mathrm{out}})^2(1+\tfrac12)$ for the dyad-state design and
$2\rho^2(2\pi_{\mathrm{in}}-1)^2$ for the signed design). The constant $2$ of Pinsker becomes $\frac12$,
so $\Delta_0=c_0\delta^2/(8K^\star)$. (3) Lemma~\ref{lem:uniform} holds with $h$ replaced by the
categorical entropy term $\sum_k\hat p_k\log\hat p_k$, the maximisation over $q\in[\eta/2,1-\eta/2]$
replaced by a maximisation over probability vectors with entries in $[\eta/2,1]$ (the maximiser of the
linear-in-counts objective $\sum_kc_k\log q_k$ over the simplex is the empirical frequency vector, which
lies in this set on the good event) on a grid of spacing $n^{-4}$ in each coordinate, intersected with
the simplex, and Hoeffding's inequality applied coordinate-wise. With these changes the
proof of Proposition~\ref{thm:concentration} goes through verbatim and gives (iii).
\end{proof}
The count (hurdle) construction satisfies (i) and (ii). A separation statement of the form (iii) for it
requires the Poisson lemmas of Section~\ref{app:poisson} in place of the Bernoulli ones, and we do not
write it out. The theorem is elementary. It identifies the problems that a method restricted to $\Phi(Y)$ cannot
pose.

\paragraph{Experiments.} The following experiment tests the theorem on four designs in which the community
structure lives only in reciprocity, sign, intensity or the combination of two layers
(Tables~\ref{tab:R7} to~\ref{tab:R6b} in Section~\ref{app:tables}, eight paired seeds each), with the
same protocol as Section~\ref{sec:simulations} of the main text and the standards applied to the reduced
data together with signed spectral clustering, signed Leiden, spectral clustering of the reciprocated
subnetwork, the directed Bernoulli and Gamma--Poisson blocks, and the best single layer as an oracle
reference. In each design the
reduced data, which are the symmetrised network $A\vee A^\top$, the unsigned network $|W|$ of the signed
matrix $W$, the presence graph $\ind{Y>0}$ of the counts and the sum of the two layers, form an Erd\H{o}s--R\'enyi
network whatever the partition, and the Bethe Hessian and Leiden applied to the reduced data return ARI
$0.00$ in every cell, as they must. The comparison of interest is with full-data methods that are not
block models. The collapsed dyad-state evidence recovers the partition once the within-community
reciprocated share reaches $0.09$ ($0.97$, and about $1$ from $0.12$ on), and so does spectral clustering
of the reciprocated subnetwork $A\circ A^\top$, the entrywise product ($0.94$ and $1.00$), while the directed Bernoulli block
model, which treats the two arcs of a dyad as independent, recovers it only at the largest reciprocity
($0.75$ at $0.18$). The collapsed Dirichlet--multinomial evidence recovers the partition at sign purity
$0.8$ and $0.9$ ($0.94$ and $1.00$), and so does signed Leiden ($0.86$ and $1.00$), with signed spectral
clustering behind ($0.22$ and $0.77$), and at purity $0.7$ no method does ($\le0.08$). The collapsed
zero-inflated evidence recovers the partition once the within intensity is four times the between one
($0.82$ at $\lambda_{\mathrm{in}}=4$ and $0.99$ at $6$), whereas the flat Gamma--Poisson evidence on the
same counts is diluted by the $92\%$ of inactive dyads ($0.12$ and $0.28$) and Leiden, the Bethe
Hessian and Infomap on the count-weighted network reach at most $0.26$. The collapsed multiplex
evidence recovers the partition from $\mathrm{SNR}_{\mathrm{KS}}=1.5$ in the first layer ($0.88$, and
$0.94$ at $2$), while the best single layer, an oracle choice, reaches $0.41$ and $0.62$.
Theorem~\ref{thm:zero} therefore rules out reduced-data methods and not full-data heuristics. In two of
the four designs a purpose-built heuristic matches the collapsed evidence, and the value of the latter
is that it needs no design-specific heuristic and covers the intensity and layer cases, for which we
know of none. These results depend on the model-based initialisation of Section~\ref{sec:inference} of the main text,
whose candidate set for these data types contains the spectral clustering of each layer, of each
category-specific subnetwork (signed, dyad-state), of the reciprocated subnetwork for directed data and
of the weighted network for counts: from a random start the single-site chain stays in a one-community
state, and the count-valued designs need the weighted spectral candidate.

\subsection{Proof of Theorem~\ref{thm:split}}\label{app:split}
Write $N=\binom \nu 2$, $E=\sum_{i<j}y_{ij}$, $d_i=\sum_{j\ne i}y_{ij}$ (so $\sum_id_i=2E$), and let
$F$ be supported in $[\theta_-,\theta_+]\subset(0,\infty)$, so that all rates $\lambda_{ij}=\theta_i\theta_j\omega_0$
then lie in $[\lambda_-,\lambda_+]=[\theta_-^2\omega_0,\theta_+^2\omega_0]$. Both type evidences share the
Gamma--Poisson factor $P(E)$ of the total, and $\log P(E)-\log p_{\ER}=-\log\Gamma(E+1)+E\log N$ (the
multinomial placement of $E$ edges among $N$ dyads), so by \eqref{eq:dc},
\begin{equation}\label{eq:dcer}
\log p_{\DC}-\log p_{\ER}=\sum_i\log d_i!+E\log2+E\log N-\log(2E)!-\log\binom{\nu +2E-1}{2E}.
\end{equation}
(The identity was also checked numerically against the definitions, to machine precision, on random
multigraphs.)

\paragraph{(a), degenerate $F$.} Put $\bar d_E=2E/\nu $ (the realised mean degree, with $\bar d_E=\bar d(1+O_p(\nu ^{-1}))$
since $E\sim\Pois(N\omega_0)$). Conditionally on $E$, the $E$ edges are independent and uniform on the
$N$ dyads, so each $d_i$ is marginally $\mathrm{Bin}(E,2/\nu )$ with mean $\bar d_E$ and variance
$\bar d_E(1-2/\nu )$, and $\sum_id_i=2E$ exactly. Write $d_i=\bar d_E+\xi_i$, so that $\sum_i\xi_i=0$. By a
Chernoff bound and a union bound over the $\nu $ nodes, $|\xi_i|\le\bar d_E/2$ for all $i$ with probability
$1-\nu e^{-c\bar d_E}$, and on that event Stirling's formula
$\log d!=d\log d-d+\tfrac12\log(2\pi d)+O(1/d)$ expanded around $\bar d_E$ gives
\begin{align*}
\sum_i\log d_i!={}&\nu \bar d_E\log\bar d_E-\nu \bar d_E+\tfrac \nu 2\log(2\pi\bar d_E)+\sum_i\frac{\xi_i^2}{2\bar d_E}\\
&+\frac12\sum_i\log\Big(1+\frac{\xi_i}{\bar d_E}\Big)+O\Big(\sum_i\frac{|\xi_i|^3}{\bar d_E^2}+\frac \nu {\bar d_E}\Big),
\end{align*}
the linear terms cancelling because $\sum_i\xi_i=0$. Conditionally on $E$, $\E\sum_i\xi_i^2=\nu \bar d_E(1-2/\nu )$
and $\Var(\sum_i\xi_i^2)=O(\nu \bar d_E^2)$ (a moment computation for the multinomial allocation), so
$\sum_i\xi_i^2/(2\bar d_E)=\nu /2+O_p(\sqrt \nu )$. Similarly $\frac12\sum_i\log(1+\xi_i/\bar d_E)=
-\sum_i\xi_i^2/(4\bar d_E^2)+O_p(\nu \bar d_E^{-3/2})=O_p(\nu /\bar d_E)$ and
$\sum_i|\xi_i|^3/\bar d_E^2=O_p(\nu \bar d_E^{-1/2})$.
Hence
\[
\sum_i\log d_i!=\nu \bar d_E\log\bar d_E-\nu \bar d_E+\tfrac \nu 2\log(2\pi\bar d_E)+\tfrac \nu 2+O_p(\sqrt \nu +\nu /\bar d_E).
\]
Next, with $N=\nu (\nu -1)/2$ and $2E=\nu \bar d_E$, Stirling's formula gives
\begin{align*}
E\log2+E\log N-\log(2E)!&=\tfrac{\nu \bar d_E}2\log\big(\nu (\nu -1)\big)-\big[\nu \bar d_E\log(\nu \bar d_E)-\nu \bar d_E+O(\log \nu )\big]\\
&=-\nu \bar d_E\log\bar d_E+\nu \bar d_E-\tfrac{\bar d_E}2+O(\log \nu +\bar d_E/\nu ),
\end{align*}
because $\frac{\nu \bar d_E}2\log(1-1/\nu )=-\frac{\bar d_E}2+O(\bar d_E/\nu )$, and
\begin{align*}
\log\binom{\nu +2E-1}{2E}&=\log\Gamma(\nu (1+\bar d_E))-\log\Gamma(\nu \bar d_E+1)-\log\Gamma(\nu )\\
&=\nu \log(1+\bar d_E)+\nu \bar d_E\log(1+1/\bar d_E)+O(\log \nu ).
\end{align*}
Collecting the three displays in \eqref{eq:dcer},
\begin{align*}
\log p_{\DC}-\log p_{\ER}&=\tfrac \nu 2\log(2\pi\bar d_E)+\tfrac \nu 2-\tfrac{\bar d_E}2-\nu \log(1+\bar d_E)-\nu \bar d_E\log(1+1/\bar d_E)\\
&\qquad+O_p(\sqrt \nu +\nu /\bar d_E+\log \nu ).
\end{align*}
Since $\nu \log(1+\bar d_E)=\nu \log\bar d_E+\nu /\bar d_E+O(\nu /\bar d_E^2)$ and $\nu \bar d_E\log(1+1/\bar d_E)=\nu -\nu /(2\bar d_E)+O(\nu /\bar d_E^2)$,
\[
\log p_{\DC}-\log p_{\ER}=-\tfrac \nu 2\log\bar d_E+\nu \,\frac{\log(2\pi)-1}2-\frac{\bar d_E}2-\frac \nu {2\bar d_E}+O_p(\sqrt \nu +\nu /\bar d_E^2+\log \nu ).
\]
In the dense regime $\bar d_E=\bar d(1+O_p(\nu ^{-1}))$ with $\bar d=\omega_0(\nu -1)$, so replacing $\bar d_E$ by
$\bar d$ changes the right side by $O_p(1)$, $\nu /\bar d_E=O(1)$, and the display is the first statement
of (a) with $c_1=(\log(2\pi)-1)/2-\bar d/(2\nu )+O(1/\bar d)$. (The script \texttt{theory\_checks/check\_split\_a.py} of the released code, with $\omega_0\in\{0.05,0.2\}$
and $\nu \in\{200,800,3200\}$, reproduces the per-node value of the right side to within $0.01$ to $0.06$, the
residual being the $O_p(\sqrt \nu )$ term. With the degree prior dropped, the binomial term is absent and the
difference is $+\tfrac \nu 2\log\bar d+O(\nu )$, so without a degree prior the DC type would be
selected on every community, which is why the typed library uses the uniform prior.)

\paragraph{(a), non-degenerate $F$.} Condition on $\theta$ and put $T=\sum_j\theta_j=\nu +O_p(\sqrt \nu )$
and $\mu_i=\E[d_i\mid\theta]=\theta_i\omega_0(T-\theta_i)$. Conditionally on $\theta$, $d_i\sim\Pois(\mu_i)$
(a sum of independent Poissons), with $\mu_i\in[\theta_-\omega_0(T-\theta_+),\theta_+\omega_0T]$, and
$\E[E\mid\theta]=\bar E:=\omega_0\sum_{i<j}\theta_i\theta_j=\frac12\sum_i\mu_i$. We show that
\begin{equation}\label{eq:a2}
\log p_{\DC}-\log p_{\ER}=\nu \bar d\,\E_F[\theta\log\theta]+O_p\big(\nu \log\bar d+\bar d\sqrt \nu \big),
\end{equation}
which is the second statement of (a) because $\bar d\to\infty$. First, $\E[\log d!]=\mu\log\mu-\mu+\frac12\log(2\pi\mu)+\frac12+O(\mu^{-1})$
for $d\sim\Pois(\mu)$, and $\sum_i[\log d_i!-\E\log d_i!]=\sum_i\log\mu_i\,(d_i-\mu_i)+O_p(\sqrt \nu )$:
the second-order term of the expansion of $\log d_i!$ around $\mu_i$ is
$\sum_i[(d_i-\mu_i)^2-\mu_i]/(2\mu_i)$, which is centred and has variance $\nu /2+O(1)$ (two summands are
correlated only through the shared edge $y_{ij}$, with covariance $(2\lambda_{ij}^2+\lambda_{ij})/(4\mu_i\mu_j)=O(\nu ^{-2})$),
and the third-order terms are $O_p(\sum_i\mu_i^{-1/2})=O_p(\nu \bar d^{-1/2})=O_p(\sqrt \nu )$ in the dense
regime. The linear term has conditional variance $\sum_{i<j}\lambda_{ij}(\log\mu_i+\log\mu_j)^2=O(\nu \bar d\log^2\bar d)$
(write $d_i-\mu_i=\sum_{j\neq i}(y_{ij}-\lambda_{ij})$), hence
$\sum_i\log d_i!=\sum_i[\mu_i\log\mu_i-\mu_i]+\frac12\sum_i\log(2\pi\mu_i)+\frac \nu 2+O_p(\sqrt{\nu \bar d}\log\bar d)$.
Second, $\log\mu_i=\log\theta_i+\log(\omega_0T)+\log(1-\theta_i/T)$, so
$\sum_i\mu_i\log\mu_i=\omega_0\sum_i\theta_i(T-\theta_i)\log\theta_i+2\bar E\log(\omega_0T)+O(\bar d)
=\omega_0T\sum_i\theta_i\log\theta_i+2\bar E\log(\omega_0T)+O(\nu )$, and by the law of large numbers
$\sum_i\theta_i\log\theta_i=\nu \E_F[\theta\log\theta]+O_p(\sqrt \nu )$, so that
$\omega_0T\sum_i\theta_i\log\theta_i=\nu \bar d\,\E_F[\theta\log\theta]+O_p(\bar d\sqrt \nu )$. Third, the
remaining terms of \eqref{eq:dcer} are the same functions of the realised $E$ as in the degenerate case,
$E\log2+E\log N-\log(2E)!-\log\binom{\nu +2E-1}{2E}=-\nu \bar d_E\log\bar d_E+\nu \bar d_E-\nu \log(1+\bar d_E)-\nu \bar d_E\log(1+1/\bar d_E)+O(\bar d_E+\log \nu )$
with $\bar d_E=2E/\nu $, and they combine with $2\bar E\log(\omega_0T)-2\bar E$ as follows: writing
$\bar d_{\bar E}=2\bar E/\nu $, one has $\omega_0T/\bar d_{\bar E}=T\nu/(T^2-\sum_i\theta_i^2)=1+O_p(\nu ^{-1/2})$, so
$2\bar E\log(\omega_0T)=\nu \bar d_{\bar E}\log\bar d_{\bar E}+O_p(\bar d\sqrt \nu )$, and
$E-\bar E=O_p(\sqrt{\nu \bar d})$ conditionally on $\theta$, so
$\nu \bar d_E\log\bar d_E-\nu \bar d_{\bar E}\log\bar d_{\bar E}=O_p(\sqrt{\nu \bar d}\log\bar d)$ and
$\nu \bar d_E-2\bar E=O_p(\sqrt{\nu \bar d})$. All other terms ($\frac12\sum_i\log(2\pi\mu_i)$, $\frac \nu 2$,
$\nu \log(1+\bar d_E)$, $\nu \bar d_E\log(1+1/\bar d_E)$, $\bar d_E$) are $O(\nu \log\bar d)$. Adding up gives
\eqref{eq:a2}. Finally $\E_F[\theta\log\theta]\ge\E_F[\theta]\log\E_F[\theta]=0$ by Jensen's inequality
applied to the strictly convex function $x\mapsto x\log x$, with equality if and only if $F$ is
degenerate.

\paragraph{(b).} Let $z$ be a labelling of $C$ into two blocks with parts $S_1,S_2$ and blocks
$u\in\{11,22,12\}$. A block with no dyads (a singleton part) has evidence $1$ and profile $0$ and is
omitted from all sums. By Lemma~\ref{lem:poisson} applied to the three blocks, on the event
$\{\max_{i<j}y_{ij}\le y^\star\}$ of Lemma~\ref{lem:poisunif} (probability $\to1$),
$\log p_{\ER\text{-}2}(Y_C\mid z)\le\ell_2(z)+C\log \nu $ uniformly in $z$ (each block has $S_u/N_u\le y^\star$
and $N_u+S_u\le \nu ^2y^\star$). By Lemma~\ref{lem:poisunif}, for any $\tau>0$, with probability $\to1$,
$\ell_2(z)\le\bar\ell_2(z)+\tau \nu ^2$ uniformly in $z$, and by \eqref{eq:poisid},
$\bar\ell_2(z)=\sum_{i<j}h_{\Pois}(\lambda_{ij})-\sum_u\sum_{(i,j)\in u}\kl_{\Pois}(\lambda_{ij}\|\bar\lambda_u)$.
On the other side, inequality~\eqref{eq:dclower} of Lemma~\ref{lem:dcprofile} at the true $(\theta,\omega_0)$ gives
$\log p_{\DC}(Y_C)\ge\ell_{\DC}(Y_C;\theta,\omega_0)-O(\nu \log \nu )$ (here $\frac{\omega_0}2\sum_i\theta_i^2\le\frac{\omega_0}2\nu \theta_+^2$,
$\bar d=O(\nu )$ and $L_P=O(\log \nu )$ because $E/N\le\lambda_+(1+o_p(1))$), and
$\ell_{\DC}(Y_C;\theta,\omega_0)=\sum_{i<j}[y_{ij}\log\lambda_{ij}-\lambda_{ij}]=\sum_{i<j}h_{\Pois}(\lambda_{ij})+O_p(\nu )$,
the fluctuation being a sum of independent centred Poisson variables with bounded coefficients and total
variance $O(\nu ^2)$. Therefore, uniformly in $z$ and with probability $\to1$,
\begin{align}
\log p_{\DC}(Y_C)-\log p_{\ER\text{-}2}(Y_C\mid z)\ &\ge\ \sum_u\sum_{(i,j)\in u}\kl_{\Pois}(\lambda_{ij}\|\bar\lambda_u)-\tau \nu ^2-O_p(\nu \log \nu )\nonumber\\
&\ge\ \frac{\omega_0^2}{2\lambda_+}\sum_u\sum_{(i,j)\in u}(\theta_i\theta_j-c_u)^2-\tau \nu ^2-O_p(\nu \log \nu ),\label{eq:bmargin}
\end{align}
where $c_u$ is the block average of $\theta_i\theta_j$ and the last inequality is \eqref{eq:poisid}.
It remains to bound the double sum below by a constant times $\nu ^2$, uniformly in $z$. By hypothesis
the support of $F$ contains three points $q_1<q_2<q_3$. Choose disjoint closed
intervals $I_t\ni q_t$ ($t=1,2,3$) of positive $F$-mass $f_t$, separated by gaps of length at least
$g>0$, and let $f_{\min}=\min_tf_t$. By the law of large numbers, eventually almost surely each $I_t$
contains at least $f_t\nu/2$ of the $\theta_i$. Fix a labelling $z$. Each interval has a majority part
$s(t)\in\{1,2\}$ containing at least half of its nodes, and by the pigeonhole principle two intervals
$t\neq t'$ share a majority part $s$. Let $A=S_s\cap\{i:\theta_i\in I_t\}$ and
$B=S_s\cap\{i:\theta_i\in I_{t'}\}$, so that $|A|,|B|\ge f_{\min}\nu /4$. For a fixed $i\in A$, the values
$\theta_i\theta_j$ for $j\in A\setminus\{i\}$ lie below $\theta_i\max I_t$ and those for $j\in B$ lie above
$\theta_i\min I_{t'}$, or the other way round, and the two ranges are separated by a gap of length at
least $\theta_ig\ge\theta_-g$. For two groups of values separated by a gap $G$ and any constant $c$,
$\sum_{j\in A\setminus\{i\}}(v_j-c)^2+\sum_{j\in B}(v_j-c)^2\ge\frac{(|A|-1)|B|}{|A|+|B|-1}G^2\ge\frac{\min(|A|-1,|B|)}2G^2$.
Summing over $i\in A$ and noting that every unordered pair is counted at most twice,
\[
\sum_{(i,j)\in u_{ss}}(\theta_i\theta_j-c_{ss})^2\ \ge\ \frac12\,|A|\cdot\frac{\min(|A|-1,|B|)}2\,\theta_-^2g^2\ \ge\ \frac{f_{\min}^2\theta_-^2g^2}{64}\,\nu ^2-O(\nu )
\]
uniformly in $z$, where $u_{ss}$ is the diagonal block of the part $S_s$. Inserting this in \eqref{eq:bmargin} with
$\kappa_0=\omega_0^2\theta_-^2g^2f_{\min}^2/(128\lambda_+)$ and $\tau=\kappa_0/2$ gives
$\min_z[\log p_{\DC}-\log p_{\ER\text{-}2}(\cdot\mid z)]\ge\frac{\kappa_0}2\nu ^2-O_p(\nu \log \nu )\ge\frac{\kappa_0}4\nu ^2$
with probability $\to1$, which is (b) with $\kappa=\kappa_0/4$.
Table~\ref{tab:splitb} reports the finite-size behaviour of this bound: for $F$ uniform
on $\{0.5,1,1.5\}$ and $\omega_0=0.2$ the split by atom beats the degree-corrected evidence at $\nu\le600$, by a
margin that falls below $140$ nats at $\nu=600$, and loses at $\nu\ge1200$.
\begin{table}[tbp]\centering\small
\caption{Finite-size behaviour of Theorem~\ref{thm:split}(b) for $F$ uniform on $\{0.5,1,1.5\}$ and $\omega_0=0.2$, with Gamma$(1,1)$ priors: exact values of $\log p_{\DC}(Y_C)-\log p_{\ER\text{-}2}(Y_C\mid z)$ for the two by-atom splits, three replicates per community size $\nu$ (script \texttt{theory\_checks/check\_split\_b.py}). A positive value means that the degree-corrected type is preferred to the split.}
\label{tab:splitb}
\begin{tabular}{lrrrrrr}\toprule
$\nu$ & \multicolumn{3}{c}{split $\{\theta_i=0.5\}\mid\{\theta_i\ge1\}$} & \multicolumn{3}{c}{split $\{\theta_i\le1\}\mid\{\theta_i=1.5\}$}\\\midrule
150 & $-134$ & $-147$ & $-130$ & $-71$ & $-106$ & $-97$\\
300 & $-184$ & $-245$ & $-172$ & $-69$ & $-57$ & $-27$\\
600 & $-136$ & $-33$ & $-72$ & $809$ & $843$ & $822$\\
1200 & $1838$ & $1677$ & $1874$ & $5406$ & $4910$ & $5136$\\
2400 & $12103$ & $12612$ & $12927$ & $25451$ & $24958$ & $26057$\\
\bottomrule\end{tabular}\end{table}

\paragraph{(c).} Let the two communities have $\nu /2$ nodes each. Every node has the same expected degree
$\mu=(\tfrac \nu 2-1)\lambda_{\mathrm{in}}+\tfrac \nu 2\lambda_{\mathrm{out}}$, and
$\bar E:=\E E=\nu \mu/2$. \emph{Lower bound for the two-block evidence.} By Lemma~\ref{lem:poisson} (other
direction) and the definition of the profile, $\log p_{\ER\text{-}2}(Y_C\mid z^\star)\ge\ell_2(z^\star)-O(\log \nu )
\ge\sum_{i<j}[y_{ij}\log\lambda^\star_{ij}-\lambda^\star_{ij}]-O(\log \nu )=\sum_{i<j}h_{\Pois}(\lambda^\star_{ij})-O_p(\nu )$,
where $\lambda^\star_{ij}\in\{\lambda_{\mathrm{in}},\lambda_{\mathrm{out}}\}$ are the true rates and the last
step is the same variance computation as in (b). \emph{Upper bound for the DC evidence.} By
inequality~\eqref{eq:dcupper} of Lemma~\ref{lem:dcprofile} and the first display of its proof,
$\log p_{\DC}(Y_C)\le\sum_id_i\log d_i-E\log(2E)-E+O(\nu \log \nu )$, because $\max_id_i=O_p(\nu )$. Now
$d_i-\mu=O_p(\sqrt \nu )$ for each $i$, $2E-\nu \mu=\sum_i(d_i-\mu)=O_p(\nu )$ (its variance is
$4\Var(E)=4\bar E$), and $\sum_id_i\log d_i=\nu \mu\log\mu+(\log\mu+1)\sum_i(d_i-\mu)+O_p(\nu )=\nu \mu\log\mu+O_p(\nu \log \nu )$.
Likewise $E\log(2E)=\frac{\nu \mu}2\log(\nu \mu)+O_p(\nu \log \nu )$ and $E=\frac{\nu \mu}2+O_p(\nu )$. Hence
\[
\log p_{\DC}(Y_C)\le\frac{\nu \mu}2\Big[\log\frac\mu \nu -1\Big]+O_p(\nu \log \nu )=\sum_{i<j}\big[\lambda^\star_{ij}\log\bar\mu-\bar\mu\big]+O_p(\nu \log \nu ),
\]
where $\bar\mu=(\lambda_{\mathrm{in}}+\lambda_{\mathrm{out}})/2$, the last equality because $\sum_{i<j}\lambda^\star_{ij}=\bar E=\nu \mu/2$, $\binom \nu 2\bar\mu=\frac{\nu \mu}2+O(\nu )$
and $\mu/\nu =\bar\mu+O(1/\nu )$. (The right side is the expected log-likelihood of the constant-rate model
$\lambda_{ij}\equiv\bar\mu$, which by the concavity of $\E\ell_{\DC}$ in $(\log\theta,\log\omega)$ and the
symmetry of the truth is the best rank-one approximation, and the bound above reaches it without using
that fact.) \emph{Difference.} Subtracting,
\begin{align*}
\log p_{\ER\text{-}2}(Y_C\mid z^\star)-\log p_{\DC}(Y_C)\ &\ge\ \sum_{i<j}\kl_{\Pois}(\lambda^\star_{ij}\|\bar\mu)-O_p(\nu \log \nu )\\
&=2\binom{\nu /2}2\kl_{\Pois}(\lambda_{\mathrm{in}}\|\bar\mu)+\frac{\nu ^2}4\kl_{\Pois}(\lambda_{\mathrm{out}}\|\bar\mu)-O_p(\nu \log \nu ),
\end{align*}
and since $2\bar\mu=\lambda_{\mathrm{in}}+\lambda_{\mathrm{out}}$,
\[
\kl_{\Pois}(\lambda_{\mathrm{in}}\|\bar\mu)+\kl_{\Pois}(\lambda_{\mathrm{out}}\|\bar\mu)=\lambda_{\mathrm{in}}\log\frac{\lambda_{\mathrm{in}}}{\bar\mu}+\lambda_{\mathrm{out}}\log\frac{\lambda_{\mathrm{out}}}{\bar\mu}=J(\lambda_{\mathrm{in}},\lambda_{\mathrm{out}}),
\]
so the bracket equals $\frac{\nu ^2}4J(1+O(1/\nu ))$, which is (c). The bound $J(u,v)\ge(u-v)^2/(4\max(u,v))$
follows from $\kl_{\Pois}(\lambda\|\mu)\ge(\lambda-\mu)^2/(2\max(\lambda,\mu))$ applied to both terms with
$|\lambda-\bar\mu|=|u-v|/2$. \qed

\subsection{Proof of Proposition~\ref{thm:cap}, the over-splitting statement and shrinkage of fixed strength}\label{app:cap}
\begin{figure}[tbp]
\centering
\includegraphics[width=0.95\linewidth]{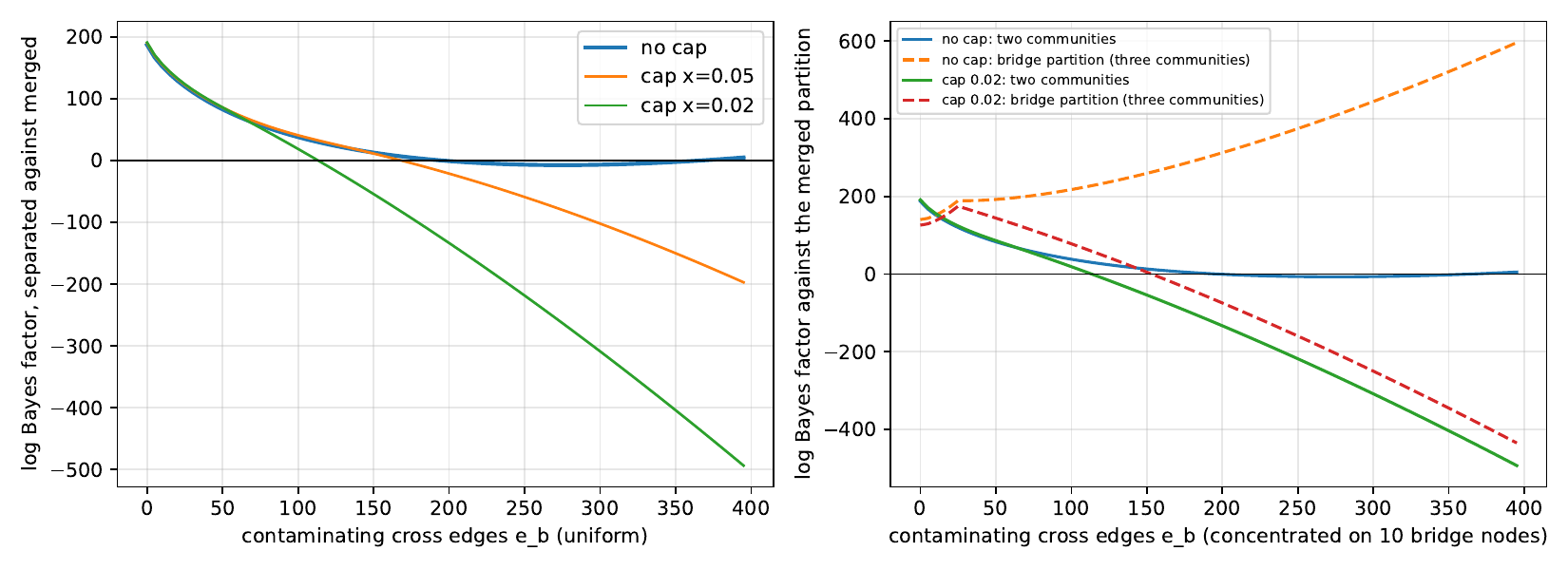}
\caption{Exact log Bayes factors under Proposition~\ref{thm:cap}. Left: log Bayes factor of the separated
against the merged partition for two 50-node communities with $p_{\mathrm{in}}=0.11$ as a function of
uniformly placed contaminating cross edges $e_b$, without a cap and with $x\in\{0.05,0.02\}$. A
positive value means that separation is preferred, and the cap lowers the number of tolerated cross
edges from 197 to 169 and 114. Right: contamination concentrated on ten bridge nodes. The three-community
bridge partition of Proposition~\ref{prop:bridge} beats the two-community partition from 15 cross
edges on with the cap and from 11 without it, and under the cap it is preferred to the merged
partition until 154 cross edges against 114 for the two-community partition.}
\label{fig:cap}
\end{figure}
Throughout, $\mathrm{BF}_x$ is a ratio of evidences (the partition prior is common to both sides and
cancels), $I_x(\alpha,\beta)=\Pbb(\mathrm{Beta}(\alpha,\beta)\le x)$, and
$\Delta_x(e,N)=\log I_x(a+e,b+N-e)-\log I_x(a,b)$. Only off-diagonal blocks carry a cap factor.

\emph{(i), first part.} By \eqref{eq:cap}, the evidence of an off-diagonal block with $e$ edges among $N$ dyads under
the truncated prior is its untruncated evidence times $I_x(a+e,b+N-e)/I_x(a,b)$. The merged partition has
no off-diagonal block and the separated partition has exactly one, with counts $(e_b,N_b)$. All other
factors coincide, so $\log\mathrm{BF}_x=\log\mathrm{BF}_1+\Delta_x(e_b,N_b)$. For fixed $N$ and $0\le e<N$
the ratio of the $\mathrm{Beta}(a+e+1,b+N-e-1)$ density to the $\mathrm{Beta}(a+e,b+N-e)$ density is
proportional to $p/(1-p)$, which is increasing on $(0,1)$. The family is therefore increasing in $e$ in
the likelihood-ratio order, hence in the usual stochastic order, and $I_x(a+e,b+N-e)$ is non-increasing
in $e$.

\emph{(i), bound.} $I_x\le1$ gives $\Delta_x(e,N)\le\log(1/I_x(a,b))$, and $I_x(1,1)=x$.

\emph{(ii)} For positive integers $\alpha,\beta$ the $\alpha$-th smallest of $\alpha+\beta-1$ independent
uniform variables has the $\mathrm{Beta}(\alpha,\beta)$ law, and it is at most $x$ if and only if at least
$\alpha$ of the uniforms are at most $x$. Hence $I_x(\alpha,\beta)=\Pbb(\mathrm{Bin}(\alpha+\beta-1,x)\ge\alpha)$.
With $\alpha=a+e$, $\beta=b+N-e$, $n'=a+b+N-1$ and $k=a+e$, the Chernoff--Hoeffding bound
$\Pbb(\mathrm{Bin}(n',x)\ge k)\le\exp\{-n'\kl(k/n'\,\|\,x)\}$ for $k/n'\ge x$ gives
$\log I_x(a+e,b+N-e)\le-n'\kl(\hat\mu\|x)$, and (ii) follows from the definition of $\Delta_x$. The
posterior mean of the block density is $(a+e)/(a+b+N)=\hat\mu\,n'/(n'+1)$, so
$0\le\hat\mu-(a+e)/(a+b+N)=\hat\mu/(n'+1)\le1/N$.

\emph{(iii)} By (i), $\log\mathrm{BF}_x(e)\le\log\mathrm{BF}_1(e)+\log(1/I_x(a,b))$ for every $e$.
If $j_x=0$ there is nothing to prove. Otherwise the slope hypothesis gives
$\log\mathrm{BF}_1(e^\star_1+j_x)\le\log\mathrm{BF}_1(e^\star_1)-j_x\zeta<-\log(1/I_x(a,b))$, because
$\log\mathrm{BF}_1(e^\star_1)<0$ and $j_x\zeta\ge\log(1/I_x(a,b))$. Hence
$\log\mathrm{BF}_x(e^\star_1+j_x)<0$ and $e^\star_x\le e^\star_1+j_x$. For the second claim,
$\log\mathrm{BF}_x(e_b)\le\log\mathrm{BF}_1(e_b)-n'\kl(\hat\mu_b\|x)+\log(1/I_x(a,b))<0$ by (i), (ii) and
the hypothesis, while $\log\mathrm{BF}_1(e_b)$ may be positive (in the example of
Section~\ref{sec:theory} it is positive for $114\le e_b\le196$).

\begin{proposition}[The cap and over-splitting]\label{prop:capsplit}
Let $\tilde z$ refine a partition $z$ with $K$ communities by splitting one community $C$ into two
parts of sizes $m_1$ and $n_C-m_1$, creating an off-diagonal block with $N'=m_1(n_C-m_1)$ dyads and
$e'$ edges, and let $\hat\mu'=(a+e')/(a+b+N'-1)$. Then
$\log\mathrm{BF}_x(\tilde z\!:\!z)=\log\mathrm{BF}_1(\tilde z\!:\!z)+\Delta_x(e',N')+R$, where $R$
collects the cap terms of the blocks joining $C$, or its two parts, to the other $K-1$ communities
and satisfies $R\le(K-1)\log(1/x)$ when $a=b=1$ and, for positive integers $a,b$,
$R\le(K-1)(2\log(1/I_x(a,b))+\log2)$ whenever every block joining $C$ to another community has
$a+e\le\lfloor(a+b+N-1)x\rfloor$. If $\hat\mu'>x$ then, by Proposition~\ref{thm:cap}(ii),
$\Delta_x(e',N')\le-(a+b+N'-1)\kl(\hat\mu'\|x)+\log(1/I_x(a,b))$. If the dyads inside $C$ are
independent Bernoulli$(p_C)$ with $p_C>x$ and $N'\to\infty$, then $\hat\mu'=p_C(1+o_p(1))$ and the cap
penalises the split at rate $N'\kl(p_C\|x)$, whereas without the cap the evidence cost of any such
split is at most $(3K+1)C_{a,b}\log(n^2+a+b+2)=O(K\log n)$ nats plus the partition prior, with
$C_{a,b}=3(a+b)+7$ the constant of Lemma~\ref{lem:stirling}.
\end{proposition}
\begin{proof}
\emph{Decomposition.} The partition $z$ has the off-diagonal blocks $(C,C')$, $C'\ne C$, and
$\tilde z$ has $(C_1,C_2)$, $(C_1,C')$ and $(C_2,C')$. Every other block is common to $z$ and $\tilde z$, and diagonal
blocks carry no cap factor. Hence
\begin{align*}
\log\mathrm{BF}_x(\tilde z\!:\!z)-\log\mathrm{BF}_1(\tilde z\!:\!z)&=\Delta_x(e',N')+R,\\
R&=\sum_{C'\neq C}\big[\Delta_x(e_{1C'},N_{1C'})+\Delta_x(e_{2C'},N_{2C'})\\
&\qquad\qquad-\Delta_x(e_{CC'},N_{CC'})\big],
\end{align*}
with $e_{CC'}=e_{1C'}+e_{2C'}$ and $N_{CC'}=N_{1C'}+N_{2C'}$.
\emph{Bound on $R$ for $a=b=1$.} Fix $C'$, write $(e_i,N_i)$ for the two sub-blocks and $(e,N)$ for their
union, and let $X_i\sim\mathrm{Bin}(N_i+1,x)$ be independent. By the identity of (c),
$I_x(1+e_i,1+N_i-e_i)=\Pbb(X_i\ge e_i+1)$, and $\{X_1\ge e_1+1,\,X_2\ge e_2+1\}\subset\{X_1+X_2\ge e+2\}$.
Now $X_1+X_2\sim\mathrm{Bin}(N+2,x)$ has the law of $X+\xi$ with $X\sim\mathrm{Bin}(N+1,x)$ and
$\xi\sim\mathrm{Bernoulli}(x)$ independent, and $\{X+\xi\ge e+2\}\subset\{X\ge e+1\}$. Therefore
$I_x(1+e_1,1+N_1-e_1)\,I_x(1+e_2,1+N_2-e_2)\le I_x(1+e,1+N-e)$, i.e.\ the summand of $R$ for $C'$ is at
most $-\log I_x(1,1)=\log(1/x)$, and $R\le(K-1)\log(1/x)$.
\emph{Bound on $R$ for positive integers $a,b$.} The two positive terms of each summand are at most
$\log(1/I_x(a,b))$ by (b). For the negative term, $-\Delta_x(e,N)\le-\log I_x(a+e,b+N-e)=
-\log\Pbb(\mathrm{Bin}(n',x)\ge a+e)$ with $n'=a+b+N-1$. Every median of $\mathrm{Bin}(n',x)$ lies in
$[\lfloor n'x\rfloor,\lceil n'x\rceil]$ \citep{Kaas1980}, so if $a+e\le\lfloor n'x\rfloor$ then
$\Pbb(\mathrm{Bin}(n',x)\ge a+e)\ge\Pbb(\mathrm{Bin}(n',x)\ge\text{median})\ge\tfrac12$ and
$-\Delta_x(e,N)\le\log2$. Under the stated condition on the blocks joining $C$ to the other communities
this gives $R\le(K-1)\big(2\log(1/I_x(a,b))+\log2\big)$.
\emph{Rate.} If $\hat\mu'>x$, (c) applied to the created block gives the stated bound on $\Delta_x(e',N')$.
If the dyads inside $C$ are independent Bernoulli$(p_C)$ with $p_C>x$, then $e'\sim\mathrm{Bin}(N',p_C)$,
$\hat\mu'\to p_C$ in probability as $N'\to\infty$, $\kl(\hat\mu'\|x)\to\kl(p_C\|x)>0$ by continuity, and
$(a+b+N'-1)\kl(\hat\mu'\|x)=N'\kl(p_C\|x)(1+o_p(1))$.
\emph{Uncapped cost.} The profile log-likelihood $\ell$ of Lemma~\ref{lem:stirling} does not decrease under
refinement: if a block with $N$ dyads and density $\hat p$ is split into sub-blocks with $(N_i,\hat p_i)$,
then $\hat p=\sum_i(N_i/N)\hat p_i$ and $\sum_iN_ih(\hat p_i)\ge Nh(\hat p)$ by convexity of $h$. The blocks
that differ between $z$ and $\tilde z$ are $1+(K-1)$ blocks of $z$ and $3+2(K-1)$ blocks of $\tilde z$, each within
$C_{a,b}\log(N_{rs}+a+b+2)\le C_{a,b}\log(n^2+a+b+2)$ of its profile term by Lemma~\ref{lem:stirling}
(this needs only $a,b\ge\frac12$). Hence
the uncapped log evidences satisfy $\log p(A\mid \tilde z)-\log p(A\mid z)\ge\ell(\tilde z)-\ell(z)-(3K+1)C_{a,b}\log(n^2+a+b+2)\ge-(3K+1)C_{a,b}\log(n^2+a+b+2)$,
deterministically. \qed

\end{proof}

\begin{proposition}[Soft shrinkage]\label{prop:soft}
Let the off-diagonal prior be Beta$(a_s,b_s)$ with $a_s=mx$ and $b_s=m(1-x)$, $a_s,b_s\ge1$, in place
of Beta$(1,1)$, and for a block with $e$ edges among $N$ dyads let
$\Delta_m(e,N)=\log\Bfun(a_s+e,b_s+N-e)-\log\Bfun(a_s,b_s)-\log\Bfun(1+e,1+N-e)$, so that the log Bayes
factor of the separated against the merged partition under shrinkage is $\log\mathrm{BF}_1$ plus the
sum of $\Delta_m$ over the off-diagonal blocks that the separation creates. Let $\pi_{a_s,b_s}$ be the
Beta$(a_s,b_s)$ density. For $1\le e\le N-1$ and $\hat p=e/N$,
\begin{align*}
\Delta_m(e,N)\ &\ge\ \log\pi_{a_s,b_s}(\hat p)-(m-2)\log\Big(1+\frac mN\Big),\\
\Delta_m(e,N)\ &\le\ \log\pi_{a_s,b_s}(\hat p)+(a_s-1)\log\Big(1+\frac{a_s}e\Big)+(b_s-1)\log\Big(1+\frac{b_s}{N-e}\Big).
\end{align*}
In particular the benefit of shrinkage for separation is at most
$\sup_p\log\pi_{a_s,b_s}(p)+(a_s-1)\log(1+a_s)+(b_s-1)\log(1+b_s)$ and its cost at a block of density
$\hat p$ is at least $\log\pi_{a_s,b_s}(\hat p)-(m-2)\log(1+m/N)$, both bounded uniformly in $N$,
against the cost $-n'\kl(\hat\mu\|x)$ of Proposition~\ref{thm:cap}(ii), which grows linearly in $N$.
\end{proposition}
The proof applies the gamma-ratio bound \eqref{eq:gammaratio} to the three ratios
$\Gamma(a_s+e)/\Gamma(1+e)$, $\Gamma(b_s+N-e)/\Gamma(1+N-e)$ and $\Gamma(a_s+b_s+N)/\Gamma(N+2)$, and
is given in Section~\ref{app:props}. For the prior of Section~\ref{sec:cap}, $m=50$ and $x=0.02$,
$\Delta_m$ is within $0.16$ nats of $\log\pi_{a_s,b_s}(\hat p)$ for $N\ge1600$ and $\hat p\le0.2$. The
exact maximum of $\Delta_m$ over $e$ is $3.83$ nats at $N=1600$, against $\log(1/x)=3.91$ for the hard
cap, so the two devices offer the same benefit for separation to one decimal, and the uniform bound of
the proposition, $192$ nats for this prior, is loose. The cost at a between density of $0.04$, $0.08$
and $0.20$ is $+1.9$, $-0.1$ and $-6.8$ nats whatever the block size, where the truncation costs $-12$,
$-85$ and $-479$ nats at $N=1600$ (script \texttt{theory\_checks/check\_soft\_bound.py}).

\subsection{Proof of Theorem~\ref{thm:resolution}}\label{app:resolution}
\begin{figure}[tbp]
\centering
\includegraphics[width=\linewidth]{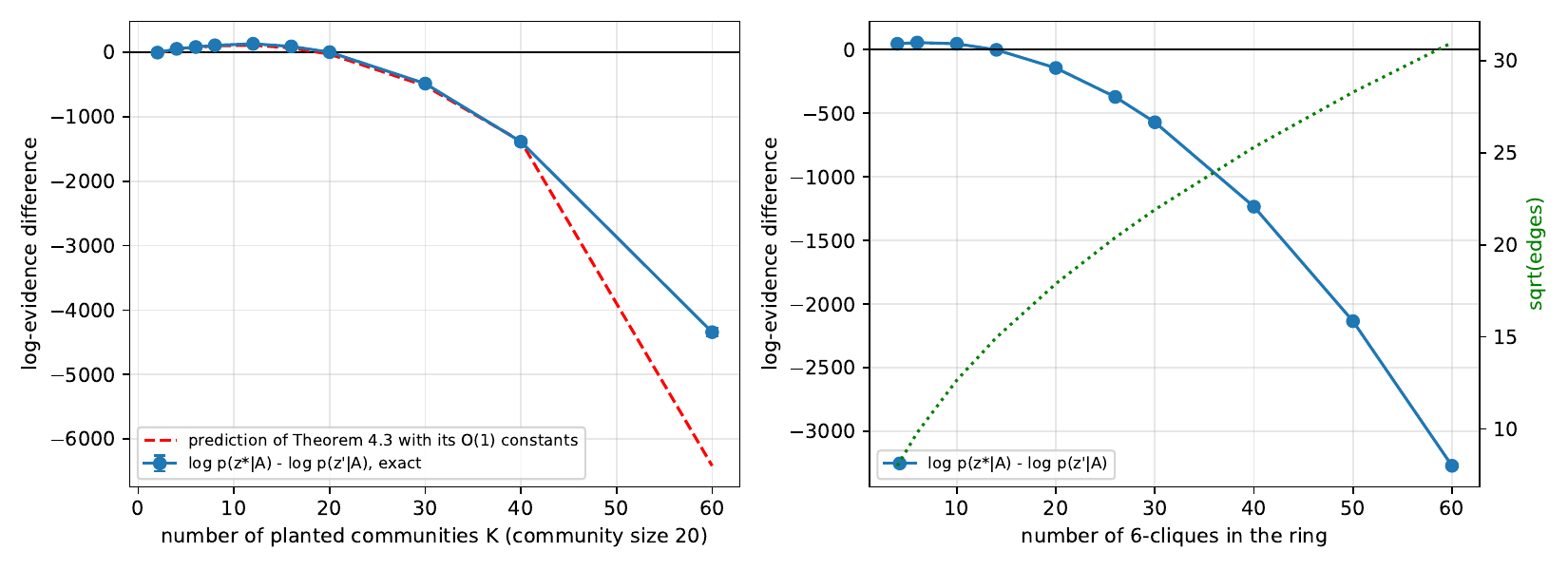}
\caption{Resolution limit under Theorem~\ref{thm:resolution}. Left: exact evidence difference between
the planted $K$-community partition, with communities of 20 nodes, within degree 6 and total between
degree 2, and its pairwise merger, over ten networks per $K$, with the prediction of the theorem
including its $O(1)$ constants. Right: rings of 6-cliques, where the planted partition loses to the
pairwise merger beyond about 14 cliques and the curve labelled $\sqrt{\mathrm{edges}}$ is drawn for reference. A positive
value means that the planted partition is preferred.}
\label{fig:resolution}
\end{figure}
Within this proof, $M$ denotes the edge count of a block, a random variable, and $\lambda$ its expectation, so that
$\lambda_{\mathrm{in}}$ and $\lambda_b$ are expected counts and not the Poisson rates of Theorem~\ref{thm:split}(c).
Write $p_{\mathrm{in}}=d_{\mathrm{in}}/(\nu -1)$ and $p_b=d_b/\nu $ for the within- and between-block edge
probabilities, so that every node has expected degree $d_{\mathrm{in}}$ inside its block and $d_b$ to each
other block. Here $N_{\mathrm{in}}=\binom \nu 2$, $\lambda_{\mathrm{in}}=N_{\mathrm{in}}p_{\mathrm{in}}=\nu d_{\mathrm{in}}/2$,
$N_b=\nu ^2$, $\lambda_b=\nu d_b$. For a block with $N$ dyads and $M$ edges put
$f_N(M)=\log\Bfun(a+M,b+N-M)-\log\Bfun(a,b)$, so that $\log p(A\mid z)=\sum_{r\le s}f_{N_{rs}(z)}(M_{rs}(z))$,
and extend $f_N$ to real $M\in[0,N]$ by the same formula. All $O(\cdot)$ constants below depend only on
$(a,b,d_{\mathrm{in}},d_b)$.

\emph{Step 1: one block.} Let $M$ be a sum of $N$ independent Bernoulli variables with success
probabilities $q_i=O(1/\nu )$, mean $\lambda=\sum_iq_i=\Theta(\nu )$ and variance $\sigma^2=\sum_iq_i(1-q_i)
=\lambda(1+O(1/\nu ))$, where $N=\Theta(\nu ^2)$. Then
\begin{equation}\label{eq:blockexp}
\E f_N(M)=\lambda\log\frac\lambda N-\lambda+\frac{\lambda^2}{2N}+\Big(a-\frac12\Big)\log\lambda-a\log N
+\frac12\log(2\pi)+\frac12-\log\Bfun(a,b)+O(\nu ^{-1/2}) ,
\end{equation}
and $(a-\frac12)\log\lambda-a\log N=\frac12\log\lambda-\log N+(a-1)\log(\lambda/N)=-(a+\frac12)\log \nu +O(1)$.
To see this, let $G=\{\lambda/2\le M\le3\lambda/2\}$. The multiplicative Chernoff bound gives
$\Pbb(G^c)\le2e^{-\lambda/12}$, and on $G^c$, $|f_N(M)|\le N\log2+C_{a,b}\log(N+a+b+2)$ by
Lemma~\ref{lem:stirling}, so $G^c$ contributes $O(\nu ^2e^{-c\nu})$ to every expectation below. On $G$ one has
$M\ge1$ and $N-M\ge N/2$ for $\nu $ large, and Stirling's formula together with \eqref{eq:gammaratio} gives
$\log\Gamma(a+M)=M\log M-M+(a-\frac12)\log M+\frac12\log(2\pi)+O(1/M)$,
$\log\Gamma(b+N-M)=(N-M)\log(N-M)-(N-M)+(b-\frac12)\log(N-M)+\frac12\log(2\pi)+O(1/N)$ and
$\log\Gamma(a+b+N)=N\log N-N+(a+b-\frac12)\log N+\frac12\log(2\pi)+O(1/N)$. Expanding
$(N-M)\log(N-M)=(N-M)\log N-M+M^2/(2N)+O(M^3/N^2)$ and $(b-\frac12)\log(N-M)=(b-\frac12)\log N+O(M/N)$,
\[
f_N(M)=M\log\frac MN-M+\frac{M^2}{2N}+\Big(a-\frac12\Big)\log M-a\log N+\frac12\log(2\pi)-\log\Bfun(a,b)+O(1/\nu )
\]
on $G$.
Now take expectations. For $g(t)=t\log t$, Taylor's formula around $\lambda$ with $|g'''|\le4/\lambda^2$ on
$G$ gives $\E[g(M)\ind G]=\lambda\log\lambda+\sigma^2/(2\lambda)+O(\lambda^{-2}\E|M-\lambda|^3)+O(\nu ^2e^{-c\nu})$,
and $\E|M-\lambda|^3\le(\E(M-\lambda)^4)^{3/4}=O(\lambda^{3/2})$ because
$\E(M-\lambda)^4\le3\sigma^4+\sum_i\E(\xi_i-q_i)^4=O(\lambda^2)$. Hence
$\E[M\log M]=\lambda\log\lambda+\frac12+O(\nu ^{-1/2})$. Likewise $\E[\log M\,\ind G]=\log\lambda+O(1/\lambda)$,
$\E M^2=\lambda^2+\sigma^2$ and $\E[M\ind G]=\lambda+O(\nu ^2e^{-c\nu})$. Collecting terms gives
\eqref{eq:blockexp}.

\emph{Step 2: expectation of the difference.} Under the planted partition every block count of $z^\star$
and of $z'$ is a sum of independent Bernoulli variables as in Step~1: $z^\star$ has $K$ diagonal blocks
$(N_{\mathrm{in}},\lambda_{\mathrm{in}})$ and $\binom K2$ off-diagonal blocks $(N_b,\lambda_b)$, and $z'$ has
$K/2$ diagonal blocks with $N'=\binom{2\nu }2=2N_{\mathrm{in}}+N_b$ dyads and mean count
$\lambda'=2\lambda_{\mathrm{in}}+\lambda_b$, and $\binom{K/2}2$ off-diagonal blocks with $4N_b$ dyads and
mean $4\lambda_b$. Apply \eqref{eq:blockexp} to every block and subtract. The terms $-\lambda$ cancel
(both partitions cover all $\binom n2$ dyads), the terms $\lambda^2/(2N)$ are $O(1)$ per block, and for a
merged off-diagonal block $4\lambda_b\log\frac{4\lambda_b}{4N_b}$ equals the sum of the four constituent
terms $\lambda_b\log\frac{\lambda_b}{N_b}$. For each of the $K/2$ merged pairs,
\begin{align*}
&2\lambda_{\mathrm{in}}\log\frac{\lambda_{\mathrm{in}}}{N_{\mathrm{in}}}+\lambda_b\log\frac{\lambda_b}{N_b}
-\lambda'\log\frac{\lambda'}{N'}\\
&\qquad=\nu d_{\mathrm{in}}\log\frac{2d_{\mathrm{in}}}{d_{\mathrm{in}}+d_b}+\nu d_b\log\frac{2d_b}{d_{\mathrm{in}}+d_b}
+\nu d_{\mathrm{in}}\log\frac{2\nu -1}{2\nu -2}+\nu d_b\log\frac{2\nu -1}{2\nu }\\
&\qquad=\nu J(d_{\mathrm{in}},d_b)+\frac{d_{\mathrm{in}}-d_b}2+O(\nu ^{-1}),
\end{align*}
so these terms contribute $\frac n2J+O(K)$. Every block of either partition contributes
$-(a+\frac12)\log \nu +O(1)$ through the logarithmic terms of \eqref{eq:blockexp}, and $z^\star$ has
$\frac{K(K+1)}2-\frac{K(K+2)}8=\frac{3K^2+2K}8$ more blocks than $z'$. Hence
\[
\E\big[\log p(A\mid z^\star)-\log p(A\mid z')\big]=\frac n2J-\frac{3K^2+2K}8\Big(a+\frac12\Big)\log \nu +O(K^2).
\]
The constants used for the prediction in Figure~\ref{fig:resolution} are exactly the $O(1)$ terms of
\eqref{eq:blockexp} with $a=b=1$ (with $\log(N+1)$ in place of $\log N$ and, for a block expected to be
empty, the exact value $-\log(N+1)$ of $f_N(0)$).

\emph{Step 3: the partition prior.} Under the CRP prior with concentration $\alpha$ the prior
is $p(z)\propto\alpha^{K(z)}\prod_k\Gamma(n_k)$, so by Stirling's formula
\begin{align*}
\log p(z^\star)-\log p(z')&=\frac K2\big[\log\alpha+2\log\Gamma(\nu )-\log\Gamma(2\nu )\big]\\
&=\frac K2\big[-2\nu \log2-\tfrac12\log \nu +O(1)\big]=-n\log2+O(K\log \nu ).
\end{align*}
Under the
Dirichlet--multinomial prior with $K_{\max}\ge K$ labels and concentration $\alpha$, $c=\alpha/K_{\max}$,
$p(z)\propto\prod_k\Gamma(c+n_k)/\Gamma(c)$ over the occupied labels, and
$\log\Gamma(c+\nu )=\log\Gamma(\nu )+c\log \nu +O(1/\nu )$ by \eqref{eq:gammaratio}, so the same computation gives
$-n\log2+O(K\log \nu )$.

\emph{Step 4: concentration.} $D:=\log p(A\mid z^\star)-\log p(A\mid z')$ is a function of the
$K(K+1)/2$ independent counts $M=(M_{rs})_{r\le s}$ of $z^\star$, each block of $z'$ being a sum of some
of them. By the Efron--Stein inequality, $\Var(D)\le\frac12\sum_{r\le s}\E[(D(M)-D(M^{(rs)}))^2]$, where
$M^{(rs)}$ replaces $M_{rs}$ by an independent copy $M'_{rs}$. Changing $M_{rs}$ changes
$f_{N_{rs}}(M_{rs})$ and the single block of $z'$ that contains it, with count
$\widetilde M=M_{rs}+(\text{other constituents})$ and $\widetilde N$ dyads, so
$\partial D/\partial M_{rs}=[\psi(a+M_{rs})-\psi(b+N_{rs}-M_{rs})]-[\psi(a+\widetilde M)-\psi(b+\widetilde N-\widetilde M)]$.
Let $G_{rs}$ be the event that $M'_{rs}$ and all constituents of $\widetilde M$ lie within a factor
$[\frac12,\frac32]$ of their means, so that $\Pbb(G_{rs}^c)\le10e^{-c\nu}$ by the Chernoff bound. On $G_{rs}$, every
point of the segment between $M$ and $M^{(rs)}$ satisfies the same bounds, and $\psi(t)=\log t+O(1/t)$
shows that $\partial D/\partial M_{rs}$ is the difference of the logarithms of two densities, both in
$[c_1/\nu ,c_2/\nu ]$, up to $O(1/\nu )$: $|\partial D/\partial M_{rs}|\le L$ for a constant $L$. Hence
$\E[(D-D^{(rs)})^2\ind{G_{rs}}]\le L^2\E(M_{rs}-M'_{rs})^2=2L^2\sigma_{rs}^2=O(\nu )$, while on $G_{rs}^c$,
$|D-D^{(rs)}|\le4\sup_M|f_{4N_b}(M)|=O(\nu ^2)$ by Lemma~\ref{lem:stirling}, so
$\E[(D-D^{(rs)})^2\ind{G_{rs}^c}]=O(\nu ^4e^{-c\nu})$. Summing, $\Var(D)=O(K^2\nu )$, and Chebyshev's inequality
gives $\Pbb(|D-\E D|\ge\varepsilon n)\le O(1/(\varepsilon^2\nu ))\to0$ for every $\varepsilon>0$, i.e.\
$D=\E D+o_p(n)$, with no restriction on $K$ beyond $\nu \to\infty$. Combining Steps~2--4 proves the expansion.
The statements about $K_{\mathrm{res}}$ follow by solving $\frac n2(J-2\log2)=\frac{3K^2+2K}8(a+\frac12)\log \nu $
for $K$ (to leading order $\frac{3K^2}{16}(2a+1)\log \nu $), and $J(d_{\mathrm{in}},d_b)\to d_{\mathrm{in}}\log2$
as $d_b\to0$ is immediate from the definition of $J$. At fixed $d_{\mathrm{out}}$ one has
$\lambda_b=\nu d_{\mathrm{out}}/(K-1)$: Step~1 requires $\lambda_b/\log \nu \to\infty$ (so that the bad-event
terms stay negligible), the constant $L$ of Step~4 becomes $O(\log K)$ because the diagonal and
off-diagonal densities now differ by a factor of order $K$, so Chebyshev's inequality needs
$\log K=o(\sqrt \nu )$ (implied by $\lambda_b/\log \nu \to\infty$, which forces $K=o(\nu /\log \nu )$), and the term $(a-\frac12)\log\lambda_b$ of each of the $(3K^2-2K)/8$ off-diagonal
blocks removed by the merger contributes $-(a-\frac12)\frac{3K^2-2K}8\log(K-1)$ to $\E D$. The remainder
is then $O(K^2\log K+K\log \nu )$. \qed

\subsection{Proof of Theorem~\ref{thm:threshold}}\label{app:threshold}
\begin{table}[tbp]\centering\small
\caption{Exact evidence threshold of Theorem~\ref{thm:threshold} in Kesten--Stigum units, the value of $\mathrm{SNR}_{\mathrm{KS}}$ at which $J_K=2\log K$, by number of communities $K$ and mean degree $d$, with the weak-signal limit $4\log K/(K-1)$. The Kesten--Stigum threshold is $1$. Values above $1$ mean that the evidence threshold lies above the Kesten--Stigum threshold.}
\label{tab:iclthr}
\begin{tabular}{rccccccc}\toprule
$K$ & $d=10$ & $d=20$ & $d=50$ & $d=200$ & $d=1000$ & $d=10^4$ & limit\\\midrule
2 & 2.64 & 2.71 & 2.75 & 2.77 & 2.77 & 2.77 & 2.77\\
3 & 2.38 & 2.36 & 2.32 & 2.27 & 2.23 & 2.21 & 2.20\\
4 & 2.21 & 2.14 & 2.05 & 1.96 & 1.90 & 1.86 & 1.85\\
5 & 2.08 & 1.98 & 1.86 & 1.75 & 1.67 & 1.63 & 1.61\\
8 & 1.85 & 1.69 & 1.52 & 1.36 & 1.27 & 1.21 & 1.19\\
10 & 1.75 & 1.57 & 1.39 & 1.21 & 1.11 & 1.05 & 1.02\\
11 & 1.71 & 1.52 & 1.33 & 1.15 & 1.05 & 0.99 & 0.96\\
12 & 1.68 & 1.48 & 1.28 & 1.10 & 0.99 & 0.93 & 0.90\\
15 & 1.59 & 1.38 & 1.17 & 0.98 & 0.87 & 0.80 & 0.77\\
20 & 1.49 & 1.26 & 1.04 & 0.84 & 0.72 & 0.66 & 0.63\\
\bottomrule\end{tabular}\end{table}

\begin{figure}[tbp]
\centering
\includegraphics[width=\linewidth]{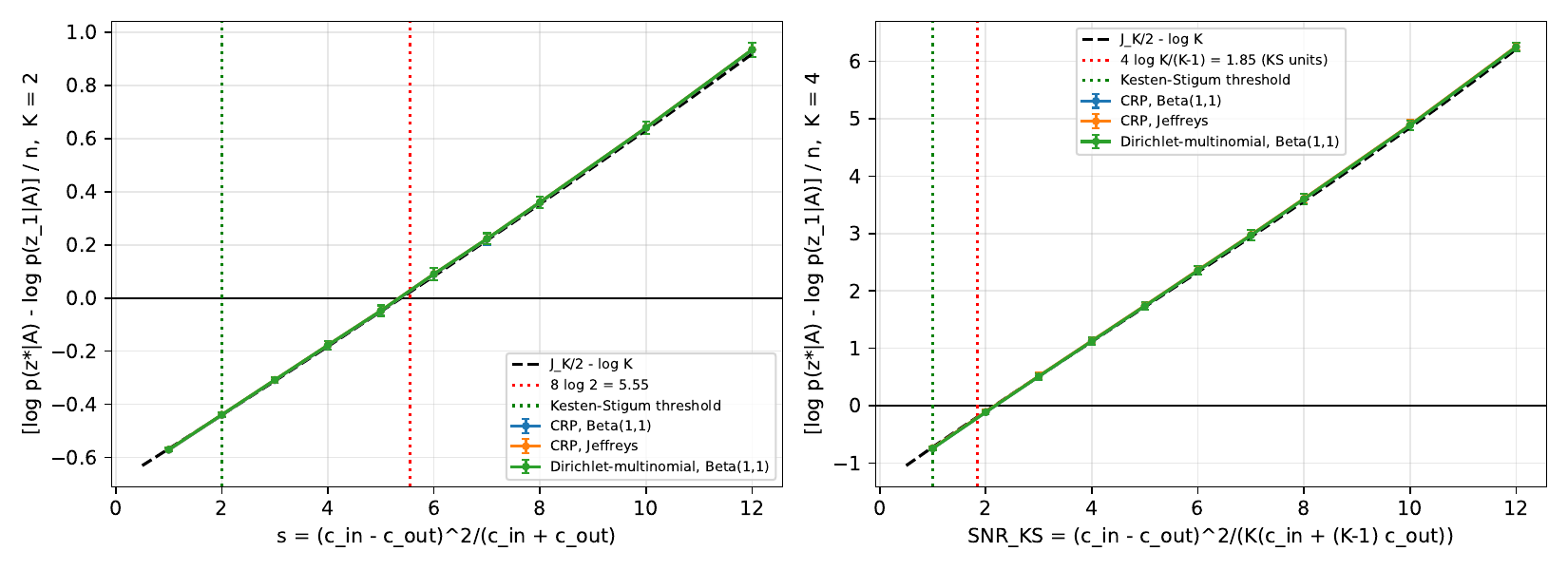}
\caption{Exact per-node difference between the collapsed log posterior of the planted partition and
of the one-community partition under Theorem~\ref{thm:threshold}, for $n=2000$, mean degree 16 and
eight networks per point, with CRP and Dirichlet--multinomial partition priors and Beta$(1,1)$ or
Jeffreys block priors, whose curves coincide. The horizontal axis is
$s=(c_{\mathrm{in}}-c_{\mathrm{out}})^2/(c_{\mathrm{in}}+c_{\mathrm{out}})$ for $K=2$ and
$\mathrm{SNR}_{\mathrm{KS}}$ for $K=4$. The dashed line is $J_K/2-\log K$, and the dotted lines mark
the Kesten--Stigum threshold and the large-degree limit of the evidence threshold, $8\log2$ for $K=2$ and $4\log4/3$
for $K=4$. A positive value means that the planted partition is preferred.}
\label{fig:threshold}
\end{figure}
\emph{Prior.} Under the CRP prior,
$\log p(z^\star)-\log p(z_1)=(K-1)\log\alpha+K\log\Gamma(n/K)-\log\Gamma(n)=-n\log K+O(K\log n)$ by
Stirling's formula. Under the Dirichlet--multinomial prior with $K_{\max}\ge K$ labels and $c=\alpha/K_{\max}$,
$\log p(z^\star)-\log p(z_1)=K\log\Gamma(c+n/K)-\log\Gamma(c+n)-(K-1)\log\Gamma(c)=-n\log K+O(K\log n)$ by
Stirling's formula and \eqref{eq:gammaratio}.

\emph{Evidence.} By Lemma~\ref{lem:stirling}, which needs $a,b\ge\frac12$,
\[
|\log p(A\mid z)-\ell(z)|\le\tfrac{K(K+1)}2C_{a,b}\log(n^2+a+b+2)=O(K^2\log n)
\]
for $z=z^\star$ and for $z=z_1$. For a block with $e\le N/2$ edges,
$Nh(e/N)=e\log(e/N)-e+e^2/(2N)+r$ with $0\le r\le\frac23e^3/N^2$ (write $u=e/N$, so that
$r=N\varphi(u)$ with $\varphi(u)=(1-u)\log(1-u)+u-u^2/2$, $\varphi(0)=\varphi'(0)=\varphi''(0)=0$ and
$0\le\varphi'''\le4$ on $[0,\frac12]$). Each block count $e_{rs}$ of $z^\star$ is Binomial with mean
$\lambda_{rs}=N_{rs}p_{rs}\asymp n$ and $N_{rs}\asymp n^2$, so $e_{rs}=\lambda_{rs}+O_p(\sqrt n)$,
$e_{rs}^2/N_{rs}=O_p(1)$ and $r=O_p(1/n)$. The same holds for the single block of $z_1$, whose count is
$\sum_{r\le s}e_{rs}$. The function $t\mapsto t\log(t/N)$ has derivative $\log(t/N)+1=O(\log n)$ on
$[\lambda_{rs}/2,2\lambda_{rs}]$, hence $e_{rs}\log(e_{rs}/N_{rs})=\lambda_{rs}\log(\lambda_{rs}/N_{rs})
+O_p(\sqrt n\log n)$, and the terms $-e_{rs}$ have the same sum in both partitions. Therefore
\[
\ell(z^\star)-\ell(z_1)=\sum_{r\le s}\lambda_{rs}\log\frac{\lambda_{rs}}{N_{rs}}-\lambda\log\frac\lambda N
+O_p(K^2\sqrt n\log n)=\sum_{r\le s}\lambda_{rs}\log\frac{p_{rs}}{\bar p}+o_p(n),
\]
where $\bar p=\lambda/N$, $\lambda=\sum_{r\le s}\lambda_{rs}$ and $N=\binom n2$.
With $K$ equal blocks of size $n/K$, $p_{rr}=c_{\mathrm{in}}/n$ on $K$ diagonal blocks of $\binom{n/K}2$ dyads and
$p_{rs}=c_{\mathrm{out}}/n$ on $\binom K2$ off-diagonal blocks of $n^2/K^2$ dyads, so
$\bar p=[c_{\mathrm{in}}(1-K/n)+(K-1)c_{\mathrm{out}}]/(K(n-1))=(d/n)(1+O(1/n))$ and
\[
\sum_{r\le s}\lambda_{rs}\log\frac{p_{rs}}{\bar p}=\frac n{2K}\Big(1-\frac Kn\Big)c_{\mathrm{in}}\log\frac{c_{\mathrm{in}}}d
+\frac{(K-1)n}{2K}\,c_{\mathrm{out}}\log\frac{c_{\mathrm{out}}}d+O(1)=\frac n2J_K+O(1).
\]
Adding the prior term and dividing by $n$ gives $J_K/2-\log K+o_p(1)$. Since $\frac1K\cdot c_{\mathrm{in}}+\frac{K-1}K\cdot c_{\mathrm{out}}=d$,
$J_K$ is the Jensen gap of $t\mapsto t\log t$ at the weights $(\frac1K,\frac{K-1}K)$ and is non-negative,
vanishing only at $c_{\mathrm{in}}=c_{\mathrm{out}}$.

\emph{Weak-signal expansion.} Write $c_{\mathrm{in}}=d+K(K-1)\delta$ and $c_{\mathrm{out}}=d-K\delta$, which leaves $d$ unchanged
and gives $c_{\mathrm{in}}-c_{\mathrm{out}}=K^2\delta$, and put $u=c_{\mathrm{in}}/d-1=K(K-1)\delta/d$
and $v=c_{\mathrm{out}}/d-1=-K\delta/d$. Then
$J_K=\frac dK\big[(1+u)\log(1+u)+(K-1)(1+v)\log(1+v)\big]$, and $(1+t)\log(1+t)=t+t^2/2+O(t^3)$ gives
$J_K=\frac dK\big[u+(K-1)v+\tfrac12u^2+\tfrac{K-1}2v^2+O(u^3+Kv^3)\big]$. The linear terms cancel,
$\frac12u^2+\frac{K-1}2v^2=\frac{K^3(K-1)\delta^2}{2d^2}$, and so
$J_K=\frac{K^2(K-1)\delta^2}{2d}\big(1+O(u)\big)=\frac{(K-1)(c_{\mathrm{in}}-c_{\mathrm{out}})^2}{2K^2d}\big(1+O((c_{\mathrm{in}}-c_{\mathrm{out}})/d)\big)
=\frac{K-1}2\,\mathrm{SNR}_{\mathrm{KS}}\big(1+O((c_{\mathrm{in}}-c_{\mathrm{out}})/d)\big)$, using
$\mathrm{SNR}_{\mathrm{KS}}=(c_{\mathrm{in}}-c_{\mathrm{out}})^2/(K(c_{\mathrm{in}}+(K-1)c_{\mathrm{out}}))=(c_{\mathrm{in}}-c_{\mathrm{out}})^2/(K^2d)$.
The large-degree criterion $\mathrm{SNR}_{\mathrm{KS}}>4\log K/(K-1)$ is therefore accurate only when
$(c_{\mathrm{in}}-c_{\mathrm{out}})/d$ is small. The finite-$d$
values quoted in the theorem solve $J_K=2\log K$ exactly at the stated mean degree: for $K=2$, in the units
$(c_{\mathrm{in}}-c_{\mathrm{out}})^2/(c_{\mathrm{in}}+c_{\mathrm{out}})$, $5.01$ ($d=5$), $5.38$ ($d=16$), $5.42$ ($d=20$), $5.54$ ($d=1000$), limit
$8\log2=5.545$, and for $K=4$, in Kesten--Stigum units, $2.24$ ($d=5$), $2.16$ ($d=16$), $2.05$ ($d=50$), $1.87$
($d=5000$), limit $4\log4/3=1.848$. The hyperparameters $(a,b)$ enter only through the
$O(K^2\log n)$ term of Lemma~\ref{lem:stirling}. \qed

\subsection{Proofs of Theorems~\ref{thm:estimator} and~\ref{thm:two} and of Corollary~\ref{cor:icl}(i)}\label{app:estimator}

\emph{Corollary~\ref{cor:icl}(i).} With $c_{\mathrm{in}}/d=1+(K-1)\lambda$ and $c_{\mathrm{out}}/d=1-\lambda$,
\[
J_K=\frac dK\Big[(1+(K-1)\lambda)\log(1+(K-1)\lambda)+(K-1)(1-\lambda)\log(1-\lambda)\Big],
\]
and $J_K>2\log K$ is $d>2K\log K/[\cdots]$, the displayed bound, which is $d^{\mathrm{upper}}_{\mathrm c}$
of \citet{Banks2016}, Theorem~1, in their parametrisation $c_{\mathrm{in}}=d(1+(K-1)\lambda)$,
$c_{\mathrm{out}}=d(1-\lambda)$. The script \texttt{theory\_checks/check\_banks\_identity.py} of the
released code evaluates both sides at the roots of Table~\ref{tab:iclthr}.

\medskip
\noindent\emph{Notation for Theorem~\ref{thm:estimator}.} For a labelling $z$ with at most
$K_{\max}$ labels, $\omega=\omega(z)$, $\pi=\pi(z)$ and the blocks $(b,b')$, $b\le b'$, are as in
Section~\ref{sec:thm-threshold}. Let $m_{bb'}(z)$ be the number of edges of $A$ inside block $(b,b')$,
$N_{bb'}(z)$ the number of dyads there, $m=\sum_{b\le b'}m_{bb'}(z)$ the total number of edges, which
does not depend on $z$, and $x=x(z)=m(z)/n$ the vector of per-node counts. Write
$\tilde\ell(z)=\sum_{b\le b'}[m_{bb'}\log(m_{bb'}/N_{bb'})-m_{bb'}]$ for the Poisson form of the
profile log-likelihood.

\begin{lemma}[Counting, prior, profile]\label{lem:count}
\begin{enumerate}[label=(\roman*),leftmargin=*]
\item The number of labellings with overlap matrix $\omega$ is at most $\exp\{nH(\omega)\}$, and the
number of distinct overlap matrices is at most $(n+1)^{KK_{\max}}$.
\item Under the CRP prior and under the Dirichlet--multinomial prior with $K_{\max}$ labels,
$\log p(z)-\log p(z^\star)=n[\log K-H(\pi(z))]+O(K_{\max}\log n)$ uniformly in $z$.
\item For every $z$, $|\log p(A\mid z)-\tilde\ell(z)|\le C\,K_{\max}^2\log n+C'\sum_{b\le b'}m_{bb'}^2/N_{bb'}$
with $C=C_{a,b}$ of Lemma~\ref{lem:stirling} and an absolute $C'$, and
$\tilde\ell(z)=n\Psi_{\omega(z)}(x(z))-m\log n-m$ up to $O(K_{\max}^2)$ when $N_{bb'}=n^2\nu_{bb'}(1+O(1/n))$.
\end{enumerate}
\end{lemma}
\begin{proof}
(i) The labellings with overlap $\omega$ are obtained by choosing, inside each planted community $a$ of
size $n/K$, which $n\omega_{ab}$ nodes receive label $b$, so their number is
$\prod_a\binom{n/K}{n\omega_{a1},\dots,n\omega_{aK_{\max}}}\le\prod_a\exp\{(n/K)H(K\omega_{a\cdot})\}=\exp\{nH(\omega)\}$
by the multinomial entropy bound. Each entry $n\omega_{ab}$ takes at most $n+1$ values.
(ii) By Stirling's formula, $\sum_b\log\Gamma(n_b+c)-\log\Gamma(n+c')=\sum_bn_b\log n_b-n\log n+O(K_{\max}\log n)
=-nH(\pi)+O(K_{\max}\log n)$ for fixed $c,c'>0$, both for the CRP ($c=0$ with $\log\Gamma(n_b)$ and a factor
$\alpha$ per occupied label) and for the Dirichlet--multinomial prior, and the value at $z^\star$ is
$-n\log K+O(K\log n)$, as in the proof of Theorem~\ref{thm:threshold}.
(iii) The first bound is Lemma~\ref{lem:stirling} applied to the at most $K_{\max}(K_{\max}+1)/2$ blocks
together with $N_{bb'}h(m/N)=m\log(m/N)-m+O(m^2/N)$, and the second is the identity
$m\log(m/N)=m\log(x/\nu)-m\log n$ summed over blocks, with $N_{bb'}=n^2\nu_{bb'}(1+O(1/n))$
contributing $O(m_{bb'}/n)$ per block.
\end{proof}

\begin{lemma}[Chernoff bound for the edge profile]\label{lem:chernoff}
Fix a labelling $z$ with overlap matrix $\omega$. The counts $m_{bb'}(z)$, $b\le b'$, are independent,
$E\,m_{bb'}(z)=n\mu_{bb'}(\omega)(1+O(1/n))$, and for every vector $k$ of non-negative integers
\[
P\big(m_{bb'}(z)=k_{bb'}\ \text{for all }b\le b'\big)\le\exp\Big\{-n\sum_{b\le b'}\kl_{\Pois}\big(k_{bb'}/n\,\|\,E\,m_{bb'}(z)/n\big)\Big\}.
\]
\end{lemma}
\begin{proof}
The blocks of $z$ partition the dyads, and the dyads are independent Bernoulli variables, so the counts
are independent, and $E\,m_{bb'}$ is the sum over the dyads of the block of $c_{\mathrm{in}}/n$ or
$c_{\mathrm{out}}/n$ according to the planted labels of the pair, which gives $n\mu_{bb'}(\omega)$ up to
the diagonal correction. For a sum $S$ of independent Bernoulli variables with mean $\mu$,
$E\,e^{\theta S}=\prod_i(1+p_i(e^\theta-1))\le\exp\{\mu(e^\theta-1)\}$ for every real $\theta$, so the
Chernoff bound gives $P(S\ge k)\le\exp\{-\kl_{\Pois}(k\|\mu)\}$ for $k\ge\mu$ and
$P(S\le k)\le\exp\{-\kl_{\Pois}(k\|\mu)\}$ for $k\le\mu$, hence $P(S=k)\le\exp\{-\kl_{\Pois}(k\|\mu)\}$
for every $k$. Independence across blocks gives the product, and $\kl_{\Pois}(k\|\mu)=n\,\kl_{\Pois}(k/n\|\mu/n)$.
\end{proof}

\begin{proof}[Proof of Theorem~\ref{thm:estimator}]
(a) Let $v<\beta^\star$, so that $\sup_{\mathrm{ov}(\omega)\le v}F(\omega)=-\eta<0$. Let $E_n$ be the
event that some labelling $z$ with $\mathrm{ov}(z)\le v$ has $\log p(z\mid A)\ge\log p(z^\star\mid A)$.
Since $\hat z$ has the largest posterior, $\{\mathrm{ov}(\hat z)\le v\}\subseteq E_n$, so it suffices to
show $P(E_n)\le e^{-n\eta/2}$ for large $n$. Let $G_n$ be the event that $m=nd/2+O(n^{2/3})$ and that
$\tilde\ell(z^\star)=n\Psi^\star-m\log n-m+O(n^{2/3})$, which has probability tending to one
by Chebyshev's inequality applied to the binomial block counts of $z^\star$ (their means are of order $n$
and $t\mapsto t\log(t/N)$ has derivative $O(\log n)$ near the means, as in the proof of
Theorem~\ref{thm:threshold}). On $G_n$, by Lemma~\ref{lem:count}(ii) and (iii), the event
$\log p(z\mid A)\ge\log p(z^\star\mid A)$ implies
\[
\Psi_{\omega(z)}(x(z))\ \ge\ \Psi^\star-\log K+H(\pi(z))-\varepsilon_n,\qquad \varepsilon_n\to0,
\]
uniformly in $z$, once the correction $\sum_{b\le b'}m_{bb'}^2/N_{bb'}$ of Lemma~\ref{lem:count}(iii) is
shown to be $o(n)$ uniformly in $z$. Let $\Delta$ be the maximum degree of $A$, which is $O(\log n)$
with probability tending to one for a bounded mean degree, and fix $\varepsilon>0$. A block with
$m_{bb'}\ge\varepsilon N_{bb'}$ spans $s$ nodes with $\varepsilon s^2/4\le m_{bb'}\le s\Delta/2$, so
$s\le2\Delta/\varepsilon$ and $m_{bb'}\le\Delta^2/\varepsilon$, and there are at most $K_{\max}^2$ blocks,
so the dense blocks contribute at most $K_{\max}^2\Delta^2/\varepsilon=O(\log^2n)$ to the correction,
while the other blocks contribute at most $\varepsilon\sum_{b\le b'}m_{bb'}=\varepsilon m$. Letting
$\varepsilon\to0$ slowly gives $o(n)$, and $G_n$ is enlarged by the event $\Delta\le C\log n$. Group the
labellings by their overlap matrix $\omega$ and their count vector $k=nx$. By Lemma~\ref{lem:count}(i)
and Lemma~\ref{lem:chernoff},
\[
P(E_n\cap G_n)\ \le\ \sum_{\omega:\mathrm{ov}(\omega)\le v}\ \sum_{k}\ e^{nH(\omega)}
\exp\Big\{-n\sum_{b\le b'}\kl_{\Pois}\big(x_{bb'}\|\mu_{bb'}(\omega)(1+O(1/n))\big)\Big\},
\]
where $k$ ranges over the at most $(nd)^{K_{\max}^2}$ integer vectors with $\sum k_{bb'}=m$ and
$\Psi_\omega(k/n)\ge\Psi^\star-\log K+H(\pi)-\varepsilon_n$. Let $F_\varepsilon(\omega)$ be the value
of~\eqref{eq:F} with the constraint $\sum x_{bb'}=d/2$ replaced by $|\sum x_{bb'}-d/2|\le\varepsilon$ and the
right side of the inequality lowered by $\varepsilon$. The infimum in~\eqref{eq:F} is a continuous
function of $(\omega,\varepsilon)$ on the compact set of overlap matrices, because the objective is
continuous and strictly convex in $x$ and the constraint set is the closure of the set where a
continuous function of $(\omega,x)$ is strictly positive, whose boundary has empty interior in the
slice $\sum x=d/2$. Hence $\sup_{\mathrm{ov}(\omega)\le v}F_{\varepsilon}(\omega)\to\sup_{\mathrm{ov}(\omega)\le v}F(\omega)=-\eta$
as $\varepsilon\to0$, and for large $n$ every term of the double sum is at most $\exp\{-n\eta/2\}$. The
number of terms is polynomial in $n$, which gives $P(E_n\cap G_n)\le n^{C}e^{-n\eta/2}$ and (a).

(b) Let $\omega$ have columns proportional to the uniform row, $\omega_{ab}=\pi_b/K$. Then every block
of $z$ has expected density $d/n$, so $\mu_{bb'}(\omega)=d\nu_{bb'}$ with $\sum_{b\le b'}\nu_{bb'}=\tfrac12$,
and $H(\omega)=H(\pi)$. Substitute $x_{bb'}=d\nu_{bb'}y_{bb'}$. The total constraint becomes
$\sum\nu_{bb'}y_{bb'}=\tfrac12$, the objective becomes
$d\sum\nu_{bb'}[y\log y-y+1]=d\sum\nu_{bb'}y_{bb'}\log y_{bb'}$, and
$\Psi_\omega(x)=\sum d\nu_{bb'}y_{bb'}\log(dy_{bb'})=\tfrac d2\log d+d\sum\nu_{bb'}y_{bb'}\log y_{bb'}$.
Since $\Psi^\star-\tfrac d2\log d=J_K/2$ (the identity $\sum_{r\le s}\lambda_{rs}\log(p_{rs}/\bar p)=\tfrac n2J_K$
of the proof of Theorem~\ref{thm:threshold}), the constraint reads
$d\sum\nu_{bb'}y_{bb'}\log y_{bb'}\ge J_K/2-\log K+H(\pi)$, and the objective is the same function. Its
range on the simplex slice is $[0,\infty)$ ($y\equiv1$ gives $0$ and concentrating the mass gives
arbitrarily large values), so the infimum is $\max\{0,J_K/2-\log K+H(\pi)\}$ and
$F(\omega)=H(\pi)-\max\{0,J_K/2-\log K+H(\pi)\}=\min\{H(\pi),\log K-J_K/2\}$. The supremum over
this set is $\log K-J_K/2$, attained at the balanced $\pi$ because $J_K\ge0$. Continuity of $F$ was
shown in (a). If $J_K>2\log K$, then $F<0$ on the set of uncorrelated matrices, which is the set
$\{\mathrm{ov}(\omega)=1/K\}$, and by continuity and compactness $F<0$ on $\{\mathrm{ov}(\omega)\le1/K+\delta\}$
for some $\delta>0$, so $\beta^\star\ge1/K+\delta$. If $J_K<2\log K$, then $F>0$ at the balanced
uncorrelated matrix, which belongs to every set $\{\mathrm{ov}\le v\}$, so the defining set of
$\beta^\star$ is empty.
\end{proof}

\begin{proof}[Proof of Theorem~\ref{thm:two}]
(a) The log likelihood ratio $\log L(A\mid z^\star)-\log L_0(A)$ is a sum over the $\binom n2$ dyads
of independent terms $A_{ij}\log(p_{ij}/p_0)+(1-A_{ij})\log((1-p_{ij})/(1-p_0))$ with
$p_{ij}\in\{c_{\mathrm{in}}/n,c_{\mathrm{out}}/n\}$ and $p_0=d/n$. Its mean is the sum of the Bernoulli
divergences $\kl(p_{ij}\|p_0)=\frac1n\kl_{\Pois}(c\|d)+O(n^{-2})$, which gives
$\frac n4[\kl_{\Pois}(c_{\mathrm{in}}\|d)+\kl_{\Pois}(c_{\mathrm{out}}\|d)]+O(1)=\frac n2J_2+O(1)$
because $c_{\mathrm{in}}+c_{\mathrm{out}}=2d$, and its variance is $O(n)$. Subtracting $n\log2$ and
dividing by $n$ gives (a) by Chebyshev's inequality.

(b) Below the Kesten--Stigum threshold the planted model $P_n$ and the Erd\H{o}s--R\'enyi model $Q_n$
with edge probability $d/n$ are mutually contiguous \citep{Mossel2015}, and $Z_n=dP_n/dQ_n$
has $E_{Q_n}Z_n=1$, so $Q_n(Z_n>M_n)\le1/M_n\to0$ for every $M_n\to\infty$, hence $P_n(Z_n>M_n)\to0$
and $Z_n=O_p(1)$ under $P_n$. Above the threshold \citet{Mossel2015} give a consistent estimator of
$(c_{\mathrm{in}},c_{\mathrm{out}})$, whose value at $c_{\mathrm{in}}=c_{\mathrm{out}}$ is excluded under
$P_n$ with probability tending to one, so there are events $E_n$ with $P_n(E_n)\to1$ and $Q_n(E_n)\to0$,
and for every $M$, $P_n(Z_n\le M)\le P_n(E_n^c)+M\,Q_n(E_n)\to0$.

(c) Write $\bar P_n(A)=\sum_zp(z)p(A\mid z)$ for the numerator of $\mathrm{BF}_n$ and $\bar Q_n(A)=p_1(A)$
for its denominator. \emph{Upper bound.} $\bar P_n$ is a probability distribution on graphs, so
$E_{Q_n}[\bar P_n/Q_n]=1$ and $\bar P_n/Q_n=O_p(1)$ under $Q_n$, hence under $P_n$ by contiguity below the
threshold. For the denominator, $\bar Q_n(A)=\int_0^1L(A\mid p)\,dp$ with $L(A\mid p)=p^m(1-p)^{N-m}$,
$N=\binom n2$. Under $P_n$, $m=Np_0+O_p(\sqrt n)$, and on the interval $|p-p_0|\le n^{-3/2}$ the
derivative of $\log L(A\mid p)$ is $(m-Np)/(p(1-p))=O_p(n^{3/2})$ and its second derivative is
$-m/p^2-(N-m)/(1-p)^2=-O(n^3)$, so $L(A\mid p)\ge L_0(A)e^{-O_p(1)}$ there and
$\bar Q_n(A)\ge2n^{-3/2}e^{-O_p(1)}L_0(A)$. Therefore
$\mathrm{BF}_n=(\bar P_n/Q_n)/(\bar Q_n/Q_n)\le O_p(1)\,n^{3/2}e^{O_p(1)}$.
\emph{Lower bound.} For every labelling $z$ with two labels, the Dirichlet--multinomial or CRP prior gives
$p(z)\ge2^{-n}/(cn)$ for a constant $c$, by $\binom n{n_1}\le2^n$ and Stirling's formula. For every $z$
the conjugate evidence is a product over its three blocks of $\int_0^1p^{e}(1-p)^{N_b-e}\,dp$, and on
each block the log integrand, a concave function of $p$, changes by at most $1/3$ on an interval around
the true rate of the block of length at least $c'n^{-3}$, by the bounds $|m_b/p-(N_b-m_b)/(1-p)|\le n^3/c_{\mathrm{out}}$
and $m_b/p^2+(N_b-m_b)/(1-p)^2\le n^4/c_{\mathrm{out}}^2$ on the relevant range of $p$. Hence
$p(A\mid z)\ge e^{-1}(c'n^{-3})^3L(A\mid z)$ with $L(A\mid z)$ the likelihood at the true rates, and
$\bar P_n(A)\ge\frac{e^{-1}c'^3}{c}\,n^{-10}\sum_z2^{-n}L(A\mid z)=c''n^{-10}Z_n(A)L_0(A)$. Finally
$\bar Q_n(A)\le\max_pL(A\mid p)=L_0(A)\exp\{N\kl(m/N\|p_0)\}=L_0(A)e^{O_p(1)}$ under $P_n$, because the
mean degree of the planted model is $d$. Combining, $\log\mathrm{BF}_n\ge\log Z_n-10\log n-O_p(1)$.
\end{proof}

\subsection{Proofs of the corollaries and propositions}\label{app:props}
This section proves Corollary~\ref{cor:icl}(ii) and (iii), Corollary~\ref{cor:res},
Proposition~\ref{prop:family} and Proposition~\ref{prop:soft}. Corollary~\ref{cor:icl}(i) and
Theorems~\ref{thm:estimator} and~\ref{thm:two} are proved in Section~\ref{app:estimator}, and
Corollary~\ref{cor:window} and Propositions~\ref{prop:finite}, \ref{prop:degsplit} and~\ref{prop:gtyped}
in Section~\ref{app:props4}.
\emph{Corollary~\ref{cor:icl}(ii) and (iii).} With $\delta=c_{\mathrm{in}}-c_{\mathrm{out}}$, $c_{\mathrm{in}}=d+(K-1)\delta/K$ and
$c_{\mathrm{out}}=d-\delta/K$, a second-order expansion of $u\log(u/d)$ around $d$ gives
$J_K=\frac{K-1}{2}\frac{\delta^2}{K^2d}\,(1+O(\delta/d))=\frac{K-1}2\mathrm{SNR}_{\mathrm{KS}}(1+O(\delta/d))$, so
$J_K=2\log K$ has the solution $\mathrm{SNR}^\star(d,K)=4\log K/(K-1)+O(d^{-1/2})$, and $4\log K/(K-1)>1$
holds exactly for $K\le10$ ($1.023$ at $K=10$, $0.959$ at $K=11$, and $4\log K/(K-1)$ decreases in $K$).
Part (ii) is the root of $J_K=2\log K$ in $\mathrm{SNR}_{\mathrm{KS}}$ for each $(d,K)$ of
Table~\ref{tab:iclthr}, computed by bisection (script \texttt{theory\_checks/check\_thresholds\_exact.py}),
and $d^\star(K)$ is the root of $\mathrm{SNR}^\star(d,K)=1$ in $d$. Part (iii): the exact integrated
classification likelihood of \citet{Come2015} is $\log p(A\mid z)+\log p(z)$ with Jeffreys block priors
and the Dirichlet--multinomial partition prior, which is the objective of the theorem with
$a=b=\frac12$, so its maximiser over $\{z^\star,z_1\}$ is $z_1$ below the threshold.

\emph{Corollary~\ref{cor:res}.} Write $D$ for the left side of Theorem~\ref{thm:resolution} and
$c=\frac12(J-2\log2)>0$. By the theorem, $D=n c-\frac{3K^2}{8}(a+\frac12)\log\nu+O(K^2+K\log\nu)+o_p(n)$,
where $-\frac{2K}8(a+\frac12)\log\nu$ has been absorbed into $O(K\log\nu)$. With $n=K\nu$ the two leading
terms are equal at $K=K_{\mathrm{res}}(\nu)$, and $D=nc\,(1-K/K_{\mathrm{res}}(\nu))+O(K^2+K\log\nu)+o_p(n)$.
If $K\le(1-\delta)K_{\mathrm{res}}(\nu)$ then $nc(1-K/K_{\mathrm{res}})\ge\delta nc$. When $K\ge1$ is
bounded, $K^2+K\log\nu=O(\log\nu)=o(n)$, and when $K\to\infty$, $n=K\nu\ge K\,K\log\nu\cdot3(2a+1)/(8(J-2\log2))\cdot(1+o(1))$
by $K\le K_{\mathrm{res}}$, so $K^2=O(n/\log\nu)=o(n)$. Hence $D\ge\delta nc-o(n)-o_p(n)>0$ with
probability tending to one. The case $K\ge(1+\delta)K_{\mathrm{res}}(\nu)$ is symmetric, with
$K=O(\nu/\log\nu)$ ensuring $K^2\log\nu=O(n)$ and $K^2=o(n)$, so that $D\le-\delta nc\,K/K_{\mathrm{res}}+o(n)+o_p(n)<0$.
When $J<2\log2$ the leading term $nc$ is negative and every other term is negative or $o_p(n)$.

\emph{Proposition~\ref{prop:family}.} For a block with $e$ edges among $N$ dyads and $\hat p=e/N$, the
Bernoulli profile is $N h(\hat p)=e\log\hat p+(N-e)\log(1-\hat p)$ and the Poisson profile is
$\ell_{\Pois}(e,N)=e\log\hat p-e$, so $N h(\hat p)-\ell_{\Pois}(e,N)=(N-e)\log(1-\hat p)+e=N g(\hat p)$.
Lemma~\ref{lem:stirling} bounds $|\log p_{\mathrm{BB}}-Nh(\hat p)|$ by $C_{a,b}\log(N+a+b+2)$ and
Lemma~\ref{lem:poisson} bounds $|\log p_{\mathrm{GP}}-\ell_{\Pois}|$ by
$(2a'+2)\log(N+e+a'+b'+2)+b'e/N+1+|\log(b'^{a'}/\Gamma(a'))|$. With $c=a+b+a'+b'+4$, $e\le N\le n^2$
gives $\log(N+a+b+2)\le\log(n^2+c)$, $\log(N+e+a'+b'+2)\le\log2+\log(n^2+c)\le2\log(n^2+c)$ and
$b'e/N+1\le(b'+1)\log(n^2+c)$, since $\log(n^2+c)\ge\log4>1$. Each block term is therefore at most
$[C_{a,b}+4a'+5+b'+|\log(b'^{a'}/\Gamma(a'))|]\log(n^2+c)$, and summing over the $K(K+1)/2$ blocks
gives (i) with $C=\frac12[C_{a,b}+4a'+5+b'+|\log(b'^{a'}/\Gamma(a'))|]$.
For (ii), $g'(p)=-\log(1-p)\ge0$ and $g''(p)=1/(1-p)\ge1$ on $[0,1)$, so $g$ is increasing and
$g(p)\ge g(0)+g'(0)p+p^2/2=p^2/2$. If a block of $z$ with $N$ dyads and density $\hat p$ is the union of
sub-blocks with $(N_i,\hat p_i)$, then $\hat p=\sum_i(N_i/N)\hat p_i$, and the second-order Taylor bound
$g(\hat p_i)\ge g(\hat p)+g'(\hat p)(\hat p_i-\hat p)+\frac12(\hat p_i-\hat p)^2$ summed with weights
$N_i$ gives $\sum_iN_ig(\hat p_i)-Ng(\hat p)\ge\frac12\sum_iN_i(\hat p_i-\hat p)^2$, since the linear
terms cancel. Summing over the blocks of $z$ gives (ii).

\emph{Proposition~\ref{prop:soft}.} By definition,
$\Delta_m(e,N)=\log\frac{\Gamma(a_s+e)}{\Gamma(1+e)}+\log\frac{\Gamma(b_s+N-e)}{\Gamma(1+N-e)}
-\log\frac{\Gamma(a_s+b_s+N)}{\Gamma(N+2)}-\log\Bfun(a_s,b_s)$. Apply \eqref{eq:gammaratio}, whose arguments are written $(x,s)$ there, with the pairs $(e,a_s)$,
$(N-e,b_s)$ and $(N+1,a_s+b_s-1)$, all admissible since $e,N-e\ge1$ and $a_s,b_s\ge1$: the first ratio lies between $(a_s-1)\log e$ and $(a_s-1)\log(e+a_s)$, the second
between $(b_s-1)\log(N-e)$ and $(b_s-1)\log(N-e+b_s)$, and the third between $(m-2)\log(N+1)$ and
$(m-2)\log(N+m)$. Since $(a_s-1)\log e+(b_s-1)\log(N-e)=\log\pi_{a_s,b_s}(\hat p)+\log\Bfun(a_s,b_s)+(m-2)\log N$,
the lower bound is $\log\pi_{a_s,b_s}(\hat p)+(m-2)\log N-(m-2)\log(N+m)$ and the upper bound is
$\log\pi_{a_s,b_s}(\hat p)+(a_s-1)\log(1+a_s/e)+(b_s-1)\log(1+b_s/(N-e))+(m-2)\log N-(m-2)\log(N+1)$,
and $(m-2)\log(N/(N+1))\le0$. The two consequences follow by taking the supremum over $\hat p$ of the
upper bound, with $\log(1+a_s/e)\le\log(1+a_s)$ and $\log(1+b_s/(N-e))\le\log(1+b_s)$, and the lower
bound at the given $\hat p$. The numbers quoted after the proposition are exact evaluations of
$\Delta_m$ (script \texttt{theory\_checks/check\_soft\_bound.py}).

\subsection{Proofs of Corollary~\ref{cor:window} and Propositions~\ref{prop:finite}, \ref{prop:degsplit} and~\ref{prop:gtyped}}\label{app:props4}

\begin{proof}[Proof of Corollary~\ref{cor:window}]
Let $z_s$ keep the classes apart and $z_m$ merge them, $n=2\nu$. By Lemma~\ref{lem:stirling} the
Bernoulli evidences are the profile log-likelihoods up to $O(\log n)$, and
$\ell(z_s)=\nu(\nu-1)h(\hat p)+\nu^2h(\hat q)$, $\ell(z_m)=\nu(2\nu-1)h(\hat{\bar p})$ with
$\hat{\bar p}=[\nu(\nu-1)\hat p+\nu^2\hat q]/[\nu(2\nu-1)]$. With $\hat p\to p$, $\hat q\to q$ and
$\hat{\bar p}\to(p+q)/2$, $\ell(z_s)-\ell(z_m)=2\nu^2\mathrm{JG}_h(p,q)+O_p(\nu)$. The partition prior
gives $\log p(z_s)-\log p(z_m)=-2\nu\log2+O(\log\nu)$ under the CRP prior
($\log\alpha+2\log\Gamma(\nu)-\log\Gamma(2\nu)$) and under the Dirichlet--multinomial prior, by
Stirling's formula. Hence $\log p(z_s\mid A)-\log p(z_m\mid A)=2\nu[\nu\mathrm{JG}_h(p,q)-\log2]+O_p(\nu)$,
which is positive for large $\nu$ if $\nu\mathrm{JG}_h>\log2$ and negative if $\nu\mathrm{JG}_h<\log2$
(the $O_p(\nu)$ term is $\nu$ times a quantity of order $|\hat p-p|+|\hat q-q|=O_p(1/\nu)$, so it is
$O_p(1)$). The Poisson evidence replaces $h$ by $\varphi$ (Lemma~\ref{lem:poisson}), and
$h=\varphi+g$ with $g$ convex (Proposition~\ref{prop:family}(ii)), so
$\mathrm{JG}_h=\mathrm{JG}_\varphi+\mathrm{JG}_g$ and $\mathrm{JG}_g\ge\tfrac18(p-q)^2\min g''\ge\tfrac18(p-q)^2$
by Taylor's formula with $g''\ge1$. The ratio of the second derivatives $h''(p)/\varphi''(p)=1/(1-p)$
gives the width of the window.
\end{proof}

\begin{proof}[Proof of Proposition~\ref{prop:finite}]
Write $\log p_{\DC}(Y_C)-\log p_{\ER\text{-}2}(Y_C\mid z_s)=[\hat\ell_{\DC}-\hat\ell_{\ER\text{-}2}]-R_\nu+[\hat\ell_{\ER\text{-}2}-\log p_{\ER\text{-}2}(Y_C\mid z_s)]$,
where $\hat\ell_{\ER\text{-}2}$ is the profile log-likelihood of the three Poisson blocks of $z_s$. The
last bracket is $O(\log\nu)$ by Lemma~\ref{lem:poisson}. For the first bracket,
$\hat\ell_{\DC}\ge\ell_{\DC}(Y_C;\theta,\omega_0)$ at the true propensities, and
$\hat\ell_{\ER\text{-}2}=\sum_b[e_b\log(e_b/N_b)-e_b]$ exceeds the log-likelihood of the three blocks at
their population rates $\bar\lambda_b$ by $\sum_bN_b\,\kl_{\Pois}(e_b/N_b\|\bar\lambda_b)=O_p(1)$. The
difference between the Poisson log-likelihoods at $(\theta,\omega_0)$ and at the block rates is
$\sum_{i<j}[y_{ij}\log(\lambda_{ij}/\bar\lambda_{b(ij)})-\lambda_{ij}+\bar\lambda_{b(ij)}]$ with
$\lambda_{ij}=\theta_i\theta_j\omega_0$, a sum of independent terms with mean
$\sum_{i<j}\kl_{\Pois}(\lambda_{ij}\|\bar\lambda_{b(ij)})$ and variance $O(\nu^2)$. Conditionally on the
atoms, the number of pairs with atom types $(t,t')$ in block $b$ is $\tfrac{\nu^2}2F(t)F(t')(1+O_p(\nu^{-1/2}))$,
so the mean is $\kappa_F(z_s)\nu^2(1+O_p(\nu^{-1/2}))$, and $\kappa_F>0$ when some block contains pairs with
two distinct rates, which happens for every split of three or more atoms into two blocks. The profile
maximisation over $(\theta,\omega)$ adds at most $O_p(\nu)$ to $\ell_{\DC}(Y_C;\theta,\omega_0)$ (the
maximum of a log-likelihood with $\nu+1$ parameters over a compact neighbourhood of the truth exceeds
its value at the truth by a quantity whose mean is of order the number of parameters). This gives the
first display, with $R_\nu\ge0$ replaced by the bound of Lemma~\ref{lem:dcprofile}: by \eqref{eq:dclower}
at $(\theta,\omega_0)$, $\log p_{\DC}(Y_C)\ge\ell_{\DC}(Y_C;\theta,\omega_0)-\tfrac{\omega_0}2\sum_i\theta_i^2-\nu(1+\log(1+2\bar d))-L_P-1$,
and $\hat\ell_{\DC}-\ell_{\DC}(Y_C;\theta,\omega_0)=O_p(\nu)$ is absorbed. The numbers quoted for the
design of Table~\ref{tab:splitb} are computed by the script named in the statement, which evaluates
$\kappa_F$ from its definition, the exact evidences of Table~\ref{tab:splitb} on the finer grid, and
$R_\nu$ from the closed form of the profile maximiser $\hat\theta_i\hat\theta_j\hat\omega=k_ik_j/(2E)$.
\end{proof}

\begin{proof}[Proof of Proposition~\ref{prop:degsplit}]
(i) Under the uniform degree prior the degree factor of a block with $\nu_r$ nodes and $e_r$
half-edges is $-\log\binom{e_r+\nu_r-1}{e_r}=-\log\Gamma(e_r+\nu_r)+\log\Gamma(e_r+1)+\log\Gamma(\nu_r)$,
and Stirling's formula gives $\log\binom{e_r+\nu_r-1}{e_r}=\nu_r\phi(e_r/\nu_r)+O(\log(e_r+\nu_r))$,
because $(e+\nu)\log(e+\nu)-e\log e-\nu\log\nu=\nu[(1+x)\log(1+x)-x\log x]$ at $x=e/\nu$. The
difference of the degree factors of $z_s$ and $z_1$ is therefore
$\nu\phi(\bar d)-\nu_1\phi(\bar d_1)-\nu_2\phi(\bar d_2)+O(\log e)$, which is non-negative by Jensen's
inequality for the concave $\phi$ ($\phi''(x)=-1/(x(1+x))$) at the weights $\nu_r/\nu$, since
$\sum_r(\nu_r/\nu)\bar d_r=\bar d$, and strictly positive of order $\nu$ when $\bar d_1\ne\bar d_2$. The
configuration count of \citet{Peixoto2017} contributes
$\log[e!\,e_{12}!\,e_{11}!!\,e_{22}!!/(e!!\,e_1!\,e_2!)]$ to the same difference, where $e_{rs}$ are the
block edge counts of $z_s$, and the uniform prior on the block edge-count matrix contributes
$-\log\binom{E+2}{2}$; at block counts within $O(\sqrt e)$ of their configuration-model values
$e_{rs}=e_re_s/(2e)$ both are $O(\log e)$ by Stirling's formula. The script
\texttt{theory\_checks/check\_degree\_split.py} evaluates the full difference on simulated communities
(Table~\ref{tab:degsplit}): the column without a degree prior, which isolates these two contributions,
stays between $-13$ and $-7$ nats for $\nu$ from $100$ to $3200$, and the uniform-prior column follows
the displayed Jensen gap within $6$ nats.
(ii) The factor $\prod_ip(k_i)$ is the same for every labelling, so it cancels from the comparison.
\end{proof}

\begin{proof}[Proof of Proposition~\ref{prop:gtyped}]
The two row factors are related by the multinomial identity
$\nu^{-e}=P_{\mathrm{mult}}(\mathbf k)\prod_ik_i!/e!$, where $P_{\mathrm{mult}}(\mathbf k)=e!\,\nu^{-e}/\prod_ik_i!$
is the probability of the degree sequence under the uniform allocation of $e$ half-edges, so the log
ratio of the display equals $-\log P_{\mathrm{mult}}(\mathbf k)-\log\binom{e+\nu-1}{e}$. If the
allocation has weights $\theta_i/\nu$ then $-\log P_{\mathrm{mult}}(\mathbf k)=-\log P_\theta(\mathbf k)+\sum_ik_i\log\theta_i$
exactly, with $P_\theta$ the multinomial probability with those weights. For degenerate $F$,
Stirling's formula applied to each $k_i!$ and to $e!$ gives
$-\log P_{\mathrm{mult}}(\mathbf k)=\sum_i[k_i\log(k_i/\bar d)+\tfrac12\log(2\pi k_i)]+O(\log e)+O(\sum_i1/k_i)$
after the terms $e\log\bar d-e\log e+e\log\nu=0$ cancel, and with $k_i=\bar d+O_p(\sqrt{\bar d})$,
$\sum_ik_i\log(k_i/\bar d)=\tfrac12\sum_i(k_i-\bar d)^2/\bar d+O_p(\nu/\sqrt{\bar d})=\tfrac\nu2(1+o_p(1))$
and $\sum_i\tfrac12\log(2\pi k_i)=\tfrac\nu2\log(2\pi\bar d)+O_p(\sqrt\nu)$, so
$-\log P_{\mathrm{mult}}(\mathbf k)=\tfrac\nu2[\log(2\pi\bar d)+1]+O_p(\sqrt\nu)$. The uniform-prior term
is $\log\binom{e+\nu-1}{e}=\nu\phi(\bar d)+O(\log e)=\nu[\log\bar d+1+O(1/\bar d)]$, and the difference is
the first case of the display. For non-degenerate $F$ the same expansion holds for $-\log P_\theta(\mathbf k)$
with $\bar d\theta_i$ in place of $\bar d$, so it is $O(\nu\log\bar d)$, while
$\sum_ik_i\log\theta_i=\bar d\sum_i\theta_i\log\theta_i(1+O_p(\bar d^{-1/2}))=\nu\bar d\,\E_F[\theta\log\theta](1+o_p(1))$
by the law of large numbers, which dominates because $\bar d\to\infty$. With the global degree prior the
term $\log\binom{e+\nu-1}{e}$ is replaced by $-\sum_i\log p(k_i)$, and the comparison of
$-\log P_{\mathrm{mult}}(\mathbf k)$ with $-\sum_i\log p(k_i)$ is the comparison of the two laws named in the
statement, whose difference per node is the divergence between the community's degree law and $p$
when the two differ and $o_p(1)$ when they coincide.
\end{proof}

\subsection{Bridge configurations}\label{app:bridge}
\begin{proposition}[Bridge configurations]\label{prop:bridge}
Let two communities $C_1,C_2$ of $\nu=50$ nodes and within density $p_{\mathrm{in}}=0.11$ receive $e_b$
contaminating cross edges, all incident to a set $B=B_1\cup B_2$ of $|B|=10$ nodes ($5$ from each
community), placed first on $B_1\times B_2$ and then in equal numbers on $B_1\times(C_2\setminus B_2)$ and
$B_2\times(C_1\setminus B_1)$, with no other cross edges. Let $z_2=\{C_1,C_2\}$,
$z_3=\{C_1\setminus B_1,\,C_2\setminus B_2,\,B\}$ and $z_1$ the merged partition, compare evidences under
Beta$(1,1)$ block priors (without partition prior) at the expected block counts (rounded, with a
contaminating edge shared between two blocks counting one half in each). Then: (i) $z_3$ has higher
evidence than $z_2$ for all $11\le e_b\le400$ without a cap, with a margin that grows roughly linearly
in $e_b$ (from $180$ nats at $e_b=100$ to $600$ at $e_b=400$), and for all $15\le e_b\le400$ with the cap
$x=0.02$, with a margin of about $59$ nats that is nearly constant in $e_b$. (ii) Without a cap, $z_3$ is
preferred to $z_1$ for every $e_b\le400$ while $z_2$ loses to $z_1$ at $e_b=197$. (iii) With the cap
$x=0.02$, $z_2$ loses to $z_1$ at $e_b=114$ while $z_3$ loses to $z_1$ only at $e_b=154$.
\end{proposition}
\begin{proof}
All block counts are determined by $(\nu,|B|,p_{\mathrm{in}},e_b)$: the core blocks $C_i\setminus B_i$ have
$0.11\cdot990$ rounded, i.e.\ $109$, edges among $990$ dyads. The bridge block $B$ has
$2+\min(e_b,25)$ edges among $45$ dyads ($2$ genuine ties within $B_1$ and $B_2$). Each block
$(B,C_i\setminus B_i)$ has $25$ genuine ties plus half of $(e_b-25)^+$ among $450$ dyads. The block
$(C_1\setminus B_1,C_2\setminus B_2)$ is empty among $2025$ dyads. The partition $z_2$ has $135$ edges in each community
among $1225$ dyads and $e_b$ among $2500$, and $z_1$ has $270+e_b$ edges among $4950$ dyads. The three
evidences are products of the corresponding Beta functions, with the factor \eqref{eq:cap} on
off-diagonal blocks when $x<1$. They are evaluated in \texttt{experiments/theory\_checks.py}
(Figure~\ref{fig:cap}, right), and the claims are read off. The mechanism of (i) is that the bridge
nodes have a connectivity profile of their own (density $0.60$ inside $B$ once $e_b\ge25$, against
$0.11$ in the cores), which the evidence detects with or without a cap. Under the cap the two
$450$-dyad blocks $(B,C_i\setminus B_i)$ of $z_3$ have density $(e_b+25)/900>x$ and pay the quadratic
penalty of Proposition~\ref{thm:cap}(ii) at a rate per contaminating edge ($2.2$ nats at $e_b=114$) that
exceeds the rate paid by the $2500$-dyad cross block of $z_2$ ($0.86$ nats), which cancels the linear
growth of the uncapped margin and leaves the nearly constant margin of $59$ nats in (i). Item (iii) then
follows: $z_3$ loses to $z_1$ where $\log\mathrm{BF}_x(z_2\!:\!z_1)$, which is $-0.9$ at $e_b=114$, reaches
$-59$, i.e.\ at $e_b=154$.
\end{proof}

\subsection{Collapsing}\label{app:collapse}
\begin{remark}[Collapsing]\label{prop:collapse}
Let $\theta$ be the block parameters and $p(z\mid Y)=\int p(z,\theta\mid Y)\,d\theta$. Three facts about
collapsing are used in the paper. (i) The collapsed posterior on $z$ is the $z$-marginal of the joint
posterior $p(z,\theta\mid Y)$, by the definition of a marginal, so the Occam behaviour of $p(Y\mid z)$ is
a property of the prior and the likelihood and not of the algorithm. (ii) Consider the systematic-scan
Gibbs sampler that updates the components $z_i$, $i=1,\dots,n$, in turn from their full conditionals
given $\theta$ and then $\theta$ from $p(\theta\mid z,Y)$, and the collapsed systematic-scan sampler
that updates the components $z_i$, $i=1,\dots,n$, in turn from $p(z_i\mid z_{-i},Y)$. Let $F_s$ and
$F_c$ be their forward operators on the mean-zero subspaces of $L^2$ of the respective stationary
distributions $p(z,\theta\mid Y)$ and $p(z\mid Y)$. Then $\|F_c\|\le\|F_s\|$. This is Theorem~1 of
\citet{Liu1994}, which orders the standard, grouped and collapsed systematic-scan samplers for any
number of components when the collapsed component is the last one visited (the three-component case is
Theorem~5.1 of \citet{LiuWongKong1994}). Here the components are $z_i$, $i=1,\dots,n$, and $\theta$,
and $\theta$ is integrated out. The norm of the forward operator is the maximal correlation between
successive states. As \citet{Liu1994} notes, an ordering of these norms does not imply the same
ordering of the spectral radii, and the result says nothing about the split and merge moves
interleaved with the sweeps or about wall-clock time. (iii) One collapsed single-site sweep costs
$O(\sum_i\deg(i)+nK^2)$ evaluations of block evidences, and one uncollapsed single-site sweep costs
$O(\sum_i\deg(i)+K\sum_i\min(\deg(i),K)+nK+K^2)$ arithmetic operations. The counts follow from
Algorithm~\ref{alg:sweep}: for node $i$ the edge counts $e_{il}$ from $i$ to every community cost
$O(\deg(i))$, and scoring a candidate community $k$ requires the $K$ block-evidence ratios of row $k$,
so $K^2$ evaluations per node. In the uncollapsed sweep, with $\theta$ drawn once per sweep ($O(K^2)$),
the log-likelihood of candidate $k$ is $\sum_l[e_{il}\log\theta_{kl}+(n_l-e_{il})\log(1-\theta_{kl})]$.
The terms $\sum_ln_l\log(1-\theta_{kl})$ are maintained for all $k$ at $O(K)$ cost per move, and the
correction $\sum_le_{il}\log[\theta_{kl}/(1-\theta_{kl})]$ has at most $\min(\deg(i),K)$ nonzero terms,
for each of $K$ candidates.
\end{remark}

\subsection{Laplace form for non-conjugate types}\label{app:occam}
\begin{theorem}[Laplace form of an approximate type evidence]\label{thm:occam}
Let the within-community evidence of type $t$ be $p_t(Y_C)=\int_\Omega e^{\ell_N(\vartheta)}\pi_t(\vartheta)\,d\vartheta$,
where $N=\binom{n_k}2$, $\Omega\subset\bb R^{d_t}$ is open, $\pi_t$ is a probability density on
$\Omega$, $\ell_N$ is twice continuously differentiable with a maximiser $\hat\vartheta_N$ such that the
closed ball $\bar B(\hat\vartheta_N,\rho)$ lies in $\Omega$, and for constants $0<c\le C$, $\gamma>0$,
$0<\pi_-\le\pi_+$ not depending on $N$:
(A1) $cN\,I\preceq-\nabla^2\ell_N(\vartheta)\preceq CN\,I$ and $\pi_-\le\pi_t(\vartheta)\le\pi_+$ on
$\bar B(\hat\vartheta_N,\rho)$, and (A2) $\ell_N(\vartheta)\le\ell_N(\hat\vartheta_N)-\gamma N$ outside
$\bar B(\hat\vartheta_N,\rho)$. Then there are $N_0$ and $M$ depending only on
$(d_t,c,C,\rho,\gamma,\pi_\pm)$ such that for all $N\ge N_0$
\[
\Big|\log p_t(Y_C)-\ell_N(\hat\vartheta_N)+\frac{d_t}2\log N\Big|\le M ,
\]
so the penalty of a $d_t$-parameter type is $d_t\log n_k+O(1)$, and if every type in $\cT$ satisfies
(A1)--(A2), $\log p(Y_C)=\max_{t\in\cT}\big\{\ell^{(t)}_N(\hat\vartheta^{(t)}_N)-\frac{d_t}2\log N\big\}+O(1)$
with the $O(1)$ bounded by $M+\log|\cT|+\max_t\log(1/\pi(t))$.
\end{theorem}
\begin{proof}
On the ball, Taylor's formula with $\nabla\ell_N(\hat\vartheta_N)=0$ and (A1) gives
$-\frac{CN}2\|u\|^2\le\ell_N(\hat\vartheta_N+u)-\ell_N(\hat\vartheta_N)\le-\frac{cN}2\|u\|^2$ for
$\|u\|\le\rho$. Hence the integral over the ball, divided by $e^{\ell_N(\hat\vartheta_N)}$, is at most
$\pi_+\int_{\bb R^{d_t}}e^{-cN\|u\|^2/2}\,du=\pi_+(2\pi/(cN))^{d_t/2}$ and at least
\[
\pi_-\int_{\|u\|\le\rho}e^{-CN\|u\|^2/2}\,du=\pi_-\Big(\frac{2\pi}{CN}\Big)^{d_t/2}\Pbb\big(\|Z\|\le\rho\sqrt{CN}\big)
\ \ge\ \frac{\pi_-}2\Big(\frac{2\pi}{CN}\Big)^{d_t/2}
\]
for $N\ge N_0$, where $Z$ is standard normal in $\bb R^{d_t}$. By (A2)
the integral outside the ball, divided by $e^{\ell_N(\hat\vartheta_N)}$, is at most $e^{-\gamma N}$,
which is at most the lower bound just obtained for $N\ge N_0$. Taking logarithms,
$\log p_t(Y_C)-\ell_N(\hat\vartheta_N)+\frac{d_t}2\log N$ lies between
$\log\frac{\pi_-}2+\frac{d_t}2\log\frac{2\pi}C$ and $\log(2\pi_+)+\frac{d_t}2\log\frac{2\pi}c$, and
$\frac{d_t}2\log N=d_t\log n_k+O(1)$. The last claim follows from
$\max_ta_t\le\log\sum_te^{a_t}\le\max_ta_t+\log|\cT|$ applied to the type-collapsed evidence
$p(Y_C)=\sum_t\pi(t)p_t(Y_C)$ of \eqref{eq:csbm} with $a_t=\log\pi(t)+\log p_t(Y_C)$.
\end{proof}
The theorem is included for the approximate types of Section~\ref{sec:types} only, which are not used
in the experiments.

\section{Catalogue, sampler and validation}\label{app:catalogue}

\subsection{Catalogue}
Every off-diagonal block, and every diagonal block of the Erd\H{o}s--R\'enyi type, is a block of a
conjugate exponential family: $p(y\mid\theta)=h_0(y)\exp\{\theta^\top T(y)-A(\theta)\}$ with prior
$\pi(\theta)\propto\exp\{\theta^\top\eta_0-\tau_0A(\theta)\}\ind{\theta\in\cS}$ on a (possibly restricted)
set $\cS$, so that with $Z_\cS(\eta,\tau)=\int_\cS\exp\{\theta^\top\eta-\tau A(\theta)\}d\theta$,
\[
p(Y_{rs}\mid z)=\Big(\prod_{\cD_{rs}(z)}h(y_{ij})\Big)\frac{Z_\cS(\eta_0+S_{rs}(z),\ \tau_0+N_{rs}(z))}{Z_\cS(\eta_0,\tau_0)}.
\]
Table~\ref{tab:catalogue} lists the instances implemented in the released code. Constants that do not
depend on $z$ (such as $\prod y_{ij}!$) are dropped consistently. When $K(z)$ varies the per-block
normalising constants are retained because the number of factors changes.

\begin{table}[H]
\centering\small
\caption{Collapsed block evidences implemented in the released code. Each row gives the data type, the block likelihood, the conjugate prior and the log block evidence up to constants that do not depend on the partition, with $e$ the edge count, $S$ the block sum, $N$ the number of dyads and $e_{\mathrm{act}}$ the number of active dyads.}
\label{tab:catalogue}
\resizebox{\linewidth}{!}{\begin{tabular}{@{}llll@{}}
\toprule
Data & Block likelihood & Prior & $\log p(Y_{rs}\mid z)$ (up to $z$-free constants)\\
\midrule
binary & $\Bern(p)$ & Beta$(a,b)$ & $\log\Bfun(a+e,b+N-e)-\log\Bfun(a,b)$\\
binary, cap $x$ & $\Bern(p)$ & Beta$(a,b)$ on $[0,x]$ & above $+\log I_x(a+e,b+N-e)-\log I_x(a,b)$\\
counts & $\Pois(\lambda)$ & Gamma$(a,b)$ & $\log\Gamma(a+S)-(a+S)\log(b+N)+a\log b-\log\Gamma(a)$\\
counts, cap $x$ & $\Pois(\lambda)$ & Gamma$(a,b)$ on $[0,x]$ & above $+\log P(a+S,(b+N)x)-\log P(a,bx)$\\
real weights & $\mathcal N(\mu,\sigma^2)$ & Normal--Inverse-Gamma$(\mu_0,\kappa_0,\alpha_0,\beta_0)$ & $-\frac N2\log2\pi+\frac12\log\frac{\kappa_0}{\kappa_N}+\log\frac{\Gamma(\alpha_N)}{\Gamma(\alpha_0)}+\alpha_0\log\beta_0-\alpha_N\log\beta_N$\\
signed / dyad states & $\Cat(\pi)$, $C\le4$ & Dirichlet$(\alpha)$ & $\log\Bfun(\alpha+\mathbf c)-\log\Bfun(\alpha)$\\
zero-inflated counts & $\Bern(p)$, $\Pois(\lambda)$ on active dyads & Beta, Gamma & Beta--Bernoulli$(e_{\mathrm{act}},N)$ + Gamma--Poisson$(S,e_{\mathrm{act}})$, with $Z_{ij}$ augmentation\\
multiplex & $\prod_\ell\Bern(p^{(\ell)})$ & Beta per layer & $\sum_\ell\log p^{(\ell)}_{rs}$\\
DC type (diagonal) & Poisson multigraph & uniform degree prior & \eqref{eq:dc}\\
global-DC (Peixoto) & Poisson multigraph & uniform edge-matrix prior & block, row and edge-matrix terms of \citet{Peixoto2017}\\
\bottomrule
\end{tabular}}
\end{table}
Here $P(\cdot,\cdot)$ is the regularised lower incomplete gamma function, $\mathbf c$ the vector of category
counts (category $0$ = no tie), and $(\kappa_N,\alpha_N,\beta_N)$ the usual Normal--Inverse-Gamma updates from
$(N,\sum y,\sum y^2)$. The zero-inflated block, with Beta$(a_p,b_p)$ prior on the activity probability and
Gamma$(a_\lambda,b_\lambda)$ prior on the intensity, samples the activity indicator $Z_{ij}$ of every zero
dyad from its collapsed conditional,
\[
\Pbb(Z_{ij}=1\mid\cdot)\propto(a_p+e^{-}_{\mathrm{act}})\Big(\frac{b_\lambda+e^{-}_{\mathrm{act}}}{b_\lambda+e^{-}_{\mathrm{act}}+1}\Big)^{a_\lambda+S},
\qquad
\Pbb(Z_{ij}=0\mid\cdot)\propto b_p+N-1-e^{-}_{\mathrm{act}},
\]
where $e^{-}_{\mathrm{act}}$ is the number of active dyads in the block excluding $(i,j)$ and $S$ the block
sum.

\subsection{One sweep}
\begin{algorithm}[tbp]
\caption{One collapsed Gibbs sweep with type-collapsed diagonal blocks.}\label{alg:sweep}
\begin{algorithmic}[1]
\Require $Y$, partition $z$, block statistics $\{S_{rs},N_{rs}\}$, within-community degrees $d_i$ and
$\sigma_k=\sum_{i\in C_k}\log d_i!$ (for the DC type), partition prior $p(z)$.
\For{$i$ in a random order}
  \State compute $c_l=\sum_{j\in C_l}T(y_{ij})$ for all $l$ and, for the DC type,
  \Statex \hspace{\algorithmicindent} $\delta_l=\log(c_l)!+\sum_{j\in C_l,\,y_{ij}\ne0}[\log(d_j+y_{ij})!-\log d_j!]$
  \State remove $i$ from $C_{z_i}$ and update the statistics of row $z_i$, the degrees $d_j$ of its
  \Statex \hspace{\algorithmicindent} neighbours in $C_{z_i}$, and $\sigma_{z_i}$
  \For{each candidate $k$ (occupied labels, and one empty label under the CRP prior)}
    \State $\log w_k\gets\log p(z_i=k\mid z_{-i})$
    \Statex \hspace{2\algorithmicindent} $+\sum_{l\ne k}\big[\log p_{kl}(S_{kl}+c_l,N_{kl}+n_l)-\log p_{kl}(S_{kl},N_{kl})\big]$
    \State $\log w_k\gets\log w_k+\log p(Y_{C_k\cup\{i\}})-\log p(Y_{C_k})$,
    \Statex \hspace{2\algorithmicindent} computed from $(S_{kk}+c_k,\ N_{kk}+n_k,\ n_k+1,\ \sigma_k+\delta_k)$
  \EndFor
  \State sample $z_i\propto w_k$ (or take the arg-max when polishing), insert $i$ and update statistics
\EndFor
\State optionally one split and merge proposal \citep{Jain2004}
\end{algorithmic}
\end{algorithm}
The per-node cost is $O(\deg(i)+K^2)$ evaluations of $\log\Bfun$, $\log\Gamma$, $\log I_x$ or
$\log P$. Incomplete beta and gamma functions are evaluated in log space by continued fractions so that
blocks far outside a cap receive large but finite penalties. Held-out dyads are removed from the dyad
counts $N_{rs}$ through per-node masked-partner counts, which is how the held-out prediction of
Section~\ref{sec:real} is computed.

\subsection{The split and merge move}
Two nodes $i\neq j$ are drawn uniformly. If $z_i=z_j=a$ (split), an empty label $b$ is drawn uniformly
among the $n_\emptyset$ empty labels, a \emph{launch state} assigns $i$ to $a$, $j$ to $b$ and the other
members of $C_a$ either at random or by the sign of the second eigenvector of the Bethe Hessian of the
subgraph induced by $C_a$ (a deterministic function of $(C_a,i,j)$), $t$ restricted Gibbs scans over
$C_a\setminus\{i,j\}$ with labels in $\{a,b\}$ are run from the launch state, and one final restricted
scan produces the proposal $z'$ with probability $q(z'\mid z)=\prod_v\Pbb(z'_v\mid\text{current state})$
(product of the restricted single-site conditionals in the order visited). The acceptance probability is
$\min\{1,\ \frac{p(z'\mid Y)}{p(z\mid Y)}\cdot\frac{n_\emptyset}{q(z'\mid z)}\}$, the factor $n_\emptyset$
being the ratio of the uniform label choice to the reverse merge (which involves no label choice). If
$z_i=a\neq z_j=b$ (merge), the same launch state is built for $C_a\cup C_b$, $t$ restricted scans are run,
the probability $q(z\mid z_{\mathrm{launch}})$ of regenerating the current state in one restricted scan is
computed, and the merge $z''$ (all of $C_b$ into $a$) is accepted with probability
$\min\{1,\ \frac{p(z''\mid Y)}{p(z\mid Y)}\cdot\frac{q(z\mid z_{\mathrm{launch}})}{n_\emptyset''}\}$ with
$n_\emptyset''$ the number of empty labels after the merge. Because the launch state and the number of
restricted scans are the same deterministic functions of $(C_a\cup C_b,i,j)$ in both directions, the
proposal densities are those of \citet{Jain2004} and the move satisfies detailed balance with respect to
the finite-label target, and the spectral launch only changes where the restricted scans start.

\subsection{Validation of the sampler}
Every closed-form marginal is tested against numerical integration. The candidate weights of
Algorithm~\ref{alg:sweep} are tested against differences of from-scratch evidences for every family,
including the type-collapsed and degree-corrected ones and the masked-dyad case. The posterior
similarity matrix of the sampler with and without split and merge moves (random or spectral launch) is
compared with the exact posterior enumerated over all $3^{12}=531{,}441$ labellings of a 12-node network
(\texttt{tests/test\_split\_merge.py}): over $1.8\cdot10^4$ retained sweeps the maximum absolute error of
the posterior similarity matrix is $0.0037$ (random launch) and $0.0047$ (spectral launch), against
$0.0055$ without split and merge moves, and the posterior probability of two occupied labels is $0.742$
against the exact $0.745$, and about $5\%$ of the proposed splits and merges are accepted.

\subsection{Correctness checks on eight planted designs}
Table~\ref{tab:correctness} runs eight strongly separated and correctly specified planted designs
(binary, counts, binary with a cap, zero-inflated counts, multiplex, Gaussian weights, dyad states, signed
ties, the multiplex design run once with three layers and once with a single layer) with the released
sampler over ten seeds, where $\mathrm{SNR}_{\mathrm{KS}}$ is the Kesten--Stigum
signal-to-noise ratio of Theorem~\ref{thm:threshold} ($\approx4.5$ for the binary design). Near-perfect recovery is
expected there and is reported only as a check of the implementation. The zero-inflated design has
$n=120$ (three groups of 40), and the capped design is the binary design with a cap of $0.02$, whose
posterior between-block densities of $0.004$ to $0.005$ never activate the cap
(Proposition~\ref{thm:cap}(i): the cap changes the evidence by at most $\log(1/0.02)=3.9$ nats here).
\begin{table}[tbp]
\centering\small
\caption{Validation of the released sampler on eight planted designs with two or three communities (ten seeds): mean ARI (sd) of the polished MAP partition against the planted partition, mean $\hat K$, and wall-clock seconds per fit. Higher ARI is better.}
\label{tab:correctness}
\resizebox{\linewidth}{!}{\begin{tabular}{lcccl}
\toprule
Design & ARI mean (sd) & $\hat K$ & seconds & note \\
\midrule
S1 Beta-Bernoulli & 0.968 (0.028) & 3.0 & 3.3 & $\mathrm{SNR}_{\rm KS}=4.45$ \\
S2 Gamma-Poisson & 1.000 (0.000) & 3.0 & 3.3 &  \\
S3 gap cap=0.02 & 0.968 (0.028) & 3.0 & 4.4 & same graphs as S1 \\
S4 ZIP (n=120) & 1.000 (0.000) & 3.0 & 0.6 & n corrected from 150 to 120 \\
S5 multiplex L=3 & 0.998 (0.006) & 3.0 & 6.6 &  \\
S5 single layer & 0.746 (0.151) & 2.8 & 3.5 &  \\
SI-B NIG Gaussian & 1.000 (0.000) & 3.0 & 3.4 &  \\
SI-C dyad-state & 0.983 (0.028) & 2.0 & 2.9 &  \\
SI-D signed & 0.975 (0.023) & 3.0 & 3.0 &  \\
\bottomrule
\end{tabular}
}
\end{table}

\section{Full tables}\label{app:tables}
\begin{figure}[t]
\centering
\includegraphics[width=0.9\linewidth]{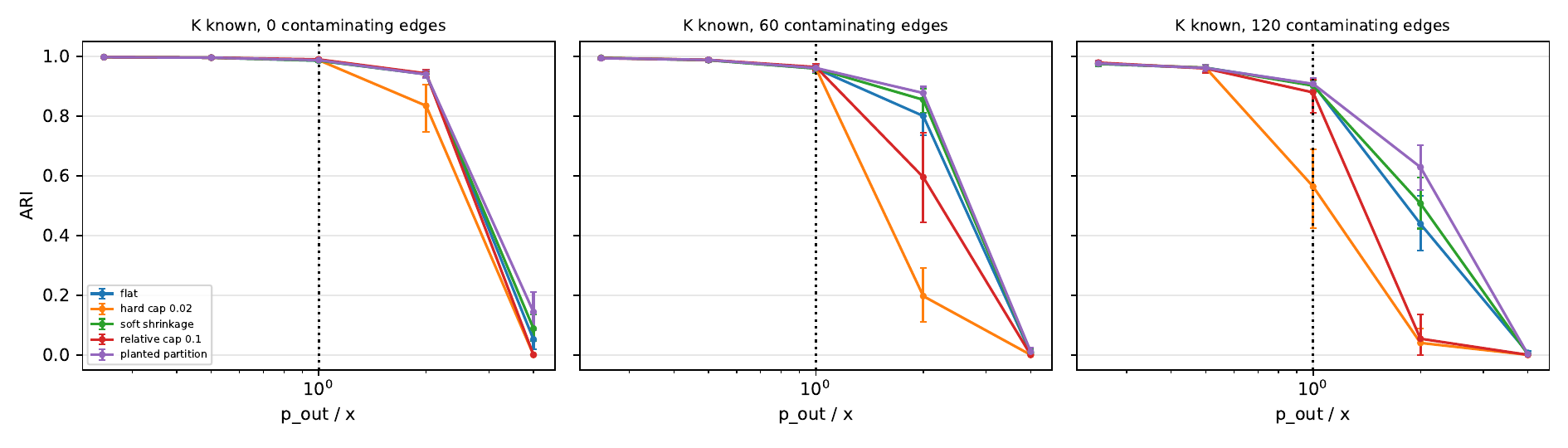}
\caption{Cap validity against contamination ($n=120$, $K=3$, $p_{\mathrm{in}}=0.2$, cap $x=0.02$,
thirty paired seeds per cell, error bars $95\%$ bootstrap intervals of the mean). The columns are $0$, $60$ and $120$ contaminating cross edges placed
uniformly at random on top of the planted network. The horizontal axis is the ratio of the planted
between density to the cap, $p_{\mathrm{out}}/x\in\{0.25,0.5,1,2,4\}$, with
$\mathrm{SNR}_{\mathrm{KS}}=7.2,6.6,5.4,3.7,1.6$, and the cap is valid to the left of the dotted line.
The curves are the ARI with $K$ known for the flat evidence without a cap, with the hard cap, with
soft shrinkage (off-diagonal Beta prior of mean $x$ and strength $50$), with the plug-in relative cap
and for the planted-partition instance, and the same five instances with a CRP prior and $K_{\max}=10$
are in Tables~\ref{tab:R16_auto_0} to~\ref{tab:R16_auto_120} of the \suppname. A higher value is
better.}
\label{fig:capexp}
\end{figure}

This section collects the numerical tables behind Section~\ref{sec:Kexp} (the overlap guarantee of
Theorem~\ref{thm:estimator} and the thermodynamic-integration values of Figure~\ref{fig:two}), the
resolution experiment against $K$ of Section~\ref{sec:Kres}, the by-degree split incentives of
Proposition~\ref{prop:degsplit}, the tables behind Figures~\ref{fig:phase} and~\ref{fig:capexp}
(Tables~\ref{tab:R15_known_0.05} to~\ref{tab:R16_auto_120}), the paired comparisons of every cell of the
mechanism-selection and constraint experiments (Tables~\ref{tab:R15pairedfull} to~\ref{tab:R16paired}),
the four reduced-data designs of Section~\ref{app:zero} (Tables~\ref{tab:R7} to~\ref{tab:R6b}), the full account
of the computation experiment of Section~\ref{sec:mixing} (Tables~\ref{tab:R11a} to~\ref{tab:R10ari}), and the
full, convergence, long-chain and initialisation tables of Section~\ref{sec:real}
(Tables~\ref{tab:real_main_full} to~\ref{tab:init_gain}) with the held-out AUC on the larger networks.

\clearpage
\subsection{Thresholds, resolution and degree priors}
\begin{table}[tbp]\centering\small
\caption{Theorem~\ref{thm:estimator} for $K=2$ and mean degree $16$ (the design of Table~\ref{tab:R17}): the exponent $F$ at the uniform overlap, which equals $\log2-J_2/2$, and the overlap guarantee $\beta^\star$ of the maximum a posteriori labelling, computed from~\eqref{eq:F} on a grid of overlap matrices (step $0.025$ in the swap and imbalance fractions, $400\times400$ points for the inner minimisation). The Kesten--Stigum threshold is $\mathrm{SNR}_{\mathrm{KS}}=1$ and the first-moment bound is $2.70$ ($s=5.4$).}
\label{tab:beta_star}
\begin{tabular}{rrrrr}\toprule
$\mathrm{SNR}_{\mathrm{KS}}$ & $s$ & $J_2$ & $F(\text{uniform})$ & $\beta^\star$\\\midrule
1 & 2 & 0.51 & $+0.44$ & --\\
2 & 4 & 1.02 & $+0.18$ & --\\
2.5 & 5 & 1.28 & $+0.05$ & --\\
2.7 & 5.4 & 1.39 & $-0.00$ & 0.52\\
3 & 6 & 1.55 & $-0.08$ & 0.62\\
4 & 8 & 2.09 & $-0.35$ & 0.76\\
5 & 10 & 2.65 & $-0.63$ & 0.82\\
5.4 & 10.8 & 2.88 & $-0.75$ & 0.84\\
6 & 12 & 3.22 & $-0.92$ & 0.86\\
8 & 16 & 4.43 & $-1.52$ & 0.91\\
10 & 20 & 5.73 & $-2.17$ & 0.93\\
\bottomrule\end{tabular}\end{table}

\begin{table}[tbp]\centering\small
\caption{Theorem~\ref{thm:two} on the design of Table~\ref{tab:R17} ($K=2$, $n=1000$, mean degree $16$, ten seeds): the log Bayes factor of the two-community model against the one-block model, $\log\mathrm{BF}_n$, with the posterior summed over all labellings by thermodynamic integration along $26$ power posteriors with parallel tempering \citep{Friel2008,Geyer1991}, and the MAP criterion $\log p(z^\star\mid A)-\log p(z_1\mid A)$ of Theorem~\ref{thm:threshold}, exact, on the same networks (mean over seeds and range). The Kesten--Stigum threshold is $\mathrm{SNR}_{\mathrm{KS}}=1$ and the first-moment bound is $2.70$.}
\label{tab:ti}
\begin{tabular}{rrrrr}\toprule
$\mathrm{SNR}_{\mathrm{KS}}$ & $s$ & $\log\mathrm{BF}_n$ (range) & MAP criterion (range) & $\tfrac n2J_2-n\log2$\\\midrule
0.5 & 1 & $-5$ ($-8$, $-2$) & $-576$ ($-595$, $-548$) & $-567$\\
1 & 2 & $-3$ ($-6$, $1$) & $-448$ ($-478$, $-412$) & $-440$\\
1.5 & 3 & $9$ ($-5$, $34$) & $-316$ ($-357$, $-269$) & $-312$\\
2 & 4 & $61$ ($31$, $93$) & $-186$ ($-231$, $-140$) & $-182$\\
2.5 & 5 & $146$ ($112$, $191$) & $-53$ ($-103$, $-10$) & $-51$\\
2.7 & 5.4 & $190$ ($142$, $238$) & $3$ ($-60$, $51$) & $2$\\
3 & 6 & $249$ ($197$, $299$) & $86$ ($27$, $134$) & $82$\\
4 & 8 & $488$ ($422$, $541$) & $357$ ($284$, $434$) & $353$\\
5 & 10 & $764$ ($681$, $848$) & $643$ ($561$, $724$) & $632$\\
\bottomrule\end{tabular}\end{table}

\begin{table}[tbp]\centering\scriptsize
\caption{Choosing $K$: planted partition with $K^\star=5$, $n=500$, mean degree $10$ (four columns) and the ring of $30$ six-cliques ($n=180$, $480$ edges, last column), 16 paired seeds: mean ARI ($95\%$ bootstrap interval) [mean $\hat K$]. Label budgets: $K_{\max}=10$ for the CRP, microcanonical and cross-validation selections and $K\in\{{2,\dots,11\}}$ for the evidence-grid and variational selections on the planted partition, $K_{\max}=60$ and $K\in\{{2,\dots,15\}}$ on the ring. Higher ARI is better. Seeds per row otherwise flat, Pitman--Yor 8, flat, Gnedin 8, flat, MFM 8, global typed, global degree prior, CRP 8, corrected BIC (Hu et al.) 8, Wang--Bickel criterion 8.}
\label{tab:R5}
\resizebox{\linewidth}{!}{\begin{tabular}{lccccc}\toprule
Method; $\mathrm{SNR}_{\mathrm{KS}}$ & 1.5 & 2.5 & 4 & 8 & ring of 30 cliques\\\midrule
flat, CRP & 0.18 (0.14, 0.22) [2.9] & \textbf{0.81 (0.80, 0.83) [5.0]} & 0.95 (0.94, 0.96) [5.0] & \textbf{1.00 (1.00, 1.00) [5.0]} & 0.38 (0.37, 0.39) [10.1]\\
flat (Jeffreys) & 0.21 (0.17, 0.26) [3.4] & \textbf{0.81 (0.79, 0.83) [5.0]} & 0.95 (0.95, 0.96) [5.0] & \textbf{1.00 (1.00, 1.00) [5.0]} & 0.47 (0.46, 0.48) [14.5]\\
flat, evidence grid & 0.08 (0.06, 0.10) [2.1] & 0.74 (0.66, 0.81) [4.6] & 0.95 (0.94, 0.96) [5.0] & \textbf{1.00 (1.00, 1.00) [5.0]} & 0.34 (0.34, 0.34) [8.0]\\
flat, Pitman--Yor & 0.14 (0.10, 0.19) [2.9] & 0.77 (0.71, 0.82) [4.9] & 0.95 (0.93, 0.96) [5.0] & \textbf{1.00 (1.00, 1.00) [5.0]} & 0.39 (0.38, 0.40) [10.5]\\
flat, Gnedin & 0.17 (0.12, 0.24) [3.6] & 0.80 (0.77, 0.82) [5.0] & 0.94 (0.93, 0.96) [5.0] & \textbf{1.00 (1.00, 1.00) [5.0]} & 0.41 (0.40, 0.42) [10.6]\\
flat, MFM & 0.17 (0.14, 0.21) [3.6] & 0.80 (0.77, 0.82) [5.0] & 0.94 (0.93, 0.96) [5.0] & \textbf{1.00 (1.00, 1.00) [5.0]} & 0.37 (0.36, 0.39) [9.9]\\
typed, CRP & 0.15 (0.11, 0.19) [2.6] & 0.79 (0.76, 0.82) [4.9] & \textbf{0.96 (0.95, 0.96) [5.0]} & \textbf{1.00 (1.00, 1.00) [5.0]} & 0.33 (0.32, 0.34) [8.3]\\
global-DC, CRP & 0.21 (0.17, 0.26) [7.4] & 0.78 (0.74, 0.81) [5.5] & 0.95 (0.94, 0.96) [5.0] & \textbf{1.00 (1.00, 1.00) [5.0]} & 0.48 (0.47, 0.49) [11.8]\\
global typed, CRP & 0.19 (0.16, 0.22) [8.3] & 0.77 (0.73, 0.81) [5.6] & 0.95 (0.94, 0.96) [5.0] & \textbf{1.00 (1.00, 1.00) [5.0]} & 0.50 (0.49, 0.52) [12.5]\\
global typed, global degree prior, CRP & 0.00 (0.00, 0.00) [1.0] & 0.77 (0.73, 0.81) [5.5] & 0.94 (0.91, 0.96) [5.1] & \textbf{1.00 (1.00, 1.00) [5.0]} & 0.47 (0.45, 0.49) [11.5]\\
microcanonical, $K$ by evidence & 0.23 (0.19, 0.27) [8.6] & 0.79 (0.75, 0.82) [5.6] & 0.95 (0.94, 0.96) [5.0] & \textbf{1.00 (1.00, 1.00) [5.0]} & \textbf{1.00 (1.00, 1.00) [30.0]}\\
variational EM with ICL & 0.13 (0.08, 0.18) [2.4] & 0.73 (0.62, 0.81) [4.5] & 0.95 (0.95, 0.96) [5.0] & \textbf{1.00 (1.00, 1.00) [5.0]} & 0.34 (0.34, 0.34) [8.0]\\
corrected BIC (Hu et al.) & 0.00 (0.00, 0.00) [1.0] & 0.65 (0.46, 0.77) [4.5] & 0.93 (0.91, 0.93) [5.0] & \textbf{1.00 (1.00, 1.00) [5.0]} & 0.34 (0.34, 0.34) [8.0]\\
Wang--Bickel criterion & \textbf{0.25 (0.20, 0.31) [3.2]} & 0.74 (0.72, 0.77) [5.0] & 0.93 (0.91, 0.93) [5.0] & \textbf{1.00 (1.00, 1.00) [5.0]} & 0.30 (0.30, 0.30) [7.0]\\
edge cross-validation & 0.00 (0.00, 0.00) [1.0] & 0.73 (0.70, 0.76) [4.9] & 0.93 (0.92, 0.94) [5.0] & \textbf{1.00 (1.00, 1.00) [5.0]} & \textbf{1.00 (0.99, 1.00) [30.8]}\\
Bethe Hessian, $K$ inferred & 0.23 (0.18, 0.28) [2.9] & 0.74 (0.71, 0.76) [4.9] & 0.93 (0.92, 0.94) [5.0] & \textbf{1.00 (1.00, 1.00) [5.0]} & \textbf{1.00 (1.00, 1.00) [30.0]}\\
Leiden & 0.16 (0.13, 0.19) [9.5] & 0.77 (0.73, 0.81) [5.4] & 0.95 (0.95, 0.96) [5.1] & \textbf{1.00 (1.00, 1.00) [5.1]} & \textbf{1.00 (1.00, 1.00) [30.0]}\\
Infomap & 0.00 (0.00, 0.00) [1.0] & 0.48 (0.43, 0.53) [16.6] & 0.94 (0.93, 0.95) [6.0] & \textbf{1.00 (1.00, 1.00) [5.5]} & \textbf{1.00 (1.00, 1.00) [30.0]}\\
\bottomrule\end{tabular}}\end{table}

\begin{table}[tbp]\centering\scriptsize
\caption{Corollary~\ref{cor:res} against $K$: planted partition with $\nu=50$ nodes per community, within density $0.1$ and between density $0.0025$ (the rates of the $K=20$ design of Section~\ref{sec:phase}), $K$ from $5$ to $40$, paired seeds (10 for $K=5$, 10 for $K=10$, 10 for $K=20$, 5 for $K=30$, 3 for $K=40$ for the flat CRP evidence, and the Pitman--Yor, Gnedin, MFM and corrected-BIC rows on the first five seeds, three at $K=40$): mean $\hat K$ ($95\%$ bootstrap interval) and mean ARI. The corollary gives $K_{\mathrm{res}}(50)=17$ under the uniform block prior and $26$ under the Jeffreys prior. Higher ARI is better.}
\label{tab:R21}
\resizebox{\linewidth}{!}{\begin{tabular}{lccccc}\toprule
Method & $K=5$ & $K=10$ & $K=20$ & $K=30$ & $K=40$\\\midrule
flat, CRP & 5.0 (5.0, 5.0) [0.97] & 10.0 (10.0, 10.0) [0.95] & 14.3 (13.9, 14.7) [0.69] & 17.0 (16.4, 17.6) [0.58] & 19.3 (19.0, 20.0) [0.52]\\
flat, CRP (Jeffreys) & 5.0 (5.0, 5.0) [0.97] & 10.0 (10.0, 10.0) [0.95] & 18.4 (18.1, 18.7) [0.85] & 21.6 (20.8, 22.0) [0.68] & 24.3 (23.0, 25.0) [0.59]\\
microcanonical, $K$ by evidence & 5.1 (5.0, 5.3) [0.96] & 10.4 (10.1, 10.7) [0.94] & 20.3 (19.7, 20.9) [0.91] & 28.6 (27.4, 29.8) [0.86] & 34.0 (31.0, 36.0) [0.76]\\
flat, Pitman--Yor & 5.0 (5.0, 5.0) [0.97] & 10.0 (10.0, 10.0) [0.94] & 14.4 (14.0, 14.8) [0.70] & 16.2 (15.6, 16.8) [0.57] & 19.3 (18.0, 20.0) [0.52]\\
flat, Gnedin & 5.0 (5.0, 5.0) [0.97] & 10.0 (10.0, 10.0) [0.95] & 14.8 (14.2, 15.4) [0.72] & 16.6 (16.2, 17.0) [0.57] & 19.0 (19.0, 19.0) [0.51]\\
flat, MFM & 5.0 (5.0, 5.0) [0.97] & 10.0 (10.0, 10.0) [0.95] & 13.8 (13.4, 14.0) [0.68] & 16.4 (16.0, 16.8) [0.57] & 18.7 (18.0, 19.0) [0.50]\\
corrected BIC & 5.0 (5.0, 5.0) [0.96] & 10.0 (10.0, 10.0) [0.93] & 17.2 (15.0, 19.4) [0.79] & 20.0 (17.8, 22.2) [0.57] & 24.0 (22.0, 26.0) [0.50]\\
Bethe Hessian, $K$ inferred & 5.0 (5.0, 5.0) [0.96] & 10.0 (10.0, 10.0) [0.93] & 20.0 (20.0, 20.0) [0.84] & 29.0 (28.2, 29.8) [0.63] & 35.0 (34.0, 36.0) [0.49]\\
Leiden & 5.8 (5.2, 6.5) [0.97] & 10.7 (10.2, 11.2) [0.96] & 20.5 (19.9, 21.2) [0.92] & 25.0 (23.6, 26.4) [0.77] & 24.3 (23.0, 25.0) [0.62]\\
flat & 5.0 (5.0, 5.0) [0.97] & 9.9 (9.7, 10.0) [0.94] & 13.4 (13.1, 13.7) [0.66] & 16.0 (16.0, 16.0) [0.57] & 17.7 (16.0, 19.0) [0.48]\\
flat (Jeffreys) & 5.0 (5.0, 5.0) [0.97] & 10.0 (10.0, 10.0) [0.95] & 17.8 (17.3, 18.3) [0.83] & 20.2 (20.0, 20.6) [0.64] & 23.0 (22.0, 24.0) [0.58]\\
\bottomrule\end{tabular}}\end{table}

\begin{table}[tbp]\centering\small
\caption{The by-degree split incentive of Proposition~\ref{prop:degsplit}: evidence of the labelling that splits one degree-corrected community (Pareto$(2.5)$ propensities, expected degree $10$) into its top quarter by degree and the rest, minus the evidence of the one-block labelling, in nats, mean of three replicates, for the microcanonical degree-corrected evidence under three degree priors and without one, for the within-degree type of the typed instance, and for the flat evidence. The column labelled prediction is $\nu[\phi(\bar d)-\pi_1\phi(\bar d_1)-\pi_2\phi(\bar d_2)]$ (script \texttt{theory\_checks/check\_degree\_split.py}).}
\label{tab:degsplit}
\resizebox{\linewidth}{!}{\begin{tabular}{rrrrrrrr}\toprule
$\nu$ & uniform prior & prediction & histogram prior & global prior & no degree prior & typed instance & flat\\\midrule
100 & 11.4 & 15.4 & 24.0 & -8.9 & -8.9 & 0.6 & 176.2\\
200 & 41.2 & 45.9 & 64.4 & -10.0 & -10.0 & -154.4 & 540.6\\
400 & 84.0 & 88.8 & 150.4 & -10.5 & -10.5 & -215.9 & 1048.8\\
800 & 287.6 & 288.8 & 349.6 & -7.1 & -7.1 & -1509.4 & 2998.5\\
1600 & 443.5 & 450.1 & 729.1 & -12.9 & -12.9 & -1297.2 & 5318.3\\
3200 & 885.2 & 890.8 & 1520.5 & -12.3 & -12.3 & -2416.7 & 10497.7\\
\bottomrule\end{tabular}}\end{table}

\clearpage
\subsection{Mechanism selection}
\begin{figure}[t]
\centering
\includegraphics[width=0.42\linewidth]{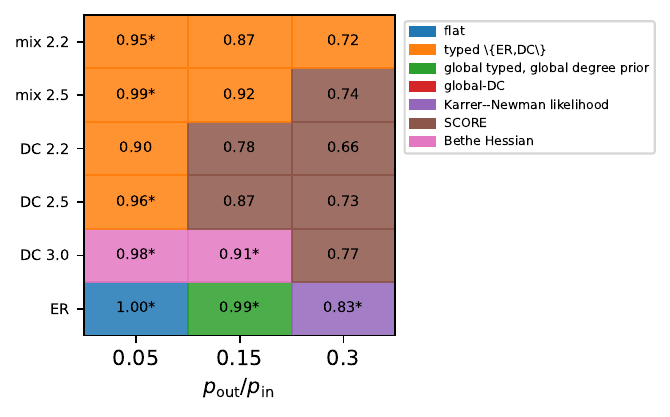}
\caption{The method with the highest mean ARI with $K$ known in each cell of Figure~\ref{fig:phase},
with its ARI, and a star when no other method is within $0.02$ of it.}
\label{fig:winner}
\end{figure}
\begin{table}[tbp]\centering\scriptsize
\caption{Mechanism selection ($n=320$, $K=4$, 20 paired seeds). Left: false-split and false-merge rates, the fraction of the 360 fits (18 cells, six designs by three ratios, times 20 seeds) with $\hat K>4$ or $\hat K<4$, with a $95\%$ bootstrap interval, for the methods that infer $K$, where lower is better. Right: type accuracy of the typed instances by design, the fraction of communities whose posterior type matches the generating mechanism, with DC denoting the degree-corrected designs, ER the homogeneous design and mix the mechanism mixtures.}
\label{tab:R15}
\noindent\resizebox{\linewidth}{!}{\begin{tabular}{lcc}\toprule Method & false split & false merge\\\midrule
flat, CRP & 0.83 (0.79, 0.87) & 0.00 (0.00, 0.00)\\
typed, CRP & 0.79 (0.75, 0.83) & 0.00 (0.00, 0.01)\\
global-DC, CRP & 0.52 (0.47, 0.57) & 0.01 (0.00, 0.02)\\
global typed, CRP & 0.84 (0.80, 0.87) & 0.00 (0.00, 0.00)\\
global typed, global degree prior, CRP & 0.64 (0.58, 0.70) & 0.00 (0.00, 0.01)\\
global typed, histogram prior, CRP & 0.74 (0.68, 0.80) & 0.00 (0.00, 0.00)\\
Bethe Hessian, $K$ inferred & 0.00 (0.00, 0.00) & 0.07 (0.04, 0.09)\\
Leiden & 0.46 (0.41, 0.51) & 0.00 (0.00, 0.00)\\
Infomap & 0.56 (0.51, 0.61) & 0.33 (0.29, 0.38)\\
\bottomrule\end{tabular}
\quad\begin{tabular}{lcccccc}\toprule type accuracy & ER & DC $\gamma$=3.0 & DC $\gamma$=2.5 & DC $\gamma$=2.2 & mix $\gamma$=2.5 & mix $\gamma$=2.2\\\midrule
typed \{ER,DC\} & 1.00 & 0.83 & 1.00 & 0.97 & 1.00 & 0.97\\
typed, CRP & 1.00 & 0.41 & 0.70 & 0.72 & 0.85 & 0.85\\
global typed & 1.00 & 0.55 & 0.81 & 0.85 & 0.95 & 0.95\\
global typed, CRP & 1.00 & 0.24 & 0.33 & 0.41 & 0.60 & 0.66\\
global typed, global degree prior & 0.78 & 1.00 & 1.00 & 0.96 & 0.98 & 0.93\\
global typed, global degree prior, CRP & 0.78 & 0.79 & 0.71 & 0.63 & 0.80 & 0.69\\
global typed, histogram prior & 1.00 & 0.70 & 0.91 & 0.90 & 0.97 & 1.00\\
global typed, histogram prior, CRP & 1.00 & 0.25 & 0.40 & 0.43 & 0.58 & 0.61\\
Yan et al.\ test, planted partition & 1.00 & 1.00 & 1.00 & 1.00 & 0.99 & 0.99\\
\bottomrule\end{tabular}}\end{table}

\begin{table}[tbp]\centering\scriptsize
\caption{Calibration of the posterior type probabilities in the mechanism-selection experiment, for the typed instances with $K$ known and inferred and for the global typed instances. Every fitted community is matched to the planted community holding the majority of its nodes, and the communities are binned by their posterior probability of the degree-corrected type. Each cell gives the number of communities in the bin and the fraction of them whose planted mechanism is degree-corrected. A calibrated posterior has that fraction inside the bin. The communities come from the seeds on which the type posteriors were stored: typed \{ER,DC\} 14 and typed, CRP 20 and global typed 20 and global typed, CRP 20 and global typed, global degree prior 20 and global typed, global degree prior, CRP 20 and global typed, histogram prior 20 and global typed, histogram prior, CRP 20 and Yan et al.\ test, planted partition 20.}
\label{tab:R15calib}
\resizebox{\linewidth}{!}{\begin{tabular}{lcccccc}\toprule
Instance & $[0,0.05)$ & $[0.05,0.25)$ & $[0.25,0.5)$ & $[0.5,0.75)$ & $[0.75,0.95)$ & $[0.95,1]$\\\midrule
typed \{ER,DC\} & 362: 0.08 & 1: 1.00 & 5: 1.00 & 6: 1.00 & -- & 633: 1.00\\
typed, CRP & 714: 0.34 & 36: 1.00 & 62: 1.00 & 177: 1.00 & 4: 1.00 & 779: 1.00\\
global typed & 684: 0.31 & 3: 1.00 & -- & 3: 1.00 & 2: 1.00 & 747: 1.00\\
global typed, CRP & 1449: 0.68 & 34: 1.00 & 16: 1.00 & 15: 1.00 & 10: 1.00 & 474: 1.00\\
global typed, global degree prior & 186: 0.04 & 81: 0.00 & 37: 0.00 & 27: 0.00 & 20: 0.00 & 641: 0.98\\
global typed, global degree prior, CRP & 480: 0.60 & 85: 0.02 & 38: 0.05 & 28: 0.04 & 22: 0.09 & 592: 0.99\\
global typed, histogram prior & 407: 0.24 & 3: 1.00 & 1: 1.00 & 2: 1.00 & 2: 1.00 & 488: 1.00\\
global typed, histogram prior, CRP & 1047: 0.70 & 21: 1.00 & 10: 1.00 & 9: 1.00 & 14: 1.00 & 370: 1.00\\
Yan et al.\ test, planted partition & 474: 0.00 & -- & -- & -- & -- & 966: 0.99\\
\bottomrule\end{tabular}}\end{table}

\begin{table}[tbp]\centering\scriptsize
\caption{Mechanism selection at scale, $K$ given: mean ARI ($95\%$ bootstrap interval) [type accuracy of the typed instances]. First three columns: the heterogeneous designs of Figure~\ref{fig:phase} at $n=3000$ ($K=4$, mean within degree $10$, $p_{\mathrm{out}}/p_{\mathrm{in}}=0.15$, 6 paired seeds). Last three columns: $K=20$ communities of $\nu$ nodes, ten Erd\H{o}s--R\'enyi and ten degree-corrected with propensities uniform on $\{0.5,1,1.5\}$, within rate $0.1$, between rate $0.1/40$ (6 paired seeds). Higher is better.}
\label{tab:R18}
\resizebox{\linewidth}{!}{\begin{tabular}{lcccccc}\toprule
Method & DC $\gamma$=2.5 & DC $\gamma$=2.2 & mix $\gamma$=2.2 & $\nu=50$ & $\nu=100$ & $\nu=250$\\\midrule
flat & 0.11 (0.08, 0.14) & 0.01 (0.01, 0.01) & 0.49 (0.37, 0.61) & 0.63 (0.59, 0.66) & 0.97 (0.97, 0.98) & \textbf{1.00 (1.00, 1.00)}\\
typed \{ER,DC\} & \textbf{0.79 (0.77, 0.80) [1.0]} & \textbf{0.74 (0.72, 0.76) [1.0]} & \textbf{0.89 (0.88, 0.90) [1.0]} & 0.59 (0.56, 0.63) [0.6] & 0.96 (0.95, 0.96) [0.6] & \textbf{1.00 (1.00, 1.00) [1.0]}\\
global-DC & 0.67 (0.65, 0.68) & 0.34 (0.26, 0.42) & 0.66 (0.64, 0.68) & 0.82 (0.77, 0.86) & \textbf{0.98 (0.97, 0.98)} & \textbf{1.00 (1.00, 1.00)}\\
global typed & 0.67 (0.65, 0.68) [1.0] & 0.25 (0.17, 0.35) [0.9] & 0.69 (0.68, 0.70) [1.0] & \textbf{0.84 (0.81, 0.86) [0.5]} & \textbf{0.98 (0.97, 0.98) [0.5]} & \textbf{1.00 (1.00, 1.00) [0.8]}\\
Karrer--Newman likelihood & 0.68 (0.66, 0.69) & 0.47 (0.44, 0.49) & 0.46 (0.38, 0.57) & 0.81 (0.77, 0.85) & \textbf{0.98 (0.97, 0.98)} & \textbf{1.00 (1.00, 1.00)}\\
SCORE & 0.72 (0.64, 0.78) & 0.73 (0.70, 0.75) & 0.37 (0.33, 0.42) & -- & -- & --\\
Bethe Hessian & 0.26 (0.22, 0.31) & 0.07 (0.04, 0.09) & 0.64 (0.54, 0.71) & 0.50 (0.47, 0.53) & 0.78 (0.73, 0.86) & \textbf{1.00 (1.00, 1.00)}\\
\bottomrule\end{tabular}}\end{table}

\begin{table}[tbp]\centering\scriptsize
\caption{Theorem~\ref{thm:split}(b) and Proposition~\ref{prop:finite} with the sampler: two degree-corrected communities of $\nu$ nodes with propensities uniform on $\{0.5,1,1.5\}$, $\omega_0=0.2$ (the design of Table~\ref{tab:splitb}), between rate $\omega_0/40$, 5 paired seeds: mean $\hat K$ ($95\%$ bootstrap interval) [mean ARI, type accuracy] for the instances that infer $K$ ($K_{\max}=8$) and mean ARI [type accuracy] with $K=2$ given. The exact evidences of Table~\ref{tab:splitb} favour the degree-corrected type over the by-atom split from $\nu\approx600$ to $800$. The flat instance carries no type.}
\label{tab:R22}
\resizebox{\linewidth}{!}{\begin{tabular}{lccc}\toprule
Method & $\nu=300$ & $\nu=600$ & $\nu=1200$\\\midrule
typed, CRP ($K_{\max}=8$) & 2.0 (2.0, 2.0) [1.00, 1.00] & 2.0 (2.0, 2.0) [1.00, 1.00] & 2.0 (2.0, 2.0) [1.00, 1.00]\\
flat, CRP ($K_{\max}=8$) & 7.4 (7.0, 7.8) [0.34] & 7.4 (7.0, 7.8) [0.33] & 8.0 (8.0, 8.0) [0.31]\\
global typed, global degree prior, CRP ($K_{\max}=8$) & 2.0 (2.0, 2.0) [1.00, 1.00] & 2.0 (2.0, 2.0) [1.00, 1.00] & 2.4 (2.0, 2.8) [0.93, 1.00]\\
typed \{ER,DC\} & 1.00 [1.00] & 1.00 [1.00] & 1.00 [1.00]\\
\bottomrule\end{tabular}}\end{table}

\begin{table}[tbp]\centering\scriptsize
\caption{Mechanism selection, $K$ known, $p_{\mathrm{out}}/p_{\mathrm{in}}=0.05$ ($n=320$, $K=4$, 20 paired seeds): mean ARI ($95\%$ bootstrap interval). Higher is better.}
\label{tab:R15_known_0.05}
\resizebox{\linewidth}{!}{\begin{tabular}{lcccccc}\toprule
Method ($K$ known) & ER & DC $\gamma$=3.0 & DC $\gamma$=2.5 & DC $\gamma$=2.2 & mix $\gamma$=2.5 & mix $\gamma$=2.2\\\midrule
flat & \textbf{1.00 (1.00, 1.00)} & 0.98 (0.97, 0.98) & 0.88 (0.80, 0.93) & 0.59 (0.48, 0.70) & 0.96 (0.92, 0.98) & 0.88 (0.83, 0.93)\\
typed \{ER,DC\} & \textbf{1.00 (1.00, 1.00)} & 0.98 (0.97, 0.98) & 0.96 (0.95, 0.97) & 0.90 (0.86, 0.93) & 0.99 (0.98, 0.99) & 0.95 (0.92, 0.97)\\
global-DC & \textbf{1.00 (1.00, 1.00)} & 0.98 (0.98, 0.98) & 0.94 (0.92, 0.95) & 0.81 (0.77, 0.85) & 0.96 (0.95, 0.98) & 0.90 (0.86, 0.93)\\
global typed & \textbf{1.00 (1.00, 1.00)} & 0.98 (0.97, 0.98) & 0.94 (0.92, 0.95) & 0.80 (0.75, 0.85) & 0.97 (0.96, 0.98) & 0.94 (0.91, 0.96)\\
global typed, global degree prior & \textbf{1.00 (1.00, 1.00)} & 0.98 (0.98, 0.98) & 0.94 (0.92, 0.95) & 0.86 (0.82, 0.89) & 0.98 (0.97, 0.99) & 0.94 (0.92, 0.96)\\
global typed, histogram prior & \textbf{1.00 (1.00, 1.00)} & 0.98 (0.98, 0.99) & 0.94 (0.92, 0.95) & 0.86 (0.82, 0.89) & 0.98 (0.97, 0.99) & 0.95 (0.93, 0.97)\\
Karrer--Newman likelihood & \textbf{1.00 (1.00, 1.00)} & 0.98 (0.97, 0.98) & 0.94 (0.92, 0.95) & 0.82 (0.78, 0.85) & 0.97 (0.95, 0.98) & 0.91 (0.88, 0.94)\\
Yan et al.\ test, planted partition & \textbf{1.00 (1.00, 1.00)} & \textbf{1.00 (1.00, 1.00)} & \textbf{1.00 (1.00, 1.00)} & \textbf{1.00 (1.00, 1.00)} & \textbf{1.00 (1.00, 1.00)} & \textbf{1.00 (1.00, 1.00)}\\
SCORE & \textbf{1.00 (0.99, 1.00)} & 0.95 (0.93, 0.96) & 0.93 (0.91, 0.95) & 0.87 (0.84, 0.90) & 0.93 (0.91, 0.95) & 0.89 (0.87, 0.90)\\
Bethe Hessian & \textbf{1.00 (1.00, 1.00)} & 0.98 (0.98, 0.99) & 0.96 (0.94, 0.97) & 0.85 (0.78, 0.90) & 0.98 (0.97, 0.99) & 0.95 (0.93, 0.96)\\
\bottomrule\end{tabular}}\end{table}
\begin{table}[tbp]\centering\scriptsize
\caption{Mechanism selection, $K$ inferred, $p_{\mathrm{out}}/p_{\mathrm{in}}=0.05$ ($n=320$, true $K=4$, 20 paired seeds): mean ARI ($95\%$ bootstrap interval) [mean $\hat K$]. Higher ARI is better.}
\label{tab:R15_auto_0.05}
\resizebox{\linewidth}{!}{\begin{tabular}{lcccccc}\toprule
Method ($K$ inferred) & ER & DC $\gamma$=3.0 & DC $\gamma$=2.5 & DC $\gamma$=2.2 & mix $\gamma$=2.5 & mix $\gamma$=2.2\\\midrule
flat, CRP & \textbf{1.00 (1.00, 1.00) [4.0]} & 0.86 (0.84, 0.87) [9.2] & 0.81 (0.79, 0.83) [9.8] & 0.73 (0.68, 0.77) [9.9] & 0.90 (0.88, 0.91) [7.7] & 0.88 (0.86, 0.89) [7.8]\\
typed, CRP & \textbf{1.00 (1.00, 1.00) [4.0]} & 0.93 (0.92, 0.94) [6.0] & 0.93 (0.92, 0.95) [5.9] & \textbf{0.87 (0.83, 0.90) [6.4]} & 0.97 (0.96, 0.98) [5.0] & 0.94 (0.91, 0.96) [5.2]\\
global-DC, CRP & \textbf{1.00 (1.00, 1.00) [4.0]} & \textbf{0.98 (0.98, 0.98) [4.0]} & 0.92 (0.89, 0.94) [4.3] & 0.77 (0.72, 0.82) [5.3] & 0.96 (0.94, 0.97) [4.3] & 0.89 (0.86, 0.92) [4.5]\\
global typed, CRP & \textbf{1.00 (1.00, 1.00) [4.0]} & 0.82 (0.80, 0.84) [7.8] & 0.78 (0.75, 0.81) [7.8] & 0.70 (0.66, 0.74) [7.8] & 0.91 (0.89, 0.92) [5.8] & 0.87 (0.85, 0.88) [6.2]\\
global typed, global degree prior, CRP & \textbf{1.00 (1.00, 1.00) [4.0]} & 0.96 (0.94, 0.97) [4.4] & 0.88 (0.84, 0.91) [5.2] & 0.78 (0.73, 0.83) [5.9] & 0.95 (0.93, 0.97) [4.6] & 0.91 (0.88, 0.93) [5.1]\\
global typed, histogram prior, CRP & \textbf{1.00 (1.00, 1.00) [4.0]} & 0.82 (0.81, 0.84) [7.8] & 0.78 (0.76, 0.80) [7.7] & 0.71 (0.68, 0.73) [8.1] & 0.91 (0.89, 0.92) [5.8] & 0.88 (0.86, 0.88) [6.0]\\
Bethe Hessian, $K$ inferred & \textbf{1.00 (1.00, 1.00) [4.0]} & \textbf{0.98 (0.98, 0.99) [4.0]} & \textbf{0.96 (0.94, 0.97) [4.0]} & 0.85 (0.78, 0.90) [4.0] & \textbf{0.98 (0.97, 0.99) [4.0]} & \textbf{0.95 (0.93, 0.96) [4.0]}\\
Leiden & \textbf{1.00 (1.00, 1.00) [4.0]} & \textbf{0.98 (0.98, 0.99) [4.2]} & 0.94 (0.92, 0.95) [4.4] & 0.81 (0.77, 0.85) [6.0] & 0.97 (0.95, 0.98) [4.3] & 0.92 (0.89, 0.94) [5.0]\\
Infomap & \textbf{1.00 (1.00, 1.00) [4.2]} & 0.97 (0.96, 0.97) [6.1] & 0.89 (0.87, 0.91) [8.9] & 0.72 (0.67, 0.76) [15.7] & 0.95 (0.93, 0.96) [7.0] & 0.88 (0.84, 0.91) [9.4]\\
\bottomrule\end{tabular}}\end{table}
\begin{table}[tbp]\centering\scriptsize
\caption{Mechanism selection, $K$ known, $p_{\mathrm{out}}/p_{\mathrm{in}}=0.15$ ($n=320$, $K=4$, 20 paired seeds): mean ARI ($95\%$ bootstrap interval). Higher is better.}
\label{tab:R15_known_0.15}
\resizebox{\linewidth}{!}{\begin{tabular}{lcccccc}\toprule
Method ($K$ known) & ER & DC $\gamma$=3.0 & DC $\gamma$=2.5 & DC $\gamma$=2.2 & mix $\gamma$=2.5 & mix $\gamma$=2.2\\\midrule
flat & 0.99 (0.98, 0.99) & 0.86 (0.83, 0.88) & 0.48 (0.38, 0.58) & 0.16 (0.12, 0.19) & 0.81 (0.77, 0.84) & 0.68 (0.63, 0.73)\\
typed \{ER,DC\} & 0.99 (0.98, 0.99) & 0.85 (0.84, 0.87) & 0.85 (0.83, 0.87) & 0.75 (0.69, 0.80) & 0.92 (0.90, 0.93) & 0.87 (0.83, 0.90)\\
global-DC & 0.99 (0.99, 0.99) & 0.86 (0.84, 0.88) & 0.71 (0.67, 0.75) & 0.50 (0.44, 0.55) & 0.82 (0.79, 0.85) & 0.72 (0.67, 0.77)\\
global typed & 0.99 (0.98, 0.99) & 0.86 (0.84, 0.87) & 0.67 (0.58, 0.73) & 0.43 (0.35, 0.51) & 0.84 (0.80, 0.86) & 0.75 (0.73, 0.78)\\
global typed, global degree prior & 0.99 (0.99, 0.99) & 0.86 (0.84, 0.88) & 0.76 (0.72, 0.80) & 0.55 (0.49, 0.61) & 0.88 (0.86, 0.90) & 0.81 (0.78, 0.85)\\
global typed, histogram prior & 0.99 (0.98, 0.99) & 0.86 (0.85, 0.88) & 0.74 (0.69, 0.78) & 0.50 (0.40, 0.60) & 0.87 (0.85, 0.90) & --\\
Karrer--Newman likelihood & 0.99 (0.99, 0.99) & 0.87 (0.85, 0.88) & 0.74 (0.71, 0.77) & 0.56 (0.50, 0.61) & 0.84 (0.81, 0.86) & 0.75 (0.71, 0.79)\\
Yan et al.\ test, planted partition & \textbf{1.00 (1.00, 1.00)} & \textbf{1.00 (1.00, 1.00)} & \textbf{1.00 (1.00, 1.00)} & \textbf{1.00 (1.00, 1.00)} & \textbf{1.00 (1.00, 1.00)} & \textbf{1.00 (1.00, 1.00)}\\
SCORE & 0.97 (0.97, 0.98) & 0.91 (0.90, 0.92) & 0.87 (0.86, 0.89) & 0.78 (0.72, 0.82) & 0.90 (0.87, 0.92) & 0.84 (0.81, 0.87)\\
Bethe Hessian & 0.98 (0.97, 0.98) & 0.91 (0.90, 0.92) & 0.83 (0.81, 0.86) & 0.58 (0.51, 0.65) & 0.88 (0.86, 0.90) & 0.80 (0.76, 0.83)\\
\bottomrule\end{tabular}}\end{table}
\begin{table}[tbp]\centering\scriptsize
\caption{Mechanism selection, $K$ inferred, $p_{\mathrm{out}}/p_{\mathrm{in}}=0.15$ ($n=320$, true $K=4$, 20 paired seeds): mean ARI ($95\%$ bootstrap interval) [mean $\hat K$]. Higher ARI is better.}
\label{tab:R15_auto_0.15}
\resizebox{\linewidth}{!}{\begin{tabular}{lcccccc}\toprule
Method ($K$ inferred) & ER & DC $\gamma$=3.0 & DC $\gamma$=2.5 & DC $\gamma$=2.2 & mix $\gamma$=2.5 & mix $\gamma$=2.2\\\midrule
flat, CRP & \textbf{0.99 (0.98, 0.99) [4.0]} & 0.78 (0.76, 0.79) [9.2] & 0.73 (0.71, 0.76) [9.8] & 0.58 (0.53, 0.63) [9.8] & 0.85 (0.83, 0.88) [7.6] & 0.82 (0.80, 0.85) [7.7]\\
typed, CRP & \textbf{0.99 (0.98, 0.99) [4.0]} & 0.83 (0.82, 0.85) [5.9] & \textbf{0.83 (0.81, 0.85) [6.1]} & \textbf{0.77 (0.74, 0.80) [6.7]} & \textbf{0.91 (0.90, 0.93) [5.2]} & \textbf{0.90 (0.89, 0.92) [5.3]}\\
global-DC, CRP & \textbf{0.99 (0.98, 0.99) [4.0]} & 0.85 (0.83, 0.87) [4.5] & 0.66 (0.60, 0.71) [4.9] & 0.49 (0.44, 0.54) [5.8] & 0.82 (0.79, 0.85) [4.5] & 0.75 (0.72, 0.78) [4.8]\\
global typed, CRP & \textbf{0.99 (0.98, 0.99) [4.0]} & 0.75 (0.74, 0.76) [8.0] & 0.67 (0.64, 0.69) [8.0] & 0.53 (0.49, 0.57) [8.0] & 0.85 (0.83, 0.87) [6.0] & 0.79 (0.76, 0.81) [6.0]\\
global typed, global degree prior, CRP & \textbf{0.99 (0.99, 0.99) [4.0]} & 0.82 (0.80, 0.85) [5.0] & 0.67 (0.62, 0.73) [5.6] & 0.51 (0.46, 0.56) [6.2] & 0.85 (0.82, 0.88) [5.1] & 0.79 (0.76, 0.83) [5.7]\\
global typed, histogram prior, CRP & \textbf{0.99 (0.98, 0.99) [4.0]} & 0.76 (0.74, 0.77) [7.8] & 0.69 (0.67, 0.70) [8.1] & 0.59 (0.54, 0.62) [8.1] & 0.84 (0.82, 0.87) [6.2] & --\\
Bethe Hessian, $K$ inferred & 0.98 (0.97, 0.98) [4.0] & \textbf{0.91 (0.90, 0.92) [4.0]} & \textbf{0.83 (0.81, 0.86) [4.0]} & 0.55 (0.47, 0.63) [3.9] & 0.88 (0.86, 0.90) [4.0] & 0.77 (0.72, 0.83) [3.9]\\
Leiden & \textbf{0.99 (0.99, 0.99) [4.0]} & 0.87 (0.85, 0.88) [4.2] & 0.72 (0.69, 0.76) [4.6] & 0.54 (0.49, 0.59) [5.7] & 0.84 (0.81, 0.87) [4.2] & 0.75 (0.70, 0.79) [4.3]\\
Infomap & \textbf{0.99 (0.99, 0.99) [4.0]} & 0.86 (0.84, 0.88) [5.5] & 0.69 (0.65, 0.73) [7.7] & 0.51 (0.47, 0.55) [11.8] & 0.82 (0.79, 0.85) [5.3] & 0.72 (0.67, 0.77) [6.5]\\
\bottomrule\end{tabular}}\end{table}
\begin{table}[tbp]\centering\scriptsize
\caption{Mechanism selection, $K$ known, $p_{\mathrm{out}}/p_{\mathrm{in}}=0.3$ ($n=320$, $K=4$, 20 paired seeds): mean ARI ($95\%$ bootstrap interval). Higher is better.}
\label{tab:R15_known_0.3}
\resizebox{\linewidth}{!}{\begin{tabular}{lcccccc}\toprule
Method ($K$ known) & ER & DC $\gamma$=3.0 & DC $\gamma$=2.5 & DC $\gamma$=2.2 & mix $\gamma$=2.5 & mix $\gamma$=2.2\\\midrule
flat & 0.83 (0.80, 0.85) & 0.25 (0.18, 0.31) & 0.09 (0.07, 0.10) & 0.05 (0.04, 0.07) & 0.43 (0.36, 0.50) & 0.42 (0.35, 0.47)\\
typed \{ER,DC\} & 0.78 (0.73, 0.82) & 0.53 (0.50, 0.56) & 0.59 (0.56, 0.62) & 0.60 (0.56, 0.63) & 0.72 (0.70, 0.75) & 0.72 (0.69, 0.75)\\
global-DC & 0.78 (0.73, 0.82) & 0.56 (0.53, 0.59) & 0.34 (0.28, 0.40) & 0.21 (0.18, 0.26) & 0.51 (0.45, 0.56) & 0.41 (0.36, 0.46)\\
global typed & 0.82 (0.78, 0.85) & 0.31 (0.24, 0.39) & 0.12 (0.10, 0.14) & 0.09 (0.08, 0.10) & 0.42 (0.38, 0.47) & 0.42 (0.38, 0.46)\\
global typed, global degree prior & 0.81 (0.76, 0.85) & 0.57 (0.53, 0.60) & 0.38 (0.30, 0.45) & 0.20 (0.13, 0.28) & 0.53 (0.47, 0.59) & 0.49 (0.45, 0.54)\\
global typed, histogram prior & 0.82 (0.78, 0.85) & 0.35 (0.27, 0.43) & 0.17 (0.12, 0.24) & 0.12 (0.09, 0.16) & 0.46 (0.39, 0.52) & --\\
Karrer--Newman likelihood & 0.83 (0.81, 0.85) & 0.61 (0.59, 0.63) & 0.47 (0.44, 0.50) & 0.32 (0.29, 0.36) & 0.61 (0.58, 0.63) & 0.52 (0.48, 0.55)\\
Yan et al.\ test, planted partition & \textbf{1.00 (1.00, 1.00)} & \textbf{1.00 (1.00, 1.00)} & \textbf{1.00 (1.00, 1.00)} & \textbf{1.00 (1.00, 1.00)} & \textbf{1.00 (1.00, 1.00)} & \textbf{1.00 (1.00, 1.00)}\\
SCORE & 0.75 (0.72, 0.77) & 0.77 (0.75, 0.79) & 0.73 (0.71, 0.76) & 0.66 (0.62, 0.69) & 0.74 (0.73, 0.76) & 0.69 (0.66, 0.72)\\
Bethe Hessian & 0.77 (0.75, 0.79) & 0.71 (0.69, 0.73) & 0.57 (0.51, 0.62) & 0.38 (0.32, 0.44) & 0.62 (0.59, 0.65) & 0.51 (0.45, 0.56)\\
\bottomrule\end{tabular}}\end{table}
\begin{table}[tbp]\centering\scriptsize
\caption{Mechanism selection, $K$ inferred, $p_{\mathrm{out}}/p_{\mathrm{in}}=0.3$ ($n=320$, true $K=4$, 20 paired seeds): mean ARI ($95\%$ bootstrap interval) [mean $\hat K$]. Higher ARI is better.}
\label{tab:R15_auto_0.3}
\resizebox{\linewidth}{!}{\begin{tabular}{lcccccc}\toprule
Method ($K$ inferred) & ER & DC $\gamma$=3.0 & DC $\gamma$=2.5 & DC $\gamma$=2.2 & mix $\gamma$=2.5 & mix $\gamma$=2.2\\\midrule
flat, CRP & \textbf{0.83 (0.82, 0.85) [4.0]} & 0.58 (0.53, 0.62) [9.1] & 0.52 (0.45, 0.59) [9.8] & 0.43 (0.36, 0.50) [9.7] & 0.72 (0.68, 0.75) [7.8] & 0.69 (0.66, 0.72) [7.7]\\
typed, CRP & 0.81 (0.78, 0.84) [4.0] & 0.53 (0.50, 0.57) [5.8] & \textbf{0.56 (0.53, 0.59) [6.0]} & \textbf{0.60 (0.55, 0.65) [6.7]} & \textbf{0.74 (0.71, 0.77) [5.0]} & \textbf{0.74 (0.70, 0.77) [5.4]}\\
global-DC, CRP & 0.81 (0.78, 0.84) [4.0] & 0.50 (0.46, 0.54) [4.9] & 0.34 (0.27, 0.40) [5.4] & 0.28 (0.24, 0.33) [5.7] & 0.59 (0.56, 0.62) [4.8] & 0.49 (0.44, 0.53) [4.7]\\
global typed, CRP & \textbf{0.83 (0.81, 0.85) [4.0]} & 0.56 (0.53, 0.58) [7.6] & 0.47 (0.41, 0.51) [8.0] & 0.40 (0.36, 0.44) [7.8] & 0.67 (0.65, 0.69) [6.0] & 0.64 (0.60, 0.67) [5.8]\\
global typed, global degree prior, CRP & 0.80 (0.76, 0.84) [4.1] & 0.46 (0.42, 0.50) [5.3] & 0.32 (0.26, 0.38) [5.8] & 0.30 (0.25, 0.35) [6.2] & 0.62 (0.59, 0.65) [5.1] & 0.60 (0.56, 0.65) [5.2]\\
global typed, histogram prior, CRP & \textbf{0.83 (0.81, 0.85) [4.0]} & 0.57 (0.54, 0.59) [7.5] & 0.44 (0.35, 0.52) [7.8] & 0.38 (0.30, 0.46) [7.6] & 0.64 (0.59, 0.68) [6.0] & --\\
Bethe Hessian, $K$ inferred & 0.75 (0.70, 0.78) [3.9] & \textbf{0.71 (0.69, 0.73) [4.0]} & 0.55 (0.48, 0.61) [3.9] & 0.36 (0.29, 0.41) [3.6] & 0.60 (0.55, 0.64) [3.9] & 0.50 (0.43, 0.56) [3.6]\\
Leiden & 0.79 (0.74, 0.84) [4.3] & 0.58 (0.55, 0.60) [5.1] & 0.42 (0.39, 0.45) [5.6] & 0.29 (0.26, 0.32) [6.0] & 0.51 (0.45, 0.57) [4.8] & 0.43 (0.37, 0.49) [5.5]\\
Infomap & 0.00 (0.00, 0.00) [1.0] & 0.00 (0.00, 0.00) [1.0] & 0.00 (0.00, 0.00) [1.0] & 0.00 (0.00, 0.00) [1.0] & 0.00 (0.00, 0.00) [1.0] & 0.00 (0.00, 0.00) [1.0]\\
\bottomrule\end{tabular}}\end{table}

\begin{table}[tbp]\centering\scriptsize
\caption{Paired differences in ARI, typed evidence minus the named method, $K$ known, in all 18 cells of the mechanism-selection experiment (20 paired seeds per cell): mean difference ($95\%$ bootstrap interval), a star marking a Wilcoxon signed-rank test significant at level $0.05$ after Holm adjustment over the 108 tests of this table. A positive value favours the typed evidence.}
\label{tab:R15pairedfull}
\resizebox{\linewidth}{!}{\begin{tabular}{llcccccc}\toprule
Design & $p_{\mathrm{out}}/p_{\mathrm{in}}$ & $-$ flat & $-$ global-DC & $-$ global typed & $-$ Karrer--Newman likelihood & $-$ SCORE & $-$ Bethe Hessian\\\midrule
ER & 0.05 & +0.00 (+0.00, +0.00) & +0.00 (-0.00, +0.00) & +0.00 (+0.00, +0.00) & +0.00 (-0.00, +0.00) & +0.00 (+0.00, +0.01) & +0.00 (-0.00, +0.00)\\
DC $\gamma$=3.0 & 0.05 & -0.00 (-0.01, +0.01) & -0.00 (-0.02, +0.00) & -0.00 (-0.01, +0.01) & -0.00 (-0.01, +0.01) & +0.03 (+0.01, +0.04) & -0.01 (-0.02, -0.00)\\
DC $\gamma$=2.5 & 0.05 & +0.09 (+0.03, +0.15)$^{*}$ & +0.03 (+0.02, +0.03)$^{*}$ & +0.02 (+0.02, +0.03)$^{*}$ & +0.03 (+0.02, +0.03)$^{*}$ & +0.03 (+0.02, +0.05) & +0.01 (-0.00, +0.02)\\
DC $\gamma$=2.2 & 0.05 & +0.31 (+0.21, +0.40)$^{*}$ & +0.08 (+0.06, +0.10)$^{*}$ & +0.09 (+0.07, +0.12)$^{*}$ & +0.08 (+0.06, +0.10)$^{*}$ & +0.02 (-0.00, +0.04) & +0.05 (+0.01, +0.10)\\
mix $\gamma$=2.5 & 0.05 & +0.03 (+0.01, +0.06) & +0.02 (+0.01, +0.03) & +0.02 (+0.01, +0.02) & +0.02 (+0.01, +0.03) & +0.05 (+0.04, +0.07)$^{*}$ & +0.01 (-0.00, +0.01)\\
mix $\gamma$=2.2 & 0.05 & +0.07 (+0.03, +0.12)$^{*}$ & +0.06 (+0.03, +0.08)$^{*}$ & +0.02 (+0.01, +0.03) & +0.04 (+0.02, +0.06)$^{*}$ & +0.07 (+0.03, +0.09)$^{*}$ & +0.01 (-0.02, +0.03)\\
ER & 0.15 & -0.00 (-0.00, +0.00) & -0.00 (-0.00, +0.00) & +0.00 (+0.00, +0.00) & +0.00 (-0.00, +0.00) & +0.02 (+0.01, +0.02)$^{*}$ & +0.01 (+0.01, +0.02)$^{*}$\\
DC $\gamma$=3.0 & 0.15 & -0.00 (-0.03, +0.03) & -0.01 (-0.03, +0.02) & -0.00 (-0.02, +0.02) & -0.01 (-0.03, +0.01) & -0.05 (-0.07, -0.04)$^{*}$ & -0.05 (-0.07, -0.04)$^{*}$\\
DC $\gamma$=2.5 & 0.15 & +0.37 (+0.27, +0.47)$^{*}$ & +0.14 (+0.10, +0.18)$^{*}$ & +0.18 (+0.12, +0.26)$^{*}$ & +0.11 (+0.08, +0.14)$^{*}$ & -0.02 (-0.04, -0.01) & +0.02 (-0.00, +0.04)\\
DC $\gamma$=2.2 & 0.15 & +0.59 (+0.53, +0.65)$^{*}$ & +0.25 (+0.21, +0.30)$^{*}$ & +0.32 (+0.24, +0.39)$^{*}$ & +0.19 (+0.15, +0.23)$^{*}$ & -0.03 (-0.07, +0.01) & +0.17 (+0.11, +0.22)$^{*}$\\
mix $\gamma$=2.5 & 0.15 & +0.11 (+0.07, +0.15)$^{*}$ & +0.09 (+0.06, +0.13)$^{*}$ & +0.08 (+0.05, +0.12)$^{*}$ & +0.08 (+0.05, +0.11)$^{*}$ & +0.02 (+0.00, +0.04) & +0.04 (+0.02, +0.06)$^{*}$\\
mix $\gamma$=2.2 & 0.15 & +0.19 (+0.13, +0.24)$^{*}$ & +0.15 (+0.11, +0.19)$^{*}$ & +0.11 (+0.08, +0.15)$^{*}$ & +0.12 (+0.08, +0.15)$^{*}$ & +0.03 (-0.02, +0.07) & +0.07 (+0.04, +0.11)$^{*}$\\
ER & 0.3 & -0.04 (-0.09, -0.01) & -0.00 (-0.01, +0.01) & -0.04 (-0.08, -0.01) & -0.05 (-0.09, -0.01) & +0.04 (-0.00, +0.07) & +0.01 (-0.04, +0.05)\\
DC $\gamma$=3.0 & 0.3 & +0.29 (+0.21, +0.35)$^{*}$ & -0.03 (-0.07, +0.01) & +0.22 (+0.13, +0.29)$^{*}$ & -0.08 (-0.11, -0.04)$^{*}$ & -0.23 (-0.26, -0.21)$^{*}$ & -0.18 (-0.21, -0.15)$^{*}$\\
DC $\gamma$=2.5 & 0.3 & +0.50 (+0.48, +0.53)$^{*}$ & +0.25 (+0.18, +0.31)$^{*}$ & +0.47 (+0.44, +0.50)$^{*}$ & +0.12 (+0.07, +0.17)$^{*}$ & -0.15 (-0.17, -0.12)$^{*}$ & +0.02 (-0.03, +0.08)\\
DC $\gamma$=2.2 & 0.3 & +0.55 (+0.51, +0.58)$^{*}$ & +0.38 (+0.33, +0.43)$^{*}$ & +0.51 (+0.47, +0.54)$^{*}$ & +0.27 (+0.24, +0.31)$^{*}$ & -0.06 (-0.10, -0.02) & +0.22 (+0.17, +0.27)$^{*}$\\
mix $\gamma$=2.5 & 0.3 & +0.30 (+0.23, +0.36)$^{*}$ & +0.21 (+0.14, +0.30)$^{*}$ & +0.30 (+0.25, +0.35)$^{*}$ & +0.12 (+0.07, +0.16)$^{*}$ & -0.02 (-0.04, +0.00) & +0.10 (+0.06, +0.15)$^{*}$\\
mix $\gamma$=2.2 & 0.3 & +0.31 (+0.23, +0.38)$^{*}$ & +0.31 (+0.26, +0.37)$^{*}$ & +0.30 (+0.24, +0.35)$^{*}$ & +0.20 (+0.16, +0.25)$^{*}$ & +0.03 (-0.01, +0.07) & +0.21 (+0.15, +0.28)$^{*}$\\
\bottomrule\end{tabular}}\end{table}

\begin{table}[tbp]\centering\scriptsize
\caption{Paired differences in ARI, typed CRP evidence minus the named method, $K$ inferred, in all 18 cells of the mechanism-selection experiment (20 paired seeds per cell): mean difference ($95\%$ bootstrap interval), a star marking a Wilcoxon signed-rank test significant at level $0.05$ after Holm adjustment over the 90 tests of this table. A positive value favours the typed evidence.}
\label{tab:R15pairedauto}
\resizebox{\linewidth}{!}{\begin{tabular}{llccccc}\toprule
Design & $p_{\mathrm{out}}/p_{\mathrm{in}}$ & $-$ flat, CRP & $-$ global-DC, CRP & $-$ global typed, CRP & $-$ Bethe Hessian, $K$ inferred & $-$ Leiden\\\midrule
ER & 0.05 & +0.00 (+0.00, +0.00) & +0.00 (-0.00, +0.00) & +0.00 (+0.00, +0.00) & +0.00 (-0.00, +0.00) & +0.00 (+0.00, +0.00)\\
DC $\gamma$=3.0 & 0.05 & +0.07 (+0.05, +0.09)$^{*}$ & -0.05 (-0.07, -0.03)$^{*}$ & +0.11 (+0.08, +0.13)$^{*}$ & -0.05 (-0.07, -0.04)$^{*}$ & -0.05 (-0.07, -0.04)$^{*}$\\
DC $\gamma$=2.5 & 0.05 & +0.13 (+0.10, +0.15)$^{*}$ & +0.02 (-0.01, +0.04) & +0.15 (+0.13, +0.18)$^{*}$ & -0.02 (-0.04, -0.00) & -0.00 (-0.02, +0.01)\\
DC $\gamma$=2.2 & 0.05 & +0.14 (+0.10, +0.19)$^{*}$ & +0.09 (+0.07, +0.12)$^{*}$ & +0.17 (+0.15, +0.20)$^{*}$ & +0.02 (-0.02, +0.08) & +0.05 (+0.03, +0.08)$^{*}$\\
mix $\gamma$=2.5 & 0.05 & +0.07 (+0.06, +0.09)$^{*}$ & +0.01 (-0.01, +0.03) & +0.06 (+0.04, +0.08)$^{*}$ & -0.01 (-0.02, -0.00) & +0.00 (-0.01, +0.02)\\
mix $\gamma$=2.2 & 0.05 & +0.06 (+0.03, +0.09)$^{*}$ & +0.05 (+0.03, +0.07)$^{*}$ & +0.07 (+0.04, +0.10)$^{*}$ & -0.00 (-0.03, +0.02) & +0.03 (+0.01, +0.05)\\
ER & 0.15 & -0.00 (-0.00, +0.00) & +0.00 (-0.00, +0.00) & +0.00 (+0.00, +0.00) & +0.01 (+0.01, +0.02)$^{*}$ & -0.00 (-0.00, +0.00)\\
DC $\gamma$=3.0 & 0.15 & +0.06 (+0.04, +0.07)$^{*}$ & -0.01 (-0.04, +0.01) & +0.08 (+0.07, +0.10)$^{*}$ & -0.08 (-0.09, -0.06)$^{*}$ & -0.03 (-0.05, -0.01)\\
DC $\gamma$=2.5 & 0.15 & +0.10 (+0.08, +0.12)$^{*}$ & +0.17 (+0.12, +0.22)$^{*}$ & +0.16 (+0.13, +0.19)$^{*}$ & -0.00 (-0.03, +0.03) & +0.11 (+0.07, +0.15)$^{*}$\\
DC $\gamma$=2.2 & 0.15 & +0.19 (+0.14, +0.24)$^{*}$ & +0.28 (+0.23, +0.34)$^{*}$ & +0.24 (+0.19, +0.28)$^{*}$ & +0.22 (+0.15, +0.29)$^{*}$ & +0.23 (+0.18, +0.28)$^{*}$\\
mix $\gamma$=2.5 & 0.15 & +0.06 (+0.04, +0.07)$^{*}$ & +0.09 (+0.05, +0.13)$^{*}$ & +0.06 (+0.05, +0.08)$^{*}$ & +0.03 (+0.01, +0.05) & +0.07 (+0.04, +0.11)$^{*}$\\
mix $\gamma$=2.2 & 0.15 & +0.08 (+0.06, +0.10)$^{*}$ & +0.15 (+0.11, +0.19)$^{*}$ & +0.11 (+0.09, +0.14)$^{*}$ & +0.13 (+0.08, +0.19)$^{*}$ & +0.16 (+0.11, +0.20)$^{*}$\\
ER & 0.3 & -0.02 (-0.05, +0.00) & +0.00 (-0.01, +0.01) & -0.01 (-0.04, +0.01) & +0.06 (+0.02, +0.12) & +0.02 (-0.01, +0.07)\\
DC $\gamma$=3.0 & 0.3 & -0.05 (-0.09, +0.00) & +0.04 (-0.01, +0.08) & -0.03 (-0.06, +0.01) & -0.18 (-0.21, -0.14)$^{*}$ & -0.04 (-0.09, -0.00)\\
DC $\gamma$=2.5 & 0.3 & +0.04 (-0.02, +0.10) & +0.23 (+0.16, +0.29)$^{*}$ & +0.10 (+0.04, +0.15) & +0.01 (-0.05, +0.09) & +0.14 (+0.09, +0.19)$^{*}$\\
DC $\gamma$=2.2 & 0.3 & +0.17 (+0.11, +0.24)$^{*}$ & +0.32 (+0.27, +0.36)$^{*}$ & +0.20 (+0.15, +0.25)$^{*}$ & +0.25 (+0.18, +0.30)$^{*}$ & +0.31 (+0.26, +0.35)$^{*}$\\
mix $\gamma$=2.5 & 0.3 & +0.02 (-0.01, +0.06) & +0.15 (+0.11, +0.20)$^{*}$ & +0.07 (+0.03, +0.11)$^{*}$ & +0.14 (+0.08, +0.21)$^{*}$ & +0.23 (+0.15, +0.31)$^{*}$\\
mix $\gamma$=2.2 & 0.3 & +0.05 (+0.01, +0.09) & +0.25 (+0.19, +0.31)$^{*}$ & +0.10 (+0.06, +0.14)$^{*}$ & +0.24 (+0.16, +0.32)$^{*}$ & +0.31 (+0.24, +0.38)$^{*}$\\
\bottomrule\end{tabular}}\end{table}

\clearpage
\subsection{Constraints between communities}
\begin{table}[tbp]\centering\scriptsize
\caption{Constraints between communities, 0 contaminating cross edges, $K$ known ($n=120$, $K=3$, $p_{\mathrm{in}}=0.2$, cap $x=0.02$, 30 paired seeds): mean ARI ($95\%$ bootstrap interval). Columns are the ratio of the planted between density to the cap. Higher is better.}
\label{tab:R16_known_0}
\resizebox{\linewidth}{!}{\begin{tabular}{lccccc}\toprule
$K$ known; $p_{\mathrm{out}}/x$ & 0.25 & 0.5 & 1 & 2 & 4\\\midrule
flat & \textbf{1.00 (0.99, 1.00)} & \textbf{1.00 (0.99, 1.00)} & \textbf{0.99 (0.98, 0.99)} & \textbf{0.94 (0.93, 0.95)} & 0.05 (0.02, 0.09)\\
hard cap 0.02 & \textbf{1.00 (0.99, 1.00)} & \textbf{1.00 (0.99, 1.00)} & \textbf{0.99 (0.98, 0.99)} & 0.83 (0.75, 0.91) & 0.00 (0.00, 0.00)\\
soft shrinkage & \textbf{1.00 (0.99, 1.00)} & \textbf{1.00 (0.99, 1.00)} & \textbf{0.99 (0.98, 0.99)} & \textbf{0.94 (0.93, 0.95)} & 0.09 (0.05, 0.13)\\
relative cap 0.1 & \textbf{1.00 (0.99, 1.00)} & \textbf{1.00 (0.99, 1.00)} & \textbf{0.99 (0.98, 0.99)} & \textbf{0.94 (0.93, 0.96)} & 0.00 (0.00, 0.00)\\
planted partition & \textbf{1.00 (0.99, 1.00)} & \textbf{1.00 (0.99, 1.00)} & \textbf{0.99 (0.98, 0.99)} & \textbf{0.94 (0.93, 0.95)} & 0.14 (0.08, 0.21)\\
variational EM & \textbf{1.00 (0.99, 1.00)} & \textbf{1.00 (0.99, 1.00)} & \textbf{0.99 (0.98, 0.99)} & \textbf{0.94 (0.92, 0.95)} & 0.31 (0.22, 0.40)\\
Bethe Hessian & \textbf{1.00 (0.99, 1.00)} & 0.99 (0.99, 1.00) & 0.98 (0.97, 0.99) & 0.91 (0.89, 0.93) & \textbf{0.43 (0.38, 0.49)}\\
Leiden & \textbf{1.00 (0.99, 1.00)} & \textbf{1.00 (0.99, 1.00)} & \textbf{0.99 (0.98, 0.99)} & 0.92 (0.89, 0.94) & 0.25 (0.22, 0.28)\\
\bottomrule\end{tabular}}\end{table}
\begin{table}[tbp]\centering\scriptsize
\caption{Constraints between communities, 0 contaminating cross edges, $K$ inferred with a CRP prior and $10$ labels (true $K=3$, 30 paired seeds): mean ARI ($95\%$ bootstrap interval) [mean $\hat K$]. Higher ARI is better.}
\label{tab:R16_auto_0}
\resizebox{\linewidth}{!}{\begin{tabular}{lccccc}\toprule
$K$ inferred; $p_{\mathrm{out}}/x$ & 0.25 & 0.5 & 1 & 2 & 4\\\midrule
flat, CRP & \textbf{1.00 (0.99, 1.00) [3.0]} & \textbf{1.00 (0.99, 1.00) [3.0]} & \textbf{0.99 (0.98, 0.99) [3.0]} & \textbf{0.94 (0.93, 0.95) [3.0]} & 0.05 (0.02, 0.08) [1.4]\\
hard cap 0.02, CRP & \textbf{1.00 (0.99, 1.00) [3.1]} & \textbf{1.00 (0.99, 1.00) [3.1]} & \textbf{0.99 (0.98, 0.99) [3.1]} & 0.82 (0.73, 0.90) [2.7] & 0.00 (0.00, 0.00) [1.0]\\
soft shrinkage, CRP & \textbf{1.00 (0.99, 1.00) [3.2]} & \textbf{1.00 (0.99, 1.00) [3.2]} & 0.98 (0.98, 0.99) [3.1] & \textbf{0.94 (0.93, 0.95) [3.0]} & 0.09 (0.05, 0.14) [1.5]\\
relative cap 0.1, CRP & \textbf{1.00 (0.99, 1.00) [3.1]} & \textbf{1.00 (0.99, 1.00) [3.1]} & \textbf{0.99 (0.99, 0.99) [3.1]} & 0.93 (0.92, 0.95) [3.1] & 0.00 (0.00, 0.00) [1.0]\\
planted partition, CRP & \textbf{1.00 (0.99, 1.00) [3.1]} & \textbf{1.00 (0.99, 1.00) [3.1]} & \textbf{0.99 (0.98, 0.99) [3.1]} & \textbf{0.94 (0.93, 0.95) [3.2]} & \textbf{0.21 (0.14, 0.28) [2.5]}\\
\bottomrule\end{tabular}}\end{table}
\begin{table}[tbp]\centering\scriptsize
\caption{Constraints between communities, 60 contaminating cross edges, $K$ known ($n=120$, $K=3$, $p_{\mathrm{in}}=0.2$, cap $x=0.02$, 30 paired seeds): mean ARI ($95\%$ bootstrap interval). Columns are the ratio of the planted between density to the cap. Higher is better.}
\label{tab:R16_known_60}
\resizebox{\linewidth}{!}{\begin{tabular}{lccccc}\toprule
$K$ known; $p_{\mathrm{out}}/x$ & 0.25 & 0.5 & 1 & 2 & 4\\\midrule
flat & 0.99 (0.99, 1.00) & \textbf{0.99 (0.98, 0.99)} & \textbf{0.96 (0.94, 0.97)} & 0.80 (0.74, 0.86) & 0.01 (0.00, 0.02)\\
hard cap 0.02 & 0.99 (0.99, 1.00) & \textbf{0.99 (0.98, 0.99)} & \textbf{0.96 (0.95, 0.97)} & 0.20 (0.11, 0.29) & 0.00 (0.00, 0.00)\\
soft shrinkage & 0.99 (0.99, 1.00) & \textbf{0.99 (0.98, 0.99)} & \textbf{0.96 (0.94, 0.97)} & 0.85 (0.81, 0.89) & 0.01 (0.00, 0.02)\\
relative cap 0.1 & 0.99 (0.99, 1.00) & \textbf{0.99 (0.98, 0.99)} & \textbf{0.96 (0.95, 0.97)} & 0.60 (0.44, 0.74) & 0.00 (0.00, 0.00)\\
planted partition & 0.99 (0.99, 1.00) & \textbf{0.99 (0.98, 0.99)} & \textbf{0.96 (0.95, 0.97)} & \textbf{0.88 (0.86, 0.90)} & 0.01 (0.00, 0.02)\\
variational EM & 0.99 (0.99, 1.00) & \textbf{0.99 (0.98, 0.99)} & \textbf{0.96 (0.94, 0.97)} & \textbf{0.88 (0.86, 0.90)} & 0.07 (0.03, 0.12)\\
Bethe Hessian & 0.99 (0.98, 0.99) & 0.98 (0.97, 0.99) & 0.95 (0.93, 0.96) & 0.83 (0.81, 0.86) & \textbf{0.27 (0.22, 0.32)}\\
Leiden & \textbf{1.00 (0.99, 1.00)} & \textbf{0.99 (0.99, 1.00)} & \textbf{0.96 (0.94, 0.97)} & 0.73 (0.67, 0.79) & 0.14 (0.12, 0.17)\\
\bottomrule\end{tabular}}\end{table}
\begin{table}[tbp]\centering\scriptsize
\caption{Constraints between communities, 60 contaminating cross edges, $K$ inferred with a CRP prior and $10$ labels (true $K=3$, 30 paired seeds): mean ARI ($95\%$ bootstrap interval) [mean $\hat K$]. Higher ARI is better.}
\label{tab:R16_auto_60}
\resizebox{\linewidth}{!}{\begin{tabular}{lccccc}\toprule
$K$ inferred; $p_{\mathrm{out}}/x$ & 0.25 & 0.5 & 1 & 2 & 4\\\midrule
flat, CRP & \textbf{0.99 (0.99, 1.00) [3.0]} & \textbf{0.99 (0.98, 0.99) [3.0]} & 0.96 (0.94, 0.97) [3.0] & 0.81 (0.75, 0.87) [2.9] & 0.01 (0.00, 0.02) [1.1]\\
hard cap 0.02, CRP & \textbf{0.99 (0.99, 1.00) [3.1]} & \textbf{0.99 (0.98, 0.99) [3.1]} & 0.96 (0.95, 0.97) [3.1] & 0.12 (0.04, 0.21) [1.3] & 0.00 (0.00, 0.00) [1.0]\\
soft shrinkage, CRP & \textbf{0.99 (0.99, 1.00) [3.2]} & \textbf{0.99 (0.98, 0.99) [3.1]} & 0.96 (0.94, 0.97) [3.2] & \textbf{0.85 (0.80, 0.89) [3.1]} & 0.01 (0.00, 0.02) [1.2]\\
relative cap 0.1, CRP & \textbf{0.99 (0.99, 1.00) [3.2]} & \textbf{0.99 (0.98, 0.99) [3.1]} & \textbf{0.97 (0.96, 0.98) [3.1]} & 0.66 (0.51, 0.79) [2.5] & 0.00 (0.00, 0.00) [1.0]\\
planted partition, CRP & \textbf{0.99 (0.99, 1.00) [3.1]} & \textbf{0.99 (0.98, 0.99) [3.1]} & 0.96 (0.94, 0.97) [3.2] & \textbf{0.85 (0.80, 0.88) [3.1]} & \textbf{0.05 (0.02, 0.08) [1.6]}\\
\bottomrule\end{tabular}}\end{table}
\begin{table}[tbp]\centering\scriptsize
\caption{Constraints between communities, 120 contaminating cross edges, $K$ known ($n=120$, $K=3$, $p_{\mathrm{in}}=0.2$, cap $x=0.02$, 30 paired seeds): mean ARI ($95\%$ bootstrap interval). Columns are the ratio of the planted between density to the cap. Higher is better.}
\label{tab:R16_known_120}
\resizebox{\linewidth}{!}{\begin{tabular}{lccccc}\toprule
$K$ known; $p_{\mathrm{out}}/x$ & 0.25 & 0.5 & 1 & 2 & 4\\\midrule
flat & \textbf{0.98 (0.97, 0.98)} & \textbf{0.96 (0.95, 0.97)} & \textbf{0.91 (0.89, 0.93)} & 0.44 (0.35, 0.53) & 0.01 (0.00, 0.01)\\
hard cap 0.02 & \textbf{0.98 (0.97, 0.98)} & \textbf{0.96 (0.95, 0.97)} & 0.56 (0.42, 0.69) & 0.04 (0.00, 0.09) & 0.00 (0.00, 0.00)\\
soft shrinkage & \textbf{0.98 (0.97, 0.98)} & \textbf{0.96 (0.95, 0.97)} & 0.90 (0.88, 0.92) & 0.51 (0.42, 0.59) & 0.00 (0.00, 0.01)\\
relative cap 0.1 & \textbf{0.98 (0.97, 0.99)} & \textbf{0.96 (0.95, 0.97)} & 0.88 (0.81, 0.92) & 0.05 (0.00, 0.14) & 0.00 (0.00, 0.00)\\
planted partition & \textbf{0.98 (0.97, 0.98)} & \textbf{0.96 (0.95, 0.97)} & \textbf{0.91 (0.89, 0.93)} & 0.63 (0.55, 0.70) & 0.00 (0.00, 0.00)\\
variational EM & \textbf{0.98 (0.97, 0.98)} & \textbf{0.96 (0.95, 0.97)} & \textbf{0.91 (0.89, 0.93)} & \textbf{0.73 (0.68, 0.78)} & 0.00 (-0.00, 0.01)\\
Bethe Hessian & 0.96 (0.95, 0.97) & 0.94 (0.92, 0.95) & 0.88 (0.86, 0.89) & 0.69 (0.65, 0.72) & \textbf{0.14 (0.12, 0.16)}\\
Leiden & \textbf{0.98 (0.97, 0.99)} & 0.95 (0.94, 0.97) & 0.87 (0.82, 0.91) & 0.48 (0.42, 0.54) & 0.08 (0.06, 0.09)\\
\bottomrule\end{tabular}}\end{table}
\begin{table}[tbp]\centering\scriptsize
\caption{Constraints between communities, 120 contaminating cross edges, $K$ inferred with a CRP prior and $10$ labels (true $K=3$, 30 paired seeds): mean ARI ($95\%$ bootstrap interval) [mean $\hat K$]. Higher ARI is better.}
\label{tab:R16_auto_120}
\resizebox{\linewidth}{!}{\begin{tabular}{lccccc}\toprule
$K$ inferred; $p_{\mathrm{out}}/x$ & 0.25 & 0.5 & 1 & 2 & 4\\\midrule
flat, CRP & \textbf{0.98 (0.97, 0.98) [3.0]} & \textbf{0.96 (0.95, 0.97) [3.0]} & \textbf{0.91 (0.89, 0.93) [3.0]} & 0.51 (0.42, 0.60) [2.4] & \textbf{0.01 (0.00, 0.02) [1.1]}\\
hard cap 0.02, CRP & \textbf{0.98 (0.97, 0.98) [3.1]} & \textbf{0.96 (0.95, 0.97) [3.1]} & 0.60 (0.48, 0.73) [2.3] & 0.01 (0.00, 0.04) [1.1] & 0.00 (0.00, 0.00) [1.0]\\
soft shrinkage, CRP & 0.97 (0.96, 0.98) [3.2] & \textbf{0.96 (0.94, 0.97) [3.2]} & \textbf{0.91 (0.89, 0.92) [3.1]} & 0.55 (0.46, 0.64) [2.5] & 0.00 (0.00, 0.00) [1.1]\\
relative cap 0.1, CRP & \textbf{0.98 (0.97, 0.99) [3.1]} & \textbf{0.96 (0.94, 0.97) [3.1]} & \textbf{0.91 (0.90, 0.93) [3.1]} & 0.05 (0.00, 0.14) [1.1] & 0.00 (0.00, 0.00) [1.0]\\
planted partition, CRP & 0.97 (0.97, 0.98) [3.3] & \textbf{0.96 (0.94, 0.97) [3.2]} & 0.90 (0.88, 0.92) [3.3] & \textbf{0.63 (0.55, 0.71) [3.0]} & 0.00 (0.00, 0.00) [1.3]\\
\bottomrule\end{tabular}}\end{table}

\begin{table}[tbp]\centering\scriptsize
\caption{Paired differences in ARI between each constrained or alternative instance and the flat evidence with the same prior on $K$, in every cell of the constraint experiment ($n=120$, $K=3$, $p_{\mathrm{in}}=0.2$, cap $x=0.02$, 30 paired seeds): mean difference ($95\%$ bootstrap interval), a star marking a Wilcoxon signed-rank test significant at level $0.05$ after Holm adjustment over the 135 tests of the table. The first five columns are differences to the flat evidence with $K=3$ given, the last four to the flat CRP evidence. A positive value favours the named instance.}
\label{tab:R16paired}
\resizebox{\linewidth}{!}{\begin{tabular}{llccccccccc}\toprule
$p_{\mathrm{out}}/x$ & edges & hard cap & soft shrinkage & relative cap & planted partition & variational EM & hard cap, CRP & soft shrinkage, CRP & relative cap, CRP & planted partition, CRP\\\midrule
0.25 & 0 & +0.00 (+0.00, +0.00) & +0.00 (+0.00, +0.00) & +0.00 (+0.00, +0.00) & +0.00 (+0.00, +0.00) & +0.00 (+0.00, +0.00) & -0.00 (-0.00, +0.00) & -0.00 (-0.00, +0.00) & -0.00 (-0.00, +0.00) & -0.00 (-0.00, +0.00)\\
0.5 & 0 & +0.00 (+0.00, +0.00) & +0.00 (+0.00, +0.00) & +0.00 (+0.00, +0.00) & +0.00 (+0.00, +0.00) & +0.00 (+0.00, +0.00) & +0.00 (-0.00, +0.00) & -0.00 (-0.00, +0.00) & +0.00 (-0.00, +0.00) & +0.00 (-0.00, +0.00)\\
1 & 0 & +0.00 (+0.00, +0.00) & +0.00 (+0.00, +0.00) & +0.00 (+0.00, +0.01) & +0.00 (+0.00, +0.00) & +0.00 (+0.00, +0.00) & +0.00 (-0.00, +0.00) & -0.00 (-0.01, +0.00) & +0.00 (+0.00, +0.01) & +0.00 (-0.00, +0.00)\\
2 & 0 & -0.11 (-0.19, -0.04) & +0.00 (+0.00, +0.00) & +0.00 (-0.01, +0.02) & -0.00 (-0.01, +0.00) & -0.00 (-0.01, +0.00) & -0.12 (-0.21, -0.05) & +0.00 (+0.00, +0.00) & -0.01 (-0.02, +0.01) & +0.00 (-0.00, +0.01)\\
4 & 0 & -0.05 (-0.09, -0.02) & +0.04 (+0.01, +0.08) & -0.05 (-0.09, -0.02) & +0.09 (+0.04, +0.16) & +0.26 (+0.16, +0.35)$^{*}$ & -0.05 (-0.08, -0.02) & +0.04 (+0.00, +0.09) & -0.05 (-0.08, -0.02) & +0.16 (+0.10, +0.23)$^{*}$\\
0.25 & 60 & +0.00 (+0.00, +0.00) & +0.00 (+0.00, +0.00) & -0.00 (-0.00, +0.00) & -0.00 (-0.00, +0.00) & +0.00 (+0.00, +0.00) & +0.00 (-0.00, +0.00) & -0.00 (-0.01, +0.00) & -0.00 (-0.01, +0.00) & +0.00 (-0.00, +0.00)\\
0.5 & 60 & +0.00 (+0.00, +0.00) & +0.00 (+0.00, +0.00) & +0.00 (-0.00, +0.00) & +0.00 (-0.00, +0.00) & +0.00 (+0.00, +0.00) & +0.00 (+0.00, +0.00) & -0.00 (-0.00, +0.00) & -0.00 (-0.01, +0.00) & +0.00 (-0.00, +0.00)\\
1 & 60 & +0.00 (+0.00, +0.00) & +0.00 (+0.00, +0.00) & +0.00 (-0.00, +0.02) & +0.00 (+0.00, +0.01) & -0.00 (-0.01, +0.00) & +0.00 (+0.00, +0.00) & -0.00 (-0.01, +0.00) & +0.01 (-0.00, +0.02) & -0.00 (-0.01, +0.00)\\
2 & 60 & -0.60 (-0.72, -0.49)$^{*}$ & +0.05 (+0.01, +0.10) & -0.21 (-0.32, -0.10) & +0.08 (+0.02, +0.14) & +0.08 (+0.03, +0.14) & -0.69 (-0.78, -0.60)$^{*}$ & +0.04 (+0.01, +0.07) & -0.15 (-0.25, -0.06) & +0.03 (-0.00, +0.08)\\
4 & 60 & -0.01 (-0.02, +0.00) & +0.00 (+0.00, +0.01) & -0.01 (-0.02, +0.00) & +0.01 (-0.01, +0.02) & +0.07 (+0.03, +0.11) & -0.01 (-0.02, +0.00) & -0.00 (-0.01, +0.00) & -0.01 (-0.02, +0.00) & +0.04 (+0.02, +0.07)\\
0.25 & 120 & +0.00 (+0.00, +0.00) & +0.00 (+0.00, +0.00) & +0.00 (-0.00, +0.01) & +0.00 (+0.00, +0.00) & +0.00 (+0.00, +0.00) & +0.00 (-0.00, +0.00) & -0.00 (-0.01, +0.00) & +0.00 (-0.00, +0.01) & -0.00 (-0.01, +0.00)\\
0.5 & 120 & -0.00 (-0.00, +0.00) & +0.00 (+0.00, +0.00) & -0.00 (-0.01, +0.01) & -0.00 (-0.00, +0.00) & -0.00 (-0.00, +0.00) & -0.00 (-0.00, +0.00) & -0.01 (-0.01, +0.00) & -0.00 (-0.01, +0.01) & -0.00 (-0.01, -0.00)\\
1 & 120 & -0.34 (-0.46, -0.21)$^{*}$ & -0.00 (-0.01, +0.00) & -0.03 (-0.08, +0.01) & +0.00 (-0.00, +0.01) & +0.00 (-0.01, +0.01) & -0.30 (-0.44, -0.18) & -0.00 (-0.01, +0.01) & +0.01 (-0.01, +0.02) & -0.01 (-0.02, +0.00)\\
2 & 120 & -0.40 (-0.49, -0.32)$^{*}$ & +0.07 (+0.02, +0.12) & -0.38 (-0.48, -0.30)$^{*}$ & +0.19 (+0.12, +0.27)$^{*}$ & +0.30 (+0.22, +0.38)$^{*}$ & -0.50 (-0.59, -0.41)$^{*}$ & +0.04 (-0.02, +0.10) & -0.46 (-0.55, -0.37)$^{*}$ & +0.12 (+0.05, +0.19)\\
4 & 120 & -0.01 (-0.01, -0.00) & -0.00 (-0.01, -0.00) & -0.01 (-0.01, -0.00) & -0.00 (-0.01, +0.00) & -0.00 (-0.01, +0.01) & -0.01 (-0.02, +0.00) & -0.01 (-0.02, +0.00) & -0.01 (-0.02, +0.00) & -0.00 (-0.02, +0.00)\\
\bottomrule\end{tabular}}\end{table}

\clearpage
\subsection{Reciprocity, sign, intensity and layer structure}
\begin{table}[tbp]\centering\small
\caption{Reciprocity-only structure ($n=200$, $K=2$, probability of no tie $0.7$ on every dyad, between-community reciprocated share $0.02$). Columns give the within-community reciprocated share. Entries are the mean ARI (sd) over eight paired seeds, and higher is better.}
\label{tab:R7}
\resizebox{\linewidth}{!}{\begin{tabular}{lccccc}\toprule
Method & 0.04 & 0.06 & 0.09 & 0.12 & 0.18\\\midrule
collapsed dyad-state & 0.00 (0.00) & 0.45 (0.28) & \textbf{0.97 (0.02)} & 0.99 (0.01) & \textbf{1.00 (0.00)}\\
collapsed directed Bernoulli & 0.00 (0.00) & -0.00 (0.00) & 0.00 (0.01) & 0.10 (0.12) & 0.75 (0.18)\\
spectral, reciprocated subnetwork & \textbf{0.04 (0.05)} & \textbf{0.52 (0.10)} & 0.94 (0.04) & \textbf{1.00 (0.01)} & \textbf{1.00 (0.00)}\\
Leiden (symmetrised) & -0.00 (0.01) & -0.00 (0.01) & -0.00 (0.01) & -0.00 (0.01) & -0.00 (0.01)\\
Bethe Hessian (symmetrised) & 0.00 (0.01) & 0.00 (0.01) & 0.00 (0.01) & 0.00 (0.01) & 0.00 (0.01)\\
\bottomrule\end{tabular}}\end{table}
\begin{table}[tbp]\centering\small
\caption{Sign-only structure ($n=200$, $K=3$, tie density $0.10$ on every dyad). Columns give the sign purity. Entries are the mean ARI (sd) over eight paired seeds, and higher is better.}
\label{tab:R8}
\resizebox{\linewidth}{!}{\begin{tabular}{lccccc}\toprule
Method & 0.55 & 0.6 & 0.7 & 0.8 & 0.9\\\midrule
collapsed signed & 0.00 (0.00) & 0.00 (0.00) & 0.01 (0.03) & \textbf{0.94 (0.04)} & \textbf{1.00 (0.00)}\\
signed spectral & 0.00 (0.00) & 0.00 (0.00) & 0.01 (0.02) & 0.22 (0.24) & 0.77 (0.23)\\
signed Leiden & \textbf{0.01 (0.01)} & \textbf{0.02 (0.01)} & \textbf{0.08 (0.03)} & 0.86 (0.19) & \textbf{1.00 (0.01)}\\
Leiden ($|W|$) & 0.00 (0.00) & 0.00 (0.00) & 0.00 (0.00) & 0.00 (0.00) & 0.00 (0.00)\\
Bethe Hessian ($|W|$) & 0.00 (0.01) & 0.00 (0.01) & 0.00 (0.01) & 0.00 (0.01) & 0.00 (0.01)\\
\bottomrule\end{tabular}}\end{table}
\begin{table}[tbp]\centering\small
\caption{Intensity-only structure in zero-inflated counts ($n=150$, $K=3$, activity $0.08$ on every dyad, $\lambda_{\mathrm{out}}=1$). Columns give $\lambda_{\mathrm{in}}$. Entries are the mean ARI (sd) over eight paired seeds, and higher is better.}
\label{tab:R13}
\resizebox{\linewidth}{!}{\begin{tabular}{lccccc}\toprule
Method & 1.5 & 2.0 & 3.0 & 4.0 & 6.0\\\midrule
collapsed zero-inflated Poisson & -0.00 (0.00) & 0.00 (0.00) & \textbf{0.07 (0.16)} & \textbf{0.82 (0.21)} & \textbf{0.99 (0.01)}\\
collapsed Gamma--Poisson & 0.00 (0.01) & 0.01 (0.02) & 0.03 (0.03) & 0.12 (0.10) & 0.28 (0.11)\\
Leiden & \textbf{0.01 (0.01)} & \textbf{0.02 (0.02)} & \textbf{0.07 (0.03)} & 0.14 (0.04) & 0.26 (0.06)\\
Bethe Hessian & \textbf{0.01 (0.01)} & 0.01 (0.02) & 0.02 (0.02) & 0.10 (0.08) & 0.16 (0.07)\\
Infomap & 0.00 (0.00) & 0.00 (0.00) & 0.01 (0.02) & 0.07 (0.08) & 0.21 (0.03)\\
flat (binarised) & 0.00 (0.00) & 0.00 (0.00) & 0.00 (0.00) & 0.00 (0.00) & 0.00 (0.00)\\
Leiden (binarised) & -0.00 (0.01) & -0.00 (0.01) & -0.00 (0.01) & -0.00 (0.01) & -0.00 (0.01)\\
Bethe Hessian (binarised) & 0.00 (0.01) & 0.00 (0.01) & 0.00 (0.01) & 0.00 (0.01) & 0.00 (0.01)\\
\bottomrule\end{tabular}}\end{table}
\begin{table}[tbp]\centering\small
\caption{Two layers with mirrored block matrices ($n=200$, $K=3$). Columns give the $\mathrm{SNR}_{\mathrm{KS}}$ of the first layer. Entries are the mean ARI (sd) over eight paired seeds, and higher is better.}
\label{tab:R6b}
\resizebox{\linewidth}{!}{\begin{tabular}{lccccc}\toprule
Method & 1.0 & 1.5 & 2.0 & 3.0 & 5.0\\\midrule
collapsed multiplex & \textbf{0.26 (0.30)} & \textbf{0.88 (0.05)} & \textbf{0.94 (0.03)} & \textbf{0.99 (0.01)} & \textbf{1.00 (0.00)}\\
best single layer (oracle) & 0.13 (0.06) & 0.41 (0.10) & 0.62 (0.11) & 0.82 (0.07) & 0.98 (0.01)\\
Leiden (summed layers) & 0.00 (0.00) & 0.00 (0.01) & 0.00 (0.01) & -0.00 (0.01) & 0.00 (0.00)\\
Bethe Hessian (summed) & 0.00 (0.01) & -0.00 (0.01) & 0.00 (0.01) & -0.00 (0.01) & -0.00 (0.01)\\
\bottomrule\end{tabular}}\end{table}

\clearpage
\subsection{Computation}\label{app:computation}
This experiment measures what collapsing, the split and merge move and the initialisation contribute,
and how the cost scales with $n$ (Tables~\ref{tab:R11a} to~\ref{tab:R10ari}). On a
planted partition with $n=300$, $K=3$ and mean degree $10$ (eight paired seeds, $400$ sweeps), the
single-site collapsed chain from a random start reaches ARI $0.95$ within $400$ sweeps on half of the
seeds at $\mathrm{SNR}_{\mathrm{KS}}=4$ and its MAP averages $0.78$. The split and merge move raises this
to seven of eight seeds and $0.94$, the same as a spectral start, with a gain in log posterior over the
random start of $350$ to $1230$ nats. An informed start accounts for most of the result at moderate and
strong signal: the default configuration improves on its initial partition by $2$ nats at
$\mathrm{SNR}_{\mathrm{KS}}=4$, by nothing at $8$, and by $34$ nats at $2$, below the evidence threshold
for $K=3$ at mean degree $10$ ($2.38$, large-degree limit $2.20$), where no configuration reaches ARI
$0.8$ and the evidence is the limitation. The uncollapsed sampler, which draws the block parameters,
has the highest effective sample size of the log posterior per second \citep{Geyer1992} and the worst
partitions ($0.26$, $0.62$ and $0.83$), because it moves quickly inside the wrong mode, so the
effective sample size is not reported as a measure of quality, and Remark~\ref{prop:collapse}(ii),
which orders the systematic-scan samplers by the maximal correlation of successive states and not by
their convergence rate, says nothing about the quality of these chains. The split and merge move as
implemented costs about twenty times the sweep it accompanies ($12$ against $0.6$ seconds for $400$
sweeps), because the restricted scans and the spectral launch state run outside the compiled kernel.
At $\mathrm{SNR}_{\mathrm{KS}}=3$, $K=5$ and mean degree $10$ with $n$ from $300$ to $10^4$ (two
seeds), the three collapsed instances scale linearly in $n$ at fixed degree ($3.6$, $10.6$, $31.6$ and
$118$ seconds for the flat evidence with $200$ sweeps, the typed instance $3$ to $13\%$ more and the
global-DC instance slightly less) and reach the same ARI as variational EM, Karrer--Newman and Leiden
($0.87$ to $0.88$ at $n=10^4$, against $0.82$ for the Bethe Hessian, while Infomap returns hundreds of
communities). The greedy optimiser of the same objective is two orders of magnitude faster ($1.4$
seconds at $n=10^4$) at no measurable cost in ARI ($0.87$ against $0.88$). On these well-separated
instances the sampler adds nothing that the greedy optimiser does not, and its advantage is confined to
the multimodal cases of Section~\ref{sec:real}, where chains from over-dispersed starts find modes that
no deterministic start reaches.
\begin{table}[tbp]\centering\small
\caption{Mixing ($n=300$, $K=3$, mean degree 10, eight paired seeds, 400 sweeps): mean ARI (sd) of the MAP partition. Higher is better.}
\label{tab:R11a}
\resizebox{\linewidth}{!}{\begin{tabular}{lccc}\toprule
Sampler by $\mathrm{SNR}_{\mathrm{KS}}$ & 2.0 & 4.0 & 8.0\\\midrule
collapsed Gibbs, random start & 0.38 (0.26) & 0.78 (0.22) & 0.94 (0.15)\\
collapsed with split and merge, random start & 0.53 (0.16) & \textbf{0.94 (0.03)} & \textbf{1.00 (0.01)}\\
collapsed Gibbs, spectral start & \textbf{0.62 (0.08)} & \textbf{0.94 (0.03)} & \textbf{1.00 (0.01)}\\
collapsed with split and merge, model-based start (default) & 0.59 (0.12) & \textbf{0.94 (0.03)} & \textbf{1.00 (0.01)}\\
uncollapsed Gibbs, random start & 0.26 (0.19) & 0.62 (0.20) & 0.83 (0.22)\\
\bottomrule\end{tabular}}\end{table}
\begin{table}[tbp]\centering\small
\caption{Mixing: fraction of seeds in which the chain reaches ARI $\ge0.95$ within 400 sweeps, with the standard deviation of the indicator over seeds in parentheses and the wall-clock seconds for the 400 sweeps in brackets.}
\label{tab:R11b}
\resizebox{\linewidth}{!}{\begin{tabular}{lccc}\toprule
Sampler by $\mathrm{SNR}_{\mathrm{KS}}$ & 2.0 & 4.0 & 8.0\\\midrule
collapsed Gibbs, random start & \textbf{0.00 (0.00) [0.6]} & 0.50 (0.53) [0.6] & 0.88 (0.35) [0.6]\\
collapsed with split and merge, random start & \textbf{0.00 (0.00) [11.4]} & \textbf{0.88 (0.35) [12.7]} & \textbf{1.00 (0.00) [12.0]}\\
collapsed Gibbs, spectral start & \textbf{0.00 (0.00) [0.6]} & \textbf{0.88 (0.35) [0.7]} & \textbf{1.00 (0.00) [0.7]}\\
collapsed with split and merge, model-based start (default) & \textbf{0.00 (0.00) [12.2]} & \textbf{0.88 (0.35) [13.0]} & \textbf{1.00 (0.00) [11.6]}\\
uncollapsed Gibbs, random start & \textbf{0.00 (0.00) [0.5]} & 0.25 (0.46) [0.5] & 0.62 (0.52) [0.4]\\
\bottomrule\end{tabular}}\end{table}
\begin{table}[tbp]\centering\small
\caption{Mixing: effective sample size of the log-posterior trace per second (Geyer's initial positive sequence estimator). A high value for a chain with a low ARI in Table~\ref{tab:R11a} means fast mixing within the wrong mode.}
\label{tab:R11c}
\resizebox{\linewidth}{!}{\begin{tabular}{lccc}\toprule
Sampler by $\mathrm{SNR}_{\mathrm{KS}}$ & 2.0 & 4.0 & 8.0\\\midrule
collapsed Gibbs, random start & 124.16 (61.90) & 268.08 (98.28) & 346.58 (167.48)\\
collapsed with split and merge, random start & 5.79 (3.44) & 11.26 (7.00) & 19.52 (8.78)\\
collapsed Gibbs, spectral start & 154.74 (84.65) & 284.70 (57.37) & 370.82 (181.67)\\
collapsed with split and merge, model-based start (default) & 7.15 (2.88) & 13.68 (2.76) & 19.36 (10.52)\\
uncollapsed Gibbs, random start & \textbf{157.49 (156.06)} & \textbf{367.26 (90.83)} & \textbf{595.02 (123.36)}\\
\bottomrule\end{tabular}}\end{table}
\begin{table}[tbp]\centering\small
\caption{Mixing: gain in log posterior (nats) of the MAP over the initial partition, with the ARI of the initial partition in brackets.}
\label{tab:R11d}
\resizebox{\linewidth}{!}{\begin{tabular}{lccc}\toprule
Sampler by $\mathrm{SNR}_{\mathrm{KS}}$ & 2.0 & 4.0 & 8.0\\\midrule
collapsed Gibbs, random start & 348.46 (11.96) [0.0] & 554.79 (56.39) [0.0] & 1178.64 (137.82) [0.0]\\
collapsed with split and merge, random start & \textbf{348.59 (11.96) [0.0]} & \textbf{597.51 (19.40) [0.0]} & \textbf{1228.49 (40.52) [0.0]}\\
collapsed Gibbs, spectral start & 33.06 (9.56) [0.6] & 5.56 (3.43) [0.9] & 0.10 (0.29) [1.0]\\
collapsed with split and merge, model-based start (default) & 33.80 (8.29) [0.6] & 1.95 (3.01) [0.9] & 0.00 (0.00) [1.0]\\
\bottomrule\end{tabular}}\end{table}
\begin{table}[tbp]\centering\small
\caption{Wall-clock seconds on one core ($\mathrm{SNR}_{\mathrm{KS}}=3$, $K=5$, mean degree 10, two seeds), mean (sd). The collapsed samplers run $200$ sweeps with model-based initialisation and split and merge moves. Lower is better.}
\label{tab:R10}
\begin{tabular}{lcccc}\toprule
Method by $n$ & 300 & 1000 & 3000 & 10000\\\midrule
flat & 3.58 (0.14) & 10.55 (0.11) & 31.59 (0.24) & 117.50 (0.98)\\
typed \{ER,DC\} & 3.69 (0.14) & 11.89 (0.25) & 34.30 (1.75) & 132.45 (2.02)\\
global-DC & 3.52 (0.22) & 9.34 (0.04) & 30.20 (1.05) & 109.27 (2.07)\\
greedy optimiser & 0.05 (0.00) & 0.11 (0.02) & 0.49 (0.05) & 1.39 (0.33)\\
variational EM & 0.03 (0.01) & 0.23 (0.02) & 13.40 (15.63) & 27.67 (4.18)\\
Karrer--Newman likelihood & 0.74 (0.14) & 3.08 (0.68) & 10.05 (0.58) & 58.92 (18.05)\\
Bethe Hessian & \textbf{0.01 (0.00)} & \textbf{0.02 (0.00)} & \textbf{0.10 (0.07)} & \textbf{0.20 (0.01)}\\
Leiden & 0.03 (0.01) & 0.12 (0.00) & 0.59 (0.18) & 4.27 (2.70)\\
Infomap & 0.04 (0.00) & 0.20 (0.00) & 0.86 (0.02) & 4.30 (0.13)\\
\bottomrule\end{tabular}\end{table}
\begin{table}[tbp]\centering\small
\caption{Accuracy in the scaling experiment of Table~\ref{tab:R10} ($\mathrm{SNR}_{\mathrm{KS}}=3$, $K=5$, mean degree 10, two seeds): mean ARI (sd) [mean $\hat K$]. Higher ARI is better.}
\label{tab:R10ari}
\resizebox{\linewidth}{!}{\begin{tabular}{lcccc}\toprule
Method by $n$ & 300 & 1000 & 3000 & 10000\\\midrule
flat & 0.87 (0.02) [5.0] & \textbf{0.91 (0.01) [5.0]} & \textbf{0.88 (0.01) [5.0]} & \textbf{0.88 (0.01) [5.0]}\\
typed \{ER,DC\} & 0.87 (0.02) [5.0] & \textbf{0.91 (0.01) [5.0]} & 0.87 (0.00) [5.0] & \textbf{0.88 (0.01) [5.0]}\\
global-DC & \textbf{0.88 (0.01) [5.0]} & 0.90 (0.01) [5.0] & 0.87 (0.01) [5.0] & 0.87 (0.00) [5.0]\\
greedy optimiser & 0.86 (0.03) [5.0] & \textbf{0.91 (0.01) [5.0]} & 0.86 (0.00) [5.0] & 0.87 (0.01) [5.0]\\
variational EM & 0.87 (0.03) [5.0] & \textbf{0.91 (0.00) [5.0]} & \textbf{0.88 (0.00) [5.0]} & \textbf{0.88 (0.01) [5.0]}\\
Karrer--Newman likelihood & 0.87 (0.03) [5.0] & \textbf{0.91 (0.02) [5.0]} & 0.87 (0.01) [5.0] & 0.87 (0.01) [5.0]\\
Bethe Hessian & 0.85 (0.03) [5.0] & 0.87 (0.01) [5.0] & 0.81 (0.01) [5.0] & 0.82 (0.01) [5.0]\\
Leiden & 0.87 (0.01) [5.0] & \textbf{0.91 (0.01) [5.0]} & 0.87 (0.02) [5.0] & 0.87 (0.01) [6.0]\\
Infomap & 0.76 (0.11) [7.0] & 0.45 (0.13) [35.0] & 0.01 (0.00) [187.0] & 0.00 (0.00) [580.0]\\
\bottomrule\end{tabular}}\end{table}

\clearpage
\subsection{Real networks}
\begin{table}[tbp]\centering\scriptsize
\caption{Family effect on the football and school networks, computed with exact evidences on the full networks and no sampling. For each network the rows are the metadata partition and the pooled MAP partitions of the flat, typed and global-DC instances with $K$ fixed to the number of metadata classes, with their $\hat K$ and ARI against the metadata. The columns give the log posterior $\log p(A\mid z)+\log p(z)$ under the homogeneous Beta--Bernoulli model, under the homogeneous Gamma--Poisson model, which is the representation on which the typed and degree-corrected instances are built, the family gap $G(z)=\sum_{\text{blocks}}N\,g(\hat p)$ of Proposition~\ref{prop:family} in nats, which the difference of the two preceding columns reproduces up to $O(K^2\log n)$, and the log posterior under the Bernoulli type library, with the Erd\H{o}s--R\'enyi type collapsed exactly and the degree-corrected type through the Laplace form with $\alpha_i\sim N(0,4)$. Bold marks the highest value of each evidence column within a network, and a higher value means a preferred partition. The within and between densities of the metadata partitions are $0.75$ and $0.036$ (football) and $0.92$ and $0.143$ (school).}
\label{tab:familygap}
\resizebox{\linewidth}{!}{\begin{tabular}{llrrrrrr}\toprule
Network & partition & $\hat K$ & ARI & Beta--Bernoulli & Gamma--Poisson & $G(z)$ & Bernoulli typed \{ER,DC\}\\\midrule
football & metadata & 12 & 1.00 & -1638.3 & -1892.1 & 257.0 & -1629.1\\
 & flat MAP & 11 & 0.86 & \textbf{-1541.4} & -1813.4 & 277.2 & \textbf{-1546.6}\\
 & typed MAP & 7 & 0.58 & -1620.8 & \textbf{-1791.2} & 171.2 & -1625.6\\
 & global-DC MAP & 10 & 0.81 & -1555.7 & -1798.3 & 244.8 & -1560.3\\
\midrule
school & metadata & 11 & 1.00 & -10868.2 & -13145.1 & 2278.4 & -10875.8\\
 & flat MAP & 11 & 0.44 & -10041.0 & \textbf{-12110.0} & 2068.9 & -10048.6\\
 & typed MAP & 11 & 0.43 & -10138.7 & -12140.4 & 2001.3 & -10138.2\\
 & global-DC MAP & 11 & 0.70 & \textbf{-10026.9} & -12354.2 & 2330.0 & \textbf{-10031.2}\\
\bottomrule\end{tabular}}\end{table}

\begin{table}[tbp]\centering\scriptsize
\caption{Semi-synthetic political blogs: the within-camp ties of one camp are replaced by an Erd\H{o}s--R\'enyi network with the same number of edges, the other camp and the between-camp ties are kept, $K=2$ given, 20 rewirings: mean ARI against the two camps ($95\%$ bootstrap interval) [type accuracy of the typed instances, the fraction of the two communities whose posterior type is Erd\H{o}s--R\'enyi for the rewired camp and degree-corrected for the other]. Higher is better.}
\label{tab:R19}
\resizebox{\linewidth}{!}{\begin{tabular}{lcc}\toprule
Method; rewired camp & liberal & conservative\\\midrule
flat & 0.27 (0.26, 0.27) & 0.14 (0.13, 0.14)\\
typed \{ER,DC\} & 0.93 (0.93, 0.93) [1.0] & \textbf{0.94 (0.94, 0.94) [1.0]}\\
global typed & 0.94 (0.94, 0.95) [1.0] & \textbf{0.94 (0.94, 0.94) [1.0]}\\
global-DC & 0.91 (0.90, 0.91) & 0.86 (0.85, 0.86)\\
Karrer--Newman likelihood & 0.91 (0.90, 0.91) & 0.88 (0.88, 0.89)\\
SCORE & 0.82 (0.81, 0.82) & 0.81 (0.80, 0.81)\\
regularised spectral & 0.90 (0.90, 0.90) & 0.85 (0.84, 0.85)\\
Bethe Hessian & \textbf{0.96 (0.96, 0.96)} & 0.12 (0.12, 0.12)\\
\bottomrule\end{tabular}}\end{table}

\begin{table}[tbp]\centering\scriptsize
\caption{Semi-synthetic football and school networks (Theorem~\ref{thm:split}(c) on two more networks): the within ties of one of the three largest classes are replaced by an Erd\H{o}s--R\'enyi network (ER) or by a degree-corrected network with Pareto$(2.5)$ propensities (DC) with the same number of edges, 10 rewirings per class and mode, $K$ given. Cells: mean ARI against the classes over the three classes and rewirings, and for the typed instances the fraction of rewired classes whose type is recovered [mean posterior probability of the degree-corrected type of the block holding the class].}
\label{tab:R23}
\resizebox{\linewidth}{!}{\begin{tabular}{lcccc}\toprule
Method & football, ER & football, DC & school, ER & school, DC\\\midrule
flat & 0.85 & 0.85 & 0.83 & 0.82\\
typed \{ER,DC\} & 0.71, 1.00 [0.00] & 0.71, 0.00 [0.00] & 0.55, 1.00 [0.00] & 0.55, 0.00 [0.00]\\
global typed & 0.82, 1.00 [0.00] & 0.82, 0.00 [0.00] & 0.70, 0.60 [0.40] & 0.69, 0.37 [0.37]\\
global typed, global degree prior & 0.82, 0.67 [0.32] & 0.82, 0.00 [0.05] & 0.70, 0.57 [0.43] & 0.70, 0.43 [0.43]\\
global-DC & 0.82 & 0.82 & 0.70 & 0.70\\
Karrer--Newman likelihood & 0.89 & 0.89 & 0.61 & 0.61\\
regularised spectral & 0.90 & 0.90 & 0.90 & 0.90\\
Bethe Hessian & 0.90 & 0.90 & 0.89 & 0.89\\
\bottomrule\end{tabular}}\end{table}

\begin{table}[tbp]\centering\small
\caption{The catalogue on count-valued ties: the primary school contact network with the number of $20$-second contacts as the tie weight, $K=11$ classes given, 4 chains per instance (mean ARI against the classes with a $95\%$ bootstrap interval over chains, mean $\hat K$). The Poisson instances see the counts and the Bernoulli instances and the spectral methods the binarised network.}
\label{tab:R24}
\begin{tabular}{lcc}\toprule Method & ARI & $\hat K$\\\midrule
flat, Gamma--Poisson (counts) & 0.91 (0.90, 0.91) & 11.0\\
typed \{ER,DC\} & 0.82 (0.81, 0.82) & 11.0\\
global-DC & 0.91 (0.91, 0.91) & 11.0\\
flat, binarised & 0.86 (0.85, 0.89) & 11.0\\
typed \{ER,DC\}, binarised & 0.54 (0.54, 0.54) & 11.0\\
global-DC, binarised & 0.71 (0.69, 0.74) & 11.0\\
Karrer--Newman, binarised & 0.68 (0.63, 0.72) & 11.0\\
regularised spectral, binarised & 0.90 (0.90, 0.90) & 11.0\\
Bethe Hessian, binarised & 0.90 (0.89, 0.90) & 11.0\\
Leiden & 0.68 (0.68, 0.68) & 7.0\\
\bottomrule\end{tabular}\end{table}

\begin{table}[tbp]\centering\small
\caption{Real networks, all eight networks and all methods: ARI against the metadata partition, with the number of metadata classes in parentheses and $K$ fixed to it for the methods that need $K$, and $\hat K$ in brackets for the methods that infer it. Collapsed samplers: pooled MAP over four chains of $300$ sweeps. Higher is better.}
\label{tab:real_main_full}
\resizebox{\linewidth}{!}{\begin{tabular}{lcccccccc}\toprule
Method & karate (2) & dolphins (2) & football (12) & polbooks (3) & polblogs (2) & school (11) & eu-core (42) & Cora (7)\\\midrule
flat & -0.01 & 0.87 & 0.86 & 0.54 & 0.00 & 0.44 & 0.16 & 0.16\\
flat (Jeffreys) & -0.01 & 0.87 & 0.86 & 0.54 & 0.00 & 0.44 & 0.20 & 0.09\\
typed \{ER,DC\} & -0.01 & 0.87 & 0.58 & 0.63 & 0.40 & 0.43 & 0.18 & 0.22\\
global typed & -0.01 & 0.87 & 0.81 & 0.63 & 0.77 & 0.64 & 0.48 & 0.19\\
global typed, global degree prior & 0.00 & 0.93 & 0.81 & 0.63 & 0.77 & 0.68 & 0.49 & 0.26\\
global-DC & 0.00 & 0.93 & 0.81 & 0.63 & 0.77 & 0.70 & 0.46 & 0.26\\
flat, CRP & 0.78 [5] & 0.87 [2] & 0.86 [11] & 0.40 [5] & 0.20 [20] & 0.59 [21] & 0.16 [42] & 0.28 [19]\\
typed, CRP & -0.01 [2] & 0.87 [2] & 0.63 [8] & 0.41 [5] & 0.19 [20] & 0.52 [17] & 0.18 [35] & 0.29 [18]\\
global typed, CRP & -0.01 [2] & 0.87 [2] & 0.81 [10] & 0.30 [6] & 0.15 [20] & 0.52 [17] & 0.47 [25] & 0.25 [21]\\
global typed, global degree prior, CRP & 0.00 [1] & 0.93 [2] & 0.81 [10] & 0.63 [3] & 0.14 [20] & 0.52 [17] & 0.43 [23] & 0.24 [21]\\
global-DC, CRP & 0.00 [1] & 0.93 [2] & 0.81 [10] & 0.63 [3] & 0.19 [20] & 0.53 [17] & 0.46 [25] & 0.25 [21]\\
flat, Pitman--Yor & 0.78 [6] & 0.87 [2] & 0.86 [11] & 0.40 [5] & 0.12 [20] & 0.59 [22] & 0.16 [41] & 0.30 [18]\\
flat, Gnedin & 0.78 [5] & 0.53 [3] & 0.86 [11] & 0.40 [5] & 0.20 [20] & 0.59 [21] & 0.16 [39] & 0.27 [17]\\
flat, MFM & 0.78 [5] & 0.87 [2] & 0.86 [11] & 0.40 [5] & 0.20 [20] & 0.59 [21] & 0.17 [39] & 0.28 [17]\\
greedy optimiser & -0.01 & 0.87 & 0.86 & 0.63 & 0.00 & 0.89 & 0.16 & 0.12\\
variational EM & -0.01 & 0.87 & 0.90 & 0.61 & 0.00 & 0.47 & 0.17 & 0.05\\
Karrer--Newman likelihood & 0.77 & 0.93 & 0.89 & 0.63 & 0.79 & 0.74 & 0.42 & 0.27\\
Bickel--Chen profile likelihood & -0.01 & 0.87 & 0.90 & 0.64 & -0.00 & 0.47 & 0.25 & 0.15\\
Yan et al.\ test, planted partition & \textbf{1.00} & \textbf{1.00} & \textbf{1.00} & \textbf{1.00} & \textbf{1.00} & \textbf{1.00} & \textbf{1.00} & \textbf{1.00}\\
corrected BIC (Hu et al.) & 0.00 [1] & 0.87 [2] & 0.73 [9] & 0.60 [3] & 0.32 [5] & 0.90 [11] & 0.52 [17] & 0.33 [8]\\
Bethe Hessian & 0.88 & 0.87 & 0.90 & 0.64 & 0.11 & 0.89 & 0.14 & 0.10\\
regularised spectral & 0.88 & 0.87 & 0.91 & 0.60 & 0.81 & 0.90 & 0.39 & 0.29\\
SCORE & 0.67 & 0.93 & 0.89 & 0.64 & 0.80 & 0.25 & 0.07 & 0.22\\
Bethe Hessian, $K$ inferred & 0.00 [1] & 0.87 [2] & 0.72 [9] & 0.61 [4] & 0.08 [7] & 0.76 [8] & 0.14 [17] & 0.09 [50]\\
Leiden & 0.46 [4] & 0.36 [6] & 0.71 [9] & 0.66 [4] & 0.76 [23] & 0.56 [6] & 0.29 [17] & 0.25 [85]\\
Louvain & 0.45 [4] & 0.36 [6] & 0.71 [9] & 0.53 [5] & 0.76 [23] & 0.56 [6] & 0.33 [18] & 0.23 [81]\\
Infomap & 0.46 [4] & 0.31 [8] & 0.90 [12] & 0.66 [5] & 0.67 [49] & 0.56 [6] & 0.27 [27] & 0.06 [257]\\
walktrap & 0.67 [2] & -0.02 [2] & 0.89 [12] & 0.67 [3] & 0.00 [12] & 0.50 [11] & 0.04 [42] & -0.00 [53]\\
\bottomrule\end{tabular}}\end{table}

\begin{table}[tbp]\centering\small
\caption{Real networks, all eight networks and all methods: held-out link prediction on $10\%$ of the dyads, AUC. Higher is better.}
\label{tab:real_pred_full}
\resizebox{\linewidth}{!}{\begin{tabular}{lcccccccc}\toprule
Method & karate (2) & dolphins (2) & football (12) & polbooks (3) & polblogs (2) & school (11) & eu-core (42) & Cora (7)\\\midrule
flat & 0.80 & 0.75 & 0.92 & 0.87 & 0.85 & 0.88 & 0.93 & 0.81\\
flat (Jeffreys) & 0.80 & 0.75 & \textbf{0.93} & 0.86 & 0.85 & 0.88 & 0.93 & 0.84\\
typed \{ER,DC\} & 0.80 & 0.77 & 0.91 & 0.86 & 0.80 & 0.88 & 0.93 & 0.82\\
global typed & 0.85 & 0.76 & 0.91 & 0.89 & 0.94 & 0.90 & \textbf{0.95} & 0.80\\
global typed, global degree prior & 0.84 & 0.78 & 0.90 & 0.89 & 0.94 & 0.90 & 0.94 & 0.81\\
global-DC & 0.84 & 0.77 & 0.91 & 0.90 & 0.94 & 0.90 & 0.94 & 0.81\\
flat, CRP & 0.87 & 0.80 & 0.92 & 0.91 & \textbf{0.95} & 0.91 & 0.93 & 0.89\\
typed, CRP & 0.84 & 0.81 & 0.92 & 0.90 & \textbf{0.95} & 0.91 & 0.93 & 0.89\\
global typed, CRP & \textbf{0.88} & 0.79 & 0.91 & 0.90 & \textbf{0.95} & 0.91 & 0.94 & 0.83\\
global typed, global degree prior, CRP & 0.84 & 0.78 & 0.91 & 0.90 & \textbf{0.95} & 0.91 & 0.94 & 0.84\\
global-DC, CRP & 0.85 & 0.77 & 0.91 & 0.90 & \textbf{0.95} & \textbf{0.92} & 0.94 & 0.84\\
flat, Pitman--Yor & \textbf{0.88} & 0.80 & \textbf{0.93} & \textbf{0.92} & 0.94 & 0.91 & 0.93 & 0.89\\
flat, Gnedin & \textbf{0.88} & 0.80 & \textbf{0.93} & 0.91 & \textbf{0.95} & 0.91 & 0.93 & 0.89\\
flat, MFM & \textbf{0.88} & \textbf{0.83} & \textbf{0.93} & 0.91 & \textbf{0.95} & 0.91 & 0.93 & \textbf{0.90}\\
greedy optimiser & 0.69 & 0.69 & \textbf{0.93} & 0.83 & 0.84 & 0.85 & 0.93 & 0.72\\
variational EM & 0.69 & 0.69 & 0.92 & 0.83 & 0.84 & 0.88 & 0.93 & 0.77\\
Karrer--Newman likelihood & 0.83 & 0.77 & 0.90 & 0.90 & 0.94 & 0.90 & 0.94 & 0.78\\
Bickel--Chen profile likelihood & 0.69 & 0.69 & 0.92 & 0.83 & 0.84 & 0.88 & 0.93 & 0.73\\
Yan et al.\ test, planted partition & 0.70 & 0.73 & 0.92 & 0.81 & 0.71 & 0.83 & 0.85 & 0.87\\
corrected BIC (Hu et al.) & 0.50 & 0.69 & 0.86 & 0.80 & 0.80 & 0.85 & 0.80 & 0.81\\
Bethe Hessian & 0.65 & 0.71 & 0.92 & 0.82 & 0.64 & 0.84 & 0.86 & 0.67\\
regularised spectral & 0.65 & 0.69 & 0.92 & 0.80 & 0.73 & 0.85 & 0.81 & 0.80\\
SCORE & 0.80 & 0.77 & 0.90 & 0.88 & 0.94 & 0.79 & 0.90 & 0.77\\
Bethe Hessian, $K$ inferred & 0.50 & 0.71 & 0.86 & 0.84 & 0.86 & 0.83 & 0.82 & 0.73\\
Leiden & 0.74 & 0.79 & 0.91 & 0.85 & 0.75 & 0.82 & 0.76 & 0.89\\
Louvain & 0.67 & 0.79 & 0.91 & 0.86 & 0.75 & 0.82 & 0.77 & 0.88\\
Infomap & 0.74 & 0.76 & 0.92 & 0.85 & 0.78 & 0.82 & 0.79 & 0.80\\
walktrap & 0.63 & 0.47 & 0.92 & 0.83 & 0.51 & 0.87 & 0.67 & 0.47\\
density & 0.50 & 0.50 & 0.50 & 0.50 & 0.50 & 0.50 & 0.50 & 0.50\\
Chung-Lu & 0.84 & 0.69 & 0.22 & 0.72 & 0.91 & 0.65 & 0.86 & 0.66\\
\bottomrule\end{tabular}}\end{table}

\begin{table}[tbp]\centering\scriptsize
\caption{Larger networks: held-out AUC (balanced held-out set of edges and non-edges, $5\%$ of the edges up to $20{,}000$). Node counts and metadata classes in the header, $K$ fixed to the number of classes for the methods that need it, $\hat K$ in brackets for the methods that infer it (label budget $\min(3K,60)$). Collapsed samplers: pooled MAP over four chains of $300$ sweeps, two chains of $200$ sweeps on the Amazon network and one chain on PubMed, whose size allowed only the instances shown. Higher is better.}
\label{tab:realv6_auc}
\resizebox{\linewidth}{!}{\begin{tabular}{lccccc}\toprule
Method & school day 2 (238, 11) & CiteSeer (1331, 6) & EuroSiS (1272, 13) & Amazon photo (7487, 8) & PubMed (19717, 3)\\\midrule
flat & 0.91 & 0.75 & 0.91 & 0.93 & 0.88\\
typed \{ER,DC\} & 0.90 & 0.77 & 0.93 & -- & --\\
global typed & 0.92 & 0.71 & 0.93 & -- & --\\
global typed, global degree prior & 0.92 & 0.72 & 0.93 & 0.96 & 0.74\\
global-DC & 0.92 & 0.72 & 0.94 & 0.96 & 0.74\\
Karrer--Newman likelihood & 0.92 & 0.70 & 0.93 & 0.95 & 0.75\\
regularised spectral & 0.87 & 0.69 & 0.90 & 0.91 & 0.72\\
Bethe Hessian & 0.87 & 0.66 & 0.81 & -- & --\\
SCORE & 0.88 & 0.69 & 0.90 & -- & --\\
flat, CRP & 0.94 [20] & 0.81 [14] & 0.94 [23] & 0.96 [24] & 0.93 [20]\\
typed, CRP & 0.94 [18] & 0.82 [11] & 0.93 [17] & -- & --\\
global typed, CRP & 0.94 [18] & 0.75 [18] & 0.94 [33] & -- & --\\
global typed, global degree prior, CRP & 0.94 [18] & 0.75 [18] & 0.94 [30] & 0.98 [24] & --\\
global-DC, CRP & 0.93 [17] & 0.74 [15] & 0.94 [29] & 0.98 [24] & --\\
flat, Pitman--Yor & 0.94 [20] & 0.80 [12] & 0.94 [23] & -- & --\\
flat, Gnedin & 0.94 [20] & 0.82 [12] & 0.95 [22] & -- & --\\
flat, MFM & 0.94 [20] & 0.81 [13] & 0.95 [21] & -- & --\\
corrected BIC & 0.88 [13] & 0.69 [6] & 0.90 [15] & 0.95 [24] & 0.77 [20]\\
Bethe Hessian, $K$ inferred & 0.85 [9] & 0.69 [28] & 0.85 [17] & 0.90 [50] & 0.75 [50]\\
Leiden & 0.82 [5] & 0.88 [85] & 0.91 [36] & 0.94 [22] & 0.90 [527]\\
\bottomrule\end{tabular}}\end{table}

\begin{table}[tbp]\centering\scriptsize
\caption{Family effect on dense binary blocks, computed with exact evidences and no sampling. For each network the table gives the within-class and between-class densities of the metadata partition, the number of communities and the ARI against the metadata of each partition, and the log posterior $\log p(A\mid z)+\log p(z)$ of the metadata partition and of the MAP partitions found by three collapsed instances (best chain of four), evaluated on the full network under the homogeneous Beta--Bernoulli model, the homogeneous Gamma--Poisson model, which is the representation on which the typed and degree-corrected instances are built, and the type-collapsed model, with $K$ fixed to the number of metadata classes. Bold marks the highest value in each column. A higher value means a preferred partition.}
\label{tab:family}
\resizebox{\linewidth}{!}{\begin{tabular}{llrrrrr}\toprule
Network & partition & $\hat K$ & ARI & Beta--Bernoulli & Gamma--Poisson & typed \{ER,DC\}\\\midrule
football ($\hat p_{\mathrm{in}}=0.75$, $\hat p_{\mathrm{out}}=0.036$) & metadata & 12 & 1.00 & -1638.3 & -1892.1 & -1899.9\\
 & flat MAP & 11 & 0.86 & \textbf{-1541.4} & -1813.4 & -1821.0\\
 & typed MAP & 7 & 0.58 & -1620.8 & \textbf{-1791.2} & \textbf{-1796.0}\\
 & global-DC MAP & 10 & 0.81 & -1555.7 & -1798.3 & -1805.3\\
\midrule
school ($\hat p_{\mathrm{in}}=0.92$, $\hat p_{\mathrm{out}}=0.143$) & metadata & 11 & 1.00 & -10868.2 & -13145.1 & -13152.8\\
 & flat MAP & 11 & 0.44 & -10041.0 & \textbf{-12110.0} & \textbf{-12117.6}\\
 & typed MAP & 11 & 0.43 & -10138.7 & -12140.4 & -12148.0\\
 & global-DC MAP & 11 & 0.70 & \textbf{-10026.9} & -12354.2 & -12361.9\\
\midrule
polbooks ($\hat p_{\mathrm{in}}=0.17$, $\hat p_{\mathrm{out}}=0.021$) & metadata & 3 & 1.00 & -1424.0 & -1460.1 & -1439.9\\
 & flat MAP & 3 & 0.54 & -1331.9 & -1381.0 & -1381.0\\
 & typed MAP & 3 & 0.63 & \textbf{-1330.7} & \textbf{-1375.9} & \textbf{-1358.8}\\
 & global-DC MAP & 3 & 0.63 & \textbf{-1330.7} & \textbf{-1375.9} & \textbf{-1358.8}\\
\midrule
polblogs ($\hat p_{\mathrm{in}}=0.04$, $\hat p_{\mathrm{out}}=0.004$) & metadata & 2 & 1.00 & -74387.5 & -74702.7 & -56164.2\\
 & flat MAP & 2 & 0.00 & \textbf{-64010.9} & \textbf{-64867.9} & -63141.0\\
 & typed MAP & 2 & 0.40 & -72763.3 & -73196.5 & \textbf{-55038.5}\\
 & global-DC MAP & 2 & 0.77 & -73488.0 & -73822.2 & -55301.1\\
\bottomrule\end{tabular}}\end{table}
\begin{table}[tbp]\centering\scriptsize
\caption{Convergence diagnostics for the collapsed samplers on the real networks (four chains of $300$ sweeps from over-dispersed starts: the model-based candidate, Bethe Hessian spectral, uniformly random and Leiden): split-$\hat R$ of the log-posterior trace / mean pairwise ARI between the MAP partitions of the chains / range of the chains' ARIs against metadata, a single value meaning that the four ARIs coincide at two decimals.}
\label{tab:real_conv}
\resizebox{\linewidth}{!}{\begin{tabular}{lcccc}\toprule
Method & karate (2) & dolphins (2) & football (12) & polbooks (3)\\\midrule
flat & 1.00 / 1.00 / -0.01 & 1.01 / 1.00 / 0.87 & 2.01 / 0.93 / 0.82--0.86 & 1.03 / 0.81 / 0.54--0.63\\
flat (Jeffreys) & 1.00 / 1.00 / -0.01 & 1.00 / 1.00 / 0.87 & 1.04 / 1.00 / 0.86 & 1.01 / 1.00 / 0.54\\
typed \{ER,DC\} & 1.00 / 1.00 / -0.01 & 1.01 / 1.00 / 0.87 & 1.17 / 0.89 / 0.53--0.60 & 1.01 / 1.00 / 0.63\\
global typed & 1.01 / 1.00 / -0.01 & 1.02 / 1.00 / 0.87 & 1.00 / 1.00 / 0.81 & 2.35 / 0.81 / 0.47--0.63\\
global typed, global degree prior & 0.99 / 1.00 / 0.00 & 1.03 / 1.00 / 0.93 & 1.00 / 1.00 / 0.81 & 1.00 / 1.00 / 0.63\\
global-DC & 1.00 / 1.00 / 0.00 & 1.01 / 1.00 / 0.93 & 1.28 / 1.00 / 0.81 & 1.00 / 1.00 / 0.63\\
flat, CRP & 1.04 / 1.00 / 0.78 & 1.07 / 0.79 / 0.53--0.87 & 1.18 / 0.93 / 0.82--0.90 & 1.28 / 0.95 / 0.40--0.43\\
typed, CRP & 1.04 / 1.00 / -0.01 & 1.19 / 1.00 / 0.87 & 1.38 / 0.71 / 0.63--0.64 & 1.55 / 0.76 / 0.37--0.41\\
global typed, CRP & 1.12 / 0.73 / -0.01--0.06 & 1.26 / 0.62 / 0.56--0.87 & 2.78 / 0.92 / 0.81--0.86 & 1.52 / 0.68 / 0.24--0.34\\
global typed, global degree prior, CRP & 1.26 / 1.00 / 0.00 & 1.34 / 1.00 / 0.93 & 2.15 / 0.97 / 0.81--0.86 & 3.01 / 0.61 / 0.37--0.63\\
global-DC, CRP & 1.03 / 1.00 / 0.00 & 1.01 / 1.00 / 0.93 & 4.11 / 0.94 / 0.81 & 1.10 / 1.00 / 0.63\\
flat, Pitman--Yor & 1.01 / 1.00 / 0.78 & 1.00 / 0.72 / 0.53--0.87 & 1.91 / 0.95 / 0.86--0.90 & 1.26 / 0.89 / 0.34--0.43\\
flat, Gnedin & 1.01 / 1.00 / 0.78 & 1.07 / 1.00 / 0.53 & 2.13 / 0.97 / 0.86--0.90 & 1.29 / 0.94 / 0.40--0.43\\
flat, MFM & 1.02 / 1.00 / 0.78 & 1.29 / 0.79 / 0.53--0.87 & 1.77 / 0.96 / 0.86--0.90 & 1.04 / 1.00 / 0.40\\
\bottomrule\end{tabular}}\\[4pt]
\resizebox{\linewidth}{!}{\begin{tabular}{lcccc}\toprule
Method & polblogs (2) & school (11) & eu-core (42) & Cora (7)\\\midrule
flat & 1.00 / 0.99 / 0.00 & 131.05 / 0.61 / 0.44--0.92 & 5.96 / 0.80 / 0.16--0.18 & 7.20 / 0.22 / 0.05--0.28\\
flat (Jeffreys) & 1.01 / 0.99 / 0.00 & 112.28 / 0.59 / 0.44--0.92 & 3.69 / 0.73 / 0.17--0.20 & 7.61 / 0.25 / 0.09--0.24\\
typed \{ER,DC\} & 1.01 / 0.96 / 0.39 & 7.34 / 0.76 / 0.43--0.57 & 1.98 / 0.79 / 0.17--0.20 & 5.26 / 0.27 / 0.20--0.29\\
global typed & 189.15 / 0.49 / -0.00--0.78 & 27.19 / 0.85 / 0.64--0.75 & 1.09 / 0.68 / 0.43--0.51 & 6.75 / 0.30 / 0.14--0.19\\
global typed, global degree prior & 199.51 / 0.49 / 0.00--0.78 & 5.22 / 0.77 / 0.46--0.75 & 2.22 / 0.67 / 0.42--0.50 & 6.30 / 0.31 / 0.18--0.26\\
global-DC & 199.51 / 0.49 / 0.00--0.78 & 2.57 / 0.86 / 0.70--0.71 & 3.92 / 0.66 / 0.44--0.49 & 6.81 / 0.32 / 0.19--0.26\\
flat, CRP & 5.13 / 0.70 / 0.10--0.20 & 3.30 / 0.94 / 0.55--0.59 & 5.81 / 0.88 / 0.16--0.17 & 7.23 / 0.30 / 0.17--0.28\\
typed, CRP & 1.32 / 0.73 / 0.12--0.21 & 2.11 / 0.92 / 0.48--0.52 & 4.10 / 0.82 / 0.17--0.19 & 5.50 / 0.29 / 0.12--0.29\\
global typed, CRP & 1.80 / 0.52 / 0.14--0.15 & 8.55 / 0.80 / 0.48--0.52 & 4.31 / 0.60 / 0.42--0.48 & 5.30 / 0.35 / 0.15--0.30\\
global typed, global degree prior, CRP & 3.18 / 0.40 / 0.11--0.19 & 4.12 / 0.81 / 0.46--0.52 & 3.50 / 0.63 / 0.42--0.46 & 6.25 / 0.44 / 0.17--0.29\\
global-DC, CRP & 1.30 / 0.60 / 0.17--0.19 & 7.77 / 0.78 / 0.50--0.57 & 8.96 / 0.60 / 0.43--0.46 & 6.31 / 0.51 / 0.19--0.30\\
flat, Pitman--Yor & 4.16 / 0.73 / 0.09--0.19 & 5.26 / 0.95 / 0.55--0.59 & 6.89 / 0.78 / 0.16--0.17 & 4.92 / 0.32 / 0.14--0.30\\
flat, Gnedin & 2.53 / 0.70 / 0.09--0.20 & 10.69 / 0.95 / 0.55--0.59 & 4.48 / 0.87 / 0.16--0.17 & 4.73 / 0.34 / 0.12--0.29\\
flat, MFM & 3.97 / 0.70 / 0.09--0.20 & 1.72 / 1.00 / 0.59 & 19.46 / 0.86 / 0.17--0.18 & 4.22 / 0.37 / 0.20--0.29\\
\bottomrule\end{tabular}}\\[4pt]
\end{table}
\begin{table}[tbp]\centering\scriptsize
\caption{Start of the best chain, the chain with the highest log posterior among the four of Table~\ref{tab:real_conv}, for each collapsed instance and network.}
\label{tab:best_start}
\resizebox{\linewidth}{!}{\begin{tabular}{lcccccccc}\toprule
Instance & karate & dolphins & football & polbooks & polblogs & school & eu-core & Cora\\\midrule
flat & model-based & model-based & model-based & spectral & spectral & random & random & spectral\\
flat (Jeffreys) & model-based & model-based & model-based & model-based & spectral & random & random & random\\
typed \{ER,DC\} & model-based & model-based & model-based & spectral & Leiden & random & random & model-based\\
global typed & random & model-based & Leiden & model-based & model-based & random & random & Leiden\\
global typed, global degree prior & model-based & model-based & random & model-based & Leiden & Leiden & spectral & Leiden\\
global-DC & model-based & model-based & random & model-based & Leiden & random & Leiden & Leiden\\
flat, CRP & model-based & model-based & random & model-based & Leiden & random & random & model-based\\
typed, CRP & model-based & model-based & Leiden & model-based & random & spectral & random & model-based\\
global typed, CRP & random & Leiden & Leiden & Leiden & Leiden & Leiden & Leiden & Leiden\\
global typed, global degree prior, CRP & model-based & model-based & random & model-based & model-based & Leiden & Leiden & model-based\\
global-DC, CRP & model-based & model-based & model-based & model-based & random & model-based & Leiden & model-based\\
flat, Pitman--Yor & model-based & model-based & Leiden & model-based & Leiden & random & random & model-based\\
flat, Gnedin & model-based & model-based & Leiden & random & Leiden & random & random & model-based\\
flat, MFM & model-based & spectral & random & random & Leiden & model-based & random & model-based\\
\bottomrule\end{tabular}}\end{table}

\begin{table}[tbp]\centering\scriptsize
\caption{Long chains on the four larger networks: four chains of $1500$ sweeps from the same over-dispersed starts as Table~\ref{tab:real_conv}, $K$ fixed to the number of metadata classes. For each instance and network: ARI of the pooled MAP against the metadata, its log posterior on the training network, the gain in log posterior over the pooled MAP of the $300$-sweep chains, the split-$\hat R$ of the log-posterior traces and the mean pairwise ARI between the chains' MAP partitions.}
\label{tab:real_long}
\resizebox{\linewidth}{!}{\begin{tabular}{llccccc}\toprule
Network & instance & ARI & log posterior & gain over $300$ sweeps & split-$\hat R$ & chain agreement\\\midrule
polblogs & flat & 0.00 & -57722.3 & +0.0 & 1.00 & 1.00\\
 & typed \{ER,DC\} & 0.39 & -51234.6 & +0.5 & 1.00 & 0.95\\
 & global-DC & 0.77 & -134624.6 & -0.2 & 197.12 & 0.49\\
 & global typed & 0.76 & -49943.6 & -0.5 & 182.92 & 0.49\\
\midrule
school & flat & 0.44 & -9074.4 & +0.0 & 132.92 & 0.61\\
 & typed \{ER,DC\} & 0.50 & -10886.9 & +78.7 & 4.57 & 0.87\\
 & global-DC & 0.70 & -41699.4 & +0.0 & 8.23 & 0.86\\
 & global typed & 0.70 & -11579.6 & +10.2 & 1.63 & 0.85\\
\midrule
eu-core & flat & 0.16 & -40991.3 & +210.5 & 3.82 & 0.83\\
 & typed \{ER,DC\} & 0.18 & -44233.6 & +88.5 & 3.17 & 0.84\\
 & global-DC & 0.46 & -125952.8 & +16.1 & 3.30 & 0.71\\
 & global typed & 0.46 & -43710.2 & +43.0 & 1.48 & 0.70\\
\midrule
Cora & flat & 0.10 & -31476.4 & +170.9 & 7.69 & 0.34\\
 & typed \{ER,DC\} & 0.31 & -31369.0 & +90.1 & 3.92 & 0.29\\
 & global-DC & 0.29 & -35279.7 & +93.2 & 5.09 & 0.33\\
 & global typed & 0.21 & -31628.7 & +97.0 & 5.35 & 0.27\\
\bottomrule\end{tabular}}\end{table}

\begin{table}[tbp]\centering\scriptsize
\caption{Gain of the sampler over its initialiser on the real networks: the increase in collapsed log posterior (nats) of the pooled-MAP partition over the best model-based initialisation candidate, and the ARI against the metadata of that candidate followed by that of the pooled MAP. A gain of $0$ means that the reported partition is the polished initialiser.}
\label{tab:init_gain}
\begin{tabular}{lcccc}\toprule
Instance & karate & dolphins & football & polbooks\\\midrule
flat & +25 (0.88$\to$-0.01) & +3 (0.87$\to$0.87) & +9 (0.90$\to$0.86) & +12 (0.64$\to$0.54)\\
global-DC & +9 (0.88$\to$0.00) & +3 (0.87$\to$0.93) & +8 (0.71$\to$0.81) & +12 (0.64$\to$0.63)\\
typed & +22 (0.88$\to$-0.01) & +0 (0.87$\to$0.87) & +13 (0.71$\to$0.58) & +12 (0.64$\to$0.63)\\
flat, CRP & +50 (0.48$\to$0.78) & +36 (0.36$\to$0.87) & +27 (0.71$\to$0.86) & +64 (0.66$\to$0.40)\\
global-DC, CRP & +35 (0.48$\to$0.00) & +58 (0.36$\to$0.93) & +8 (0.71$\to$0.81) & +15 (0.66$\to$0.63)\\
typed, CRP & +46 (0.48$\to$-0.01) & +49 (0.36$\to$0.87) & +7 (0.71$\to$0.63) & +38 (0.66$\to$0.41)\\
\bottomrule\end{tabular}\\[4pt]
\begin{tabular}{lcccc}\toprule
Instance & polblogs & school & eu-core & Cora\\\midrule
flat & +6220 (0.11$\to$0.00) & +296 (0.91$\to$0.44) & +11773 (0.12$\to$0.16) & +1088 (0.10$\to$0.16)\\
global-DC & +119 (0.75$\to$0.77) & +208 (0.91$\to$0.70) & +3304 (0.31$\to$0.46) & +808 (0.10$\to$0.26)\\
typed & +297 (0.81$\to$0.40) & +494 (0.90$\to$0.43) & +9839 (0.31$\to$0.18) & +637 (0.29$\to$0.22)\\
flat, CRP & +12346 (0.04$\to$0.20) & +1980 (0.56$\to$0.59) & +14203 (0.16$\to$0.16) & +2019 (0.14$\to$0.28)\\
global-DC, CRP & +2503 (0.76$\to$0.19) & +570 (0.56$\to$0.53) & +3228 (0.31$\to$0.46) & +1194 (0.24$\to$0.25)\\
typed, CRP & +4865 (0.76$\to$0.19) & +1093 (0.56$\to$0.52) & +9920 (0.31$\to$0.18) & +1474 (0.26$\to$0.29)\\
\bottomrule\end{tabular}\\[4pt]
\end{table}

\end{document}